\documentclass[a4paper, oneside, 11 pt]{book}
\usepackage[T1]{fontenc} 
\usepackage[utf8]{inputenc}                 
\usepackage[english]{babel}                            
\usepackage{amsmath}
\usepackage{amsfonts}  
\usepackage{nicefrac}
\usepackage[standard,amsmath,thref,hyperref,thmmarks]{ntheorem}
\usepackage[backref=true, style=alphabetic, maxnames=10]{biblatex}
\usepackage{hyperref}
\usepackage{bm}
\usepackage{geometry}
\usepackage{frontespizio}
\usepackage{titletoc}
\usepackage{eqparbox}
\usepackage{amssymb}
\usepackage{thmtools}
\usepackage{tikz-cd}
\usepackage{mathrsfs}
\usepackage{enumerate}
\usepackage{multicol}
\usepackage{longtable}
\usepackage{pdfpages}
\usepackage{array}
\usepackage[all]{xy}
\usepackage{listings}
\usepackage{listings}
\usepackage{xcolor}
\definecolor{maccolor}{rgb}{0.3,0.3,0.8}
\lstdefinelanguage{Macaulay2}
{
basicstyle={\ttfamily},
keywordstyle={\color{maccolor!80!black}},
commentstyle={\color{gray}},
stringstyle={\color{red!40!black}},
rulecolor=\color{maccolor},
basewidth={1.2ex}, 
sensitive=false,
morecomment=[l]{--},
morecomment=[s]{-*}{*-},
morestring=[b]",
escapechar={`},
escapebegin={\rmfamily},
morekeywords={about,abs,AbstractToricVarieties,accumulate,Acknowledgement,acos,acosh,acot,addCancelTask,addDependencyTask,addEndFunction,addHook,AdditionalPaths,addStartFunction,addStartTask,Adjacent,adjoint,AdjointIdeal,AffineVariety,AfterEval,AfterNoPrint,AfterPrint,agm,AInfinity,alarm,AlgebraicSplines,Algorithm,Alignment,all,AllCodimensions,allowableThreads,ambient,analyticSpread,Analyzer,AnalyzeSheafOnP1,ancestor,ancestors,ANCHOR,and,andP,AngleBarList,ann,annihilator,antipode,any,append,applicationDirectory,applicationDirectorySuffix,apply,applyKeys,applyPairs,applyTable,applyValues,apropos,argument,Array,arXiv,Ascending,ascii,asin,asinh,ass,assert,associatedGradedRing,associatedPrimes,AssociativeAlgebras,AssociativeExpression,atan,atan2,atEndOfFile,Authors,autoload,AuxiliaryFiles,backtrace,Bag,Bareiss,baseFilename,BaseFunction,baseName,baseRing,baseRings,BaseRow,BasicList,basis,BasisElementLimit,Bayer,BeforePrint,beginDocumentation,BeginningMacaulay2,Benchmark,benchmark,Bertini,BesselJ,BesselY,betti,BettiCharacters,BettiTally,between,BGG,BIBasis,Binary,BinaryOperation,Binomial,binomial,BinomialEdgeIdeals,Binomials,BKZ,BlockMatrix,BLOCKQUOTE,BODY,Body,BoijSoederberg,BOLD,Book3264Examples,Boolean,BooleanGB,borel,Boxes,BR,break,Browse,Bruns,cache,CacheExampleOutput,CacheFunction,CacheTable,cacheValue,CallLimit,cancelTask,capture,catch,Caveat,CC,CDATA,ceiling,Center,centerString,Certification,ChainComplex,chainComplex,ChainComplexExtras,ChainComplexMap,ChainComplexOperations,ChangeMatrix,char,CharacteristicClasses,characters,charAnalyzer,check,CheckDocumentation,chi,Chordal,class,Classic,clean,clearAll,clearEcho,clearOutput,close,closeIn,closeOut,ClosestFit,CODE,code,codim,CodimensionLimit,coefficient,CoefficientRing,coefficientRing,coefficients,Cofactor,CohenEngine,CohenTopLevel,CoherentSheaf,CohomCalg,cohomology,coimage,CoincidentRootLoci,coker,cokernel,collectGarbage,columnAdd,columnate,columnMult,columnPermute,columnRankProfile,columnSwap,combine,Command,commandInterpreter,commandLine,COMMENT,commonest,commonRing,comodule,CompactMatrix,compactMatrixForm,CompiledFunction,CompiledFunctionBody,CompiledFunctionClosure,Complement,complement,complete,CompleteIntersection,CompleteIntersectionResolutions,Complexes,ComplexField,components,compose,compositions,compress,concatenate,conductor,ConductorElement,cone,Configuration,ConformalBlocks,conjugate,connectionCount,Consequences,Constant,Constants,constParser,content,continue,contract,Contributors,ConvexInterface,conwayPolynomial,ConwayPolynomials,copy,copyDirectory,copyFile,copyright,Core,CorrespondenceScrolls,cos,cosh,cot,CotangentSchubert,cotangentSheaf,coth,cover,coverMap,cpuTime,createTask,Cremona,csc,csch,current,currentColumnNumber,currentDirectory,currentFileDirectory,currentFileName,currentLayout,currentLineNumber,currentPackage,currentString,currentTime,Cyclotomic,Database,Date,DD,dd,deadParser,debug,debugError,DebuggingMode,debuggingMode,debugLevel,DecomposableSparseSystems,Decompose,decompose,deepSplice,Default,default,defaultPrecision,Degree,degree,degreeLength,DegreeLift,DegreeLimit,DegreeMap,DegreeOrder,DegreeRank,Degrees,degrees,degreesMonoid,degreesRing,delete,demark,denominator,Dense,Density,Depth,depth,Descending,Descent,Describe,describe,Description,det,determinant,DeterminantalRepresentations,DGAlgebras,diagonalMatrix,diameter,Dictionary,dictionary,dictionaryPath,diff,DiffAlg,difference,dim,directSum,disassemble,discriminant,dismiss,Dispatch,distinguished,DIV,Divide,divideByVariable,DivideConquer,DividedPowers,Divisor,DL,Dmodules,do,doc,docExample,docTemplate,document,DocumentTag,Down,drop,DT,dual,eagonNorthcott,EagonResolution,echoOff,echoOn,EdgeIdeals,edit,EigenSolver,eigenvalues,eigenvectors,eint,EisenbudHunekeVasconcelos,elapsedTime,elapsedTiming,elements,Eliminate,eliminate,Elimination,EliminationMatrices,EllipticCurves,EllipticIntegrals,else,EM,Email,End,end,endl,endPackage,Engine,engineDebugLevel,EngineRing,EngineTests,entries,EnumerationCurves,environment,Equation,EquivariantGB,erase,erf,erfc,error,errorDepth,euler,EulerConstant,eulers,even,EXAMPLE,ExampleFiles,ExampleItem,examples,ExampleSystems,Exclude,exec,exit,exp,expectedReesIdeal,expm1,exponents,export,exportFrom,exportMutable,Expression,expression,Ext,extend,ExteriorIdeals,ExteriorModules,exteriorPower,Factor,factor,false,Fano,FastMinors,FastNonminimal,FGLM,File,fileDictionaries,fileExecutable,fileExists,fileExitHooks,fileLength,fileMode,FileName,FilePosition,fileReadable,fileTime,fileWritable,fillMatrix,findFiles,findHeft,FindOne,findProgram,findSynonyms,FiniteFittingIdeals,First,first,firstkey,FirstPackage,fittingIdeal,flagLookup,FlatMonoid,flatten,flattenRing,Flexible,flip,floor,flush,fold,FollowLinks,for,forceGB,fork,FormalGroupLaws,Format,format,formation,FourierMotzkin,FourTiTwo,fpLLL,frac,fraction,FractionField,frames,FrobeniusThresholds,from,fromDividedPowers,fromDual,Function,FunctionApplication,FunctionBody,functionBody,FunctionClosure,FunctionFieldDesingularization,fusePairs,futureParser,GaloisField,Gamma,gb,GBDegrees,gbRemove,gbSnapshot,gbTrace,gcd,gcdCoefficients,gcdLLL,GCstats,genera,GeneralOrderedMonoid,GenerateAssertions,generateAssertions,generator,generators,Generic,GenericInitialIdeal,genericMatrix,genericSkewMatrix,genericSymmetricMatrix,gens,genus,get,getc,getChangeMatrix,getenv,getGlobalSymbol,getNetFile,getNonUnit,getPrimeWithRootOfUnity,getSymbol,getWWW,GF,gfanInterface,Givens,GKMVarieties,GLex,Global,global,globalAssign,globalAssignFunction,GlobalAssignHook,globalAssignment,globalAssignmentHooks,GlobalDictionary,GlobalHookStore,globalReleaseFunction,GlobalReleaseHook,Gorenstein,GradedLieAlgebras,GradedModule,gradedModule,GradedModuleMap,gradedModuleMap,gramm,GraphicalModels,GraphicalModelsMLE,Graphics,graphIdeal,graphRing,Graphs,Grassmannian,GRevLex,GroebnerBasis,groebnerBasis,GroebnerBasisOptions,GroebnerStrata,GroebnerWalk,groupID,GroupLex,GroupRevLex,GTZ,Hadamard,handleInterrupts,HardDegreeLimit,hash,HashTable,hashTable,HEAD,HEADER1,HEADER2,HEADER3,HEADER4,HEADER5,HEADER6,HeaderType,Heading,Headline,Heft,heft,Height,height,help,Hermite,hermite,Hermitian,HH,hh,HigherCIOperators,HighestWeights,Hilbert,hilbertFunction,hilbertPolynomial,hilbertSeries,HodgeIntegrals,hold,Holder,Hom,homeDirectory,HomePage,Homogeneous,Homogeneous2,homogenize,homology,homomorphism,HomotopyLieAlgebra,hooks,horizontalJoin,HorizontalSpace,HR,HREF,HTML,html,httpHeaders,Hybrid,HyperplaneArrangements,Hypertext,hypertext,HypertextContainer,HypertextParagraph,icFracP,icFractions,icMap,icPIdeal,id,Ideal,ideal,idealizer,identity,if,IgnoreExampleErrors,ii,image,imaginaryPart,IMG,ImmutableType,importFrom,in,incomparable,Increment,independentSets,indeterminate,IndeterminateNumber,Index,index,indexComponents,IndexedVariable,IndexedVariableTable,indices,inducedMap,inducesWellDefinedMap,InexactField,InexactFieldFamily,InexactNumber,InfiniteNumber,infinity,info,InfoDirSection,infoHelp,Inhomogeneous,input,Inputs,insert,installAssignmentMethod,installedPackages,installHilbertFunction,installMethod,installMinprimes,installPackage,InstallPrefix,instance,instances,IntegralClosure,integralClosure,integrate,IntermediateMarkUpType,interpreterDepth,intersect,intersectInP,Intersection,intersection,interval,InvariantRing,inverse,InverseMethod,inversePermutation,Inverses,inverseSystem,InverseSystems,Invertible,InvolutiveBases,irreducibleCharacteristicSeries,irreducibleDecomposition,isAffineRing,isANumber,isBorel,isCanceled,isCommutative,isConstant,isDirectory,isDirectSum,isEmpty,isField,isFinite,isFinitePrimeField,isFreeModule,isGlobalSymbol,isHomogeneous,isIdeal,isInfinite,isInjective,isInputFile,isIsomorphism,isLinearType,isListener,isLLL,isMember,isModule,isMonomialIdeal,isNormal,isOpen,isOutputFile,isPolynomialRing,isPrimary,isPrime,isPrimitive,isPseudoprime,isQuotientModule,isQuotientOf,isQuotientRing,isReady,isReal,isReduction,isRegularFile,isRing,isSkewCommutative,isSorted,isSquareFree,isStandardGradedPolynomialRing,isSubmodule,isSubquotient,isSubset,isSupportedInZeroLocus,isSurjective,isTable,isUnit,isWellDefined,isWeylAlgebra,ITALIC,Iterate,Jacobian,jacobian,jacobianDual,Jets,Join,join,Jupyter,K3Carpets,K3Surfaces,Keep,KeepFiles,KeepZeroes,ker,kernel,kernelLLL,kernelOfLocalization,Key,keys,Keyword,Keywords,kill,koszul,Kronecker,KustinMiller,LABEL,last,lastMatch,LATER,LatticePolytopes,Layout,lcm,leadCoefficient,leadComponent,leadMonomial,leadTerm,Left,left,length,LengthLimit,letterParser,Lex,LexIdeals,LI,Licenses,LieTypes,lift,liftable,Limit,limitFiles,limitProcesses,Linear,LinearAlgebra,LinearTruncations,lineNumber,lines,LINK,linkFile,List,list,listForm,listLocalSymbols,listSymbols,listUserSymbols,LITERAL,LLL,LLLBases,lngamma,load,loadDepth,LoadDocumentation,loadedFiles,loadedPackages,loadPackage,Local,local,localDictionaries,LocalDictionary,localize,LocalRings,locate,log,log1p,LongPolynomial,lookup,lookupCount,LowerBound,LUdecomposition,M0nbar,M2CODE,Macaulay2Doc,makeDirectory,MakeDocumentation,makeDocumentTag,MakeHTML,MakeInfo,MakeLinks,makePackageIndex,MakePDF,makeS2,Manipulator,map,MapExpression,MapleInterface,markedGB,Markov,MarkUpType,match,mathML,Matrix,matrix,MatrixExpression,Matroids,max,maxAllowableThreads,maxExponent,MaximalRank,maxPosition,MaxReductionCount,MCMApproximations,member,memoize,memoizeClear,memoizeValues,MENU,merge,mergePairs,META,method,MethodFunction,MethodFunctionBinary,MethodFunctionSingle,MethodFunctionWithOptions,methodOptions,methods,midpoint,min,minExponent,mingens,mingle,minimalBetti,MinimalGenerators,MinimalMatrix,minimalPresentation,minimalPresentationMap,minimalPresentationMapInv,MinimalPrimes,minimalPrimes,minimalReduction,Minimize,minimizeFilename,MinimumVersion,minors,minPosition,minPres,minprimes,Minus,minus,Miura,MixedMultiplicity,mkdir,mod,Module,module,ModuleDeformations,modulo,MonodromySolver,Monoid,monoid,MonoidElement,Monomial,MonomialAlgebras,monomialCurveIdeal,MonomialIdeal,monomialIdeal,MonomialIntegerPrograms,MonomialOrbits,MonomialOrder,Monomials,monomials,MonomialSize,monomialSubideal,moveFile,multidegree,multidoc,multigraded,MultigradedBettiTally,MultiGradedRationalMap,multiplicity,MultiplicitySequence,MultiplierIdeals,MultiplierIdealsDim2,MultiprojectiveVarieties,mutable,MutableHashTable,mutableIdentity,MutableList,MutableMatrix,mutableMatrix,NAGtypes,Name,nanosleep,Nauty,NautyGraphs,NCAlgebra,NCLex,needs,needsPackage,Net,net,NetFile,netList,new,newClass,newCoordinateSystem,NewFromMethod,newline,NewMethod,newNetFile,NewOfFromMethod,NewOfMethod,newPackage,newRing,nextkey,nextPrime,nil,NNParser,NoetherianOperators,NoetherNormalization,NonAssociativeProduct,NonminimalComplexes,nonspaceAnalyzer,NoPrint,norm,normalCone,Normaliz,NormalToricVarieties,not,Nothing,notify,notImplemented,NTL,null,nullaryMethods,nullhomotopy,nullParser,nullSpace,Number,number,NumberedVerticalList,numcols,numColumns,numerator,numeric,NumericalAlgebraicGeometry,NumericalCertification,NumericalImplicitization,NumericalLinearAlgebra,NumericalSchubertCalculus,numericInterval,NumericSolutions,numgens,numRows,numrows,odd,oeis,of,ofClass,OL,OldPolyhedra,OldToricVectorBundles,on,OneExpression,OnlineLookup,OO,oo,ooo,oooo,openDatabase,openDatabaseOut,openFiles,openIn,openInOut,openListener,OpenMath,openOut,openOutAppend,operatorAttributes,Option,OptionalComponentsPresent,optionalSignParser,Options,options,OptionTable,optP,or,Order,order,OrderedMonoid,orP,OutputDictionary,Outputs,override,pack,Package,package,PackageCitations,PackageDictionary,PackageExports,PackageImports,PackageTemplate,packageTemplate,pad,pager,PairLimit,pairs,PairsRemaining,PARA,Parametrization,parent,Parenthesize,Parser,Parsing,part,Partition,partition,partitions,parts,path,pdim,peek,PencilsOfQuadrics,Permanents,permanents,permutations,pfaffians,PHCpack,PhylogeneticTrees,pi,PieriMaps,pivots,PlaneCurveSingularities,plus,poincare,poincareN,Points,polarize,poly,Polyhedra,Polymake,PolynomialRing,Posets,Position,position,positions,PositivityToricBundles,POSIX,Postfix,Power,power,powermod,PRE,Precision,precision,Prefix,prefixDirectory,prefixPath,preimage,prepend,presentation,pretty,primaryComponent,PrimaryDecomposition,primaryDecomposition,PrimaryTag,PrimitiveElement,Print,print,printerr,printingAccuracy,printingLeadLimit,printingPrecision,printingSeparator,printingTimeLimit,printingTrailLimit,printString,printWidth,processID,Product,product,ProductOrder,profile,profileSummary,Program,programPaths,ProgramRun,Proj,Projective,ProjectiveHilbertPolynomial,projectiveHilbertPolynomial,ProjectiveVariety,promote,protect,Prune,prune,PruneComplex,pruningMap,Pseudocode,pseudocode,pseudoRemainder,Pullback,PushForward,pushForward,Python,QQ,QQParser,QRDecomposition,QthPower,Quasidegrees,QuaternaryQuartics,QuillenSuslin,quit,Quotient,quotient,quotientRemainder,QuotientRing,Radical,radical,RadicalCodim1,radicalContainment,RaiseError,random,RandomCanonicalCurves,RandomComplexes,RandomCurves,RandomCurvesOverVerySmallFiniteFields,RandomGenus14Curves,RandomIdeals,randomKRationalPoint,RandomMonomialIdeals,randomMutableMatrix,RandomObjects,RandomPlaneCurves,RandomPoints,RandomSpaceCurves,Range,rank,RationalMaps,RationalPoints,RationalPoints2,ReactionNetworks,read,readDirectory,readlink,readPackage,RealField,RealFP,realPart,realpath,RealQP,RealQP1,RealRoots,RealRR,RealXD,recursionDepth,recursionLimit,Reduce,reducedRowEchelonForm,reduceHilbert,reductionNumber,ReesAlgebra,reesAlgebra,reesAlgebraIdeal,reesIdeal,References,ReflexivePolytopesDB,regex,regexQuote,registerFinalizer,regSeqInIdeal,Regularity,regularity,relations,RelativeCanonicalResolution,relativizeFilename,Reload,remainder,RemakeAllDocumentation,remove,removeDirectory,removeFile,removeLowestDimension,reorganize,replace,RerunExamples,res,reshape,ResidualIntersections,ResLengthThree,Resolution,resolution,ResolutionsOfStanleyReisnerRings,restart,Result,resultant,Resultants,return,returnCode,Reverse,reverse,RevLex,Right,right,Ring,ring,RingElement,RingFamily,ringFromFractions,RingMap,rootPath,roots,rootURI,rotate,round,rowAdd,RowExpression,rowMult,rowPermute,rowRankProfile,rowSwap,RR,RRi,rsort,run,RunDirectory,RunExamples,RunExternalM2,runHooks,runLengthEncode,runProgram,same,saturate,Saturation,scan,scanKeys,scanLines,scanPairs,scanValues,schedule,schreyerOrder,Schubert,Schubert2,SchurComplexes,SchurFunctors,SchurRings,SCRIPT,scriptCommandLine,ScriptedFunctor,SCSCP,searchPath,sec,sech,SectionRing,SeeAlso,seeParsing,SegreClasses,select,selectInSubring,selectVariables,SelfInitializingType,SemidefiniteProgramming,Seminormalization,separate,SeparateExec,separateRegexp,Sequence,sequence,Serialization,serialNumber,Set,set,setEcho,setGroupID,setIOExclusive,setIOSynchronized,setIOUnSynchronized,setRandomSeed,setup,setupEmacs,sheaf,SheafExpression,sheafExt,sheafHom,SheafOfRings,shield,ShimoyamaYokoyama,short,show,showClassStructure,showHtml,showStructure,showTex,showUserStructure,SimpleDoc,simpleDocFrob,SimplicialComplexes,SimplicialDecomposability,SimplicialPosets,SimplifyFractions,sin,singularLocus,sinh,size,size2,SizeLimit,SkewCommutative,SlackIdeals,sleep,SLnEquivariantMatrices,SLPexpressions,SMALL,smithNormalForm,solve,someTerms,Sort,sort,sortColumns,SortStrategy,source,SourceCode,SourceRing,SPACE,SpaceCurves,SPAN,span,SparseMonomialVectorExpression,SparseResultants,SparseVectorExpression,Spec,SpechtModule,SpecialFanoFourfolds,specialFiber,specialFiberIdeal,SpectralSequences,splice,splitWWW,sqrt,SRdeformations,stack,stacksProject,Standard,standardForm,standardPairs,StartWithOneMinor,stashValue,StatePolytope,StatGraphs,status,stderr,stdio,step,StopBeforeComputation,stopIfError,StopWithMinimalGenerators,Strategy,String,STRONG,StronglyStableIdeals,STYLE,Style,style,SUB,sub,SubalgebraBases,sublists,submatrix,submatrixByDegrees,Subnodes,subquotient,SubringLimit,Subscript,subscript,SUBSECTION,subsets,substitute,substring,subtable,Sugarless,Sum,sum,SumOfTwists,SumsOfSquares,SUP,super,SuperLinearAlgebra,Superscript,superscript,support,SVD,SVDComplexes,switch,SwitchingFields,sylvesterMatrix,Symbol,symbol,SymbolBody,symbolBody,SymbolicPowers,symlinkDirectory,symlinkFile,symmetricAlgebra,symmetricAlgebraIdeal,symmetricKernel,SymmetricPolynomials,symmetricPower,synonym,SYNOPSIS,syz,Syzygies,SyzygyLimit,SyzygyMatrix,SyzygyRows,syzygyScheme,TABLE,Table,table,take,Tally,tally,tan,TangentCone,tangentCone,tangentSheaf,tanh,target,Task,taskResult,TateOnProducts,TD,temporaryFileName,tensor,tensorAssociativity,TensorComplexes,terminalParser,terms,TEST,Test,testExample,testHunekeQuestion,TestIdeals,TestInput,tests,TEX,tex,TeXmacs,texMath,Text,TH,then,Thing,ThinSincereQuivers,ThreadedGB,threadVariable,Threshold,throw,Time,time,times,timing,TITLE,TO,to,TO2,toAbsolutePath,toCC,toDividedPowers,toDual,toExternalString,toField,TOH,toList,toLower,top,top,topCoefficients,Topcom,topComponents,topLevelMode,Tor,TorAlgebra,Toric,ToricInvariants,ToricTopology,ToricVectorBundles,toRR,toRRi,toSequence,toString,TotalPairs,toUpper,TR,trace,transpose,TriangularSets,Tries,Trim,trim,Triplets,Tropical,true,Truncate,truncate,truncateOutput,Truncations,try,TSpreadIdeals,TT,tutorial,Type,TypicalValue,typicalValues,UL,ultimate,unbag,uncurry,Undo,undocumented,uniform,uninstallAllPackages,uninstallPackage,Unique,unique,Units,Unmixed,unsequence,unstack,Up,UpdateOnly,UpperTriangular,URL,urlEncode,Usage,use,UseCachedExampleOutput,UseHilbertFunction,UserMode,userSymbols,UseSyzygies,utf8,utf8check,validate,value,values,Variable,VariableBaseName,Variables,Variety,variety,vars,Vasconcelos,Vector,vector,VectorExpression,VectorFields,VectorGraphics,Verbose,Verbosity,Verify,VersalDeformations,versalEmbedding,Version,version,VerticalList,VerticalSpace,viewHelp,VirtualResolutions,VirtualTally,VisibleList,Visualize,wait,WebApp,wedgeProduct,weightRange,Weights,WeylAlgebra,WeylGroups,when,whichGm,while,width,wikipedia,Wrap,wrap,WrapperType,XML,xor,youngest,zero,ZeroExpression,zeta,ZZ,ZZParser}
}
\lstalias{Macaulay2output}{Macaulay2}

\newtheorem{te}{Theorem}[chapter]
\newtheorem{pr}[te]{Proposition}
\newtheorem{lem}[te]{Lemma}
\newtheorem{co}[te]{Corollary}
\newtheorem{con}[te]{Conjecture}

\newtheorem{de}[te]{Definition}
\newtheorem{no}[te]{Notation}
\newtheorem{es}[te]{Example}
\newtheorem{Os}[te]{Remark}
\newtheorem{Qu}[te]{Question}

\newcommand{\N}{\mathbb{N}}
\newcommand{\Z}{\mathbb{Z}}
\newcommand{\Q}{\mathbb{Q}}

\newcommand{\C}{\mathbb{C}}

\newcommand{\p}{\mathbb{P}}
\newcommand{\os}{\mathcal{O}}
\newcommand{\osn}{\os_{\p^n}}
\newcommand{\e}{\mathcal{E}}
\newcommand{\f}{\mathcal{F}}
\newcommand{\g}{\mathcal{G}}
\newcommand{\id}{\mathcal{I}}
\newcommand{\cb}{\mathcal{B}}
\newcommand{\cc}{\mathcal{C}}
\newcommand{\cd}{\mathcal{D}}
\newcommand{\ci}{\mathcal{I}}
\newcommand{\h}{\mathcal{H}}
\newcommand{\ck}{\mathcal{K}}
\newcommand{\cl}{\mathcal{L}}
\newcommand{\cm}{\mathcal{M}}
\newcommand{\n}{\mathcal{N}}
\newcommand{\cq}{\mathcal{Q}}
\newcommand{\car}{\mathcal{R}}
\newcommand{\ct}{\mathcal{T}}
\newcommand{\cu}{\mathcal{U}}
\newcommand{\cV}{\mathcal{V}}
\newcommand{\cX}{\mathcal{X}}
\newcommand{\cY}{\mathcal{Y}}
\newcommand{\cW}{\mathcal{W}}

\newcommand{\scr}[1]{\mathscr{#1}}
\newcommand{\gel}[1]{\underline{\textbf{#1)}}}
\newcommand{\du}{^\vee}
\newcommand{\dd}{^{\vee\vee}}

\newcommand{\ra}{\rightarrow}

\newcommand{\xra}{\xrightarrow}
\newcommand{\mt}{\mapsto}

\newcommand{\inv}{^{-1}}
\newcommand{\til}{\widetilde}

\newcommand{\gk}{\mathbf{k}}
\newcommand{\gi}{\mathbf{I}}

\newcommand{\sesl}[4]{\begin{equation}\label{#4} 0\ra #1\ra #2\ra #3 \ra 0\end{equation}}
\newcommand{\ses}[3]{\begin{equation*} 0\ra #1\ra #2\ra #3 \ra 0 \end{equation*}}

\newcommand{\on}[1]{\os_{\mathbb{P}^{#1}}}

\newcommand{\ext}{\mathscr{E}xt}
\newcommand{\Hom}{\mathscr{H}om}

\newcommand{\De}{\begin{de}}
\newcommand{\Ne}{\end{de}}
\newcommand{\Te}{\begin{te}}
\newcommand{\Ma}{\end{te}}
\newcommand{\Prop}{\begin{pr}}
\newcommand{\One}{\end{pr}}
\newcommand{\Es}{\begin{es}}
\newcommand{\Io}{\end{es}}
\newcommand{\Le}{\begin{lem}}
\newcommand{\ma}{\end{lem}}
\newcommand{\Co}{\begin{co}}
\newcommand{\io}{\end{co}}
\newcommand{\Con}{\begin{con}}
\newcommand{\ura}{\end{con}}
\newcommand{\Oss}{\begin{Os}}
\newcommand{\one}{\end{Os}}
\newcommand{\Do}{\begin{Qu}}
\newcommand{\da}{\end{Qu}}
\newcommand{\No}{\begin{no}}
\newcommand{\ione}{\end{no}}

\newcommand{\dia}{ \begin{center}\begin{tikzcd}}
\newcommand{\mma}{\end{tikzcd}\end{center}}

\hypersetup{
    colorlinks,
    linkcolor={green!80!blue},
    citecolor={blue!90!black},
    urlcolor={blue!70!black}
}

\begin{document}

\includepdf[pages=-]{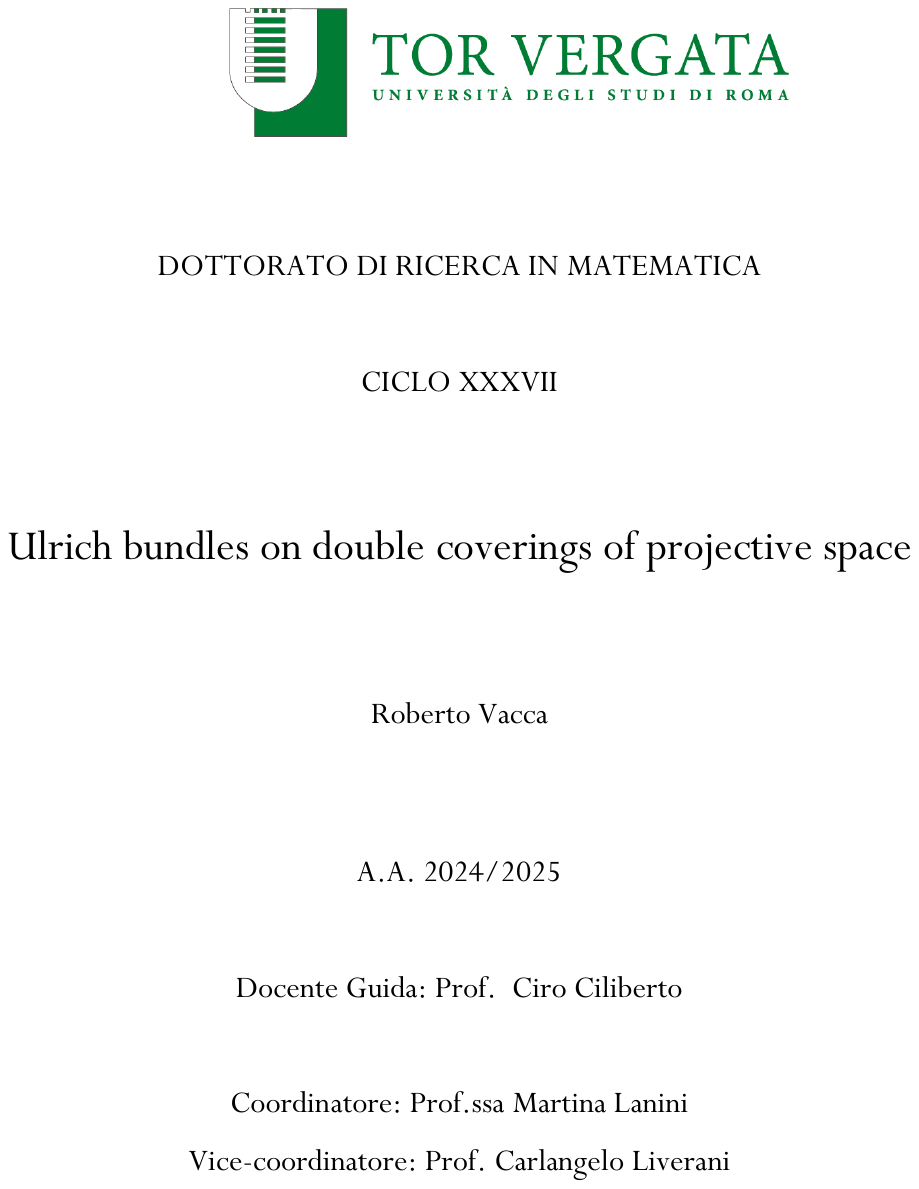}

\newpage

\chapter*{Acknowledgments}

This work is my PhD thesis.
I warmly thank Ciro Ciliberto, my advisor, for his care and patience during those years.
All his suggestions and encouragements greatly improved the mathematical content as well as the readability of this text.
\noindent Moreover, I thank Nelson Alvarado, Valerio Buttinelli, Davide Gori, Angelo Lopez and Antonio Rapagnetta for many useful conversations.

Part of this work have been done during my 3 months stay in University Paul Sabatier in Toulouse, so I thank Thomas Dedieu and Laurent Manivel for their hospitality, and Fulvio Gesmundo for suggesting how to use Macaulay2 for \thref{join}.

\chapter*{Abstract}

Fixed a polarised variety $X$, we can ask if it admits Ulrich bundles and, in case, what is their minimal possible rank. In this thesis, after recalling general properties of Ulrich sheaves, we show that any finite covering of $\mathbb{P}^n$ that embeds as a divisor in a weighted projective space with weights $(1^{n+1},m)$ admits Ulrich sheaves, by using matrix factorisations. Among these varieties, we focus on double coverings of with $n\ge3$. Through Hartshorne--Serre correspondence, which we review along the way, we prove that the general such $X$ admits a rank $2$ Ulrich sheaf if and only if $n=3$ and $m=2,3,4$, and characterise the zero loci of their sections. Moreover, we construct generically smooth components of the expected dimension of their moduli spaces, analyse the action of the natural involution on them and the restriction of those bundles to low degree hypersurfaces. For $m=2,3$, we verify the existence of slope-stable Ulrich bundles of all the possible ranks.

\tableofcontents

\newpage

\chapter*{Introduction}
\addcontentsline{toc}{chapter}{Introduction}

The study of vector bundles is of fundamental importance to understand an algebraic variety $X$, for example it encodes information on the morphism from $X$ to projective spaces or, more generally, to Grassmannians.
A cornerstone in this area is the dictionary between algebraic vector bundles and finitely generated projective modules over commutative rings created by Serre in \cite{Serre_fac}. 
Despite the attention this topic received in the following years, we still lack a good understanding already for vector bundles on projective spaces.
For example, in 1974 Hartshorne \cite{Har_small_codim}[Conj. 6.3] conjectured that any vector bundle of rank $2$ on $\p^n$ for $n\ge 7$\footnote{actually even of $\p^6$ there are no examples of indecomposable rank $2$ vector bundles} is actually split; as of today this is wide open.

In this work, we study a special class of vector bundles, or in general coherent sheaves: Ulrich sheaves.

\noindent\textbf{Definition.}
    A coherent sheaf $\e$ on an $n$-dimensional proper scheme $X$ endowed with an ample and globally generated line bundle $H$ is said to be Ulrich if 
\[h^j(X,\e(iH))=0 \qquad \text{for} \;\; 0\leq j\leq n,\, -n\leq i\leq -1.\]
This condition is quite a strong one and has many remarkable consequence on the properties of the underlying variety admitting such a sheaf: determinantal representations of Chow form, Boij-Soderberg theory of cone of cohomology tables, minimal rank conjecture for resolutions of ideals of points...
One major point in the interest in them comes from the Eisenbud and Schreyer's question, formulated in 2003 in \cite{EisSch}, regarding their existence on arbitrary projective varieties.

\section*{Origins}

Ulrich sheaves/modules are at the crossroads of algebraic geometry and commutative algebra.
They were first introduced in 1984 by Bernd Ulrich in the realm of commutative algebra in \cite{Ulrich} in order to study properties of Cohen-Macaulay rings.
Specifically, it was already known, see \cite{Sally}, a bound for the maximal number of elements of a minimal set of generators for a finitely generated Cohen-Macaulay module $M$ of positive rank over a local Cohen-Macaulay ring $R$.
In \cite{Ulrich}[Thm. 3.1], Ulrich found a criterion to check whether $R$ is Gorenstein by just checking $Ext^i_R(M,R)=0$ for $0\le i\le dim(R)$ where $M$ is a fixed $R$-module with a sufficiently high number of generators.
Therefore, he was led to ask if modules realising the optimal bound, for this reason sometimes called \textit{maximally generated maximal Cohen-Macaulay} modules, exist on any such ring $R$.

Only a few examples of Ulrich modules, on $1$ or $2$-dimensional rings, were known before matrix factorisations entered the story.
In 1980 Eisenbud showed that Cohen-Macaulay modules of maximal dimension over local rings $R/(h)$ with $R$ regular ring and $h$ not a divisor of zero (\textit{hypersurface rings}) are in correspondence with matrix factorisations of powers of $h$, see \cite{Eis-resolutions}.
A matrix factorisation for $h$ is simply the datum of some matrices whose product is $h$ times the identity matrix of some order.
The above connection carries over on graded rings, see \cite{Backelin_Herzog_Sanders} and \cite{Brenna_Herzog_Ulrich}, where, if $R$ is a polynomial ring, it is shown that the above maximality property of $M$ is equivalent to the existence of a linear resolution. 
In \cite{BHU} the existence of matrix factorizations has been shown under quite general assumptions on the ring $R$, giving the existence of Ulrich modules on complete intersections in projective space.

In 2000, Beauville proved the equivalence between the datum of an aCM vector bundle on a hypersurface in $\p^n$ and determinantal representation of its equation, see \cite{Beadeterminantal}, essentially rediscovering in a geometric context the previous results on modules. 
Then, in 2003 in \cite{EisSch} a sheaf-theoretic definition has been given for the sheaves corresponding to Ulrich modules on projective schemes and the existence problem was posed.
Their main motivation was the fact that Ulrich sheaves on embedded schemes $X\subset\p^N$ admit, analogously to the module case, linear resolutions and from those Eisenbud and Schreyer construct a determinantal presentation for the Chow form of $X$.
Moreover, they showed existence of Ulrich sheaves on any polarised curve and on Segre--Veronese varieties.

After this paper, the search for Ulrich bundles and the investigation of their properties received much attention until nowadays, especially on surfaces and $3$-folds.
For example, it has been shown that the following classes of varieties admit Ulrich bundles: complete intersections in projective space, Grassmannians, surfaces of Kodaira dimension $\le 0$, Fano $3$-folds of index $\le 2$ or of index $1$ and cyclic Picard group.  
Among many contributions, we cite the fact that those sheaves are always Gieseker-semistable, proved in \cite{CH}, and the fact that on any smooth complex surface there are polarisations admitting rank $2$ Ulrich bundles, see \cite{CoskunHuizenga}.
Nevertheless, this theory is still quite mysterious already on $3$-folds and not much is known in higher dimension except for the very special varieties cited above.
Furthermore, in many cases we only know existence of such sheaves and have very raw estimates fort their minimal rank.

\section*{Our work}

We describe in some detail the ideas and results contained in this work.
For the sake of presentation, we will give some statements under non-optimal assumptions, but leave references to more precise statements in the full text.

In this thesis we study Ulrich bundles with respect to an ample and globally generated polarisation therefore, rather then seeing our proper, $n$-dimensional scheme $X$ as a subscheme of $\p^N$ we will see it as a finite covering $f:X\ra\p^n$, in the spirit of Noether normalisation.
In particular, in this setting we always set $H$ to be the pullback of some hyperplane in $\p^n$.
Aside from a few exceptions, see, for example, \cite{MRPL}, only recently the case of ample and globally generated polarisations was considered systematically, as in the survey \cite{AC} or in \cite{Valerio+}.
In particular, there is a lack of results linking the presence of Ulrich sheaves with the geometry of the cover $f$, which we will try to fill in some very special cases like cyclic coverings. 

On one side our weaker assumption grants a cleaner functoriality under finite pushforwards, see \thref{proiezione}, since finite pullabcks preserve ampleness but no very ampleness in general, but on the other hand we only get a weakening of the linear resolution property, see \thref{definizioneUlrich}, in particular we lose the determinantal presentation given in \cite{EisSch}.
Our point of view is that this is not much of a loss, since it is well known that from Ulrich sheaves for some polarisation $H$ we can construct Ulrich sheaves for all multiples $mH$ as shown in \thref{polarizzazionimultiple}, see also \autoref{veronese}, but actually presents a new geometric picture to study.
In our setting, instead of looking for properties of the embedding $X\subset \p^N$ we should focus on properties of the finite map $f:X\ra\p^n$, that is the structure of $f_*\os_X$.
For example, we observe in \thref{restrizionesezioni} that zero loci of an Ulrich sheaf cannot contain fibers of $f$, a restriction that becomes quite severe when this map has a small degree.

A remarkable case is that of cyclic coverings, for which the existence of Ulrich sheaves has already been shown in \cite{KuNaPa} for the degree $2$ case and in \cite{ParPin} for arbitrary degrees.
Actually, this result can also be deduced from \cite{HanselkaKummer}[Thm. 8.1] using the existence of matrix factorisations, established in \cite{BHU}.
Following this approach, we get a further generalisation.
We call \textit{divisorial} a finite covering $f:X\ra\p^n$ such that $f$ factors as $X\subset\p(\osn\oplus\osn(m))$ followed by the standard projection to $\p^n$ or, equivalently, through a weighted projective space with weights $(1^{n+1},m)$.
Such varieties are divisors, even though in some modification of the usual projective space, and therefore can be essentially treated as zeroes of a single polynomial.

\noindent\textbf{Theorem}[\thref{ulrichmatrix''}]
\textit{Let $f:X\ra\p^n$ be a divisorial covering then $X$ admits an Ulrich sheaf respect to $f^*\osn(1)$.}

\medskip

Note that this result holds over an arbitrary field.
In fact, it is also more accurate, giving an upper bound for the minimal rank depending on the equation of $X$ and some properties of the field.

As a corollary, we recover many previous results.
More interestingly, combining these results with our generalisation of an older fact proved in \cite{Casnati_wild}, see \thref{pullback}, we also get some new instances of existence.

\noindent\textbf{Corollary}[\thref{divci}, \thref{horikawa'}]
\textit{Anytime a finite surjective morphism $X\ra X_1$, with $X_1$ a complete intersection of $\p^n$, can be written as a pullback $X_1\times_{\p^n}X_2$ where $X_2\ra\p^n$ is a divisorial covering, then $X$ admits Ulrich sheaves.
As a special example, we get Ulrich sheaves on Horikawa surfaces which are double coverings of the quadric surface.}

\medskip

Using a result \cite{Laz} on triple coverings we deduce the next corollary.

\noindent\textbf{Corollary}[\thref{triple}]
\textit{If $X$ is a smooth variety of dimension $\ge4$ over $\C$ then any morphism $f:X\ra\p^n$ of degree $3$ gives an Ulrich bundle on $X$ for $f^*\osn(1)$.}

We have no hope to generalise this result to higher degrees, even assuming high enough dimension, without some further insight because Lazarsfeld's result does not generalise.
At this point, we don't even have a clear strategy to treat the low degree triple covers which are not covered by this result.

\medskip

Another way to characterise divisorial coverings is to ask $f_*\os_X=\oplus_{i=0}^{d-1}\osn(mi)$.
Therefore, it seems natural to seek for coverings in which this sheaf has some special structure.
We observe that, for a finite surjective morphism $f:X\ra\p^n$, by Horrock's splitting criterion $f_*\os_X$ is completely split if and only if $(X,\os_X(1))$ is aCM, that is $\os_X$ has no intermediate cohomology.
This seems quite a natural class of varieties on which pursue existence of Ulrich sheaves, especially in view of the fact that they admit aCM sheaves as showed in \cite{FaenziPons_acm}.

\smallskip

From this point on, our work is focused on double coverings, with the aim of studying rank $2$ Ulrich sheaves.

If we can construct such a rank $2$ Ulrich bundle then we have a morphism from $X$ to some Grassmannian but we also have a new variety: $\p(\e)$ endowed with the projection $\pi:\p(\e)\ra X$.
We studied its other contraction.

\noindent\textbf{Theorem}[\thref{coni}, \thref{nonampio}]
\textit{For $n\ge 3$, let $f:X\ra\p^n$ be a smooth double covering over $\C$ and $\e$ a rank $2$ Ulrich bundle on it.
Then both Nef and Pseudo-effective cone of $\p(\e)$ are generated by the pullback of $H$ and the class of $\os_{\p(\e)}(1)$.
Moreover, the latter line bundle induces a surjective contraction to $\p^3$ whose fibers are projected isomorphically by $\pi$ onto zeroes of sections of $\iota^*\e$, where $\iota$ is the involution associated to $f$.}

A fundamental result in our work is the following characterisation that links Ulrich bundles with both geometric and algebraic properties of our varieties.

\noindent\textbf{Theorem}[\thref{section-rk2}, \thref{rk2-section}, \thref{branchrango2}, \thref{fisso'}]
\textit{Suppose $\gk=\overline{\gk}$, $char(\gk)\neq 2$ and $n\ge 3$. 
Let $f:X\ra\p^n$ be a smooth double covering with branch locus $B\subset \p^n$ of degree $2m$ and equation $b=0$.
If $Pic(X)\cong \Z H$ then the following are equivalent:
\begin{itemize}
    \item $b=p_0^2+p_1p_2+p_3p_4$ for $p_i\in H^0(\p^n,\osn(m))$
    \item there is $Y\subset X$ mapped by $f$ isomorphically on a complete intersection of two hypersurfaces of degree $m$
    \item there is a rank $2$ Ulrich sheaf $\e$ on $X$.
\end{itemize}
Moreover, under the above assumptions, we get an exact sequence
    \sesl{\os_X}{\e}{\id_Y(m)}{seq1'}
Furthermore, if $\iota$ is the involution of the covering then $\iota^*\e=\e$ if and only if we can choose $Y$ entirely contained in the ramification divisor of $f$.    
}

We do not know a similar link between higher rank Ulrich bundles and properties of coverings even in the degree two case.

\medskip
The first result can be used as an existence criterion for rank $2$ Ulrich bundles on those double coverings.
In particular, we can find $X$ of arbitrary large dimension which support them, even though for $n\ge 5$ those varieties will always be singular.
On the other hand, the second result is an application of Hartshorne-Serre correspondence, see \thref{scorrispondenza}, which constructs rank $2$ sheaves from the zero loci of their sections.
The strong assumptions we have on Ulrich sheaves impose strong conditions on the corresponding $Y$.

We are mostly interested in the case there exist smooth double coverings admitting rank $2$ Ulrich bundles, even better, to understand when the general one does.
Some dimensional estimates on polynomials easily show that it can happen only for $n=2$, in which case existence has already been shown in \cite{SebastianTripathi}, and for $n=3$ but only with branch locus of degree $2m=4,6,8$.
Note that, in the first two cases $X$ is Fano, in particular existence of rank $2$ Ulrich bundles is highly expected since it holds in the very ample case as shown in \cite{ArrondoCosta}, \cite{BF}, \cite{Bea}, \cite{CFK3} and \cite{CFK3.1},
while in the third is Calabi-Yau, where we have only examples which are in some way complete intersections.
In this last case we also show that those bundles are \textit{spherical}.

\noindent\textbf{Theorem}[\thref{ulrichr2doublecoverexistence}, \thref{double3folds}]
\textit{Suppose $\gk=\C$.
The general double covering $f:X\ra\p^3$ branched along a divisor of degree $2m=4,6,8$ admits stable rank $2$ Ulrich bundles.
Furthermore, there are reduced components in the moduli spaces of those sheaves of dimension, respectively, $5,6,0$.
}

\medskip

The existence part is based on understanding some multiplication maps between spaces of polynomials and showing that the general $b\in H^0(\p^3,\on{3}(2m))$ can be written as $p_0^2+p_1p_2+p_3p_4$.
Stability is clear since those sheaves are always semistable and their destabilising subsheaves have to be Ulrich, but there are no Ulrich line bundles on such varieties.

Moreover, considering Hartshorne-Serre's correspondence in families, see \thref{modularcorrispondenza} and \thref{mappamodulare}, we get a morphism, actually a $\p^3$-bundle, from an open subset of the Hilbert scheme parametrising the $Y$-s as in \eqref{seq1'} to the moduli space containing stable Ulrich bundles.
Therefore, we reduce the study of their moduli space to that of the Hilbert scheme, and standard deformation theory tells us to investigate the normal bundles $\n_{Y/X}$.
Notably, also the smoothness of this Hilbert scheme, and hence of the moduli space of sheaves, reduces to a purely algebraic problem: whether the degree $2m$ part of the ideal generated by the $p_i$-s contains all the degree $2m$ polynomials or not.
Finally, fixing a base point, we have an Abel--Jacobi morphism from those Hilbert schemes to the intermediate Jacobian of $X$, which factors through the Hartshorne--Serre morphism and hence gives a rational morphism from the moduli spaces of Ulrich bundles to the intermediate Jacobian of $X$.
For a general quartic double solid ($m=2$) \cite{Voisinqds}[Lem. 4.3 and Prop. 4.5] shows that the aforementioned rational map is generically finite and the closure of its image is a $5$-dimensional component of the singular locus of the theta divisor, in particular giving another proof of the non-rationality of such varieties.

\smallskip

In the final part, we study the moduli spaces of such bundles and construct higher rank ones.
The following results are in accordance with what happens for other Fano $3$-folds of index $1$ and $2$, see \cite{CFK3.1},\cite{CFK3}, except for the fact that on some special $X$ non-smooth points representing stable bundles could appear, see \thref{qdsrk2}.
Such singularities come from sheaves fixed by $\iota$, a pathology already spotted, in the quite more general setting of Bridgeland moduli spaces of Enriques categories, in \cite{PPZ}.

\noindent\textbf{Theorem}[\thref{rango>sextic}, \thref{qdswild}]
\textit{Suppose $\gk=\C$.
\begin{itemize}
    \item The general double covering $f:X\ra\p^3$ branched along a divisor of degree $6$ admits stable Ulrich bundles of any even rank $2\rho$, in particular there are generically smooth components of their moduli spaces having dimension $5\rho^2+1$.
    \item The general double covering $f:X\ra\p^3$ branched along a divisor of degree $4$ admits stable Ulrich bundles of any rank $r\ge 2$, in particular there are generically smooth components of their moduli spaces having dimension $r^2+1$.
\end{itemize}
}

This theorem is proved by a strategy, introduced by \cite{CH} and \cite{CFK3}, which we reviewed and tried to generalise as much as possible in \thref{wildext}.
In the first case, we make extensions of two Ulrich bundles and then show that they deform to stable one; note that we know there cannot be odd rank Ulrich bundles.
In the second case, we extend Ulrich bundles with $\id_l(1)$ where $l\subset X$ is a line and then deform.
Our proof that the general such deformation is Ulrich is more involved than the one given in \cite{CFK3}, since we believe that there is a subtle gap in their presentation; see \thref{R} for more details.
However, we believe that the same reasoning applies to other Fano $3$-folds of index $2$ as well.

We think the same should apply for $2m=8$, replacing $\id_Y(1)$ with $\id_C(2)$ where $C$ has genus $1$ and degree $4$, but it requires some more work.

\smallskip

In conclusion, we studied how these moduli spaces behave under the restriction to some surface in $X$.
Tyurin's theorem, see \cite{Beafk3}, tells us that restriction to anti-canonical divisors is an étale map from our moduli spaces onto Lagrangian subvarieties of the moduli spaces of sheaves on these surfaces, which are Hyperk\"{a}hler varieties.
In the case $m=2$, see \thref{tyurin}, we can show that this map is injective.

We have not been able to tackle the irreducibility of those moduli spaces.
Equivalently, since the locus in the Hilbert scheme containing subvarieties $Y$ fitting in \eqref{seq1'} is a $\p^3$ bundle over those moduli spaces we could ask for the irreducibility of the latter. 
This is quite a delicate matter since there are Fano $3$-folds with non-cyclic Picard group, albeit not being products, with reducible moduli spaces of Ulrich bundles, see \cite{CasnatiFilipMalaspina}.
On the contrary, it seems that the only case in which irreducibility is known (for those Fano varieties) is the cubic $3$-fold, see \cite{Druel} for the rank $2$ case and \cite{FeyzbakhshPertusi} for the general case.

A remarkable result we have in the case $m=2$, i.e. when $X$ has index $2$, is that the general Ulrich bundle on a smooth hyperplane section of $X$ actually extends to the whole $X$.
Up to our knowledge, the only other example in which such a statement has been verified is for the cubic $3$-fold, see \cite{CH}.
Our approach is similar, but we do not rely on computer algebra computations.

\noindent\textbf{Theorem}[\thref{restrizdomi}, \thref{qdpwild}]
\textit{Let $X$ be a general complex quartic double solid and $\Sigma$ be a smooth divisor in $|\os_X(1)|$.
For any $r\ge 2$, the restriction map $\rho_\Sigma:\scr{U}_r\ra \cm^{ss}_r$, which sends Ulrich bundles of rank $r$ on $X$ to Ulrich bundles of rank $r$ on $S$ is generically étale and the target is irreducible.}
\medskip

In the rank $2$ case, the condition that $\rho_\Sigma$ is étale in the point $[\e]$ is equivalent to $\e|_Y\cong\n_{Y/X}$ being Ulrich, where $Y$ is the zero locus of a section of $\e$, see \thref{restretale}.
Note that there are no ramification points in the Ulrich locus by \thref{finitefixedsheaves}.
Finally, the above result can be seen as an interpolation statement; see \thref{interpolazione}.

\noindent\textbf{Corollary}
Let $X$ be a general quartic double solid.
For any $r\ge 2$, given $r^2$ general points on some smooth $\Sigma\in |\os_X(1)|$ there is a smooth degree $r^2$ and genus $\frac{2}{3}r(r^2-3r+2)+(r-1)^2$ projectively normal curve passing through them.

\section*{Generalisations and relations to other notions}

The definition of Ulrich sheaf on a polarised variety $(X,H)$ we gave is a straightforward one to state but encloses many geometric properties, see for example the other equivalent definitions in \thref{definizioneUlrich}.  

In this work, we will investigate thoroughly the relation between the existence of Ulrich sheaves on divisors in some special weighted projective spaces and linear determinantal representations of their equations. Moreover, we will briefly recall how, in the case of usual projective space, this correspondence can be upgraded to arbitrary subscheme if we substitute the equations with the Chow form in \autoref{section1.5}.
On the contrary, since we will not deal with Boij-S\"{o}derberg theory in the thesis we want to recap some material here, see also the survey \cite{Floy}.
To a resolution of a graded module $M$ over a polynomial ring we can associate its \textit{graded Betti numbers}, in the form of a matrix (\textit{Betti tables}) of integers. 
In \cite{Boij-Sod}, the authors consider the rational cone spanned by those vectors that come from Cohen-Macaulay modules $M$ and they conjectured that the extremal rays of this cone should come from special resolutions, called \textit{pure}.
This conjecture has been proved in \cite{Eis_Sch-bettinumbersgrademodules} where the authors also exhibit a duality between Betti tables and cohomology tables of vector bundle on $\p^n$, which encode all the integers $h^j(\p^n,\e(i))$ for $i,j\in\Z$.
Also the convex cone spanned by those is studied and their extremal rays are determined, they come from \textit{supernatural vector bundles}.
Anytime we have a finite surjective morphism $f:X\ra\p^n$, we can consider the subcone of the cone of cohomology tables on $\p^n$ obtained from the sheaves of the form $f_*\e$.
In \cite{EisSch_BStheory} it is remarked that those two cones coincide if and only if there is some $\e$ on $X$ such that $f_*\e\cong\osn^\rho$, that is if and only if there is an $f$-Ulrich sheaf.
One direction is clear while, if there is an $f$-Ulrich sheaf $\e$ then a sheaf $E$ on $\p^n$ sits in the same ray as $E^\rho\cong f_*(\e\otimes f^*E)$ by projection formula.

We can easily generalise the given definition of Ulrich sheaves: on one side we can ask for fewer vanishing conditions and on the other we could relax the assumptions on which type of objects are involved.

We have seen that the cohomology table of an Ulrich sheaf has the same zeroes as the one of $\osn$.
If we drop the requirements on the top cohomology, that is $h^n(X,\e(-n))=0$, the we get the definition of (initialised) aCM sheaf, that is $h^j(X,\e(i))=0$ for $1\le j\le n-1$ and $i\in \Z$, or for $(i,j)=(-1,0)$.
Those have become popular after Horrock's splitting criterion, see \cite{Horrocks}, which says that on $\p^n$ the only aCM vector bundles withouth intermediate cohomology are direct sums of line bundles.
There has been much work seeking conditions in order to generalise this theorem on other varieties, constructing examples of aCM sheaves and studying their moduli spaces, for example, see the survey \cite{Ottaviani_split} or \cite{madonna2001}.
In the case of hypersurfaces in $\p^n$ we still get a link with determinantal presentations in \cite{Beadeterminantal}[Thm. A].
Moreover, anytime a reduced scheme $(X,\os_X(1))$ is aCM then, except for some classified pairs, we can find infinitely many aCM sheaves on it, see \cite{FaenziPons_acm}.
A similar problem for Ulrich sheaves is still open.

A far more general class of vector bundles, which specialise to Ulrich ones but was considered quite earlier, is the one of instanton bundles.
They first were defined on $\p^3$ as vector bundles coming from solutions of the Yang-Mills equation on the $4$-dimensional sphere, see \cite{AHDM}.
Many purely algebro-geometric extensions of this definition have been given for vector bundles on other varieties  besides projective spaces, see \cite{Faenzi11} and \cite{Kuznetsov_instanton} for Fano $3$-folds and \cite{AC} for arbitrary varieties, or even for objects in the bounded derived category of coherent sheaves on some variety, see \cite{Comaschi-Jardim-Martinex-Mu}.
A detailed study of instanton bundles in the Fano $3$-fold case is in \cite{ComaschiFaenzi}.
We remark that in \cite{CMR_Veronese3} it is shown that Ulrich bundles on Veronese three-folds actually coincide with some special instanton bundles on $\p^3$.

Finally, we could also ask for objects in the bounded derived category of $X$ satisfying the Ulrich vanishing, with an \textit{ample sequence} taking the role of the ample divisor.
Some investigations in this direction have been carried out in \cite{Yoshida} and \cite{Deaconu}.

\section*{Content}

This work is divided in a preliminary chapter and four main chapters.

In chapter $0$, we fix the notation used in the rest of the thesis and recall some preliminaries on scheme theory, sheaf theory and moduli spaces.

Chapter $1$ deals with the general theory of Ulrich sheaves on projective varieties, discussing properties and existence results.

Chapter $2$ presents some technical results on Hartshorne-Serre correspondence and computations with extensions and deformations of families of sheaves.

In chapter $3$ we expose the link between Ulrich modules and matrix factorisations, with particular focus on double coverings of projective space.

In chapter $4$ we investigate the moduli spaces of rank $2$ Ulrich bundles on some double coverings of $\p^3$ and their relation to moduli spaces of vector bundles on special contained in those $3$-folds.

\chapter*{Notation}
\addcontentsline{toc}{chapter}{Notation}
\begin{longtable}{cp{0.75\textwidth}}
        $\gk$ & base field \\
        $\p^n$ & projective space of dimension $n$ over $\gk$\\
        $(X,H)$ & polarised scheme \\
        $\os_X$ & structure sheaf of $X$ \\
        $\os_{X,x}$ & stalk of $\os_X$ in the point $x\in X$ \\
        $m_x$ & maximal ideal of $\os_{X,x}$\\
        $k(x):=\os_{X,x}/m_x$ & skyscraper sheaf supported on $x$ with value its residue field\\
        $\dim(X)$ & dimension of the scheme $X$ over $\gk$\\
        $\omega_X$ & dualizing sheaf of $X$ \\
        $K_X$ & canonical divisor of $X$ \\
        $\mathcal{T}_X$ & tangent bundle of $X$ \\
        $Pic(X)$ & Picard group of $X$ \\
        $|H|$ & complete linear system associated to the Cartier divisor $H$ \\
        $\id_{Y/X}$ (or simply $\id_Y$) & ideal sheaf of the subscheme $Y$ in $X$ \\
        $\n_{Y/X}$ (or simply $\n_Y$) & normal sheaf of the subscheme $Y$ in $X$ \\
        $\e\du:=\Hom_{\os_X}(\e,\os_X)$ & dual sheaf of $\e$ \\
        $\bigwedge^l\e$ & $l$-th wedge bundle of $\e$\\
        $det(\e):=\bigwedge^{rk(\e)}\e$ & determinant bundle\\
        $\ext^i(\e,\f):=(R^i\Hom(\e,-))(\f)$ & ext-sheaf \\
        $\ext^i_\tau(\e,\f):=((R^i\tau_*)\Hom(\e,-))(\f)$ & relative ext-sheaf respect to the morphism $\tau$ \\
        $\scr{T}or_i(\e,\f):=(L_i(\e\otimes-))(\f)$ & Tor-sheaf \\
        $Spec_{\os_X}(\e)$ (or simply $Spec(\e)$) & spectrum of $\e$ \\
        $\p(\e):=Proj_{\os_X}(Sym(\e))$ & projectivization of $\e$ \\
        $h^j(X,\e):=dim_\gk H^j(X,\e)$ & \\
        $ext^j(\e,\f):=dim_\gk Ext^j(\e,\f)$ & \\
        $\chi(\e):=\sum_{j\in \N}(-1)^jh^j(X,\e)$ & Euler characteristic of $\e$ \\
        $\chi(\e,\f):=\sum_{j\in \N}(-1)^jext^j(\e,\f)$ & Euler product of $\e$ and $\f$\\
        $P(\e)(t):=\chi(\e(tH))$ & Hilbert polynomial of $\e$ respect to $H$ \\
        $p(\e)(t):=\dfrac{P(\e)(t)}{rk(\e)}$ & reduced Hilbert polynomial of $\e$\\ 
        $A^i(X)$ & $i$-th Chow group of $X$ \\
        $c_i(\e)$ & $i$-th Chern class of $\e$, as an element in $A^i(X)$ \\
        $c_i(X):=c_i(\mathcal{T}_X)$ & $i$-th Chern class of $X$ \\
        $td(X):=td(\mathcal{T}_X)$ & Todd class of a non-singular variety $X$ \\
        $pd_X(\e_x)$ & projective dimension of $\e_x$ as an $\os_{X,x}$-module \\
        $depth_X(\e_x)$ & depth of $\e_x$ as an $\os_{X,x}$-module \\
        $\cong$ & isomorphism \\
        $\sim$ & linear equivalence for divisors \\
        $\gi_l$ & identity matrix of rank $l$\\
        $\p(1^{n+1},m)$ & weighted projective space with weights $(1,\dots, 1,m)$\\
        $P_m:=\p(\osn\oplus\osn(m))$ & projective closure of the total space of $\osn(m)$\\
\end{longtable}

\chapter{Conventions and preliminaries}

We will fix most of the notation here, while presenting the preliminaries.
For scheme theory and cohomology, we refer mainly to \cite{Har} or to \cite{GW1} and \cite{GW2}.
For stability and moduli spaces of sheaves we follow \cite{HuyLeh}.
Concerning Chow groups and intersection theory, we refer to \cite{Ful}.

\startcontents[chapters]
\printcontents[chapters]{}{1}{}

\section{Schemes and sheaves}

Fix an arbitrary field $\gk$.
We will write simply $\p^n$ instead of $\p^n_\gk$.
All of our schemes will be Noetherian over $\gk$.
A variety is an integral, separated, finite-type scheme over $Spec(\mathbf{k})$.
All sheaves on a scheme will be quasi-coherent, if not otherwise stated.
"Locally free sheaf" and "vector bundle" will be used interchangeably.

Given a morphism $g:X\ra Y$ between two schemes we can define two functorial operations on sheaves: pullback and pushforward.
The functor $g^*$ sends sheaves on $Y$ to sheaves on $X$, preserves quasi-coherence and is right exact; if $g$ is flat it is also left exact.

The functor $g_*$ sends sheaves on $X$ to sheaves on $Y$, preserves quasi-coherence due to our noetherian assumption, see \cite{GW1}[Cor. 10.27 and Cor. 10.24], and is left exact; if $g$ is affine it is also right exact.
As a special case, if $g$ is affine it \textbf{preserves cohomology}, meaning that for any quasi-coherent sheaf $\f$ on $X$ we have $H^j(X,\f)\cong H^j(Y,g_*\f)$ for all $j\in \Z$, see \cite{GW2}[Cor. 22.6 1)].
In particular, if $\f$ is coherent and $X,Y$ proper then the dimensions of those $\gk$ vector spaces are equal $h^j(X,\f)= h^j(Y,g_*\f)$.

Pullback and pushforward are \textbf{adjoint functors} meaning that, for every morphism $g:X\ra Y$ and any two quasi-coherent sheaves $\f,\g$ on $X,Y$ respectively there is a functorial isomorphism
\[Hom_{\os_X}(g^*\g,\f)\cong Hom_{\os_Y}(\g,g_*\f).\]
This induces natural morphisms $g^*g_*\f\ra \f$ and $\g\ra g_*g^*\g$.

\Le\label{counitsurj}
If $g:X\ra Y$ is an affine morphism then for any $\f$ quasi-coherent sheaf on $X$ the morphism $g^*g_*\f\ra \f$ is surjective.
\ma
\begin{proof}
This is a local property so, being $g$ affine, we can reduce to the case $g:Spec(A)\ra Spec(B)$.
There is an equivalence between the categories of quasi-coherent sheaves on those schemes and, respectively, the ones of modules over $A,B$, see \cite{GW1}[Cor. 7.17].
Hence, translating everything in the language of modules over rings with \cite{GW2}[Prop. 7.24], we reduce to the following easy claim: let $\varphi:B\ra A$ be a ring map and $M$ be an $A$-module then the map $M\otimes_B A\ra M$ given by $m\otimes a\ra a.m$ is surjective.
\end{proof}

Another ubiquitous identity relating $g^*$ and $g_*$ is the \textbf{projection formula}.
With the above notation, by \cite{GW2}[Prop. 22.81] we always have a morphism of sheaves 
\[g_*(\f)\otimes\g\ra g_*(\f\otimes g^*\g),\]
which is an isomorphism, for example, if $\g$ is locally free of finite rank or if $g$ is affine. 

For any closed subscheme $i:Y\hookrightarrow X$ we have an exact sequence of sheaves
\ses{\id_Y}{\os_X}{\os_Y}
where, by abuse of notation, we will write $\os_Y$ instead of $i_*\os_Y$.

We end by recalling some formulas.
We know all the cohomology of line bundles on projective space:
\begin{equation}\label{bott}
    h^j(\p^n,\osn(i))= \begin{cases}
        \binom{n+i}{n} \hfill i\geq 0,\; j=0 \\
        \binom{-i-n-1}{n} \qquad \hfill  i\leq -n-1,\; j=n\\
        0 \hfill \text{otherwise} \\        
    \end{cases}
\end{equation}

On a projective scheme $X\subset \p^N$ of pure dimension $n$ over a field $\gk$ we define $\omega_X:=\ext^{N-n}(\os_X,\osn(-n-1))$ to be the \textbf{dualizing sheaf}, in particular $\omega_{\p^n}\cong \osn(-n-1)$.
By \cite{AltKle_duality}[I Prop. 2.3], if $Y\subset X\subset \p^N$ is a pure $s$-dimensional closed subscheme of $X$, then we have $\omega_Y\cong \ext^{n-s}(\os_Y,\omega_X)$.
If $X$ is Cohen-Macaulay then, for every coherent sheaf $\e$ Serre duality holds:
\[H^j(X,\e)\du\cong Ext^{n-j}(\e,\omega_X) \qquad \text{for all}\quad j\in\Z,\]
see \cite{AltKle_duality}[I (1.3)].
In addition, if $X$ is a locally complete intersection, or more generally Gorenstein, and reduced then $\omega_X$ is a line bundle \cite{AltKle_duality}[I Cor. 2.6, Prop. 2.8].

\section{Polarised schemes}

Let $\cl$ be a line bundle on a scheme $X$.
$\cl$ is \textbf{globally generated} if there exists some $\rho\geq 0$ and a surjective morphism $\os_X^\rho\ra\cl$.
$\cl$ is \textbf{ample} if there is some $l>0$ and an embedding $i:X\ra \p^N$ such that $i^*\osn(1)\cong \cl^{\otimes l}$.
A Cartier divisor $H$ is said ample and globally generated if the corresponding line bundle $\os_X(H)$ is such.

In the following, unless otherwise specified, by a 
\textbf{polarised scheme} $(X,H)$ (or $(X,\cl)$) we mean a proper scheme $X$ over $\mathbf{k}$ with an ample and globally generated divisor $H$ (line bundle $\cl$); in particular $X$ is projective over $\gk$.
We will reserve the letter $d$ for the \textbf{degree} of $H$, i.e. $d:=H^n$.
For technical reasons, we will assume $X$ to be equidimensional, and for simplicity we also suppose $X$ to be connected.
Given a coherent sheaf $\e$ on a polarised scheme we write $\e(i)$ instead of $\e\otimes\os_X(iH)$.

Given a polarised scheme $(X,H)$, the complete linear system $|H|$ defines a finite morphism $\phi:X\ra\p^N$ since $H$ is globally generated.
Moreover, being $X$ proper and $H$ ample, $\phi$ is a finite morphism by \cite{GW1}[Thm. 13.84 2)].
If $\gk$ is infinite we also have Noether's normalization: we can construct a finite surjective morphism $f:X\ra\p^n$ such that $f^*\osn(1)\cong\os_X(H)$, see \cite{GW1}[Thm. 13.89].

\section{Depth and projective dimension}

Here we follow \cite{BrunsHerzog}.
Consider a noetherian, commutative ring $R$ and a non-zero finitely generated $R$-module $M$.
The elements $a_1,\dots ,a_l\in R$ form a \textbf{regular sequence} if for each $1\leq i\leq l$ we have that $a_i$ is not a zero-divisor in $M/(a_1,\dots, a_{i-1})M$ and $(a_1,\dots, a_i)M\neq M$.
Given an ideal $I\subset R$, the $\mathbf{I}$-\textbf{grade} of $M$, denoted by $grade_I(M)$, is the maximal length of a regular sequence if $M\neq IM$, and $+\infty$ otherwise.
In particular, if $M\neq IM$ we have $grade(M)\leq dim(Supp(M))\leq dim(R)$.
We will usually deal with local rings $R$, then the grade computed with respect to the unique maximal ideal is called \textbf{depth} and written $depth_R(M)$, or simply $depth(M)$.

From now on we suppose for simplicity that $R$ is local.
The module $M$ is said \textbf{Cohen-Macaulay} if $depth(M)=dim(M)$, where the dimension of $M$ is, as usual, the Krull dimension of its support $R/Ann(M)$.
Note that, if $M\neq 0$ we have $depth(M)\leq dim(M)\leq dim(R)$.
In the same setting, we call \textbf{projective dimension} of $M$, denoted $pd_R(M)$, the minimum length of a resolution of $M$ by finite free $R$-modules.
This minimum can be $+\infty$ in general, but if $R$ is regular it is always a finite number.
\textbf{Auslander-Buchsbaum formula} links depth and projective dimension of a module.
If $M$ is a finitely generated module of finite projective dimension over a local ring $R$, by \cite{BrunsHerzog}[Thm. 1.3.3] we have 
\[depth_R(M)+pd_R(M)=depth_R(R).\]
Here it comes a standard consequence of this formula.

\Le\label{codimensione}
Let $R$ be a Cohen-Macaulay local ring and $M$ a finite $R$-module of finite projective dimension.
Then $pd(M)\ge codim(Supp(M))$ and we have equality if and only if $M$ is Cohen-Macaulay.
\ma
\begin{proof}
    By the above formula we have 
    \[pd(M)=depth(R)-depth(M)=dim(R)-depth(M)\geq dim(R)-dim(M)\]
    hence we derive the stated inequality.
    This is an equality if and only if $dim(M)=depth(M)$ that is $M$ is Cohen-Macaulay.
\end{proof}

If $\e$ is a sheaf on $X$ we will usually write $depth_X(\e_x)$ instead of $depth_{\os_{X,x}}(\e_x)$.
For future use, we rephrase \cite{Stacks}[\href{https://stacks.math.columbia.edu/tag/0AUK}{Tag 0AUK}].

\Le\label{depthfinite}
Let $g:X\ra Y$ be a finite morphism of Noetherian schemes and $\e$ a coherent sheaf on $X$.
Fix some closed point $y\in g(X)$ and denote by $x_i$ the points in the schematic preimage $g\inv(y)$.
We have 
\[depth_Y(g_*\e)_y=\min_i\{depth_X(\e_{x_i})\}\]
\ma
\begin{proof}
Since $g$ is finite, taking the pullback diagram 
\dia
Spec(S) \ar[r,"g'"] \ar[d] & Spec(\os_{Y,y}) \ar[d] \\
X \ar[r, "g"] & Y \\
\mma
we get a finite ring morphism $g':\os_{Y,y}\ra S$, being finiteness stable under base change \cite{GW1}[Prop. 12.11 2)].
Call $\mathfrak{m}_i$ the maximal ideals of $S$.
We can consider $\e_S:=\e\otimes_{\os_X} S$ both as an $S$-module and as a $\os_{Y,y}$-module, through the map $g'$.
Then, by \cite{Stacks}[\href{https://stacks.math.columbia.edu/tag/0AUK}{Tag 0AUK}] we have $depth_{\os_{Y,y}}(\e_S)=\min_i\{grade_{\mathfrak{m}_i}(\e_S)\}.$
The stalk $(g_*\e)_y$ is the pullback of $g_*\e$ on $Spec(\os_{Y,y})$ therefore, by base change \cite{GW2}[Lem. 22.88], it is isomorphic to $\e_S$ seen as a $\os_{Y,y}$-module; hence $depth_{\os_{Y,y}}(\e_S)=depth_Y(g_*\e)_y$.
Finally, being $x_i$ closed points and hence $\mathfrak{m}_i$ maximal ideals, by \cite{BrunsHerzog}[Prop. 1.2.10 a)] we have $grade_{\mathfrak{m}_i}(\e_S)=depth_{S_{\mathfrak{m}_i}}((\e_S)_{\mathfrak{m}_i})=depth_{\os_{X,x_i}}(\e_{x_i})$.
Putting all together, we conclude.
\end{proof}

Next, we introduce Serre's conditions, following \cite{EGA4.2}[Definition 5.7.2] and \cite{BrunsHerzog}\footnote{It seems that there is no agreement on this definition in the literature, see \url{https://mathoverflow.net/questions/22228/what-is-serres-condition-s-n-for-sheaves}}.
\De
Let $\e$ be a coherent sheaf on $X$.
Given $k\in\N$, we say that $\e$ satisfies \textbf{condition} $\mathbf{S_k}$ if $depth(\e_x)\geq min\{k,dim(\e_x)\}$ for all $x\in X$.
We call $\e$ Cohen-Macaulay if $\e_x$ is Cohen-Macaulay for all $x\in X$.
In particular, a sheaf is Cohen-Macaulay if and only if it satisfies $S_k$ for all $k$.
\Ne

Note that the property of being Cohen-Macaulay can be checked just on closed points $x\in X$ thanks to \cite{Stacks}[\href{https://stacks.math.columbia.edu/tag/0AAG}{Tag 0AAG}].
The main point of those are the following characterisations, see \cite[\href{https://stacks.math.columbia.edu/tag/0AXY}{Tag 0AXY}]{Stacks} and  \cite[\href{https://stacks.math.columbia.edu/tag/0AY6}{Tag 0AY6}]{Stacks}
\Le\label{s12}
Let $\e$ be a coherent sheaf on a variety $X$.
\begin{itemize}
    \item $\e$ is torsion-free if and only if it satisfies $S_1$ and its support is $X$.
    \item If $X$ is normal then, $\e$ is reflexive if and only if it satisfies $S_2$ and is torsion-free.
\end{itemize}
\ma

Finally, we recall this well-known fact for future reference.
\Le\label{extpd}
If $X$ is smooth, then for a coherent sheaf $\e$ the following are equivalent:
\begin{enumerate}
    \item $\e$ has a locally free resolution of length $\leq l$
     \item $\ext^j(\e,\cl)=0$ for all $j>l$ and for any locally free sheaf $\cl$
     \item $\ext^j(\e,\cl)=0$ for all $j>l$ and for some locally free sheaf $\cl$
    \item $pd(\e_x)\leq l$ for all $x\in X$.
\end{enumerate} 
\ma
\begin{proof}
    $\mathbf{1.\Rightarrow 2.}$
    The sheaf $\ext^j(\e,\cl)$ can be computed as the $j$-th homology of the complex obtained by applying $\Hom(-,\cl)$ to a locally free resolution of $\e$, see \cite{Har}[III Prop. 6.5], so we get the claimed vanishing.
    
    $\mathbf{2.\Rightarrow 3.}$
    This is trivial.
    
    $\mathbf{3.\Rightarrow 4.}$
    By \cite{Har}[III Prop. 6.8] we have $\ext^j(\e,\cl)_x\cong Ext^j(\e_x,\cl_x)$.
    Then our assumption implies that for every $x\in X$ we get $Ext^j(\e_x,\os_{X,x})=0$ for all $i>l$ so that by \cite{Har}[III Exercise 6.6 b)] we conclude.
    
    $\mathbf{4.\Rightarrow 1.}$
    This follows from \cite{Har}[III Exercise 6.5 c)].
\end{proof}

\section{Stability and flatness}\label{defstability}

In this section, we will consider the stability of sheaves following \cite{HuyLeh}, and we will focus on the simple setting that we will need in the rest of this work.

Given a polarised scheme $(X,H)$ we define the \textbf{Hilbert polynomial} of a sheaf $\e$ to be $P(\e)(t):=\chi(\e(tH))$.
Consider now a polarised $n$-dimensional variety $(X,H)$, hence an integral scheme, and a torsion-free sheaf $\e$.
The rank of $\e$ will be its \textbf{generic rank}, that is the dimension of $\e\otimes k(\eta)$ as a vector space over $k(\eta)$, where $\eta$ is the generic point of $X$.
Note that $rk(\e)>0$ since $\e$ is torsion-free.

Being $Supp(\e)=X$, this polynomial must be of degree exactly equal to $n=dim(X)$.
Then we call \textbf{reduced Hilbert polynomial} $p(\e)(t):=\dfrac{P(\e)(t)}{rk(\e)}$.
We define $\e$ to be \textbf{Gieseker (semi-)stable} if it is torsion-free and for all subsheaves $0\neq\f\subsetneq \e$ we have $p(\f)<p(\e)\; (\text{respectively}\; p(\f)\leq p(\e))$, where $p(\f)<(\le)p(\e)$ if and only if $p(\f)(t)<(\le)p(\e)(t)$ for $t>>0$.
Usually we just say that $\e$ is (semi-)stable.
Note that Gieseker (semi-)stability can be checked by looking at torsion-free quotients instead of subsheaves, see \cite{HuyLeh}[Prop. 1.2.6].

Given a sheaf $\e$ we can write its Hilbert polynomial as 
\[P(\e)(t)=\sum_{i=0}^n \alpha_i(\e)\dfrac{t^i}{i!} \quad \text{and set} \qquad deg(\e):=\alpha_{n-1}(\e)-rk(\e)\cdot\alpha_{n-1}(\os_X).\]
The definition of $deg(\e)$ depends on $H$ but on a smooth variety, using Grothendieck--Riemann--Roch formula, we can show that $deg(\e)=c_1(\e)\cdot H^{n-1}$.

For a torsion-free coherent sheaf $\e$, we define the \textbf{slope} $\mu(\e):=\dfrac{deg(\e)}{rk(\e)}$.
We say that $\e$ is \textbf{slope (semi-)stable} if  for all subsheaves $\f\subset \e$ with $0<rk(\f)<rk(\e)$ we have $\mu(\f)<\mu(\e)\; (\text{respectively}\; \mu(\f)\leq \mu(\e))$.
Also slope (semi-)stability can be checked on torsion-free quotients, see \cite{OSS}[Chap. 2, Thm. 1.2.2].

We have, see \cite{HuyLeh}[Lem. 1.2.13], the following implications between these concepts:
\begin{center}
    slope stable $\Rightarrow$ stable $\Rightarrow$ semi-stable $\Rightarrow$ slope semi-stable.
\end{center}
Note that torsion-free rank $1$ sheaves are always slope stable.

The first set of items in the following is contained in \cite{HuyLeh}[Prop. 1.2.7] so we prove just the second one.
\Prop\label{stabilitàtrick}
Let $\f,\g$ be two torsion-free, semistable sheaves on a variety $X$. 
\begin{itemize}
    \item If $p(\f)>p(\g)$ then $Hom(\f,\g)=0$.
    \item If $p(\f)=p(\g)$ then any non-zero morphism $\f\ra \g$ is injective if $\f$ stable and surjective if $\g$ stable.
    \item If $p(\f)=p(\g)$, $rk(\f)=rk(\g)$ and at least one of them is stable then any non-zero map $\f\ra \g$ is an isomorphism.
\end{itemize}
Let $\f,\g$ be two torsion-free, slope-semistable sheaves on a variety $X$.
\begin{enumerate}
    \item If $\mu(\f)>\mu(\g)$ then $Hom(\f,\g)=0$. 
    \item If $\mu(\f)=\mu(\g)$ and they are both slope-stable then any non-zero map $\f\ra \g$ is injective and $rk(\f)=rk(\g)$.
\end{enumerate}
\One
\begin{proof}
\gel{1}
    Any morphism $\varphi:\f\ra \g$ can be factored in $\f\twoheadrightarrow \ci \hookrightarrow\g$ with $\ci$ the image of $\varphi$, which is torsion-free.
    If this sheaf is non-zero then we have $\mu(f)\leq \mu(\ci)$ by semistability of $\f$.
    Since $\ci\subset \g$ either $rk(\ci)<rk(\g)$ and then $\mu(\ci)\leq \mu(\g)$ by semistability or $rk(\ci)=rk(\g)$ and then is clear that $\mu(\ci)\leq \mu(\g)$.

    \gel{2}
    With the above notation we get either $\f\cong \ci$, meaning that $\varphi$ is injective, or $\mu(ker(\varphi)<\mu(\f)<\mu(\ci)$ by stability of $\f$.
    But the second case cannot happen since either $rk(\ci)<rk(\g)$ and then $\mu(\ci)<\mu(\g)=\mu(\f)$ by slope-stability or $rk(\ci)=rk(\g)$ and then is clear that $\mu(\ci)\leq \mu(\g)$, hence in both cases we have a contradiction.
    It follows that $\varphi:\f\hookrightarrow\g$ but now, arguing as above, by slope-stability of $\g$ the only possibility is that $rk(\f)=rk(\g)$.
\end{proof}

Given a semi-stable sheaf $\e$ with reduced Hilbert polynomial $p$, we can always construct a \textbf{Jordan-Holder filtration}
\[0=\e_0\subset \e_1\subset \dots \subset\e_l=\e\]
with the property that $\e_i/\e_{i-1}$ are torsion-free, stable and with reduced Hilbert polynomial equal to $p$.
Even if the filtration is not unique, the sheaf $gr(\e):=\oplus_{i=0}^l\e_i/\e_{i-1}$ is well defined up to isomorphism.
Two sheaves $\e_1,\e_2$ are called \textbf{S-equivalent} if $gr(\e_1)\cong gr(\e_2)$.

Recall that a \textit{family of sheaves} $\scr{F}$ on $X$ flat over a base scheme $B$ is a sheaf on $X\times B$ flat over $B$.
For any $b\in B$ we will usually denote by $\scr{F}_b$ the pullback of $\scr{F}$ for the morphism $b\hookrightarrow B$.
Flatness implies that, if $B$ is connected, the Hilbert polynomial $P(\scr{F}_b)$ is the same for each $b\in B$.
In particular, also $\chi(\scr{F}_b)$ does not depend on $b$.
Even though the single cohomology groups can jump we always have a semicontinuity result.
Actually, it holds more generally for ext groups: given any two flat families $\scr{E},\scr{F}$ of coherent sheaves on a smooth scheme $X$, the function from $B$ to $\N$ given by 
\[b\mt ext^j(\scr{E}_b,\scr{F}_b)=dim_\gk Ext^j(\scr{E}_b,\scr{F}_b)\]
is upper semi-continuous, see \cite{BanicaPutinarSchumacher}[Satz 3 i)].
Moreover, the \textbf{Euler product}
\[\chi(\scr{E}_b,\scr{F}_b):=\sum_{i\in\N}ext^j(\scr{E}_b,\scr{F}_b)\]
is locally constant on $B$ by \cite{BanicaPutinarSchumacher}[Satz 3 iii)].

We end with a useful technical lemma that we will need in the following, the strategy of proof is the same used in \cite{HuyLeh}[Lem. 2.3.1].
\Le\label{crit_stab}
Let $(X,H)$ be a polarised, $n$-dimensional variety. 
Suppose $\scr{F}$ is a family of torsion-free (slope-)semistable sheaves on $X$ flat over an irreducible base $\cc$.
Then exactly one of the following holds:
\begin{itemize}
    \item $\scr{F}_c$ is (slope-)stable for $c\in\cc$ general 
    \item there is an irreducible scheme $\cq$ with a projective morphism $\cq\ra \cc$ and two families $\scr{E},\scr{Q}$ of sheaves on $X$ flat over $\cq$ such that for general $q\in\cq$ both $\scr{E}_q,\scr{Q}_q$ are torsion-free and, calling $\scr{F}_\cq$ the base change of $\scr{F}$ to $\cq$, we have an exact sequence of families
    \ses{\scr{E}}{\scr{F}_\cq}{\scr{Q}}
   whose restrictions to $X_q$ is (slope-)destabilising for all $q\in\cq$.
\end{itemize}
Moreover, in the Gieseker-stable case we can suppose that $\scr{E}_q,\scr{Q}_q$ are torsion-free for all $q\in\cq$.
While, in the slope-stable case, if $L_q$ is the quotient of $\scr{Q}_q$ by its torsion, then the surjection $(\scr{F}_\cq)_q\twoheadrightarrow L_q$ is stille slope-desatbilising.
\ma
\begin{proof}
If the first holds then we are done, so we suppose that it does not hold.
Since being (slope-)stable is an open property, see \cite{HuyLeh}[Lem. 2.3.1], then we deduce that for all $c\in\cc$ the sheaf $\scr{F}_c$ is not (slope-)stable.
Since $\scr{F}_c$ is already (slope-)semistable, this implies that it admits a quotient $\scr{F}_c\ra Q_c$ such that $(\mu(\scr{F}_c)=\mu(Q_c)) \; p(\scr{F}_c)=p(Q_c)$. 
Moreover, we can always assume $Q_c$ torsion-free by (\cite{OSS}[Chap. 2, Thm. 1.2.2]) \cite{HuyLeh}[Lem. 1.2.6].

Consider the set $\mathcal{P}$ of degree $n$ polynomials $P$ such that there is $ c\in\cc$ with $\scr{F}_c\twoheadrightarrow Q$, the latter is torsion-free and satisfies $P(Q)=P,\; p(Q)=p(\scr{F}_c)$
(for the slope-stability case we just replace the last condition with $\mu(Q)=\mu(\scr{F}_c)$).
This set is clearly finite since it is just made of the polynomials $l\cdot p(\scr{F}_c)$ for $1\le l< rk(\scr{F}_c)$ (in the slope case this is still a finite set by \cite{HuyLeh}[Lem. 1.7.9, Remark 1.7.10]).
For a fixed polynomial $P$, by \cite{HuyLeh}[Thm. 2.2.4] we can construct the relative Quot-scheme $\mathcal{Q}(P):=\mathcal{Q}_{X\times\cc/\cc}(\scr{F},P)$ representing the functor of quotients of $\scr{F}$ which are flat over $\cc$ and having fiberwise Hilbert polynomial $P$.
Moreover, each of them has a projective morphisms $\phi_P:\mathcal{Q}(P)\ra\cc$.
Consider the subscheme $\mathcal{Q}(P)'$ of $\mathcal{Q}(P)$ made of those components whose general sheaf is torsion-free with the restricted map $\phi_P':\mathcal{Q}(P)'\ra\cc$.
We are supposing that no sheaf in $\scr{F}$ is (slope-)stable, hence the union of the images of $\phi_P'$ for $P\in\mathcal{P}$ covers all $\cc$.
Being $\cc$ irreducible, we deduce that there is $P\in\mathcal{P}$ and some irreducible component $\mathcal{Q}(P)^\circ$ in $\mathcal{Q}(P)'$ such that $\phi_P^\circ:\mathcal{Q}(P)^\circ\ra\cc$ is surjective.

We define $\cq:=\mathcal{Q}(P)^\circ$ which has a surjective and projective morphism to $\cc$.
Moreover define $\scr{Q}$ to be the universal quotient sheaf on $\cq$, which is flat over $\cq$, so that we have a surjection $(\phi_P^\circ)^*\scr{F}\twoheadrightarrow \scr{Q}$.
Finally, define $\scr{K}$ to be the kernel of this surjection, so it follows that it is flat over $\cq$ and we have a sequence
\sesl{\scr{E}}{(id_X\times\phi_P)^*\scr{F}}{\scr{Q}}{univ}
The general sheaf in $\scr{Q}$ is torsion-free since by construction there exists at least one $q\in\cq$ such that $\scr{Q}_q$ is torsion-free and this is an open property, see \cite{HuyLeh}[Lem. 2.3.1].
Since the kernel of torsion-free sheaves is still torsion-free, the same applies to $\scr{E}$.
Finally, for all $q\in\cq$ the restriction of \eqref{univ} to $X_q$ is (slope-)destabilising since, by construction, there exists at least one such point and (slope) Hilbert polynomial are constant in flat families. 

Furthermore, in the Gieseker-semistable case, actually for all $q\in\cq$ we have $\scr{Q}_q$ torsion-free. 
Indeed, consider the sequence
\ses{T_Q}{\scr{Q}_q}{F_Q}
where $T_Q$ is the torsion subsheaf of $\scr{Q}_q$ hence $\scr{Q}_q,F_Q$ share the same rank.
If $T_Q\neq 0$ then $P(T_Q)> 0$ and it follows that 
\[p(F_Q)=\dfrac{P(F_Q)}{rk(F_Q)}=\dfrac{P(\scr{Q}_q)-P(T_Q)}{rk(F_Q)}<\dfrac{P(\scr{Q}_q)}{rk(F_Q)}=\dfrac{P(\scr{Q}_q)}{rk(\scr{Q}_q)}=p(\scr{Q}_q)\]
hence the quotient $\scr{Q}_q\ra F_Q$ contradicts semi-stability.
It follows that also $\scr{K}_q$ is torsion-free for all $q\in\cq$.
Finally, in the slope-stable case, by composition we get a surjection $\f\twoheadrightarrow L_q$ hence by slope-semistability of $\scr{F}$ we have
\[\mu(L_q)\ge \mu(\f)=\mu(\scr{Q}_{q})=\mu(L_q)+\dfrac{c_1(T)}{rk(\scr{Q}_{q})}\]
hence $c_1(T)=0$ and $\mu(L)=\mu(\f)$.
\end{proof}

\section{Hilbert scheme and moduli spaces of sheaves}\label{moduli}

The Hilbert polynomial of any subvariety $Y\subset X$ is defined to be the Hilbert polynomial of its ideal sheaf $\id_Y$.
Grothendieck, \cite{Gro_Hilb}, showed that the functor which sends any scheme $T$ to the set of families of subschemes in $X\times T$ flat over $T$ with fixed Hilbert polynomial $P$ is representable by a projective scheme: the \textbf{Hilbert scheme} $\mathcal{H}_P$.
This follows from the existence of a universal family of those subschemes over $\mathcal{H}_P$.
We can study the local properties of $\mathcal{H}_P$ by deformation theory.
In particular, if $Y\subset X$ has Hilbert polynomial $P$ and normal sheaf $\n_{Y/X}:=\Hom(\id_Y,\os_Y)$ then the tangent space to $\mathcal{H}_P$ in $[Y]$ is naturally identified to $H^0(X,\n_{Y/X})$, see \cite{Sernesi}[Thm. 4.3.5].
Moreover, if $Y\subset X$ is a regular embedding and $H^1(X,\n_{Y/X})=0$ then $[Y]$ is a smooth point.

Now, let us assume that $\gk$ is algebraically closed of characteristic $0$.
Fixed a polarised scheme $(X,H)$ and some polynomial $P\in \Q[t]$, we can form the Gieseker-Maruyama moduli space $\cm$ of torsion-free, semi-stable sheaves with Hilbert polynomial $P$, as shown in \cite{HuyLeh}[Thm. 4.3.4]. 
$\cm$ is a projective scheme whose closed points parametrise $S$-equivalence classes of semi-stable sheaves on $X$ with Hilbert polynomial $P$.
$\cm$ is a coarse moduli space, meaning that any family $\scr{E}$ of torsion-free, semi-stable sheaves on $X$ with Hilbert polynomial $P$ flat over a scheme $B$ determines a morphism $B\ra \cm$, sending $b\mt [\scr{E}_b]$.
For a stable sheaf $S$-equivalence and isomorphism class coincide, since its only Jordan-H\"older factor is the sheaf itself.
Hence, we obtain a subset $\cm^s\subset \cm$ whose points actually represent (isomorphism classes of) stable sheaves.
This is an open subscheme.
The completions of the local rings of those points pro-represent the local deformation functor of the sheaf they parametrise.
Therefore, the tangent space to $\cm^s$ in a point $[\e]$ is isomorphic to the space of first order deformations of $\e$, which is $Ext^1(\e,\e)$.
In addition, $Ext^2(\e,\e)$ contains the obstructions, hence if it is $0$ then $\cm^s$ is smooth in $[\e]$, see \cite{HuyLeh}[Thm. 4.5.1 and Cor. 4.5.2].

In the following, we will also need to consider \textbf{simple sheaves}, that is, sheaves $\e$ such that $Hom(\e,\e)\cong\gk$.
Being $\gk$ algebraically closed any stable sheaf is simple, see \cite{HuyLeh}[Cor. 1.2.8], hence in any family of simple sheaves there is an open locus parametrising stable ones.

The proof in \cite{CH}[Prop. 2.10] shows that on a smooth complex variety $X$ any bounded family of simple sheaves $\e$ with Hilbert polynomial $P$ and with $ext^2(\e,\e)=0$ admits a smooth \textbf{modular family} $\scr{F}$.
This means that there exists a scheme $\scr{S}$ and a family $\scr{F}$ of simple sheaves on $X$ flat over $\scr{S}$ with Hilbert polynomial $P$ such that:
\begin{itemize}
    \item each isomorphism class of simple sheaves with Hilbert polynomial $P$ appears at least once and at most finitely many times as $\scr{F}_s$ for some $s\in\scr{S}$
    \item for each $s\in \scr{S}$ the completion $\widehat{\os_{\scr{S},s}}$ of the local ring $\os_{\scr{S},s}$ pro-represents the local deformation functor of $\scr{F}_s$; in particular this means that we can study local properties of $\scr{S}$ with deformation theory, as in the stable case
    \item for any other family $\scr{F}'$ flat over $\scr{S}'$ of such sheaves there is a scheme $\scr{S}''$ with 
    \dia
    \scr{S}'' \ar[r, "q"] \ar[d, "p"] & \scr{S}' \\
    \scr{S}\\
    \mma
    $q$ étale and $p^*\scr{F}\cong p_2^*\scr{F}'$.
\end{itemize}

Consider the functor $\scr{S}pl$ which assigns to any base scheme $B$ the flat families of torsion-free, simple sheaves on $X$ with Hilbert polynomial $P$ flat over $B$, where two families over $B$ are identified if they differ by the publlback of a line bundle from $B$.
It is known that the étale sheafification of $\scr{S}pl$ is representable by an \textit{algebraic space}, see \cite{AltKle_compact}[Thm. 7.4], even with fewer assumptions than Hartshorne.
Therefore, there exists an étale presentation for this functor: a scheme $Spl$ with an étale surjective morphism $\psi:Spl\ra \scr{S}pl$ and to this morphism there corresponds a family on $Spl$ which is modular.
Indeed, each isomorphism class of simple sheaf with the fixed Hilbert polynomial appears as a closed point in $\scr{S}pl$ hence, being $\psi$ étale and surjective, it appear at most finitely many times in $Spl$.
The second listed property follows again by etaleness of $\psi$, since such maps preserve completions of local rings, while for the third is enough to define $S''=S'\times_{\scr{S}pl} Spl$ and recall that any morphism to $S''\ra\scr{S}pl$ determines a family on $S''$. 

\chapter{Ulrich sheaves}

Ulrich vector bundles are the main object of investigation in this work.
In this chapter, we review the basic theory of Ulrich sheaves with respect to an ample and globally generated polarisation.
Along with this, we try also to give some motivation for the study of such sheaves.

Most of the results presented in this chapter are already known, frequently under restrictive hypothesis.
Therefore, we want not only to create a unified reference for the following chapters but also to state and prove everything under the weakest possible assumptions on $X$, $H$, and the field $\gk$ that we are able to work with.
Sometimes, when the proofs already existing in the literature carry over with only minor modifications, we will give only statements and references.
We mainly quote from the works \cite{EisSch}, \cite{CH}, \cite{CKM_clifford}, \cite{Cos}, \cite{Bea} and \cite{AC}, see also the book \cite{CMRPL} on this topic.

We emphasize that working over arbitrary fields has meaningful applications.
For example, in \cite{HanselkaKummer} existence of Ulrich bundles is connected to the old problem of writing a polynomial with real coefficients as a sum of squares of real polynomials.

In the first section, we define Ulrich sheaves, see \thref{defUlrich} and \thref{definizioneUlrich}, and give some of their basic properties, among which regularity, see \thref{corollario}.
As a side remark, we will note that considering ample and globally generated polarizations instead of very ample ones is often more natural and comes at a low cost.
For example, we will show how to reduce the existence problem for Ulrich sheaves on arbitrary schemes to the case of normal varieties.

The second one deals with ways to construct new Ulrich sheaves from already existing ones.
Most of the main results are contained in the literature, but we will extract from them some interesting properties of Ulrich sheaves.
For example, we prove \thref{restrizione} which tells that the restriction of Ulrich sheaves to divisors in the linear system $|H|$ is again Ulrich, from which we deduce geometrical properties of the zero loci of sections of Ulrich bundles \thref{restrizionesezioni}.

In the third section we collect the known existence results on Ulrich sheaves that can be found in the literature, for example we make an overview on Ulrich bundles on Fano $3$-folds, \thref{Fano3i2}, and blow-ups of $\p^3$, \thref{blp3}.
So in particular we focus on curves, surfaces, some special varieties, and low rank Ulrich bundles.

In the last two sections, we study numerical properties like Hilbert polynomial, \thref{hilbertpol}, Chern classes, \thref{c12}, and stability, \thref{stability}, of Ulrich sheaves, which are the prerequisites to, respectively, the search for new examples of such sheaves and the study of their moduli spaces.
Finally, we spend some pages on positivity properties of Ulrich bundles.

\startcontents[chapters]
\printcontents[chapters]{}{1}{}

\section{Introducing Ulrich sheaves}

Fix a field $\gk$.
Let $(X,H)$ be a polarised scheme, that is a proper scheme with an ample and globally generated divisor and denote $\os_X(iH)$ simply by $\os_X(i)$.
We assume $X$ to be an equidimensional scheme of dimension $n$; this will be necessary to state the definition of Ulrich sheaf.
Moreover, without loss of generality, we will always assume that $X$ is connected.

We will start by recalling various equivalent definitions of Ulrich sheaves.
Then in the second, third and fourth sections, we will deduce some properties of such objects.
Finally we discuss a bit the existence conjectures.

\subsection{Definition}

\De\label{defUlrich}
A coherent sheaf $\e$ on a polarised scheme $(X,H)$ of dimension $n$ is said \textbf{Ulrich} if 
\[h^j(X,\e(iH))=0 \qquad \text{for} \;\; 0\leq j\leq n,\, -n\leq i\leq -1.\]
\Ne

If $n=0$ then any sheaf is Ulrich so we will ignore this case.
In addition, we always tacitly assume $\e\neq 0$, since the zero sheaf trivially satisfies the definition.
Note that, contrary to most of the literature, we do not assume $H$ to be very ample.

\No
When the scheme $X$ is clear from the context, we will frequently call $\e$ an $H$-Ulrich sheaf, or $\os_X(H)$-Ulrich sheaf.
The datum of an ample and globally generated line bundle $\os_X(H)$ is equivalent to that of a finite map $\phi:X\ra\p^N$ such that $\phi^*\osn(1)\cong \os_X(H)$.
Therefore, we will simply call $\phi$-Ulrich an Ulrich sheaf on $(X,\phi^*\osn(1))$. 
\ione

We will start by noting that the property of being Ulrich has a nice functorial behaviour for finite morphisms of polarised pairs.
Note that, for a finite morphism $g$, $g^*$ preserve the property of being ample and globally generated, but not of being very ample in general.

\Le\label{proiezione}
Let $(X,H), (X',H')$ be two polarised schemes of dimension $n$.
Suppose that $g:X'\ra X$ is a finite morphism such that $g^*H\sim H'$ and that $\e$ is a coherent sheaf on $X'$.
Then, $\e$ is $H'$-Ulrich if and only if $g_*\e$ is $H$-Ulrich.
\ma
\begin{proof}
    Since the map $g$ is finite it preserves cohomology, i.e. we have $h^j(X',\e)=h^j(X',g_*\e)$ for all $j\in \Z$.
    Moreover, by \textit{projection formula}, see \cite{GW2}[Thm. 22.81], we have $g_*(\e\otimes g^*\os_X(iH))\cong g_*(\e)\otimes\os_X(iH)$ for all $i\in \Z$.
    Putting those two observations together we conclude that 
    \[h^j(X',\e(iH'))=h^j(X',\e\otimes g^*\os_X(iH))=h^j(X,g_*(\e\otimes g^*\os_X(iH)))=\]
    \[=h^j(X,g_*(\e)\otimes\os_X(iH)),\]
    from which the thesis follows.
\end{proof}

Next, we will give some characterisations of Ulrich sheaves.
If in \thref{proiezione} we choose $g$ to be the morphism given by $|H|$ then the above result allows us to reduce to the very ample case.
Most part of the proofs are taken from \cite{EisSch}[Thm. 2.1], \cite{Bea}[Thm. 2.3] or \cite{AC}[Thm. 1.4].

\Te\label{definizioneUlrich}
Let $(X,H)$ be a polarised, $n$-dimensional scheme over $\gk$ and $\e$ a coherent sheaf on $X$.
The following are equivalent:
\begin{enumerate}[i)]
    \item $\e$ is an Ulrich sheaf for $(X,\os_X(H))$
    \item $h^j(X,\e(-j))=0$ for $1\leq j\leq n$ and $h^{j}(X,\e(-j-1))=0$ for $0\leq j\leq n-1$
    \item Denote $\phi:X\ra\p^N$ the finite morphism given by the complete linear system $|H|$ and define $c:=N-n$, then we have a \textbf{linear resolution} 
    \[0\ra \os_{\p^N}(-c)^{r_c}\ra \dots \ra \os_{\p^N}(-1)^{r_1}\ra \os_{\p^N}^{r_0}\ra \phi_*\e\ra 0\]
    \item $h^j(X,\e(i))=0\;$  if  $\; \begin{cases}
        j=0, \;i<0 \\
        1\leq j\leq n-1,\; i\in \Z \\
        j=n,\; i\geq -n.\\
    \end{cases}$
   
\end{enumerate}
Moreover, if there exists a finite morphism $f:X\ra\p^n$ such that $f^*\os_{\p^n}(1)=\os_X(H)$\footnote{for example, it holds if $\gk$ is infinite by \cite{GW1}[Thm. 13.89]} then the above conditions are equivalent to:
\begin{enumerate}
    \item[v)] $f_*\e\cong \osn^\rho$.
\end{enumerate}
\Ma
\begin{proof}
    \underline{$\mathbf{i)\Rightarrow ii)}$} 
    Is trivial.
    
    \underline{$\mathbf{ii)\Rightarrow iii)}$} 
    Since, by definition, we have $\phi^*\osn(1)=\os_X(H)$ then by projection formula and the fact that $\phi_*$ preserves cohomology, we deduce that $h^j(\p^N,\phi_*\e(-j))=0=h^{j-1}(\p^N,\phi_*\e(-j))$ for $1\leq j\leq n$.
    Moreover, being the support of $\phi_*(\e)$ contained in $\phi(X)$ hence $n$-dimensional, by Grothendieck's vanishing theorem, see \cite{GW2}[Thm. 21.57]. we have $h^j(\p^N,\phi_*\e(-j))=0$ for all $j>n$.
    
    We define inductively sheaves $\ck_l$ for $l=0,\dots, c+1$ such that 
    \begin{itemize}
        \item $\ck_0=\phi_*\e$
        \item $h^j(\p^N,\ck_l(-j))=0$ for $1\leq j\leq N$ and $h^{j}(\p^N,\ck_l(-1-j))=0$ for $0\leq j\leq n+l-1$
        \item evaluation of global sections of $\ck_l$ gives an exact sequence 
        \sesl{\ck_{l+1}(-1)}{H^0(\p^N,\ck_l)\otimes\os_{\p^N}}{\ck_l}{gg''}
    \end{itemize}
    The plan is as follows.
    We already proved that $\phi_*\e$ verifies the second requirement.
    Next, we show that if $\ck_l$ satisfies the second condition than it is globally generated, so that $\ck_{l+1}$ is well defined.
    From the sequence obtained in this way, we prove that $\ck_{l+1}$ satisfies the second requirement, so that we can argue inductively.
    Finally, we see that $\ck_{c+1}=0$ so that $\ck_c\cong \os_{\p^N}^{r_c}$ and glueing the sequences \eqref{gg''} for all $l$ we get the claimed resolution, with $r_l=h^0(\p^N,\ck_l)$.
    
    If $h^i(\p^N,\ck_l(-i))=0$ for $1\leq i\leq N$ then $\ck_l$ is $0$ regular for Castelnuovo-Mumford, see \cite{Mum}[Lecture 14] so, by the proposition in the above reference, we know that $\ck_l$ is globally generated.

    We will suppose that the second condition is true for $\ck_l$ for some $0\leq l\leq c$ and deduce it holds also for $\ck_{l+1}$.   
    Twisting \eqref{gg''} by $\os_{\p^N}(-p)$ and taking the long exact sequence in cohomology we deduce that $h^q(\p^N,\ck_l(-p))=h^{q+1}(\p^N,\ck_{l+1}(-p-1))$ for $1\leq p\leq N$ and any $q$.
    In particular, as soon as $1\leq j-1\leq N-1$ we can put $p=q=j-1$ and get $h^{j}(\p^N,\ck_{l+1}(-j))=0$ for $2\leq j\leq N$.
    Similarly, putting $q=j-1$ and $p=j$ we get $h^{j}(\p^N,\ck_{l+1}(-1-j))=0$ for $1\leq j\leq n+(l+1)-1$.
    Being the left morphism in \eqref{gg''} the evaluation of global sections of $\ck_l$, we know that $h^0(\p^N,\ck_{l+1}(-1))=0$. 
    Moreover, being $h^1(\p^N,\os_{\p^N})=0$, by the above cohomology sequence we also get $h^1(\p^N,\ck_l(-1))=0$, which was the last vanishing we desired.
    
    Finally, by the sequence \eqref{gg''} with $l=c$ we get $h^0(\p^N,\ck_{c+1}(-1))=0$.
    But, by the second property $\ck_{c+1}(-1)$ is $0$-regular hence globally generated, so we conclude that $\ck_{c+1}(-1)=0$ therefore $\ck_c\cong H^0(\p^N,\ck_c)\otimes\os_{\p^N}$.
    
    \underline{$\mathbf{iii)\Rightarrow i)}$}
    As in the previous part, we can break the above resolution in short exact sequences of the form
    \ses{\ck_{l+1}(-l-1)}{H^0(\p^N,\ck_l)\otimes\os_{\p^N}(-l)}{\ck_l(-l)}
    where $\ck_0\cong \phi_*\e$ and $\ck_c\cong \os_{\p^N}^{r_c}$.
    Recalling \eqref{bott}, if we are in one of the following three cases
    \begin{itemize}
        \item $j=0$ and $i<0 $
        \item $1\leq j\leq n-1$ and $i\in \Z$ 
        \item $j=n$ and $i\geq -n$
    \end{itemize} 
    from the cohomology of the above short sequences tensored by $\os_{\p^N}(i)$ we get
    \[h^j(\p^N,\phi_*\e(i))= h^{j+c}(\p^N,\ck_c(i-c))=h^{j+c}(\p^N,\os_{\p^N}(i-c))=0.\]
    Recall that by projection formula and the fact that finite morphisms preserve cohomology we get that $h^j(\p^N,\phi_*\e(i))=h^j(X,\e(i))$, so we conclude.
    
    \underline{$\mathbf{iv)\Rightarrow i)}$}
    Obvious.
    
    \underline{$\mathbf{i)\Leftrightarrow v)}$}
    From $iii)$, we know that $\osn^\rho$ are the only Ulrich sheaves on $(\p^n,\osn(1))$ ($\phi=id$), therefore the claim is an application of \thref{proiezione}.
\end{proof}

\subsection{First properties}

We present some simple properties which are easy consequences of the definition(s) of Ulrich sheaf.
All of these are well-known.
This appears in \cite{CKM_clifford}[Prop. 2.14].
\Le\label{2-3}
Let $(X,H)$ be as always our polarised scheme.
Suppose we have an exact sequence of coherent sheaves
\[0\ra\e\ra\f\ra\g\ra 0.\]
If any two among $\e,\f,\g$ are Ulrich sheaves then also the last one is.
Moreover, if the sequence is split, that is $\f\cong \e\oplus\g$, and $\f$ is Ulrich then both $\e,\g$ must be Ulrich.
\ma
\begin{proof}
    Fix $-n\leq i\leq -1$ and consider the long exact sequence in cohomology given by
    \[0\ra\e(i)\ra\f(i)\ra\g(i)\ra 0.\]
    If two among $\e,\f,\g$ are Ulrich then each cohomology group of the third one is surrounded by zeroes, hence must vanish.
    In the split case we know that $H^j(X,\f(i))\cong H^j(X,\e(i))\oplus H^j(X,\g(i))$ for all $i,j$, so the conclusion is clear.
\end{proof}

As already noted, we can easily characterize Ulrich sheaves on $(\p^n,\osn(1))$ see, for example, \cite{Bea}[§4]
\Le\label{oulrich}
As always, let $(X,H)$ be our polarised scheme.
\begin{itemize}
    \item If a line bundle of the form $\os_X(i)$ is Ulrich then $i=0$ and $(X,H)\cong (\p^n,\osn(1))$.
    \item The only Ulrich sheaves on $(X,H)\cong (\p^n,\osn(1))$ are direct sums of $\osn$.
\end{itemize}
\ma
\begin{proof}
    If $\os_X(i)$ is Ulrich then $h^0(X,\os_X(i-1))=0$, hence $i\leq 0$, but from \thref{definizioneUlrich} $iii)$ an Ulrich sheaf has sections, hence $i\geq 0$.
    So we conclude that $i=0$.
    Call $\phi:X\ra\p^N$ the map induced by the linear system $|H|$.
    Since $h^0(\p^n,\phi_*\os_X)=h^0(X,\os_X)=1$, from \thref{definizioneUlrich} $iii)$ we get a surjective $\os_{\p^N}\twoheadrightarrow\phi_*\os_X$ hence, being the map $\phi$ finite, by \cite{GW1}[Cor. 12.2] $\phi$ is the closed embedding $Spec(f_*\os_X)\subset\p^N$.
    Moreover, if $X\neq \p^N$, then its codimension is positive hence, again by \thref{definizioneUlrich} $iii)$, we get a surjective morphism $\os_{\p^N}(-1)^{r_1}\ra\id_X$.
    But only the ideals of linear subspaces are generated by linear polynomials hence we have our first claim.

Similarly, if $(X,H)\cong (\p^n,\osn(1))$ then $\phi$ is an isomorphism and by \thref{definizioneUlrich} $iii)$ the only Ulrich sheaves are $\osn^{r_0}$ for some $r_0\geq 0$.
\end{proof}

\Co\label{pic1'}
If $Pic(X)\cong \Z H$ then there are no Ulrich line bundles on $(X,H)$ except that for $(\p^n,\osn(1))$.
\io

\subsection{Local properties of Ulrich sheaves}

In this subsection, we will prove regularity properties of Ulrich sheaves, clearly depending on regularity of $X$ itself.
We start recalling some terminology.

\De
Suppose $X$ is a scheme and $\e$ a coherent sheaf.
We call $\e$ \textbf{maximally Cohen-Macaulay} (mCM) if $depth_X(\e_x)=dim(\os_{X,x})$ for all $x\in X$.
\Ne
This property can be checked just by looking at closed points $x\in X$: indeed, once it is guaranteed on those, then $Supp(\e)=X$ and the other local rings are obtained by localization, hence we just apply \cite{Stacks}[\href{https://stacks.math.columbia.edu/tag/0AAG}{Tag 0AAG}].
in particular, the stalks of a mCM sheaf are Cohen-Macaulay, in particular satisfy $S_k$ for all $k$, and of the highest possible dimension.
Let us show a couple of properties of mCM sheaves before continuing.

\Le\label{lcmlf}$\qquad$ 

\begin{itemize}
    \item A mCM sheaf $\e$ on a smooth scheme $X$ is locally free.
    \item If $g:X\ra Y$ is finite and surjective then $\e$ is mCM on $X$ if and only if $g_*\e$ is mCM on $Y$. 
\end{itemize}
\ma
\begin{proof} 
    For the first claim, since $X$ is smooth, by Auslander-Buchsbaum formula we have $pd(\e_x)=dim(\os_{X,x})-depth(\e_x)=0$, hence we conclude by \thref{extpd}.
    
    For the second, we can suppose $X,Y$ equidimensional, then $dim(X)=dim(Y)$ being $g$ finite and surjective, hence we have $dim(\os_{X,x})=dim(\os_{Y,g(x)})$.
    Using \thref{depthfinite} we get $depth_Y(g_*\e)_y=min_i\{depth_X\e_{x_i}\}$, where for all $y\in Y$ we denote by $x_i$ the points in $g\inv(y)$. 
    Now it is clear that $depth_Y(g_*\e)_y=dim(\os_{Y,g(x)})$ if and only if $min_i\{depth_X\e_{x_i}\}=dim(\os_{X,x})$ if and only if $depth_X\e_{x_i}=dim(\os_{X,x})$.
\end{proof}

The following claims are well-known in case $X$ is smooth.

\Te\label{corollario}
Let $(X,H)$ be a polarized $n$-dimensional scheme.
If $\e$ is an $H$-Ulrich sheaf then:
\begin{enumerate}[i)]
\item $\e$ is globally generated.
\item $\e$ is Cohen-Macaulay on its support.
\item If $X$ is irreducible then $\e$ has full support, in particular is mCM.
\item If $X$ is integral then $\e$ is torsion-free.
\item If $X$ is a normal variety then $\e$ is reflexive.
\item If $X$ is a smooth variety then $\e$ is locally free.
\end{enumerate}
\Ma 
\begin{proof}
\gel{i}
By \thref{definizioneUlrich} $iii)$, recall that $\phi:X\ra\p^N$ is the finite morphism given by $|H|$, we have a surjective map $\os_{\p^N}^{r_0}\ra \phi_*\e$.
Since pullback is right exact we get a surjective morphism $\os_{X}^{r_0}\ra \phi^*\phi_*\e$, therefore it is enough to show that the map $\phi^*\phi_*\e\ra\e$ given by the adjunction of $\phi_*$ and $\phi^*$ is surjective.
But this is the content of \thref{counitsurj}.

\gel{ii}
Our goal is to prove that for any closed point $x\in Supp(\e_x)$ we have $depth(\e_x)= dim(\os_{X,x})$.
For any $x\in Supp(\e)$ set $y:=\phi(x)$.
From \thref{depthfinite} we know that $depth_X(\e_x)\geq depth_{\p^N}(\phi_*\e)_y$. 
Recall that by \thref{definizioneUlrich} $iii)$ we have a locally free resolution of length $c$ of $\phi_*(\e)$ hence $pd_{\p^N}((\phi_*\e)_y)\leq c$.
From Auslander-Buchsbaum formula we get 
\[depth_{\p^N}(\phi_*(\e)_y)=dim(\os_{\p^N,y})-pd_{\p^N}((\phi_*\e)_y)\geq dim(\os_{\p^N,y})-c=dim(\os_{X,x}),\]
hence $depth_X(\e_x)\geq dim(\os_{X,x})$, so that $\e$ is Cohen-Macaulay on its support.

\gel{iii}
We continue with the above notation.
By \thref{codimensione}, we know that $codim(\phi_*\e_x)\le pd(\e_x)\le c$ hence $dim(Supp(\phi_*\e))\ge N-c=n$.
Since $\phi$ is finite we have $dim(\phi(X))=n$ hence, being $X$ and $\phi(X)$ irreducible, $\phi_*\e$ has full support on $\phi(X)$ and the same applies to $\e$.

\gel{iv} 
Since $X$ is integral and it follows by $ii)$ that $\e$ satisfies $S_1$ we just apply \thref{s12}. 

\gel{v}  
Being $X$ integral, by $iii)$ the sheaf $\e$ is torsion-free.
Moreover, by $ii)$ it satisfies $S_2$ therefore, since $X$ is normal, we can conclude applying \thref{s12}.

\gel{vi}
Being $X$ smooth, this is a consequence of $ii)$ and \thref{lcmlf}.
\end{proof}

As a corollary of the previous theorem, on an integral scheme, we can identify Ulrich sheaves inside the class of initialised sheaves without intermediate cohomology as the ones having the maximal possible number of sections.
This is analogue to the original definition of Ulrich modules, given in \cite{Ulrich}, as maximal Cohen-Macaulay modules with the maximal cardinality for a minimal set of generators.
Let us start recalling the needed definitions.
\De
Let $(X,H)$ be a polarised scheme of dimension $n$.
We say that $\e$ is \textbf{initialised} if $h^0(X,\e(-1))=0$ but $h^0(X,\e)\neq 0$.
We say that $\e$ \textbf{has no intermediate cohomology} if
\[h^j(X,\e(i))=0 \qquad \forall \; i\in\Z,\; 1\leq j\leq n-1.\]
A mCM sheaf without intermediate cohomology is called \textbf{arithmetically Cohen-Macaulay} (aCM).
\Ne

In \cite{EisSch}[Prop. 2.1] it has already been shown that for very ample polarisations a sheaf is Ulrich if and only if its graded module of sections satisfies this maximality property.
More similar in spirit, but different in form, is the proof in \cite{Casnati2017}[Thm. 5.4], which covers the rank $2$ case under some regularity assumption.

\Co
Let $(X,H)$ be a polarised integral scheme over an infinite field $\gk$.
Let $\e$ be an initialized aCM sheaf. 
Then $\e$ is Ulrich if and only if it has the maximal possible $h^0(X,\e)$ among aCM initialised sheaves of the same rank $r$.
\io
\begin{proof}
    Consider a finite surjective morphism $f:X\ra\p^n$ such that $f^*\osn(1)\cong\os_X(1)$.
    Define $d:=H^n$ the degree of $f$.
    Since $\e$ has no intermediate cohomology and is initialised, by projection formula and the fact that finite morphisms preserve cohomology, the same applies to $f_*\e$.
    Moreover, $f_*\e$ is mCM by \thref{lcmlf} since $\e$ is mCM by \thref{corollario}.
    We deduce that $f_*\e$ is locally free by \thref{lcmlf}.
    By Horrocks' splitting criterion \cite{OSS}[Thm. 2.3.1]\footnote{This result is stated with $\gk=\C$ but is not necessary for the proof} we must have $f_*\e\cong \oplus_{i=1}^{rd}\osn(a_i)$, where we used $rk(f_*\e)=rd$.
    Being $f_*\e$ initialised we have $a_i\leq 0$ and hence $h^0(X,\e)\leq rd$.
    Moreover, the equality holds if and only if $a_i=0$ for all $i$, that is equivalent to $\e$ being Ulrich by \thref{definizioneUlrich}.
\end{proof}

Since most of the literature deals with smooth varieties, and hence Ulrich bundles, let us present a couple of pathologic examples.

\Es
\begin{itemize}
\item If $X$ is reducible then an Ulrich sheaf could be non-lCM simply because its support could be different from the entire $X$.
For example, suppose $X=X_1\cup X_2$ is an equidimensional scheme.
Then, any Ulrich sheaf for $(X_1,H|_{X_1})$ is Ulrich also seen as a sheaf on $X$.
Indeed, the inclusion $X_1\hookrightarrow X$ is finite and then we can apply \thref{proiezione}.

\item We show an example of non-locally free but torsion-free Ulrich sheaf.
Take a nodal integral plane cubic curve and its normalization $\nu:\p^1\ra C\subset \p^2$.
We have $\nu^*\os_C(1)\cong \os_{\p^1}(3)$ and it is easily verified that $\os_{\p^1}(2)$ is Ulrich on $(\p^1,\os_{\p^1}(3))$.
Therefore, by \thref{proiezione}, we know that $\nu_*\os_{\p^1}(2)$ is Ulrich on $(C,\os_C(1))$, we need only to show that is not locally free.

Let $P\in C$ be a smooth point.
$\nu^*\os_C(P)\cong \os_{\p^1}(1)$ so by projection formula we have 
\[\nu_*\os_{\p^1}(2)\cong \nu_*(\os_{\p^1}\otimes \nu^*\os_C(2P))\cong \nu_*(\os_{\p^1})\otimes\os_C(2P)\]
and, being $\os_C(2P)$ locally free, $\nu_*\os_{\p^1}(2)$ is locally free if and only if $\nu_*(\os_{\p^1})$ is.
But the last sheaf is not locally free as shown in \cite{Har}[III Example 9.7.1].
\end{itemize}
\Io

\subsection{U-Duality}

Recall that for a Gorenstein scheme the dualizing sheaf $\omega_X$ is a line bundle.

\Prop\label{duale}
Let $(X,H)$ be a Gorenstein polarised scheme of dimension $n$.
Then a vector bundle $\e$ is $H$-Ulrich if and only if $\e\du(n+1)\otimes\omega_X$ is.
\One
\begin{proof}
$\e$ is $H$-Ulrich if and only if
\[h^j(X,\e(i))=0 \quad \forall \; j\in\Z,\; -n\leq i\leq -1.\]
By Serre duality and reflexivity of locally free sheaves, this is equivalent to
\begin{equation}\label{1}
    h^{n-j}(X,\e\du(-i)\otimes\omega_X)=0 \quad  \; j\in\Z\; -n\leq i\leq -1.
\end{equation}
Now set
\[n+1+i'=-i, \qquad  j'=n-j,\]
then the previous equations become
\[h^{j'}(X,\e\du(n+1+i')\otimes\omega_X)=h^{j'}(X,(\e\du(n+1)\otimes\omega_X)(i'))=0 \quad  j\in\Z,\; -n\leq i'\leq -1\]
which is exactly the $H$-Ulrich condition for  $\e\du(n+1)\otimes\omega_X$.
\end{proof}

\De\label{dual}
For any vector bundle, we call $\e\du(n+1)\otimes\omega_X$ the \textbf{U-dual} of $\e$ and we say that $\e$ is \textbf{self U-dual} if $\e\cong \e\du(n+1)\otimes\omega_X$.
\Ne
\thref{duale} tells us that a vector bundle is Ulrich if and only if its U-dual is.

In case a sheaf has this kind of symmetry, it has been noted in \cite{EisSch}[Cor. 2.3] that the Ulrich condition becomes much simpler due to Serre duality and a theorem of Mumford, see \cite{Mum}[Lecture 14].
\Le
Let $(X,H)$ be a Gorenstein polarised $n$-dimensional scheme.
A self U-dual vector bundle $\e$ is Ulrich if and only if $h^i(X,\e(-i))=0$ for all $i>0$, i.e. $\e$ is $0$-regular.
\ma
See also \cite{EisSch}[§3] for a more general discussion of this property.

\subsection{The existence problem}\label{section1.5}

We start by explaining how, if we are only interested in the existence of Ulrich sheaves, then we can easily reduce to the case of normal varieties and indecomposable sheaves.

\Oss\label{riduzioneesistenza}
Given any scheme $X$ over $\gk$ we can find a normal variety $\Tilde{X}$ and a finite morphism $\nu:\Tilde{X}\ra X$.
Indeed, it is enough to choose one irreducible component of $X$, consider it with the reduced induced structure, and then normalise it: this is a candidate $\Tilde{X}$ since all those operations give rise to finite morphisms.
Suppose now $X$ is equidimensional and we fix an ample and globally generated polarisation $H$ on it.
Then, by \thref{proiezione}, $\e$ is an Ulrich sheaf on $(\Tilde{X},\nu^*H)$ if and only if $\nu_*\e$ is an Ulrich sheaf on $(X,H)$.
Furthermore, such an $\e$ is always reflexive by \thref{corollario}.
The sheaf $\e$ can be assumed indecomposable by \thref{2-3}.

In particular, note that $\nu^*$ preserves the property of being ample and globally generated but not the one of being very ample in general, so if we want to make this reduction, we have to allow ample and globally generated polarizations.
\one

The reverse argument does not hold in general.
In the next example, we will see how an Ulrich bundle on a reducible, but connected, curve may not restrict to Ulrich bundles on its components.

\Es
The tangent bundle $\ct_{\p^2}$ is known to be Ulrich on $(\p^2,\on{2}(2))$, for example in \cite{EisSch}[Thm. 5.9] for $d=2$ we get the Euler sequence. 
Then, by \thref{restrizione}, for any $Q\in |\on{2}(2)|$ also $\ct_{\p^2}|_Q$ is Ulrich on $(Q,\os_Q(2))$.
On the other hand, $Q$ can be chosen to be the union of two lines; call $L$ one of them.
But we have $\ct_{\p^2}|_L\cong \on{1}(1)\oplus\on{1}(2)$ which is not Ulrich.
\Io

Here we list some variations of the famous problem posed by Eisenbud and Schreyer in their seminal paper \cite{EisSch}.
\Do\label{congetturaesistenza}

\noindent\begin{itemize}
    \item (Strong) Does every polarised variety $(X,H)$ admit an Ulrich sheaf?
    \item (Original) Does every polarised variety $(X,H)$ admit an Ulrich sheaf if $H$ is very ample?
    \item (Weak) Given a polarised variety $(X,H)$, is there an $l>0$ such that $(X,lH)$ admits an Ulrich sheaf?
    \item (Very weak) Given a variety $X$, is there a polarisation $H$ such that $(X,H)$ admits an Ulrich sheaf?
    \end{itemize}
    \da

   The interest of the authors towards these sheaves came from the main result of the above paper.
   To introduce this, let us recall the notion of \textbf{Chow Form}.
   Let $X\subset \p^N$ be a codimension $c$ subvariety.
   Consider the incidence variety $\mathscr{I}$ parametrizing pairs $(P,L)$ where $P\in\p^N$ and $L$ is a $(c-1)$-dimensional linear subspace of $\p^N$ containing $P$.
   Considering the two projections gives the following diagram
   \dia
& \mathscr{I} \ar[dr, "\pi_2"] \ar[ld, "\pi_1"]\\
\p^N  & & Gr(c,N+1) \\
\mma
where $Gr(c,N+1)$ is the Grassmannian parametrising $(c-1)$-dimensional linear subspaces of $\p^N$.
The Chow form of $X$ is defined as the equation of the divisor $\pi_2(\pi_1\inv X)$.
If we know that $\pi_2(\pi_1\inv X)$ is a divisor, then its equation can be expressed as a polynomial since $Pic(Gr(c,N+1))\cong \Z$ and the ample generator is actually very ample and induces a linearly normal embedding, called \textit{Plücker embedding}.

Let us show that this is indeed a divisor.
A general linear variety of dimension $c-1$ intersecting $X$ will do it in just one point, hence $\pi_2|_{\pi_1\inv(X)}$ is birational.
Recall that $dim(Gr(l,k))=l(k-l)$ hence $dim(Gr(c,N+1))=c(N+1-c)$.
The fiber of $\pi_1$ over a point $P$ can be identified with $Gr(c-1,N)$ hence has dimension $(c-1)(N+1-c)$ and $\pi_1\inv X$ is irreducible. 
It follows that \[dim(\pi_2(\pi_1\inv X))=N-c+(c-1)(N+1-c)=c(N+1-c)-1\]
hence $\pi_2(\pi_1\inv X)$ is a divisor.

Still assuming that $X$ is integral, we can also define the Chow form of a sheaf $\f$ on $X$.
Since $X$ is integral, $\f$ has a well defined rank $r$ and we define its Chow form to be the $r$-th power of the one of $X$.
Therefore we can rephrase \cite{EisSch}[Thm. 0.3] as follows.

   \Te
If there is an Ulrich sheaf of rank $r$ on a variety $X$ then the $r$-th power of its Chow form can be written as the determinant of a matrix whose entries are linear forms on $Gr(c,N+1)$.
In case $r=2$ we can actually write this Chow form as the Pfaffian of an anti-symmetric matrix.
   \Ma

   Moreover, in Eisenbud and Schreyer's paper also the following problem is presented.
\Do   
     If there are Ulrich sheaves on $(X,H)$, what is the minimum possible rank of such a sheaf?
     We call this rank \textbf{Ulrich complexity} of the pair.
\da

There is a known obstruction to the existence of Ulrich line bundles, see \cite{Bea}[§ 4] and \cite{AC}[Thm. 8.2].
\Prop\label{pic1}
Let $(X,H)$ be a smooth polarized variety of dimension $n\geq 2$ over an algebraically closed field $\mathbf{k}$.
Suppose $Pic(X)\cong \Z=<\os_X(L)>$ with $h^0(X,\os_X(L))\neq 0$.
Then there are Ulrich line bundles on $(X,H)$ if and only if it is isomorphic to $(\p^n,\osn(1))$.
\One
\begin{proof}
From our assumption $\os_X(H)=\os_X(hL)$ with $h>0$, since $H$ is ample and $L$ effective.
 Suppose by contradiction there is an Ulrich line bundle $\e=\os_X(eL)$ for some $e\in\Z$.
 Being $X$ smooth the divisor $K_X$ is Cartier and hence $K_X\sim kL$.
If $\e$ is $H$-Ulrich then we must have 
 \[h^0(X,\os_X((e-h)L))=0 \qquad 0=h^n(X,\os_X((e-nh)L))=h^0(X,\os_X((nh-e+k)L)),\]
 where we used Serre duality.
Being $L$ effective, from the previous equations we get
 \[e-h<0 \qquad nh-e+k<0\]
 which is equivalent to
\begin{equation}\label{diseq}
h-1\geq e \qquad e\geq nh+k+1.
\end{equation}
Combining these two we get $h-1\geq nh+k+1$ and so
 \[  k\leq -(n-1)h-2\leq -n-1. \] 
 But \cite{KacKol}[Thm. 1] implies that this is equivalent to $X\cong \p^n$, which clearly gives $\os_X(L)\cong\osn(1)$. 

 Being $n>1$ then $k=-n-1$ and $\eqref{diseq}$ implies
 \[nh\leq e-k-1=e+n\leq h-1+n \; \Rightarrow \; h\leq 1\] so $h=1$ hence $H\sim L$.
 \end{proof}

\section{Constructing new Ulrich sheaves from older ones}

In this section, we will always assume that $\gk$ is infinite.
Recall that under this assumption, given an $n$-dimensional polarised scheme $(X,H)$, we can always construct a finite (surjective) morphism $f:X\ra\p^n$ such that $f^*\osn(1)\cong \os_X(1)$, \cite{GW1}[Thm. 13.89].
From now on the symnol $f$ will be reserved to such a morphism.

In this section we will show how, starting with Ulrich sheaves on $(X,H)$, we can construct Ulrich sheaves on some "modification" of it.
As an example we already saw how Ulrich sheaves behave for pushforward under finite morphisms, see \thref{proiezione}.
Here we will treat the case of pullbacks, extensions, products, blow-ups and change of polarization.
Only in the first three cases we will actually prove some original result, while for the other we limit ourselves in giving statements and references.

\subsection{Intersections and pullbacks}

In general, we cannot expect an Ulrich sheaf to pullback to another Ulrich sheaf under any morphism.
Indeed, $\osn$ is Ulrich on $(\p^n,\osn(1))$ and pullbacks always to the structure sheaf which can be Ulrich only on $(\p^m,\os_{\p^m}(1))$ by \thref{oulrich}.
On the other hand, in some special situations we can do something.
We start with a known result, see \cite{Casnati_wild}[Thm. 1.3] or \cite{CMRPL}[Thm. 4.3.1], which we prove in slightly more generality.
\Prop\label{intersezione}
Suppose $X_1,X_2\subset \p^N$ are equidimensional closed subschemes of codimension $c_1,c_2>0$ whose schematic intersection $X:=X_1\cap X_2$ is non-empty and equidimensional of codimension $c:=c_1+c_2$.
Moreover, assume $X_1$ integral.
If $\e_1,\e_2$ are Ulrich sheaves on $X_1,X_2$ with respect to the restrictions of $\on{N}(1)$ then, $\e_1\otimes\e_2$ is Ulrich on $(X,\os_X(1))$.
\One
\begin{proof}
    Consider the resolution of $\e_2$ 
    \[0\ra\os_{\p^N}^{r_{c_2}}(-c_2)\ra \dots \ra \os_{\p^N}^{r_{0}}\ra \e_2\ra 0\]
    given by \thref{definizioneUlrich} $iii)$.
    Note that exactness and the fact that $\e_2$ is supported on a proper subscheme implies that $\sum_{i=0}^{c_2}(-1)^ir_i=0$.
    Now tensor this resolution by $\e_1$.
    We get a sequence on the variety $X_1$
    \begin{equation}\label{twist}
        0\ra\e_1^{r_{c_2}}(-c_2)\xra{\phi} \dots \ra \e_1^{r_{0}}\ra \e_1\otimes\e_2\ra 0
    \end{equation}
    which is clearly right exact.
    Actually, it is also left exact: indeed, $\e_1\otimes\e_2$ is supported on a proper subscheme of $X_1$, but $X_1$ is integral hence the kernel of $\phi$ would be a rank $\sum_{i=0}^{c_2}(-1)^irk(\e_1)\cdot r_i =0$, hence torsion, subsheaf of $\e_1^{r_{c_2}}(-c_2)$ but $\e_1$ is torsion-free by \thref{corollario}.

Now, $\e_1\otimes\e_2$ is a sheaf on $X$, which has codimension $c=c_1+c_2$, hence dimension $n-c$, so to prove that it is Ulrich we need to show that $h^j(X,\e_1\otimes\e_2(-i))=0$ for all $j$ if $1\le i\le n-c_1-c_2$. 
Note that $\e_1$ is Ulrich on $X_1$ which has codimension $c_1$ hence dimension $n-c_1$.
It follows that $h^j(X_1,\e_1(-l))=0$ for all $j$ if $1\le l\le n-c_1$.
Therefore, if we twist \eqref{twist} by $\os_{X_1}(-i)$ with $1\le i\le n-c_1-c_2$, then all the terms except $\e_1\otimes\e_2(-i)$ are direct sums of $\e_1(-l)$ for some $1\le l\le n-c_1$ hence have no cohomology and we conclude the same for $\e_1\otimes\e_2(-i)$.
\end{proof}

If the polarisation is not very ample we will need to adapt a bit the statement.
Suppose we have two polarised schemes $(X_i,H_i)$.
We can always suppose $H_i$ give two finite maps to the same projective space such that at least one of them is also surjective.
In this setting we can apply the following.

\Prop\label{pullback}
Suppose we have two equidimensional schemes $X_1,X_2$ as in
\dia
X_3:=X_1\times_{\p^n}X_2 \ar[r, "p_2"] \ar[d, "p_1"] & X_2 \ar[d,twoheadrightarrow ,"f_2"] \\
X_1 \ar[r, "f_1"] & \p^n \\
\mma
with finite morphisms $f_1,f_2$, such that $f_2$ is also surjective.
If $\e_i$ is an Ulrich sheaf on $(X_i,f_i^*\osn(1))$ then $p_1^*\e_1\otimes p_2^*\e_2$ is Ulrich on $(X_3,p_1^*f_1^*\osn(1))$.
\One
\begin{proof}
    First of all, $p_1,p_2$ are pullbacks of finite morphisms and hence finite, so that $p_1\circ f_1=p_2\circ f_2$ is finite and $(X_3,p_1^*f_1^*\osn(1))$ is a polarised variety.
    Moreover, $p_1$ is surjective being $f_2$ such.
    By projection formula we have $(p_1)_*(p_1^*\e_1\otimes p_2^*\e_2)\cong \e_1\otimes (p_1)_*p_2^*(\e_2)$.
    Being $f_2$ finite, by base change \cite{GW1}[Thm. 12.6], we have 
    \[(p_1)_*p_2^*(\e_2)\cong (f_1)^*(f_2)_*\e_2\cong (f_1)^*\osn^N\cong \os_{X_1}^N\]
    where we used \thref{definizioneUlrich} $vi)$.
    We conclude that $(p_1)_*(p_1^*\e_1\otimes p_2^*\e_2)\cong \e_1^N$ is Ulrich on $(X_1,f_1^*\osn(1))$ hence by \thref{proiezione} we proved that $p_1^*\e_1\otimes p_2^*\e_2$ is Ulrich.    
\end{proof}

Let us give an example.
\Es\label{horikawa}
Suppose $\gk=\C$.
With the notations of $\thref{pullback}$, let $X_1\subset \p^3$ to be a smooth quadric, $X_2$ to be a smooth $K3$ surface of degree $2$ and $f_i$ the two associated degree $2$ covers of $\p^2$.
If $B\subset \p^2$ is the branch locus of $f_2$ then $p_1:X_3\ra X_1$ is a degree $2$ covering branched along the divisor $f_1^*B$.
Since by adjunction we have $\on{2}(B)\cong \on{2}(6)$, it follows that $X_3$ is a special instance of the surfaces studied by Horikawa in \cite{Horikawa1}[Thm. 1.6 iii)]\footnote{with his notation, we are considering the case $d=0$ and $n=3$}.

We always consider Ulrich bundles respect to the polarisation given by pulling back $\on{2}(1)$. 
We know $Pic(X_1)\cong \Z^2$ and it is generated by the classes of two intersecting lines; it is easily verified that they are Ulrich line bundles.
On $X_2$ there are stable Ulrich bundles $\f$ of rank $2$, it can be shown as for the other $K3$ surfaces in \cite{FaenziK3} or it follows from \cite{KuNaPa}.
We deduce that there are rank $2$ Ulrich bundles on $X_3$.
Up to our knowledge this is the only Horikawa surface, except for double coverings of $\p^2$, where existence of Ulrich bundles is known.
Actually, we will see in \thref{horikawa'} that the general such $X_3$ admit Ulrich sheaves.
\Io

\subsection{Restrictions and extension of Ulrich sheaves}

An immediate application of the theorem above is the generalization of \cite{CH}[Lem. 2.4] and \cite{AC}[Cor. 4.5].

\Co\label{restrizione}
Let $\e$ be an Ulrich sheaf on an $n$-dimensional polarised scheme $(X,H)$ over an infinite field $\gk$.
Suppose $n\geq 1$.
For any $D\in |\os_X(1)|$ define $\e_D:=\e\otimes\os_D$ and $\os_D(1):=\os_X(1)|_D$, then $\e_D$ is Ulrich on $(D,\os_D(1))$.
The same applies to any complete intersection of divisors in $|\os_X(1)|$
\io
\begin{proof}
    In \thref{pullback}, we choose $X_1\cong \p^{N-1}$ such that $f^*X_1=D$ with the Ulrich sheaf $\os_{\p^{N-1}}$.
    If we set $X_2=X$ with $\e$ as Ulrich sheaf then we immediately get that $\e|_D$ is Ulrich. 
    By arguing inductively we conclude the second assertion.
\end{proof}

\Oss
Actually, by applying \thref{intcomplete}, we can get a better result: any complete intersections in $X$ of divisors $D_1,\dots, D_l$ which are pullback of divisors in $\p^n$ admits Ulrich sheaves.
However, in principle we cannot conclude the same for all complete intersections of divisors linearly equivalent to some multiple of $H$.
\one

A better analysis of the previous setup will allow us to investigate the zero loci of sections of Ulrich sheaves; this will be especially interesting when the polarisation chosen has low degree.

\Prop\label{restrizionesezioni}
Let $\e$ be an Ulrich sheaf on an $n$-dimensional polarised scheme $(X,H)$ over an infinite field $\gk$ with $n\geq 1$.
\begin{itemize}
    \item The restriction of global sections gives an isomorphism $H^0(X,\e)\cong H^0(D,\e|_D)$ for any $D\in |\os_X(1)|$.
    \item If $Y$ is the zero locus of some non-zero section of $\e$, then $Y$ cannot contain any complete intersection of divisors in $|\os_X(1)|$.
    \item If $f:X\ra\p^n$ is a finite morphism such that $f^*\osn(1)\cong \os_X(H)$ then $Y$ is not of the form $f^*Z$ for some $Z\subset\p^n$.
\end{itemize}
\One
\begin{proof}
By definition of $D$, we have the exact sequence
\[0\ra\os_X(-1)\ra\os_X\ra\os_D\ra 0\]
which tensored with $\e$ stays exact at least on the right.
Calling $\ck$ the kernel of the resulting left morphism, we obtain
\[0\ra\ck\ra\e(-1)\ra\e\ra\e_D\ra 0.\]
By \thref{definizioneUlrich} $v)$ we have $f_*\e\cong \osn^\rho$, so that, by projection formula, we have $f_*(\e(-1))\cong f_*(\e)\otimes\osn(-1)\cong \osn(-1)^\rho$.
Applying $f_*$, which is exact since $f$ is finite, we get 
\[0\ra f_*\ck\ra\osn(-1)^\rho\ra\osn^\rho\ra f_*(\e_D)\ra 0.\]
Since $f_*(\e_D)$ is supported on $f(D)$, by the additivity of rank we get that $f_*\ck$ is a torsion sheaf so, being $\osn(-1)$ locally free, it must be $0$.
By \thref{counitsurj} we have a surjective morphism $0=f^*f_*\ck\ra\ck$ hence $\ck=0$ and we have an exact sequence
\ses{\e(-1)}{\e}{\e|_D}
In particular, since $h^0(X,\e(-1))=0=h^1(X,\e(-1))$ we get an isomorphism between the global sections of $\e$ and $\e|_D$.

    Applying iteratively this fact, we prove that for any complete intersection $W$ of divisors in $|\os_X(1)|$ the restriction of global sections gives an isomorphism $H^0(X,\e)\cong H^0(W,\e|_W)$.
    In particular, no section of $\e$ can vanish on $W$, that is $W$ is not contained in $Y$.
    For the last claim, it suffices to note that a scheme of the form $f^*Z$ contains all the fibres over the points in $Z$.
    But each such fibre is a complete intersection of divisors in $|\os_X(1)|$, since a point in $\p^n$ is a complete intersection of hyperplanes.
\end{proof}

Since no zero locus $Y$ of an Ulrich sheaf can be of the form $f^*Z$ this proposition suggests, but does not prove, that maybe the cycle $[Y]$ is not of the form $f^*[Z]$.
Considering the case of $c_1$, if $Pic(X)$ is cyclic then the first Chern class is forced to be a multiple of $H$, but we have seen that in this case we have no Ulrich line bundles, see \thref{pic1'}.
This raises the following question.

\Do
Which Chern class $c_k(\e)$ of an Ulrich bundle can be multiples of $H^{n-k}$?
\da

Contrary to the rank $1$ case, we have plenty of example of higher rank Ulrich bundles having $c_1$ proportional to $H$, this remarkable situation will be analysed in \thref{special}.

Here we prove a natural follow up to \thref{restrizione}, this is essentially a consequence of the fact that a vector bundle on $\p^n$ ($n\ge 3$) is split if and only if its restriction to an hyperplane is such.
\Prop\label{estensione}
Let $\e$ be a coherent sheaf on an $n$-dimensional polarised scheme $(X,H)$ over an infinite field $\gk$.
Suppose $n\geq 3$.
If $\e$ is mCM 
and there is $D\in |\os_X(1)|$ such that $\e_D$ is Ulrich then $\e$ is Ulrich.
\One
\begin{proof}
Consider a finite surjective morphism $f:X\ra\p^n$ such that $f^*\osn(1)\cong\os_X(H)$.
Call $D':=f(D)\cong \p^{n-1}$.
By \thref{definizioneUlrich} $v)$ to prove that $\e$ is Ulrich is enough to prove that $f_*\e$ is a direct sum of trivial bundles.
By the same result, we know that $f_*(\e|_D)\cong \on{n-1}^\rho$ since $\e|_D$ is Ulrich.
Note that, being $f$ affine, by base change, see \cite{GW2}[Lem. 22.88], we have $(f_*\e)|_{D'}\cong f_*(\e|_D)\cong \on{n-1}^\rho$.
Being $\e$ mCM we deduce, using both parts of \thref{lcmlf}, that $f_*\e$ is locally free.
Now we can apply Horrocks' splitting criterion in the form \cite{OSS}[Thm. 2.3.2] to $f_*\e$ and conclude that it splits.
Since its restriction to $D'$ is trivial, $f_*\e$ must be itself trivial.
\end{proof}

\Oss
In the converse statement we cannot let $dim(X)=2$, as it happens for Horrock's splitting.
One counterexample is as follows.
Consider a point $P\in\p^2$ then from the exact sequence
\sesl{\id_P}{\os_{\p^2}}{\os_P}{p} we get 
\[Ext^1(\id_P,\os_{\p^2})\cong Ext^2(\os_P,\os_{\p^2})\cong H^0(\p^2,\os_P)=\mathbf{k}\]
by Serre duality.
We have only one non-split extension of the form
\sesl{\os_{\p^2}}{\e}{\id_P}{e'}
$\e$ is not Ulrich, actually not even ACM, being 
\[h^1(\p^2,\e(-1))=h^1(\p^2,\id_P(-1))=h^0(X,\os_P(-1))=\mathbf{k}\]
but restricting $\e$ to a line $l$ not containing $P$ it will become $\os_{l}^2$ hence Ulrich.
So, if we prove that $\e$ is locally CM, that in this case is equivalent to locally free, we have our example.
From \eqref{p} we get 
\[\ext^1(\id_P,\os_{\p^2})\cong \ext^2(\os_P,\os_{\p^2})\cong \os_P\]
and $\ext^2(\id_P,\os_{\p^2})=0$.
By \eqref{e'} we have
$\ext^2(\e,\os_{\p^2})\cong \ext^2(\id_P,\os_{\p^2})=0$
and 
\[0\ra \os_{\p^2}\ra E\du\ra \os_{\p^2}\ra \os_P \ra \ext^1(\e,\os_{\p^2})\ra 0\]
The induced map $\os_{\p^2}\ra \os_P$ is non-zero on global sections, because the extension is not split.
Therefore, it must be surjective as a map of sheaves, implying that $\ext^1(\e,\os_{\p^2})=0$ hence $\e$ is locally free by \thref{extpd}.
\one

\subsection{Other modifications}

The following has already appeared (for $H$ very ample) in \cite{EisSch}[Thm. 2.6].
\Prop\label{prodotti}
Let $(X_i,H_i)$ be two polarised varieties of dimensions $n_i$ over a field $\gk$ with Ulrich sheaves $\e_i$, with $i=1,2$.
If we call $\pi_i$ the two projections $\pi_i:X_1\times X_2\ra X_i$ then $\pi_1^*\e_1\otimes\pi_2^*\e_2(n_1)$ and $\pi_1^*\e_1(n_2)\otimes\pi_2^*\e_2$ are Ulrich on $(X_1\times X_2,\pi_1^*H_1\otimes\pi_2^*H_2)$.
\One

The following is proved as \cite{EisSch}[Thm. 5.4].
\Prop\label{polarizzazionimultiple}
Let $f:X\ra\p^n$ be a finite surjective morphism and $\e$ a $f$-Ulrich sheaf.
For any sheaf $\f$ that is Ulrich on $(\p^n,\osn(l))$ it holds that $\e\otimes f^*\f$ is Ulrich on $(X,f^*\osn(l))$.
\One

This is analogous to \cite{Kim}[Thm. 0.1]
\Prop\label{blowup}
Consider a smooth polarised variety $(X,H)$ and the blow-up $\pi:\Tilde{X}\ra X$ in a point.
Suppose that $\pi^*H-E$ is still ample and globally generated, where $E$ is the exceptional divisor.
If $\e$ is an Ulrich sheaf on $(X,H)$ then $\pi^*\e(-E)$ is Ulrich on $(\Tilde{X},\pi^*H-E)$.
\One

Finally, the projective bundle construction can be treated as in \cite{Hochenegger}[Prop. 3.1], see \cite{Hochenegger}[Thm. A and B] for applications to curves and surfaces.
\Prop\label{projbundle}
Let $(X,H)$ be a smooth, $n$-dimensional polarised variety and $\pi: \p(E) \ra X$, where $E$ is locally free of rank $r+1$.
Suppose $H_\p:= \pi^* H + \xi$ is ample and globally generated, with $\xi$ the relative hyperplane section, and $\e$ is a locally free sheaf on $X$ with $H^j(X,\e) = 0$ for all $j$.
Then $\pi^* \e(H_\p)$ is Ulrich if and only if
\[ Ext^j\bigl(Sym^k E, \e(-c_1(E)-(r+1+k)H)\bigr) = 0
\qquad 0\le k\le n-2,\; j\ge 0.\]
\One

\section{Existence results}

We expose the main existence results of Ulrich sheaves without giving proofs.
We selected a restricted number of statements, which are by no means exhaustive, from the vast existing literature; for some other references, see the book \cite{CMRPL}.

\subsection{Curves}

We call a curve a $1$-dimensional connected scheme over $\gk$.
The existence of Ulrich sheaves on curves is well understood, and actually we will see that the strongest possible statement in \thref{congetturaesistenza} holds.

Let us start with an historical remark.
\Oss
Any irreducible curve is birational to a plane curve, but we will see in \thref{intcomplete1} that existence of Ulrich sheaves for divisors in projective space is equivalent to existence of linear determinantal representations for some power of their equations.
It seems that the first to investigate this topic was Hesse in \cite{Hesse}.
It was proven in \cite{Dixon} that general plane curves admit determinantal representations and in \cite{Dickson} that any such curve satisfies this property.
For more on this topic, see \cite{Dolgachev}[Historical remark p.205].
\one

The existence of Ulrich sheaves on curves has been proved in \cite{EisSch}.
By definition, a sheaf $\e$ on a polarised curve $(X,H)$ is Ulrich if and only if $h^j(X,\e(-1))=0$ for $j=0,1$.
In other words, Ulrich sheaves respect to $H$ are exactly sheaves with no cohomology twisted by $H$.
This implies that the theory of Ulrich sheaves on curves does not depend on the polarisation.
We summarise part of the content of \cite{EisSch}[Thm. 4.3, Thm. 4.4 and Cor. 4.5] in the following.
\Te\label{curve}
Fix a polarised curve $(X,H)$ over some field $\gk$.
\begin{itemize}
\item All Ulrich sheaves are of the form $\cl(1)$ with $h^0(C,\cl)=0=h^1(C,\cl)$.
\item If $\gk$ is algebraically closed then any $(X,H)$ admits Ulrich sheaves of rank $1$.  
\item for arbitrary $\gk$ there always are Ulrich sheaves on a $(X,H)$.
\end{itemize}
\Ma

\Oss
Over an algebraically closed field, if $X$ is smooth then line bundles $\cl$ such that $h^j(X,\cl)=0$ for $j=0,1$ must have degree $g(X)-1$ by Riemann--Roch.
For $g=0$ we get only the line bundle $\on{1}(-1)$ and for $g=1$ it is $\os_X$.
While for $g\ge 2$ they form a divisor in $Pic(X)^{g-1}$ called \textit{theta-divisor}.
\one

In \cite{Cos}[4.1.1, 4.1.2 and Thm. 4.1] the next fact is proven.
\Prop
Suppose $(X,H)$ is a smooth polarised complex curve  of genus $g$, then
\begin{itemize}
    \item $g=0$ implies all Ulrich bundles are direct sums of Ulrich line bundles
    \item $g=1$ implies that for any rank $r\ge2$ there is just a $1$-dimensional family of stable Ulrich bundle of rank $r$
    \item $g\ge2$ there are stable Ulrich bundles of any rank and their moduli space has dimension $r^2(g-1)+1$
\end{itemize}
\One

Those three behaviours are usually called, respectively, \textit{finite},\textit{tame} and \textit{wild representation type} after Gabriel's theorem, see \cite{FaenziPons_acm}[§1]. 
 
\subsection{Surfaces}

The theory of Ulrich sheaves on surfaces becomes much more interesting.
Note that in \cite{Dickson} it is proved that if a general form of degree $d>2$ in $\p^n$ is linear determinantal then either $n=2$ or $n=3$ and $d=3$.
In modern language, see the upcoming \thref{intcomplete1}, this is saying that all the hypersurfaces of degree $d$ in $\p^n$ with $n\ge 3$ admitting Ulrich line bundles are quadrics and cubic surfaces, therefore already for surfaces we need to search for higher rank ones.

Many examples have been worked out case by case and we also know the behaviour of the moduli spaces for non-positive Kodaira dimension, but still there are many open questions and is a very active area of research.
For example, at least over $\C$, we know that the weak formulation of \thref{congetturaesistenza} holds but we do not know if the stronger ones hold in general or not.
Actually, to our knowledge, the only tentative counterexample to the standard version of \thref{congetturaesistenza} is a surface with special numerical properties, see \cite{Beauville_controes}[p. 18-20], whose existence is still open.

We start by analysing some cases depending on the Kodaira dimension of our surface.
A detailed account is in \cite{CMRPL}[Chapter 5] here we give just a rough sketch.
Note that, due to \thref{blowup}, we can reduce the existence problem to minimal surfaces, sorted by Kodaira dimension.

For $(\p^2,\on{2}(m))$ we know exactly for which values of $m,r$ there are rank $r$ Ulrich bundles.
For a general discussion on this and also blow-ups of $\p^2$ see the next \autoref{veronese}.
In the case of complex Hirzebruch surfaces a complete classification of the possible Ulrich bundles and a proof of their existence has been achieved with two different methods in \cite{Antonelli}[Thm. 1.1 and Thm. 1.3] and \cite{CoskunHuizenga}[Thm. 4.6].
On general ruled surfaces we know existence of Ulrich bundles and stable rank $2$ bundles by \cite{ACM} and \cite{ACCMRM}.
This completes the case of Kodaira dimension $-\infty$.

The case of complex $K3$ surfaces has been studied in \cite{AFO}[Thm. 0.4] and \cite{FaenziK3}[Thm. 1, Thm. 2].
The following complete classification is in \cite{CoskunNuerYoshioka}[Thm. 4.4].
\Te
If $(X,mH)$ is a smooth polarised $K3$ surface with $Pic(X)\cong \Z H$ over $\gk=\overline{\gk}$ then there are stable rank $r$ Ulrich bundles as soon as $2|rm$ and their moduli space has dimension $r^2(\dfrac{m^2H^2}{4}+2)+2$.
\Ma

The existence of rank $2$ Ulrich bundles on abelian surfaces has been proved in \cite{Beabelian}.
More generally, in \cite{Bea}[Prop. 5.1] there is a uniform proof for surfaces of Kodaira dimension $0$.
The case of Enriques surfaces has been studied also in \cite{BN}, where it is proven the following.
\Te
On an unnodal Enriques surface there is a polarisation (the \textit{Fano polarisation}) admitting stable Ulrich bundles of any rank $r\ge 1$.
Moreover, their moduli space, when non-empty, has dimension $c_1^2-19r^2+1$, where $c_1$ is their  first Chern class. 
\Ma

\Es
For higher Kodaira dimension we just list some examples of surfaces that have been studied in the literature:
\begin{itemize}
    \item elliptic surfaces which are Weierstrass fibrations and with $q=0$, \cite{MRPL_elliptic}
    \item non-special polarised surfaces with $q=p_g=0$, \cite{Casnati_pgq0}
    \item non-special polarised surfaces with $q=1$ and $p_g=0$, \cite{Casnati_pg0_q1}
    \item non-special polarised surfaces with $q=0$ and $p_g\ge 1$, \cite{Casnati_reg_noneg}
    \item surfaces with $q\le 1$ or $q\ge 1$ and minimal model with cyclic Picard group, \cite{Lopez_surface}
    \item surfaces with Picard group generated by the very ample $K_X$ and such that $q<\chi(\os_X)-1$, \cite{BasuPal}.
\end{itemize}
We refer to the corresponding papers for precise statements.
\Io

Given a polarised surface $(X,H)$ in general it does not admit rank $2$ Ulrich bundles, even for $H$ vary ample.
For example, in \cite{Beadeterminantal}[Prop. 7.6] it is shown that a general surface in $\p^3$ of degree $d$ is Pfaffian, hence admits rank $2$ Ulrich bundles, if and only if $d\le 15$, . 
Nevertheless, the following remarkable asymptotic result holds, see \cite{CoskunHuizengaNuer}[Prop. 6.9], see also \cite{CoskunHuizenga}[Thm. 4.3]

\Te
Take a smooth, complex, polarised surface $(X,H)$. There exists a positive integer $m_0$ such that for all $m \geq m_0$, the polarised surface $(X,mH)$ admits an Ulrich bundle of every positive even rank. Moreover, if $K_X$ (respectively, $K_X + H$) is divisible by $2$ in $Pic(X)$ and $2m \geq m_0$ (respectively, $2m+1 \geq m_0$), then $(X, 2mH)$ (respectively, $(X, (2m+1)H)$) admits an Ulrich bundle of every rank $r \geq 2$.
\Ma

\Oss
Note that an analogous result does not hold for Ulrich line bundles, hence, even up to taking multiples of a given polarisation, in general the Chow form of a surface is not a determinant of linear forms.
One simple reason is that having cyclic Picard group (generated by an effective divisor) is an obstruction even when we are allowed to change polarisation, see \thref{pic1}.

Even worse, given a polarisation $H$ admitting Ulrich line bundles it could happen that its multiples do not.
One example is clearly $(\p^n,\osn(1))$ for $n\ge 2$.
More interestingly, this happens for ruled surfaces as well, see \cite{ACM}[Thm. 2.1].
However, note that for $X=\p^1\times\p^1$ the line bundles $\os_X(n,n-1)$ and $\os_X(n-1,n)$ are Ulrich on $(X,\os_X(n,n))$ for all $n\ge 1$.
\one

\subsection{Veronese varieties and blow-ups of $\p^n$}\label{veronese}

Let us start with an Ulrich bundle on $(X,H)$.
We have seen in \thref{polarizzazionimultiple} that in order to obtain Ulrich bundles on $(X,mH)$ for any $m>1$ is enough to construct Ulrich bundles on the $m$-th Veronese embedding of $\p^n$.
This is settled by the two different constructions given in \cite{EisSch}[Cor. 5.7] and \cite{Bea}[Thm. 3.1], we deduce the following.
\Te\label{Veronese}
There exists Ulrich bundles of rank $n!$ on $(\p^n,\osn(m))$ for any $n,m$.
Moreover, there are equivariant, for the action of $GL(n+1)$, ones of rank $d^{\binom{n}{2}}$.
\Ma
A precise classification of the triples $(n,m,r)$ for which there are rank $r$ Ulrich bundles on $(\p^n,\osn(m))$ is available only for $n=2$ in \cite{CoskunGenc}, \cite{CostaMiroroig2}[Thm. 1] and for $n=3$ in \cite{FaenziPretti}[Thm. 1].
In higher dimension the only case in which all the possible triples $(n,m,r)$ are completely understood is when $n!|m$, see \cite{FaenziPretti}[Lem. 4].
In general, the situation is definitely more complicated; see \cite{FaenziPretti}[Remark 3] and \cite{LopezRayveronese}[Thm. 1].

This result has a nice interaction with \thref{blowup}, since to get ample polarisations on blow-ups we usually need to start with much more positive ones.
It was shown already in \cite{EisSch}[Thm. 6.6] that, under the assumption $\gk$ infinite, there are polarisations on blow-ups of $\p^2$ admitting rank $2$ Ulrich sheaves, see also \cite{AproduCostaMiroroig}[Cor. 3.8].
Among those we have del Pezzo surfaces on which Ulrich line bundles were studied in \cite{PL_T}, \cite{CKM_delpezzo}[Prop. 2.19].
Even better is the following result, proved in \cite{CFK2}[Thm. 2.1 and Thm. 4.8], by directly arguing on those blown-up surfaces.
\Te
Suppose $\gk=\C$.
On the blow-up of $\p^2$ in any set of sufficiently general points there are polarizations admitting Ulrich bundles of any rank $r\ge 1$ and if we increase $r$ we get arbitrary large families of non-isomorphic sheaves.
\Ma
In the case of del Pezzo surfaces this has been proved in \cite{MRPL}[Thm. 4.10].
\Te
Suppose $\gk=\C$.
For any $r\ge 2$ there are rank $r$ simple, special Ulrich bundles on del Pezzo surfaces with respect to the ample (even if non-globally generated) $-K_X$.
They are stable with respect to some other polarisation and their moduli space has dimension $r^2+1$.
\Ma
See also \cite{CoskunHuizenga}[Thm. 4.10].

We will give an outline of what is known on Ulrich bundles on blow-ups of $\p^n$.
We will always suppose $n\geq 3$ to be the dimension of our variety $X$.
Let $\pi:X\ra\p^n$ be the blow-up in $k$ general points $P_1,\dots, P_k$.
Call $E_i:=\pi\inv(P_i)$ the exceptional divisor corresponding to $P_i$ and set $E:=\sum_{i=1}^kE_i$.
Moreover, call $L:=c_1(\pi^*\on{3}(1))$.
$Pic(X)$ is freely generated by $L$ and the $E_i$-s.
Finally, set $H_m:=mL-E$.
Fixed any $m\geq 3$, if the points $P_1,\dots P_k$ are general enough then the divisor $H_m$ is very ample for $0\leq k\leq \binom{m+n}{m}-2n-2$ by \cite{Coppens}.

The choice of $H_m$ is particularly well behaved for the study of Ulrich bundles since from an iterated application of \cite{Kim}[Thm. 0.1] we deduce the following.

\Le
If $\e$ is an Ulrich bundle on $(\p^n,\osn(m))$ then $\Tilde{\e}:=\pi^*\e(-E)$ is Ulrich on $(X,H_m)$.
\ma

Together with the existence of Ulrich bundles on $(\p^n,\osn(m))$ for any $n,m$, we deduce the following.

\Te\label{blp3}
For every $k,n\in \N$ there exists $m\geq 3$ such that the blow-up of $\p^n$ in $k$ general points has an Ulrich bundle respect to the very ample polarization $H_m$.
\Ma

Not all Ulrich bundles on $X$ are obtained as pullback through $\pi$: those of this form are characterised by their restriction to the exceptional locus $E$, see \cite{CasnatiKim}[Thm.. 4.2].

\subsection{Complete intersection}

Existence of Ulrich sheaves, or even aCM ones, on hypersurfaces in projective space have been investigated in \cite{Beadeterminantal} through the connection between sheaves on hypersurfaces and linear determinantal representations of their equation.
A consequence of \cite{Beadeterminantal}[Thm. A] is the following.
\Te\label{intcomplete1}
Let $X\subset \p^n$ be an integral hypersurfaces of degree $d$ and equation $p=0$.
A coherent sheaf $\e$ of rank $r$ on $X$ is an Ulrich sheaf on $(X,\os_X(1))$ if and only if it fits in a sequence like
\[0\ra\osn(-1)^{rd}\xra{A}\osn^{rd}\ra\e\ra 0\]
where $det(A)=p^r$.
\Ma

But actually, many works on the relating matrix factorisations and modules over rings that are quotient of regular rings by principal ideals predate this paper, see for example \cite{Eis-resolutions}, \cite{BHU}, \cite{Brenna_Herzog_Ulrich}.
More generally, there are Ulrich sheaves on all complete intersections in projective space.
We will give a proof of this fact in \thref{intcomplete} and the subsequent corollaries, essentially the same given in \cite{BHU} or \cite{CMRPL}.
For a different approach to construct sheaves on divisors using invariant theory, see \cite{Manivel_ulrich}.
Here we just discuss briefly the case of quadrics and cubics.

Ulrich bundles on quadrics actually have a long history: they were already studied in connection with Clifford algebras in \cite{BEH} and known as spinor bundles in \cite{Ottaviani_spinor}.
Moreover, also come out in order to construct a semi-orthogonal decomposition of the bounded derived category of quadrics, see \cite{Kapranov}.
A proof of the following is in \cite{Bea}[Prop. 2.5].
\Te
For a smooth quadric $Q\subset \p^{n+1}$ with $n$ even (resp. odd) there are exactly two (resp. one) indecomposable Ulrich bundle of rank equal to $2^{\left\lfloor \frac{n+1}{2}\right\rfloor}$
\Ma
In \cite{EisSchci} there is a further study of complete intersections of quadrics.

Already for cubic hypersurfaces the situation is poorly understood in dimension $\ge 5$, as we do not have a precise result for the minimal rank of Ulrich bundles.
On smooth cubic surfaces there always are Ulrich line bundles, this has been proved classically in the language of determinantal representations.
On smooth, complex cubic $3$-folds we know existence of rank $2$ Ulrich bundles since \cite{Druel}, where it is also proven that their moduli space is birational to the intermediate Jacobian of those varieties.
Actually there are Ulrich bundles of any rank $\ge 2$, see \cite{LMS} and \cite{CFK3}.

In the case of $4$-folds we know that the general such surface admits only Ulrich bundles of rank divisible by $6$, \cite{FaeKim}.
Nevertheless, there are divisors in the moduli space of such varieties whose general member parametrises $4$-dimensional cubics with Ulrich bundles of rank $2$ or $3$, see \cite{TN}.
We remark that having a rank $2$ Ulrich bundle is equivalent to giving a Pfaffian representation for the equation of this variety but Pfaffian cubic fourfolds are known to be rational.

\subsection{Grassmannians}

For homogeneous varieties it makes sense to search for Ulrich bundles which are also invariant under the action of the acting group.
For Grassmann varieties, a complete classification in terms of Schur functors has been obtained in \cite{CMR_Grassmann}[Thm. 3.6].
The following is a corollary of their more general statement.
\Te\label{Gr}
Given a complex vector space $V$, for any $1\leq k\leq n-1$ there are Ulrich vector bundles on the Grassmannian $Gr(k,V)$ invariant for the action of $GL(V)$.
\Ma

\Es
Using this and \thref{intersezione} we can construct an Ulrich bundle on the complete intersection of two copies of the Grassmannian $Gr(2,5)$ suitably embedded in $\p^9$, see \cite{Casnati_wild}[Example 4.1].
This variety, usually called $GPK^3$ threefold is Calabi-Yau.
This is one of the few Calabi-Yau-s of dimension $\ge 3$ known to admit Ulrich bundles, except for complete intersections in $\p^n$ and products of other varieties.
Up to our knowledge, the only other one is the \textit{octic double solid}, see \thref{ulrichr2doublecoverexistence}.
\Io

\subsection{Del Pezzo varieties}

Recall that a \textit{Fano variety} is a smooth projective variety $X$ such that $-K_X$ is an ample divisor.
The divisibility of $-K_X$ in $Pic(X)$, that is the maximum integer $i$ such that $\dfrac{-K_X}{i}$ is still a divisor, is called \textit{index}.

The work on Ulrich bundles on Del Pezzo varieties actually pre-dates their sheaf theoretic definition.
Already in 2000 in \cite{ArrondoCosta} the rank $2$ case on degree $3,4,5$ del Pezzo $3$-folds were considered, since they were particular aCM bundles.
Again as a special case of a more general notion, this time that of an instanton bundle, they appeared in \cite{Faenzi11} and \cite{Kuznetsov_instanton}.

The following is essentially a collection of results already in the literature.
\Te\label{Fano3i2}
Let $X$ be a smooth complex Del Pezzo variety, i.e. a Fano variety of dimension $n$ and index $n-1$.
If $H:=\dfrac{-K_X}{n-1}$ is globally generated then there are Ulrich bundles on $(X,H)$.
\Ma
\begin{proof}
    Such varieties $X$ are classified.
    If $n=2$ then $X$ is one of the classical \textit{Del Pezzo surfaces}, which are either blow-ups of $\p^2$ in $k$ general points, with $0\leq k\leq 8$, or $\p^1\times\p^1$.

    If $n\geq 3$ their classification is essentially due to T. Fujita and V. Iskovskikh, see \cite{IskPro}[Thm. 3.3.1].
    We call $H^n=d$ the degree of $X$.
    Note that we must have $d\leq 8$.
    The condition $H$ globally generated rules out only the case in which $d=1$, hence $X$ is a degree $6$ hypersurface in a weighted projective space $\p(1^n,2,3)$, see \cite{IskPro}[Thm. 3.2.4].
    
    For $d=2$ we find that $X$ is a double cover of $\p^n$ and hence the existence of Ulrich bundles was proved in \cite{KuNaPa}.
    It also follows from \thref{ulrichmatrix''}, moreover see \thref{ulrichr2doublecoverexistence} and \thref{qdswild} for the possible ranks.
    For $d=3,4$ those are cubic hypersurfaces or complete intersection of two quadrics, hence existence of Ulrich bundles was already known; see \thref{intcomplete} for a proof.
    
    The other possible varieties have a bounded dimension.
    If $d=5$ then $X$ is a linear section of $6$-dimensional Grassmannian $(Gr(2,5),\os_{Gr(2,5)}(1))$, where $\os_{Gr(2,5)}(1)$ gives the Plucker embedding.
    This means that $X$ is a complete intersection of $l$ divisors in $|\os_{Gr(2,5)}(1)|$, with $0\leq l\leq 3$, and $\os_X(H)\cong \os_{Gr(2,5)}(1)|_X$.
    Therefore, by \thref{restrizione}, it is enough to know the existence of Ulrich bundles in $(Gr(2,5),\os_{Gr(2,5)}(1))$ to conclude, but this has already been addressed in \thref{Gr}.

    If $d\geq 6$ then $X$ is either $\p^2\times\p^2$ or has dimension $3$.
    For the first one see \thref{prodotti}.
    For Fano $3$-folds of index $2$ we know a lot about Ulrich bundle, see \cite{Bea} or \cite{CFK3}.
\end{proof}

\subsection{Other $3$-folds}

The only Fano $3$-folds left out from the previous results are the ones having index $1$.
All those varieties are classified, see \cite{IskPro}.
If their Picard groups is cyclic and generated by a very ample line bundle then the existence of rank $2$ Ulrich bundles follows from the proof of Proposition 3.15 in \cite{BF}, even though it was not explicitly stated.
More generally, existence of smooth components of the expected dimension in the moduli spaces of any even rank Ulrich bundles on those varieties has been proven in \cite{CFK3.1}[Thm. 1.3].
If we remove the very ampleness condition on the generator of the Picard group we get two new families of varieties, one of them are the sextic double solids which we will treat in \thref{rango>sextic}.
When the Picard group is not cyclic the more interesting cases are the ones which are not products.
We know by \cite{Genc} that seldom we have Ulrich line bundles and it seems that the other case have not been yet attacked.

\smallskip

Another class of varieties which have received much attention is that of scrolls.
We already have seen some results in case of surfaces.
For the $3$-fold case we refer to \cite{FaniaLelliPL}, \cite{FaniaFlamini1} and \cite{FaniaFlamini2}.

\section{Stability of Ulrich sheaves}

In this section, we will consider a polarised variety $(X,H)$ over an infinite field $\gk$.
Under such assumption $\e$ will be torsion-free by \thref{corollario}.
Stability of sheaves will always be respect to $H$.
We will call $d$ the degree of $H$, that is $H^n$.
The goal of the upcoming few pages is to present the proof that any Ulrich sheaf on a variety is semistable.
This will be essential to construct moduli spaces of Ulrich sheaves.
Before doing this, we compute the Hilbert polynomial of an Ulrich sheaf.

\subsection{Hilbert polynomial}

\Prop\label{hilbertpol}
Let $(X,H)$ be an $n$-dimensional polarised variety over an infinite field $\gk$.
Set $d=H^n$.
If $\e$ is an Ulrich sheaf of rank $r$ then its Hilbert polynomial is $P(\e)(t)=rd\binom{t+n}{t}$.
In particular, $h^0(X,\e)=rd$.
\One
\begin{proof}
    By \thref{definizioneUlrich} $v)$, we know that for a finite morphism $f:X\ra\p^n$ such that $f^*\osn(1)\cong \os_X(1)$ we have $f_*\e\cong \osn^\rho$.
    Since finite morphisms preserve cohomology and using projection formula we easily deduce that $P(\e)=P(f_*\e)=\rho\cdot P(\osn)$.
    Note that $\rho=rk(f_*\e)=d\cdot rk(\e)=dr$ since $f$ is of degree $d$.
    Moreover, the Hilbert polynomial of $\osn$ is $\binom{t+n}{t}$, hence we get the claim.

    The last part follows from $\chi(\e)=p(\e)(0)=rd$ and the fact that $h^j(X,\e)=0$ for $j>0$.
\end{proof}

\Co\label{rhps}
The reduced Hilbert polynomial of an Ulrich sheaf depends only on $n$ and $d=H^n$.
In particular, fixed $n,d$ all Ulrich sheaves have the same slope.
\io
\begin{proof}
    The first claim is clear being $p(\e):=P(\e)/rk(\e)=d\binom{t+n}{n}$.
    For the second, recall that the slope is actually encoded by the reduced Hilbert polynomial, see \cite{HuyLeh}[Def. 1.2.2, 1.2.3 and 1.2.11].
\end{proof}

\Co
The family of Ulrich sheaves of fixed rank over $X$ is bounded. 
\io
\begin{proof}
    We immediately deduce the thesis from \cite{HuyLeh}[Lem. 1.7.9], recalling that Ulrich sheaves are globally generated, see \thref{corollario}.
\end{proof}

\subsection{Semistability}

The following is the main result of this section.
It already appeared, for example, in \cite{CH}[Thm. 2.9] and \cite{CKM_clifford}[Thm. 2.11, Lem. 2.15], under additional assumptions.

\Prop\label{stability}
Let $(X,H)$ be an $n$-dimensional polarised variety over an infinite field.
\begin{enumerate}[1)]
    \item An Ulrich sheaf is semistable.
    \item If $\e$ is Ulrich and we have an exact sequence of coherent sheaves 
    \[0\ra \f\ra \e \ra \g\ra 0\]
     with $\mu(\f)=\mu(\e)$ and $\g$ torsion-free then also $\f,\g$ are Ulrich. 
     \item Each factor of the Jordan-Holder filtration of an Ulrich sheaf $\e$ must be a stable Ulrich sheaf. 
    \item If $\e$ is stable then it is also slope stable. 
\end{enumerate}
\One
\begin{proof}
\gel{1}
Consider a finite surjective map $f:X\ra\p^n$ with $f^*\os(1)=\os_X(H)$.
        Let $\f\subseteq \e$ be a subsheaf then, being $f_*$ left exact, we have $f_*\f\subseteq f_*\e=\os_{\p^n}^\rho$.
        Since $f$ is finite it preserves cohomology and since has degree $d=H^n$ multiplies ranks by $d$.
        Summing up, for any sheaf $\f$ on $X$ we have $P(f_*\f)=P(\f)$ and $rk(f_*\f)=d\cdot rk(\f)$ therefore, $d\cdot p(f_*\f)=p(\f)$.
        So, we need only to prove $p(f_*\f)\leq p(\osn^\rho)$, but this is clear since $\os_{\p^n}^\rho$ is semi-stable.
        
       \gel{2}
        Since $f$ is finite, $f_*$ is exact and we get
        \[0\ra f_*\f\ra \osn^\rho \ra f_*\g\ra 0\]
        where, by the discussion in $1)$ and the fact that the degree of a sheaf can be read from its Hilbert polynomial, we have $\mu(f_*\f)=\mu(f_*\e)=0$ so that also $\mu(\g)=0$.
        In addition, $f_*\g$ is torsion free by \cite{EGA1}[Chapter 1 Thm. 7.4.5] being $f$ dominant.
        To conclude we argue that any torsion-free quotient of $\osn^\rho$ with zero slope must be trivial, this is a standard result, see for example \cite{OSS}[Cor. p.28], we will give a proof in \thref{trivialquot} below.
        Indeed, this proves that $f_*\g$ is trivial hence also $f_*\f$ is and then both $\e,\g$ are Ulrich by the last claim in \thref{definizioneUlrich}.
    
\gel{3}
Consider a Jordan-Holder filtration $0=\e_0\subset \dots \subset \e_l\subset \e_{l+1}=\e$, so that $\e_i/\e_{i-1}$ is stable and has the same reduced Hilbert polynomial as $\e$.
    We can argue by induction on the length $l$ of this filtration.
    If $l=0$ then $\e$ is stable and we are done.
    Otherwise, it will be enough to prove that $\e_l$ is Ulrich so that we can apply the inductive hypothesis.
    We have the exact sequence
    \ses{\e_l}{\e}{\e/\e_{l}}
    with $p(\e/\e_l)=p(\e)$ hence, by additivity of Euler characteristic and rank, it follows that $p(\e_l)=p(\e)$, in particular they have the same slope.
    Observe that $\e/\e_{l}$ is not torsion, otherwise $rk(\e)=rk(\e_l)$ and $p(\e_l)=p(\e)$ would imply $P(\e_l)=P(\e)$ and hence $P(\e/\e_{l})=0$.
    Then $\e/\e_{l}$ must be torsion-free by definition of stability, hence we conclude by $ii)$.
    
\gel{4}
Suppose $\e$ is Ulrich but not slope stable then, there is a non-trivial $\f\subset \e$ having the same slope as $\e$ and such that $0<rk(\f)<rk(\e)$.
If we can show that there is $\f'\subset \e$ having the slope of $\e$ and such that $\e/\f'$ is torsion-free, then by $2)$ we conclude that $\f'$ is Ulrich hence $\e$ is not stable, since $\f',\e$ share the same reduced Hilbert polynomial.

Call $\g:=\e/\f$ and define $\f'$ to be the saturation of $\f$, that is the kernel of the composition $\e\ra\g\ra F_{\g}$, where this last sheaf is the torsion-free part of $\g$.
Clearly $\mu(\f')\leq \mu(\e)$ by semistability, so we need only to show the other inequality.
Note that $\f\subseteq \f'$ and $rk(\g)=rk(F_\g)$ implies $rk(\f)=rk(\f')$.
Therefore
\[p(\f')=\dfrac{P(\f')}{rk(\f')}=\dfrac{P(\f)+P(\f'/\f)}{rk(\f)}\geq \dfrac{P(\f)}{rk(\f)}=p(\f)\]
in particular $\mu(\f')\geq \mu(\f)=\mu(\e)$ as desired.
\end{proof}

\Le\label{trivialquot}
Suppose $(X,H)$ is a Cohen-Macaulay polarised variety over some field $\gk$, then any torsion-free globally generated sheaf of slope $0$ is trivial.
\ma
\begin{proof}
Suppose we have a surjective morphism $\os_X^\rho\ra \g$ and call $\f$ its kernel.
In particular, $\f$ is torsion-free and $\mu(\f)=0$.
    We proceed by induction on $rk(\f)$, leaving $\rho$ arbitrary.
        If $rk(\f)=0$, then $\f=0$, since it is torsion-free, and we are done.
        For the inductive step, if $rk(\f)>0$ then there is a summand of $\os_X^\rho$ such that the composition $\f\hookrightarrow \os_X^\rho\twoheadrightarrow\os_X$ is non-zero; call the resulting map $\psi$.
        We get a commutative diagram
        \dia
        & 0 \ar[d] & 0 \ar[d] & 0 \ar[d] &  \\
       0 \ar[r] & Ker(\psi) \ar[r] \ar[d] &  \os_X^{\rho-1} \ar[r] \ar[d] & \ck \ar[r] \ar[d] & 0 \\
       0 \ar[r] & \f \ar[r] \ar[d] &  \os_X^{\rho} \ar[r] \ar[d] & \g \ar[r] \ar[d] & 0 \\
       0 \ar[r] & Im(\psi) \ar[r] \ar[d]  &  \os_X \ar[r] \ar[d] & \cc \ar[r] \ar[d] & 0 \\
       & 0  & 0  & 0  &  \\
        \mma
    By semistability of $\os_X^{\rho},\os_X$ we get that $deg(ker(\psi)),deg(Im(\psi))\leq 0$ but, since their sum is $deg(\f)=0$ being $\mu(\f)=0$, they have to be both zero.
    It follows that we can apply our induction hypothesis on the first line of this diagram, since $rk(ker(\psi))=rk(\f)-1$ and $\ck$ is torsion-free and of slope $0$.
    We deduce that $\ck\cong\osn^r$ for some $r\geq 0$. 
    We have $Im(\psi)\subset \os_X$ hence $Im(\psi)=\id_Y$ and $\cc\cong \os_Y$ for some subscheme $Y\subset X$.
    Being $\mu(Im(\psi))=0$, we have also $\mu(\os_Y)=0$ hence $c_1(\os_Y)=0$.
    All components of $Y$ have codimension bigger that $1$, since $c_1(\os_Y)=[Y_d]$ where $Y_d$ is the divisorial part of $Y$, it follows by excision sequence \cite{Ful}[Prop. 1.8] or directly by \cite{Ful}[Ex. 15.3.1].
    Then, $\cc\cong \os_Y$ and by Serre duality
    \[ext^1(\cc,\ck)=ext^1(\os_Y,\os_X^r)=r\cdot h^{n-1}(Y,\os_Y\otimes\omega_X)=0\]
    by dimensional reasons.
    But $\g\cong \os_X^r\oplus\os_Y$ is torsion free hence $Y=\emptyset$ and we conclude.
\end{proof}
In particular, note that the existence of an Ulrich sheaf implies that there are slope-stable ones.

As a corollary of the previous result, we can get more "positivity" properties for Ulrich sheaves different from $\os_X$.
The first result is in \cite{Lopezpositivity}[Lem. 2.1] while the second already appeared in \cite{CFK3}[Lem. 2.1].
\Co\label{grado0}
Let $(X,H)$ be an $n$-dimensional polarised variety over an infinite field and $\e$ a non-zero Ulrich sheaf on it.
\begin{enumerate}[1]
    \item $\mu(\e)\geq 0$ and we have equality if and only if $(X,H)\cong (\p^n,\osn(1))$.
    \item If $(X,H)\neq (\p^n,\osn(1))$ then for every Ulrich sheaf $\e$ we have $H^0(X,\e\du)=0$.
\end{enumerate}
\io
\begin{proof}
\gel{1}
By \thref{corollario} $i)$ an Ulrich sheaf is always globally generated hence we have a surjective morphism $\os_X^{\rho}\ra\e$.
Being $\os_X$ slope semistable, it follows by \thref{stabilitàtrick} that $0=\mu(\os_X)\leq \mu(\e)$.
Now suppose that it is an equality.
We know that $\e$ is torsion-free by \thref{corollario} and $\mu(\e)=0$ by assumption so by \thref{trivialquot} we deduce that $\e\cong \os_X$ .
Since $\os_X$ is Ulrich we conclude with \thref{oulrich} that $(X,\os_X(H))\cong (\p^n,\osn(1))$.

\gel{2}
    If $\mu(\e)>0$ then there are no morphisms $\e\ra\os_X$ by \thref{stabilitàtrick}, hence $0=hom(\e,\os_X)=h^0(X,\e\du)$.
\end{proof}

\subsection{Weak Brill-Noether property}

The fact that the reduced Hilbert polynomial of an Ulrich sheaf only depends on numerical properties of $(X,H)$ implies that all Ulrich sheaves of the same rank live in the same moduli space of semistable sheaves on $(X,H)$.
Then, we can restate the Ulrich condition as an existence condition plus a \textit{weak Brill-Noether} condition.
Let us recall the following definition from \cite{CoskunHuizengaNuer}[Definition 3.3]
\De
A coherent sheaf $\e$ on a variety $X$ satisfies \textbf{wBN} (weak Brill-Noether) if $h^j(X,\e)\neq 0$ for at most one $j\in\N$.
\Ne

Note that on a smooth complex curve is known that a general stable vector bundle in any moduli space satisfies wBN; see, for example, \cite{CoskunHuizengaNuer}[Thm. 3.8], but in general is false for higher-dimensional varieties.
Nevertheless, fixed $(X,H)$, we could ask if there are special Hilbert polynomials such that wBN holds for the general sheaf in some component of the corresponding moduli spaces.

\Le\label{wBN}
Let $(X,H)$ be a polarised variety over an infinite field.
An Ulrich sheaf $\e$ of rank $r$ on $(X,H)$ is a sheaf with Hilbert polynomial $rd\binom{t+n}{n}$ such that $\e(t)$ satisfies wBN for $-n\leq t\leq -1$.
When such an $\e$ is stable then a general sheaf in its moduli space has the same property.
\ma
\begin{proof}
If $\e$ is an Ulrich sheaf then $\e(t)$ have the required wBN property for $-n\leq t\leq -1$.
Moreover, its Hilbert polynomial has been computed in \thref{hilbertpol}.
    For the converse, note that 
    \[\sum_{j\in \N}(-1)^jh^j(X,\f(t))=\chi(\f(t))=p(\f)(t)=0 \quad \text{for}\quad -n\leq t\leq -1\]
    but at most one element in the left sum can be non-zero, by assumption,
    therefore they must all vanish, and we conclude that $\e$ is Ulrich by definition.

    For the second claim recall that Ulrich sheaves are always semi-stable by \thref{stability} and have a well determined Hilbert polynomial by \thref{corollario}.
Being Ulrich is encoded by a cohomology vanishing, hence by upper-semicontinuity of cohomology it is open in flat families.
Since the moduli space is constructed as a quotients of one such family, the universal one over the adequate $Quot$ scheme, stable Ulrich sheaves form an open subset in the stable locus of those moduli spaces by \cite{HuyLeh}[Cor. 4.3.5].
\end{proof}

We present an example in which the wBN condition is actually automatic.
We follow the proof given in \cite{CH}[Lem. 4.4] but we think that slope-stability is necessary, instead of simple stability.

\Prop
Let $X$ be a smooth del Pezzo surface with $H=-K_X$ ample and globally generated.
Any slope-stable locally free sheaf $\e$ with Hilbert polynomial $P(\e)(t)=rd\binom{t+2}{2}$ is an Ulrich bundle on $(X,H)$.
\One
\begin{proof}
    If $d:=H^2$ then by Riemann--Roch
    \[p(\os_X(1))(t)=P(\os_X(1))(t)=\chi(\os_X(t+1))=\chi(\os_X)+\dfrac{(t+1)H\cdot ((t+1)H-K_X)}{2}=\]
    \[1+d\dfrac{(t+1)(t+2)}{2}=1+d\binom{t+2}{2}>\dfrac{P(\e)(t)}{r}=p(\e)(t).\]
    Since $\e$ is slope stable, the above inequality implies that $0=hom(\os_X(1),\e)=h^0(X,\e(-1))$.
    This also implies $h^0(X,\e(-2))=0$.
    Since the dual of a slope-stable vector bundle stays slope-stable, see \cite{OSS}[Lem. 1.2.4], we have that the above reasoning applies to $\e\du(2)$, hence $h^0(X,\e\du(1))=0$.
    By Serre duality this last equation gives $h^2(X,\e(-2))=0$ hence also $h^2(X,\e(-1))=0$.
    It follows that wBN property holds for $\e(-1)$ and $\e(-2)$ and we conclude by \thref{wBN}.
\end{proof}

\section{Numerical properties}

Next, we will deal a bit with numerical properties of Ulrich bundles.
We start by computing the standard formulas for the degrees of the first and second Chern class.
Then we will discuss briefly positivity results.
After recalling a criterion to ensure ampleness for Ulrich bundles, we focus on discussing when positivity of their Chern classes fails. 

\subsection{Chern classes}

In this section, we assume $X$ to be a smooth variety and $\gk$ to be infinite.
Our strategy will be as direct as possible: applying Grothendieck's version of the Riemann--Roch theorem for a finite morphism $f:X\ra\p^n$ we get all the degrees of Chern classes of $\e$ at once. 
The formula for the first Chern class of $\e$ seems to date back to \cite{CH}[Lem. 2.4], although in a slightly different form.
The formula for $c_2$ can be found, for example, in \cite{CKM_delpezzo}[Prop. 2.19]. 

\Prop\label{c12}
For an Ulrich bundle $\e$ of rank $r$ on a smooth polarised variety of dimension $n$ over an infinite field $\gk$ we have
\[c_1(\e)\cdot H^{n-1}=\dfrac{r((n+1)H+K_X)\cdot H^{n-1}}{2}\]
\[c_2(\e)\cdot H^{n-2}=\left(\dfrac{r}{12}\left(K_X^2+c_2(X)-\dfrac{1}{2}H^2(3n^2+5n+2)\right)+\frac{c_1(\e)}{2}\left( c_1(\e)-K_X\right)\right)H^{n-2}\]
\One
\begin{proof}
    Since $\gk$ is infinite, we have a finite surjective morphism $f:X\ra\p^n$ given by a subsystem of $|H|$ and $f_*\e\cong \osn^{rd}$, with $d=H^n$.
    We will denote intersection product by a $"\cdot"$.
    Moreover, by abuse of notation, we will often identify a zero cycle with its degree.
    
    By Grothendieck--Riemann--Roch theorem \cite{Ful}[Thm. 15.2] we have 
    \begin{equation}\label{grr}
        f_*(ch(\e)\cdot td(X))=ch(f_*\e)\cdot td(\p^n)=rd\cdot td(\p^n)
    \end{equation}
    We will need to know those classes just up to degree $2$.
    From \cite{Ful}[Example 3.2.3] we know 
    \[ch(\e)=r+c_1(\e)+\dfrac{c_1(\e)^2}{2}-c_2(\e)+\dots\]
    By \cite{Ful}[Example 3.2.4 and Example 3.2.11] we have 
    \[td(X)=1-\dfrac{K_X}{2}+\dfrac{K_X^2+c_2(X)}{12}+\dots \qquad td(\p^n)=1+\dfrac{(n+1)}{2}H_\p+\dfrac{(n+1)^2+\binom{n+1}{2}}{12}H_\p^2+\dots\]
    where $H_\p$ is the hyperplane class on $\p^n$.
    Now we can rewrite the degree one (in the grading of the Chow group $A^\bullet(X)$) terms in \eqref{grr} as
    \[f_*\left(c_1(\e)-\dfrac{rK_X}{2}\right)=\dfrac{rd(n+1)}{2}H_\p.\]
    Recall that, by definition, $H\sim f^*H_{\p}$, so multiplying both sides of the previous equality by $H_{\p}^{n-1}$ we get 
    \[\dfrac{rd(n+1)}{2}=\dfrac{rd(n+1)}{2}H_\p^n=f_*\left(c_1(\e)-\dfrac{rK_X}{2}\right)H_{\p}^{n-1}=f_*\left(\left(c_1(\e)-\dfrac{rK_X}{2}\right)H^{n-1}\right)\]
    where in the last step we used projection formula, see \cite{Ful}[Thm. 2.3 c)].
    Rearranging this last formula we get the first claim.

    For the second one, we make similar computations but with the degree $2$ part of \eqref{grr}.
    \[f_*\left(\dfrac{c_1(\e)^2}{2}-c_2(\e)-\dfrac{c_1(\e)\cdot K_X}{2}+\dfrac{r(K_X^2+c_2(X))}{12}\right)=\dfrac{rd(n+1)(3n+2)}{24}H_\p^2.\]
    Multiplying by $H_{\p}^{n-2}$, using projection formula and rearranging we get the desired result.
\end{proof}

In principle, the same proof allows us to give formulas for all the degrees of the Chern classes of $\e$, but they become increasingly more complicated and we will not need them for the rest of this work.

The previous result suggests a particular value for the first Chern class of a (Ulrich) bundle.
\De\label{special}
Suppose $(X,H)$ is an $n$-dimensional Gorenstein  polarized scheme, we will call a rank $r$ vector bundle $\e$ \textbf{special} if $c_1(\e)\sim \dfrac{r((n+1)H+c_1(\omega_X))}{2}$.
\Ne

Clearly, if $Pic(X)\cong\Z$ then every (Ulrich) sheaf must be special; however, note that special Ulrich bundles are very rare, see \thref{pic1}.
In general there are plenty of non-special Ulrich bundles.

The notions of being U-dual and special are linked even for arbitrary vector bundles.
\Co\label{specialdual}
Let $(X,H)$ be a smooth projective variety and $\e$ be a vector bundle of rank $r$.
If $\e$ is self U-dual then $2c_1(\e)\sim r((n+1)H+K_X)$.
Somehow conversely, a special bundle of rank $1$ or $2$ is always self U-dual .
\io
\begin{proof}
Suppose that $\e\cong\e\du((n+1)H+K_X)$.
Then taking the first Chern class on both sides we get $c_1(\e)\sim c_1(\e\du)+r((n+1)H+K_X)$, where we used \cite{Ful}[Example 3.2.2], hence the claim.

For the converse, a special line bundle must be $\e\cong \os_X\left(\dfrac{(n+1)H+K_X}{2}\right)$
    then we see easily that $\e\cong \e\du((n+1)H+K_X)$.
    If $rk(\e)=2$ then $\e\du\otimes det(\e)\cong\e$ and by definition of speciality $det(\e)\cong\os_X((n+1)H+K_X)$ so we conclude as before.
\end{proof}

\subsection{Positivity}

Let $\e$ be an Ulrich bundle on a smooth $n$-dimensional polarized variety $(X,H)$, then $\e$ is globally generated by \thref{corollario}.
In general, we cannot expect more positivity since $\osn$ is Ulrich on $(\p^n,\osn)$ but we can try to control where positivity fails.
We start with a definition.
\De
A vector bundle $\e$ is called \textbf{ample} (\textbf{very ample}, \textbf{big}) if the line bundle $\os_{\p(\e)}(1)$ is ample (very ample, big).
\Ne

The following are obtained from \cite{LopezSierra_geom}[Thm. 1] and \cite{Valerio+}[Thm. 4].
\Te\label{}
Let $(X,H)$ be a complex, smooth polarised variety such that either $H$ is very ample or there is a linear sub-system $V\subset |H|$ inducing a morphism étale on the image.
Then the following are equivalent:
\begin{itemize}
    \item $\e$ is very ample
    \item $\e$ is ample
    \item for any curve $C\subset X$ such that $H\cdot C=1$ we have $\e|_C$ is ample
\end{itemize}
\Ma

Let us go back to the weakest positivity property that all Ulrich bundles have.
Since $\e$ is globally generated hence $det(\e)$ is globally generated, in particular $c_1(\e)$ effective.
Actually, all $c_i(\e)$ are effective since they can be constructed as degeneracy loci, see \cite{Ful}[Example 14.4.3].
Since $|c_1(\e)|$ has no base-points we also have $c_1(\e)^l$ is effective for any $l>0$, since can be represented by a complete intersection of divisors in that linear system, and then $c_1^l(\e)\cdot H^{n-l}=0$ implies $c_1(\e)^l\sim 0$ by ampleness of $H$.
Therefore, $c_1^l(\e)\cdot H^{n-l}=0$ if and only if $c_1(\e)^l\sim 0$.

We have already seen that if $c_1(\e)\sim 0$ then $(X,H)\cong (\p^n,\osn(1))$, see \thref{grado0}.

\Oss
\begin{itemize}
    \item There is a classification of the pairs $(X,H)$ admitting an Ulrich bundle such that $c_1(\e)^l\sim 0$ over $\gk=\C$.
    For $dim(X)=2,3$ it is contained in \cite{LopezMunoz21}[Thm. 1 and Thm. 2], for $dim(X)=4$ in \cite{LopezMunozSierra23}[Thm. 1] and in arbitrary dimension if $l\leq \left\lfloor\dfrac{n+1}{2}\right\rfloor$ in \cite{LopezSierra_geom}[Cor. 4].
    \item For special Ulrich bundles we expect positivity of $c_1(\e)$.
Indeed, $\e$ is globally generated hence also $det(\e)$ is such.
It follows that $det(\e)$ non-big implies that $0=c_1(\e)^n=(lH)^n=l^nH^n$ hence $l=0$ and we conclude that $(X,H)$ is $(\p^n,\osn(1))$ by \thref{grado0}.
\end{itemize}   
\one
 
The following comes from discussions with A.F. Lopez, I thank him for allowing me to reproduce it here.
\Le\label{c2=0}
If $\gk=\C$.
If $c_2(\e)\sim 0$ then $c_1(\e)^2=0$ hence $c_1(\e)$ is not big. 
It follows that either $(X,H)\cong (\p^n,\osn(1))$ or $(X,H)\cong (\p(E),\os_{\p(E)}(1))$ with $E$ vector bundle over a curve and $\e$ is the pullback of a vector bundle on $C$. 
\ma
\begin{proof}
Up to taking smooth hyperplane sections we can suppose $dim(X)=2$, indeed if $H_1,\dots H_{n-2}\in |H|$ then $0=c_1(\e|_{\cap H_i})^2=c_1(\e)^2\cdot H^{n-2}$ and by ampleness of $H$ it is equivalent to $c_1(\e)^2\sim 0$.
    Since $\e$ is globally generated and $c_2(\e)\sim 0$, if we put $D\sim c_1(\e)$ and $r=rk(\e)$, we have 
    \ses{\os_X^{r-1}}{\e}{\os_X(D)}
    If the above sequence is split then $\os_X$ is Ulrich and hence $(X,H)\cong (\p^n,\osn(1))$.
    If not then 
    \[0\neq Ext^1(\os_X(D),\os_X)\cong H^1(X,\os_X(-D))\cong H^{n-1}(X,\os_X(D+K_X))\du\]
    hence by Kawamata-Viehweg $D$ cannot be both big and Nef.
    Being $\e$ globally generated also $c_1(\e)$ is such, in particular $D$ is Nef and hence it is not big, by the above argument.
    This implies that $c_1(\e)^2=0$.

    The last claim follows by classification, precisely from \cite{LopezSierra_geom}[Cor. 4].
\end{proof}

We will see that, since $c_2(\e)$ is effective, when $c_2(\e)\neq 0$ we can always write an (Ulrich) bundle of rank $r$ as an extension of $\os_X^{r-1}$ with some twisted ideal sheaf.

\Do
What can we say if $c_2(\e)^l\sim 0$?
\da

\chapter{Constructing sheaves}

In this chapter we present the standard Hartshorne-Serre correspondence to construct torsion-free sheaves on a smooth variety.
From one side it can be used as a way of constructing sheaves exploiting the geometry of the variety, on the other, sometimes vector bundles are easier to handle than specific subvarieties, as happens for line bundles and hypersurfaces.

Hartshorne-Serre correspondence, in its original formulations, is essentially the extension of the divisor/line bundle correspondence for higher ranks.
In the rank $2$ case, this is the assertion that vector bundle with effective second Chern class can be recovered from the zero locus of any of its sections.
In general, zero loci of sections are substituted by degeneracy loci of codimension $2$, which are always controlled by the second Chern class.

In this mindset, we will derive more explicit criteria for the existence of Ulrich bundles, characterising the needed degeneracy loci.
Finally, we want to address the problem of directly constructing higher rank Ulrich bundles starting from low rank ones.
Any extension of two Ulrich bundles is Ulrich by \thref{2-3} but is not stable, hence we will need to search for deformations of such sheaves.

In the first section, we discuss Hartshorne-Serre correspondence on different levels of generality, see \thref{sheafcorrispondenza} and \thref{vectorcorrispondenza}.
In particular, we also present a relative version, \thref{modularcorrispondenza}, which will be useful when studying moduli space.
The second section is a specialization of the previous one to the Ulrich case, \thref{Ulrichcorrispondenza} and \thref{Ulrichcorrispondenza'}, and will be useful later.

The third section will be quite technical.
The goal is to build a machinery, \thref{wildext}, that, if we can verify some cohomological vanishing, will take as input families of stable vector bundles and give us similar families of higher rank.
The main technical tool will be deformation theory and the existence of modular families of simple sheaves.

\startcontents[chapters]
\printcontents[chapters]{}{1}{}

\section{Hartshorne-Serre correspondence}

The original idea of a correspondence between rank $2$ vector bundles and codimension $2$ subvarieties goes back to \cite{Serre60} in the affine case and is developed for rank $2$ bundles on $\p^3$ in \cite{Har78}[Thm. 1.1] and for reflexive sheaves in \cite{Har_stab_refl}[Thm. 4.1].
A statement for vector bundles of any rank $\geq 2$ is in \cite{Mae}[Thm. 1.7].
In \cite{Arr} there is quite a different proof with respect to the previous ones.
We are going to follow mainly the ideas contained in \cite{CFK3}[Thm. 3.1] and \cite{Mae}.

\subsection{General case}

We start recalling some facts we will use in the proof of our main theorems.

\Le\label{degenerazione'}
Let $X$ be a smooth variety over a field $\gk$.
\begin{enumerate}
    \item Suppose that we have an exact sequence
\ses{\e}{\f}{Q}
with $\e$ locally free, $\f$ torsion-free and $Q$ locally free outside a closed subset of codimension at least $2$ in $X$.
Then $Q$ is torsion-free.

\item If $Q$ is any rank $1$ torsion-free sheaf then it is an ideal sheaf twisted by $det(Q)$\footnote{Note that, being $X$ smooth, the determinant is a well defined line bundle since we can compute it from a locally free resolution.}.
\end{enumerate}
\ma
\begin{proof}
\begin{enumerate}
    \item 
    By assumption $Q$ has full support on $X$ hence by \thref{s12} we just need to verify that satisfies $S_1$. 
But, by \cite[\href{https://stacks.math.columbia.edu/tag/00LX}{Tag 00LX}]{Stacks} we have \[depth(Q_x)\geq \min\{depth(\e_x)-1, depth(\f_x)\}\]
so is clear that $Q$ could satisfies $S_1$ except, possibly, for some $x\in X$ such that $depth(\e_x)=1$.
Being $\e$ locally free we have $\e_x\cong \os_{X,x}^\rho$ hence $1=depth(\e_x)=depth(\os_{X,x})=dim(\os_{X,x})$ but in this case we already know by assumption that $Q_x$ is free.
   
   \item Since $Q$ is torsion-free, the natural map $Q\ra Q^{\vee\vee}$ is injective.
    The latter sheaf is reflexive of rank $1$ hence, being $X$ smooth, it must be locally free by \cite{Har_stab_refl}[Prop. 1.9].
    On the locus $U$ where $Q$ is locally free the inclusion in $Q^{\vee\vee}$ is an isomorphism.
    It follows that those sheaves share the same determinant bundle, being $codim(U)\geq 2$.
    We conclude that $Q^{\vee\vee}$, being locally free, must be isomorphic to $det(Q)$ so that we have an inclusion
    \[Q\otimes det(Q)\du\hookrightarrow  Q^{\vee\vee}\otimes det(Q)\du\cong \os_X.\]
\end{enumerate} 
\end{proof}

The following is a weak version of the Hartshorne-Serre correspondence.
Given two sheaves $\e,\f$ and $V\subset H^0(X,\e), V'\subset H^0(X,\f)$, we say that $(\e,V)$ is isomorphic to $(\f,V')$ is there is an isomorphism $\e\ra\f$ sending $V$ to $V'$.

\Prop\label{sheafcorrispondenza}
Let $X$ be a smooth projective variety of dimension $n\geq 2$ over some algebraically closed field $\gk$.
Fixed a divisor $D$, there is a bijective correspondence between 
\begin{itemize}
\item isomorphism classes of pairs $(\e,V)$ such that
\begin{enumerate}[1)]
    \item $\e$ is a torsion-free sheaf of rank $r>1$, $det(\e)\cong\os_X(D)$ and has no direct summand isomorphic to $\os_X$
    \item $V\subset H^0(X,\e)$ has dimension $r-1$, the morphisms $V\otimes\os_X\ra \e$ is injective and drops rank on a locus of pure codimension $2$
\end{enumerate}
\item pairs $(Y,W)$  such that
\begin{enumerate}[i)]
\item $Y\subset X$ is a closed subscheme of pure codimension $2$ 
    \item $W\subseteq Ext^1(\id_Y(D),\os_X)$ has dimension $r-1$
    \end{enumerate}
    \end{itemize}
    moreover they fit in the exact sequence (where we can put $W\du$ instead of $V$)
\begin{equation}\label{scorrispondenza}
    0\ra  V\otimes\os_X\ra \e\ra \id_Y(D)\ra 0.
\end{equation}
\One
\begin{proof}
    Let us start with $(\e,V)$.
  The vector space $V\subset H^0(X,\e)$ gives an evaluation morphism whose degeneracy locus has pure codimension $2$ by hypothesis: this will be $Y$ and hence we have the sequence \eqref{scorrispondenza} by \thref{degenerazione'}. 
  Applying $\Hom(-,\os_X)$ to
\sesl{\id_Y(D)}{\os_X(D)}{\os_Y(D)}{subd'} 
we get the isomorphisms $(\id_Y(D))\du\cong \os_X(-D)$, since being $Y$ of codimension $2$ we have $\Hom(\os_Y(D),\os_X)\cong 0\cong \ext^1(\os_Y(D),\os_X)$ by \cite{AltKle_duality}[§IV Lem. 5.1].
Hence, applying $Hom(-,\os_X)$ to \eqref{scorrispondenza} we get
\begin{equation}\label{long}
     H^0(\os_X(-D))\hookrightarrow Hom(\e,\os_X)\xra{\phi} V\du\ra Ext^1(\id_Y(D),\os_X)\ra Ext^1(\e,\os_X)\ra \dots 
\end{equation}
Define $W$ to be the image of $V\du$ inside $Ext^1(\id_Y(D),\os_X)$; we will show that they are isomorphic, in particular $r-1=dim(V\du)=dim(W)$, by showing that $\phi=0$.
Indeed, if we call $ev$ the inclusion $V\otimes\os_X\ra\e$, the map $\phi$ sends a morphism $a:\e\ra \os_X$ to the composition
\[V\xra{H^0(ev)}H^0(X,\e)\xra{H^0(a)}H^0(X,\os_X)\cong \gk.\]
If this composition is non-zero then we have a non-zero map $\os_X\ra\e\ra\os_X$ which must then be an isomorphism.
In particular, $\os_X$ would be a direct summand of $\e$ contradicting our assumption.
So we have proved that $1)+2)\Longrightarrow i)+ii)$.

Let us now start from $(Y,W)$.
We have
\[Ext^1(\id_Y(D),W\du\otimes\os_X)\cong Ext^1(\id_Y(D),\os_X)\otimes W\du\cong Hom(W,Ext^1(\id_Y(D),\os_X))\] 
where in the first step we used additivity of $Ext$ and the second follows by a standard isomorphism of vector spaces, since $W\du$ is finite dimensional.
This identification implies that the inclusion $W\subseteq  Ext^1(\id_Y(D),\os_X)$ gives the sequence \eqref{scorrispondenza}, after we set $V:=W\du$.
By \eqref{long}, being the map $V\du\ra Ext^1(\id_Y(D),\os_X)$ injective, we have $\phi=0$.
This implies that $h^0(X,\e\du)=Hom(\e,\os_X)=H^0(X,\os_X(-D))$, hence if $\e$ had a trivial summand then $-D$ would be effective.
We will derive a contradiction from this.
Indeed, consider the global section of $\e$ corresponding to a trivial direct summand of it.
This section must live outside $V$ because $\phi=0$. 
Therefore, by the cohomology sequence of \eqref{scorrispondenza}, we would have $h^0(X,\id_Y(D))>0$ and, since $Y$ is non-empty, $D$ should be effective and non-zero, a contradiction.
The rank and determinant of $\e$ are now determined by \eqref{scorrispondenza}.
In addition, $\e$ is torsion-free being an extension of torsion-free sheaves, so we have proved that $i)+ii)\Longrightarrow 1)+2)$.

The correspondence is clearly bijective on the second coordinate of our pairs by biduality for finite-dimensional vector spaces.
Starting with $(Y,W)$, we construct $\e$ as an extension \ses{W\du\otimes\os_X}{\e}{\id_Y(D)}
so if we take the degeneracy locus of the morphism $W\du\otimes\os_X\ra \e$ we get back the same $Y$ we started with.
On the other hand, beginning with $(\e,V)$, after the construction of $Y$ we can identify uniquely the isomorphism class of $\e$ as the class of the extension it gives inside $Ext^1(\id_Y(D),V\otimes\os_X)$.
\end{proof}

\Co\label{r00}
If $h^1(X,\os_X(-D))=0=h^2(X,\os_X(-D))$ then $Ext^1(\id_Y(D),\os_X)\cong H^0(Y,\omega_Y(-D-K_X))$ hence in $ii)$ we can ask for $W\subseteq H^0(Y,\omega_Y(-D-K_X))$.
\io
\begin{proof}
Applying $Hom(-,\os_X)$ to the sequence \eqref{subd'} we get
\[\dots H^1(X,\os_X(-D))\ra Ext^1(\id_Y(D),\os_X)\ra Ext^2(\os_Y(D),\os_X)\ra H^2(X,\os_X(-D)) \dots \]
From the supposed vanishing we conclude
\[Ext^1(\id_Y(D),\os_X)\cong Ext^2(\os_Y(D),\os_X)\cong H^0(Y,\omega_Y(-D-K_X))\]
where the last isomorphism is proven as follows.
By applying Serre duality twice we get
 \[Ext^2(\os_Y(D),\os_X)\cong H^{n-2}(X,\os_Y(D+K_X))\du\cong\]
 \[ \cong H^{n-2}(Y,\os_Y(D+K_X))\du\cong Hom(\os_Y(D+K_X), \omega_Y)\cong H^0(Y,\omega_Y(-D-K_X)).\]
 Note that the second to last step is possible by the first part of \cite{AltKle_duality}[I (1.3)].
\end{proof}

\subsection{Locally free case}

We continue with the above notations.
Our goal is to specialize \thref{sheafcorrispondenza} adding more regularity both on $\e$ and $Y$.
Then, under those assumptions, we will give an explicit description of the morphism $\e\ra\id_Y(D)$ in \eqref{scorrispondenza}.
For simplicity, we first state a lemma.

\Le\label{cmcodim2}
Let $X$ be a smooth projective variety and $Y$ a closed subscheme such that $dim(Y)\leq dim(X)-2$.
The following are equivalent:
\begin{enumerate}
    \item $\id_Y$ has a locally-free resolution of length $1$
    \item $\ext^i(\os_Y,\os_X)=0$ for all $i>2$
    \item $Y$ is Cohen-Macaulay of pure codimension $2$.
\end{enumerate}
\ma
\begin{proof}
By \thref{extpd}, $\id_Y$ has a locally-free resolution of length $1$ if and only if $\ext^i(\id_Y,\os_X)=0$ for all $i>1$.
From the exact sequence 
\ses{\id_Y}{\os_X}{\os_Y}
this is equivalent to $\ext^i(\os_Y,\os_X)=0$ for all $i>2$.

Assuming those, we can prove $3)$ as follows.
By \thref{extpd} the second condition is equivalent to $pd(\os_{Y,y})\leq 2$ for all $y\in Y$.
But, by \thref{codimensione}, this last condition is satisfied if and only if $codim(Y)\leq 2$ and $Y$ Cohen-Macaulay, as desired.
\end{proof}

Note that being $X$ smooth, we have $\omega_X\cong \os_X(K_X)$.
In the long exact sequence obtained applying $\Hom(-,\os_X)$ to \eqref{subd'} we get
\begin{equation}\label{contoduale'}
    \ext^1(\id_Y(D),\os_X)\cong \ext^2(\os_Y(D),\os_X)\cong \ext^2(\os_Y,\os_X)\otimes\os_X(-D)\cong\omega_Y(-D-K_X)
\end{equation} being $\omega_Y\cong\ext^2(\os_Y,\os_X(K_X))$.
We will denote by $d$ the boundary map 
\begin{equation}\label{boundary}
    d:Ext^1(\id_Y(D),\os_X)\ra Ext^2(\os_Y(D),\os_X)\cong H^0(Y,\omega_Y(-D-K_X)),
\end{equation}
obtained by applying $Hom(-,\os_X)$ to \eqref{subd'}, recall also the final part of \thref{r00}.

Here is the Hartshorne-Serre's correspondence.
\Te\label{vectorcorrispondenza}
Let $X$ be a smooth projective variety of dimension $n\geq 2$ over some algebraically closed field and $D$ a divisor.
There is a bijective correspondence between 
\begin{itemize}
\item isomorphism classes of pairs $(\e,V)$ such that
\begin{enumerate}[1)]
    \item $\e$ is a vector bundle of rank $r>1$, $det(\e)\cong\os_X(D)$ and has no direct summand isomorphic to $\os_X$
    \item $V\subset H^0(X,\e)$ has dimension $r-1$, the morphisms $V\otimes\os_X\ra \e$ is injective and drops rank on a locus of pure codimension $2$
\end{enumerate}
\item pairs $(Y,W)$  such that
\begin{enumerate}[i)]
\item $Y\subset X$ is a Cohen-Macaulay closed subscheme of pure codimension $2$ 
    \item $W\subseteq Ext^1(\id_Y(D),\os_X)$ has dimension $r-1$ and the evaluation morphism $d(W)\otimes\os_X \ra \omega_Y(-D-K_X)$ is surjective
    \end{enumerate}
    \end{itemize}
    moreover they fit in the exact sequence
\begin{equation}\label{vcorrispondenza}
    0\ra  V\otimes\os_X\ra \e\ra \id_Y(D)\ra 0.
\end{equation}

\Ma
\begin{proof}
Having the blueprint of \thref{sheafcorrispondenza}, we need only to show that $\e$ is locally free if and only if $Y$ is Cohen-Macaulay and $d(W)$ generates $\omega_Y(-D-K_X)$.
By \thref{extpd}, $\e$ is locally free if and only if $\ext^i(\e,\os_X)=0$ for all $i>0$.
Applying $\Hom(-,\os_X)$ to \eqref{vcorrispondenza}, this is equivalent to $\ext^{i}(\id_Y(D),\os_X)=0$ for all $i>1$ and the morphism $W\otimes\os_X\ra \ext^{1}(\id_Y(D),\os_X)$ is surjective.
The first property is equivalent to $Y$ being Cohen-Macaulay by \thref{cmcodim2}.
We have already seen in \eqref{contoduale'} that  $\ext^1(\id_Y(D),\os_X)\cong \omega_Y(-D-K_X)$ hence, the second property above is exactly the second condition which we have to verify.
\end{proof}

\Co\label{hsc2}
In the above setting, $Y\sim c_2(\e)$.
Suppose that we have another pair $(Y',V')$ that fits in a sequence analogous to \eqref{vcorrispondenza}, in particular $Y\sim Y'$, and call $\overline{V}:=<V,V'>\subseteq H^0(X,\e)$.
If $Y,Y'$ are integral, $dim(\overline{V})=r$ and $\overline{V}\otimes\os_X\ra\e$ drops rank along a normal integral divisor $\Delta\in|\os_X(D)|$ then $\Delta$ contains both $Y,Y'$ and they are linearly equivalent on it.
\io
\begin{proof}
    From 
    \ses{\id_{Y_b}}{\os_X}{\os_{Y_b}}
    we get $ c_2(\id_{Y})\sim -c_2(\os_{Y})\sim Y$ where in the last step we used \cite{Ful}[Example 15.3.1].
    Applying \cite{Ful}[Thm. 3.2 e)] on \eqref{vcorrispondenza} twisted by $\os_X(-D)$ we have
    \[c_2(\id_{Y})=c_2(\e(-D))=c_2(\e)\]
    where in the last step we used \cite{Ful}[Example 3.2.2 p.56] and $c_1(\e)\sim D$.
    
    For the second part we can construct the diagram
    \dia
    & & & \os_X \ar[d] & \\
    0 \ar[r] & V \otimes \os_X \ar[d] \ar[r] & \e \ar[d, "id"] \ar[r] & \id_Y(D) \ar[d] \ar[r] & 0 & \\  
        0 \ar[r] & \overline{V} \otimes \os_X \ar[r] \ar[d] & \e \ar[r] & \id_{Y/\Delta}(D) \ar[r] & 0 & \\ 
        & \os_X & & & \\
    \mma
    Being $Y$ integral and $\Delta$ normal and integral by \cite{Schwede}[Prop. 3.4] we have $\id_{Y/\Delta}\cong \os_\Delta(-Y)$, meant as a reflexive rank $1$ sheaf.
    The same diagram holds replacing $V$ with $V'$ and $Y$ with $Y'$, hence $\os_{\Delta}(-Y)\cong \os_\Delta(-Y')$ implying that $Y$ and $Y'$ are linearly equivalent by \cite{Schwede}[Prop. 3.11].    
\end{proof}

\thref{sheafcorrispondenza} and \thref{vectorcorrispondenza} are two extremal cases, meaning that one can deduce other intermediary cases linking the depth of stalks of $\e$ and regularity properties of $Y$; see, for example, \cite{Har_stab_refl}[Thm. 4.1].

\Oss\label{r01}
If $r=2$ then $\e|_Y\cong \n_Y$ is locally free hence $Y$ is a local complete intersection.
This is a well-known fact but we will recall the proof.\footnote{The fact that a Gorenstein codimension $2$ subscheme of a smooth scheme is locally complete intersection is originally due to Serre \cite{Serre60}[§7 Prop. 5]}
First, restrict \eqref{vcorrispondenza} to $Y$.
This new sequence will be right exact and the left map will be $0$.
Therefore 
\[\e|_Y\cong \id_Y(D)|_Y\cong \n_Y\du(D)\]
hence the conormal sheaf is locally free, meaning that $Y$ is a locally complete intersection.
Recalling that taking determinants commutes with restricting the bundle we have 
\[\os_Y(D)\cong det(\e|_Y)\cong det(\n_Y\du(D))\cong det(\n_Y\du)\otimes \os_Y(D)^{\otimes 2}\]
when the last isomorphism comes from \cite{Ful}[Example 3.2.2] and the fact that $\n_Y\du$ has rank $2$.
We conclude that $det(\n_Y\du)\cong \os_X(-D)$ hence, again from the rank $2$ hypothesis, we have
\[\n_Y\cong \n_Y\du\otimes det(\n_Y)\cong \n_Y\du(D)\cong \e|_Y.\]
\one

During the proof of \cite{Har78}[Thm. 1.1], the sequence \eqref{vcorrispondenza}, in the case $r=2$, was thought as the Koszul complex of a section of $\e$, hence we have an explicit description of the morphism $\e\ra\id_Y(D)$.
Here we want to look at what happens for any $r\geq 2$.
Recall that for a locally free sheaf $\e$ of rank $r$ we have a perfect pairing 
\[\e\otimes\bigwedge^{r-1}\e\ra \bigwedge^{r}\e\cong \os_X(D),\]
in particular $H^0(X,\bigwedge^{r-1}\e)\cong Hom(\e,\os_X(D))$.
\Le\label{wedge}
Let $X$ be a smooth projective variety and suppose we have an exact sequence
\begin{equation}\label{cor'''}
    0\ra  V\otimes\os_X\ra \e\ra \id_Y(D)\ra 0
\end{equation}
with $\e$ locally free of rank $r$ and $Y$ of codimension $2$.
Then the morphism $\e\ra \id_Y(D)\subset \os_X(D)$ is, up to multiplication by a scalar, the one determined by $[V]\in H^0(X,\bigwedge^{r-1}\e)$.
\ma
\begin{proof}
Call $v:\e\ra\os_X(D)$ the morphism given by $[V]\in H^0(X,\bigwedge^{r-1}\e)$, that is, if we choose any basis $v_1,\dots v_{r-1}$ of $V$ this morphism sends a section $s$ of $\e$ to $s\wedge v_1\wedge\dots \wedge v_{r-1}$.
Therefore, $V\otimes\os_X$, seen as a subsheaf of $\e$, is contained in $ker(v)$ so we have a surjective morphism $\id_Y(D)\ra Im(v)$.
From \eqref{cor'''}, the sections in $V$ are linearly independent on the general point of $X$ hence the morphism $v$ in non-zero.
This implies that its image has rank $1$, since $\os_X(D)$ is torsion-free.
But then the surjective map $\id_Y(D)\ra Im(v)$ must be an isomorphism, since its kernel would be a torsion subsheaf of $\id_Y(D)$, hence the claim.
\end{proof}

\Co\label{surjdeg}
In the setting of \thref{wedge}, if $H^1(X,\os_X)=0$ then all the sections in $H^0(X,\id_Y(D))$ are of the form $s\wedge v_1\wedge\dots \wedge v_{r-1}$ for some $s\in H^0(X,\e)$.
\io

\subsection{Relative version}

Our next aim is to prove that Hartshorne-Serre construction works in families, so that we do not have only a function between two sets but an algebraic morphism between the base of a family of subvarieties and the base of a family of vector bundles.
We believe this result is well known, see for example \cite{MarTikcubic}[Lem. 5.2], but could not find a detailed proof in the literature.

Note that for families over a smooth base $B$ essentially we can just apply \thref{sheafcorrispondenza} or \thref{vectorcorrispondenza} to $X\times B$.
For the general case, we need just a couple of addenda. 
We start with some notation.
Fixed a scheme $B$ over $\gk$, we consider a family $\Pi:\cX\ra B$.
Moreover, given $\cY\subset \cX$ and a sheaf $\scr{E}$ on it, for any $b\in B$ we set $\cX_b:= \cX\times_B\{b\}$, $\cY_b:= \cY\times_B\{b\}$ and $\scr{E}_b:=\scr{E}|_{X_b}$.

\Prop\label{familycorrispondenza}
Fix a variety $B$ over $\gk=\overline{\gk}$.
Let $\Pi:\cX\ra B$ be a smooth, proper morphism of relative dimension $n\geq 2$ and $\cd\subset\cX$ a Cartier divisor flat over $B$.
There is a bijective correspondence between 
\begin{itemize}
\item isomorphism classes of pairs $(\scr{E},\cV)$ such that
\begin{enumerate}[1)]
    \item $\scr{E}$ is a torsion-free sheaf of rank $r>1$ on $\cX$ flat over $B$, $\scr{E}_b$ is torsion-free and has no direct summand isomorphic to $\os_{\cX_b}$ and $det(\scr{E}_b)\cong\os_{\cX_b}(\cd_b)$ for all $b\in B$,
    \item $\cV\subset H^0(\cX,\scr{E})$, for all $b\in B$ we have $dim(\cV\otimes\os_{\cX_b})=r-1$, the morphism $\cV\otimes\os_{\cX_b}\ra \scr{E}_b$ is injective and drops rank on a locus of pure codimension $2$;
\end{enumerate}
\item pairs $(\cY,\cW)$  such that
\begin{enumerate}[i)]
\item $\cY\subset \cX$ is a subscheme of pure codimension $2$, flat and proper over $B$
    \item $\cW\subseteq Ext^1(\id_\cY(\cd),\os_{\cX})$ such that $dim(\cW\otimes k(b))=r-1$ for all $b\in B$.
    \end{enumerate}
    \end{itemize}
    Moreover, they fit in the exact sequence
\begin{equation}\label{fcorrispondenza}
    0\ra  \cV\otimes\os_{\cX}\ra \scr{E}\ra \id_{\cY}(\cd)\ra 0
\end{equation}
whose restrictions to $\cX_b$ are sequences of the form \eqref{scorrispondenza}.
\One
\begin{proof}
We proceed similarly to \thref{sheafcorrispondenza}.
Starting from $(\scr{E},\cV)$ we construct the sheaf $\cc$ as the cokernel of $\cV\otimes\os_{\cX}\ra \scr{E}$ and call $\cY$ the locus where this map drops rank.
By assumption $2)$ we know that the sequence
\ses{\cV\otimes\os_{\cX}}{\scr{E}}{\cc}
remains exact when restricted to $\cX_b$ for all $b\in B$ therefore, by \cite{HuyLeh}[Lem. 2.1.4] we know that $\cc$ is flat over $B$.
Moreover, by \thref{sheafcorrispondenza} we have $\cc_b\cong \id_{Y_b}(\cd_b)$ for some $Y_b\subset \cX_b$, in particular is torsion-free, hence $\cc$ itself is torsion-free.
It follows that we have an injection $\cc\hookrightarrow\cc^{\vee\vee}$ and this last sheaf is locally free by \cite{Kollar}[Lem. 6.13].
We deduce that $\cc$ is a twist of some ideal sheaf but, being $\cc_b\cong \id_{Y_b}(\cd_b)$ for all $b\in B$, we must have $\cc\cong \id_{\cY}(\cd)$ where $\cY\subset \cX$ is such that $\cY_b\cong Y_b$.
$\cY$ is flat over $B$ by \cite{HuyLeh}[Lem. 2.1.4] because $\os_\cY$ is the cokernel of the morphism 
\[\cc(-\cd)\cong\id_\cY\hookrightarrow\os_{\cX},\]
which stays injective when restricted to $\cX_b$ for all $b$.
$\cY$ is proper over $B$ since it is closed in $\cX$ and $\Pi:\cX\ra B$ is proper.
Finally, $\cY$ is of pure codimension $2$ by \cite{Har}[III Cor. 9.6] since $\cY_b$ enjoy this property for all $b\in B$. 
Taking the dual of \eqref{fcorrispondenza}, the map $\scr{E}\du\ra\cV\du\otimes\os_{\cX}$ has to be zero, since its restrictions to $\cX_b$ is zero for all $b\in B$ by \thref{sheafcorrispondenza}, hence $\cW\cong \cV\du$ injects in $Ext^1(\id_{\cY}(\cd),\os_{\cX})$.

Let us start from $(\cY,\cW)$.
First of all, by \cite{Har}[III Cor. 9.6], $\cY_b$ is also purely of codimension $2$ for all $b$.
$\cY$ flat over $B$ means that $\os_{\cY}$ is flat over $B$, in particular restricting the sequence
\ses{\id_{\cY}}{\os_{\cX}}{\os_{\cY}}
to $\cX_b$ it remains exact and then $(\id_{\cY})_b\cong\id_{\cY_b}$.
We can define $\scr{E}$ by the sequence \eqref{fcorrispondenza} and $\cV:=\cW\du$.
Being $\cY$ flat over $B$ also $\id_\cY$ is such, being kernel of a morphism between flat sheaves.
Also $\scr{E}$ is flat and torsion-free, being extension of flat torsion-free sheaves.
Since $\id_\cY(\cd)$ is flat, the restriction of $\eqref{fcorrispondenza}$ to $\cX_b$ stays exact for all $b\in B$.
This implies that $\scr{E}_b$ is torsion-free being extension of such sheaves.
Finally, $\scr{E}_b$ has the correct determinant and rank, and no trivial direct summand by \thref{sheafcorrispondenza}.
\end{proof}

The difficulty in applying such a result consist in finding $\cV,\cW$ over the all $B$.
Actually, for our aim, that is proving that the Hartshorne-Serre correspondence is algebraic, is enough that such $\cV$ or $\cW$ exist locally around any point of $B$.
In \cite{Lange}[Prop. 2.2 and Prop. 2.3] the existence of a local/global extension group $\cW$ is discussed cohomologically.

Note that, by semicontinuity, we can always assume that $ext^1(\id_{\cY_b}(\cd_b),\os_{\cX_b})$ to be constant over $B$.
We will work out one extremal case, which will be the one we will need in applications: that is the construction of vector bundle of rank equal to $ext^1(\id_{\cY_b}(\cd_b),\os_{\cX_b})+1$.
In general, to construct lower rank ones we would have to consider a Grassmannian bundle parametrising pairs $(\cY_b,W_b)$ for some $W_b\subset Ext^1(\id_{\cY_b}(\cd_b),\os_{\cX_b})$ of fixed dimension.
For any morphism $V\ra B$ we set $\cX_V:=\cX\times_B V$.

\Prop\label{modularcorrispondenza}
Let $\Pi:\cX\ra B$ be a smooth, proper morphism of relative dimension $n\ge 2$ over the quasi-projective variety $B$ over $\gk=\overline{\gk}$.
Suppose $\cY$ is a family of closed codimension $2$ subschemes of $\cX$ flat over $B$ and assume there exists some divisor $\cd$ on $\cX$ such that $ext^1(\id_{\cY_b}(\cd_b),\os_{\cX_b})=r-1$ is constant for all $b\in B$.
Then there is an affine open cover $U\ra B$ and a rank $r$ torsion-free sheaf $\scr{E}$ on $\cX_U$, such that $\scr{E}_b$ corresponds to $(\cY_b,Ext^1(\id_{\cY_b}(\cd_b),\os_{\cX_b}))$ through Hartshorne-Serre for all $b\in U$.
In particular, the sheaf $\scr{E}_b$ is locally-free for those $b\in U$ over which $\cY_b$ is Cohen-Macaulay and the evaluation morphism \[d\left(Ext^1(\id_{\cY_b}(\cd_b),\os_{\cX_b})\right)\otimes\os_{\cX_b} \ra \omega_{\cY_b}(-D-K_X)\] 
is surjective.  
\One
\begin{proof}
Consider the relative ext sheaf $\ext^1_{\Pi}(\id_\cY(\cd),\os_{\cX})$ on $B$ obtained by sheafification of the pre-sheaf $V\mt Ext^1(\id_\cY(\cd)|_{\cX_V},\os_{\cX_V})$.
By \cite{BanicaPutinarSchumacher}[Satz 3 (ii)]\footnote{Here it is proved over the complex numbers but works over any field by \cite{Lange}[Thm. 1.4]}, since $ext^1(\id_{\cY_b}(\cd_b),\os_{\cX_b})=r-1$ is constant over the reduced scheme $B$, we get that $\ext^1_{\Pi}(\id_\cY(\cd),\os_{\cX})$ is locally free of rank $r-1$ on $B$.
We define $\{U_i\}$ to be some affine open cover of $B$ trivialising this vector bundle.
Then, for all $i$ we also have
\[H^0(U_i,\os_{U_i})^{r-1}\cong H^0(U_i,\ext^1_{\Pi}(\id_\cY(\cd),\os_{\cX})|_{U_i})\cong Ext^1(\id_\cY(\cd)|_{X_{U_i}},\os_{X_{U_i}})\]
by \cite{Lange}[(2) p.102]\footnote{the proof is similar to \cite{Har}[III Prop. 8.5]}.
Define $U$ to be the disjoint union of the $U_i$-s, hence $X_U$ is the disjoint union of $X_{U_i}$ therefore, we can construct our family $\scr{E}$ separately on the $X_{U_i}$.
    We can apply \thref{familycorrispondenza} taking $\cW_i:= Ext^1(\id_\cY(\cd)|_{X_{U_i}},\os_{X_{U_i}})$.
    Since $(\cW_i)_b=ext^1(\id_{\cY_b}(\cd_b),\os_{X_b})$ by construction, \thref{sheafcorrispondenza} tells us that $\scr{E}_b$ corresponds to $(\cY_b,ext^1(\id_{\cY_b}(\cd_b),\os_{X_b}))$ and is a torsion-free sheaf of rank $r$.
    Moreover, by \thref{vectorcorrispondenza}, we deduce the criterion to have $\scr{E}_b$ locally free.
\end{proof}

Note that, by flatness of $\scr{E}$ and irreducibility of $B$, if we have a line bundle $\os_{\cX}(1)$ on $\cX$ relatively ample to $\Pi$ then all the sheaves $\scr{E}_b$ would share the same Hilbert polynomial $P$.
At this point, we can form the \textit{relative moduli space of semistable sheaves} with this Hilbert polynomial as in \cite{HuyLeh}[Thm. 4.3.7]: this is a scheme $\mathcal{M}_{\cX/B}(P)$ with a projective morphism $p:\mathcal{M}_{\cX/B}(P)\ra B$ which is a coarse moduli space for the functor of flat families of semistable sheaves.
Moreover, for a closed point $b\in B$ we get $(\mathcal{M}_{\cX/B}(P))_b=\mathcal{M}_{\cX_b}(P)$, the corresponding moduli space of semistable sheaves on $\cX_b$ with Hilbert polynomial $P$ respect to the polarisation $\os_{\cX_b}(1)$.

\Co\label{sezionemodulare}
In the setting of \thref{modularcorrispondenza}, suppose that there is some line bundle $\os_\cX(1)$ relatively very ample respect to $\Pi$. 
If $\scr{E}_b$ is semi-stable for all $b\in B$ then there is an algebraic section $HS:B\ra \mathcal{M}_{\cX/B}(P)$ of the morphism $p:\mathcal{M}_{\cX/B}(P)\ra B$.
\io
\begin{proof}
    From \thref{sheafcorrispondenza} we have a set-theoretic map $HS:B\ra \mathcal{M}_{\cX/B}$ which composed with $p$ is the identity.
    The property of being algebraic is a local one, so can be checked on an affine cover. 
    By \thref{modularcorrispondenza} we can find an open cover $U\ra B$ and a semi-stable sheaf on $\cX_U$ flat over $U$ which, by universal property of $\mathcal{M}_{\cX/B}$, induces a morphism $U\ra \mathcal{M}_{\cX/B}$.
    This last one agrees with $HS$ on the closed points since given $\cY_b$ we have a uniquely determined $W_b\subset Ext^1(\id_{\cY_b}(\cd_b),\os_{\cX_b})$ of dimension $r-1$, by assumption in \thref{modularcorrispondenza}, and from $(\cY_b,W_b)$ we have a unique sheaf by \thref{sheafcorrispondenza}.
    Therefore, $HS:B\ra \mathcal{M}_{\cX/B}$ is algebraic. 
\end{proof}

\Oss
In the above setting, fix $r\in \N$.
Since the dimension $n$ of $\cX_b$ and the degree $d$ of the polarisation $\os_{\cX_b}(1)$ are constant along $B$, Ulrich sheaves on $(\cX_b,\os_{X_b}(1))$ share the same Hilbert polynomial, see \thref{rhps}, and are semistable, see \thref{stability}, for all $b\in B$.
Therefore, they all fit in the moduli space $\mathcal{M}_{\cX/B}\left(rd\binom{t+n}{n}\right)$.
Moreover, by semicontinuity of cohomology, anytime we have a family $\scr{E}$ as in the proposition above, if there is some $b\in B$ such that $\scr{E}_b$ is Ulrich the same holds true on an open neighbourhood of $b$.
\one

\Co\label{mappamodulare}
 In the special case $\cX=X\times B$ and $\os_\cX(1)\cong \pi_1^*\os_X(1)$, where $\pi_1:X\times B\ra X$ is the projection, then we get a morphism $HS:B\ra \mathcal{M}_{X}(P)$.
Moreover, if two points $b,b'\in B$ have the same image $\e$ and it is a stable sheaf then $\cY_b$ and $\cY_b'$ are rationally equivalent on $X$.
\io
\begin{proof}
    In the case $\cX=X\times B$ and $\os_\cX(1)\cong \pi_1^*\os_X(1)$ then $\mathcal{M}_{\cX/B}\cong X\times \mathcal{M}_{X}$ hence composing $HS$ with the projection to the first factor we get the claim.
    Since the $S$-equivalence class of a stable sheaf $\e$ is the sheaf itself, if $b,b'$ hit a stable sheaf then they both come from the Hartshorne-Serre construction on $\e$, they just correspond to different subspaces $V\subset H^0(X,\e)$.
    From this and \thref{hsc2} we conclude.
\end{proof}

\section{Ulrich bundles and Hartshorne-Serre correspondence}

We apply the previous general framework to the theory of Ulrich bundles.
During the preparation of this work the preprint \cite{LopezRaysubv} came out.
Our main result in the Ulrich bundle case, see \thref{Ulrichcorrispondenza'}, is similar to Thm. 1 in the previous paper.
Let us just note that, while they focus on sorting out the "nicest" possible degeneracy loci from Ulrich bundles, our approach is aimed at giving a bijective correspondence between pairs of vector bundles with a space of sections and subvarieties with sections of a twist of their dualizing sheaf.

\subsection{Proving Ulrichness}

Our goal is to specialise Hartshorne-Serre's correspondence to Ulrich bundles.
As we have already seen how to construct vector bundles as extensions, in this section we will suppose their existence and only care about additional properties we want to put on those sheaves.
Later, we will "glue" everything in a complete statement.

The Ulrich condition on $\e$ implies restrictions on $Y$, the first of which regards the Hilbert polynomial.
The situation is more clear when $\os_X(i)$ have no intermediate cohomology for $-n\leq i\leq -1$, for example when $(X,\os_X(1))$ is aCM or $\gk=\C$, by Kodaira vanishing.
The following is the natural generalisation of \cite{CFK3}[Thm. 3.1]

\Le\label{verificaulrich}
Let $(X,H)$ be a $n$-dimensional polarised variety over some infinite field $\gk$ such that $h^j(X,\os_X(i))=0$ for all $j\leq n-1$ and $-n\leq i\leq -1$.
Given a Cartier divisor $D$ on $X$, consider a rank $r$ sheaf $\e$ fitting in 
\sesl{\os_X^{r-1}}{\e}{\id_Y(D)}{s'''}
with $Y$ of pure codimension $2$.
Then $\e$ is Ulrich if and only if 
\begin{itemize}
    \item $h^n(X,\e(-n))=0$
    \item $h^j(X,\id_Y(D+iH))=0$ for all $0\leq j\leq n-2$ and $-n\leq i\leq -1$
    \item $P(\id_Y(D))=rd\binom{t+n}{n}-(r-1)P(\os_X)$
\end{itemize}
\ma
\begin{proof}
    Clearly, if $\e$ is Ulrich then the first condition is satisfied and the second comes from \eqref{s'''} and the assumptions on $(X,H)$.
    By \thref{hilbertpol} also the one regarding the Hilbert polynomial is verified.

    On the other direction, we derive that $P(\e)=rd\binom{t+n}{n}$, which is the correct polynomial for such a bundle.
    Moreover, from the first item in the list it follows that $h^n(X,\e(i))=0$ for all $i\geq -n$.
    Using the cohomology of \eqref{s'''} and our vanishings we get $h^j(X,\e(i))=0$ for all $0\leq j\leq n-2$ and $-n\leq i\leq -1$.
    Finally, since $P(\e)(t)=\chi(\e(t))=\sum_{j=0}^n(-1)^jh^j(X,\e(t))$ vanishes for $t=-1,\dots, -n$, we conclude that $h^{n-1}(X,\e(i))=0$ for $-n\leq i\leq -1$ hence $\e$ is Ulrich.
\end{proof}

In the case $n=2$, the above criterion is even simpler.
\Co\label{verificaulrich2}
Let $(X,H)$ be a polarised integral surface over some infinite field $\gk$ such that $h^1(X,\os_X(-1))=0$.
Given a Cartier divisor $D$ on $X$, consider a rank $r$ sheaf $\e$ fitting in 
\ses{\os_X^{r-1}}{\e}{\id_Y(D)}
with $Y$ of dimension $0$.
Then $\e$ is Ulrich if and only if 
\begin{itemize}
    \item $h^2(X,\e(-2))=0$ 
    \item $h^0(X,\id_Y(D-H))=0$ 
    \item $length(Y)=P(\os_X(D))-rd\binom{t+n}{n}+(r-1)P(\os_X)$
\end{itemize}
\io
\begin{proof}
The first two requirements are left unchanged.
Note that 
    The last condition in this statement is equivalent to the last condition in the previous one, since $P(\os_Y)=length(Y)$ and we have the sequence
    \ses{\id_Y(D)}{\os_X(D)}{\os_Y}
\end{proof}

The condition $h^n(X,\e(-n))=0$ is quite annoying, so we will treat it separately taking inspiration from \cite{ArrMad}[Prop. 4.1] and \cite{CFK3}[Remark 3.6].
\Le\label{surj_totale}
Assume that $X$ is a smooth variety and we have
\sesl{W\du\otimes\os_X}{\e}{\id_Y(D)}{s''}
with $Y$ of pure codimension $2$.
Then $h^n(X,\e(-n))=0$ is equivalent to each of the following:
\begin{enumerate}[1)]
    \item the boundary map $H^{n-1}(X,\id_Y(D-nH))\ra H^n(X,W\du\otimes\os_X(-nH))$ 
    is surjective and $h^0(X,\os_X(K_X+nH-D))=0$ 
    \item the multiplication map $d(W)\otimes H^0(X,\os_X(K_X+nH))\ra H^0(Y,\omega_Y(nH-D))$\footnote{note that by \eqref{boundary} we have $d(W)\subset H^0(Y,\omega_Y(-D-K_X))$} is injective and $h^0(X,\os_X(K_X+nH-D))=0$, assuming we know $h^1(X,\os_X(K_X+nH-D))=0$ 
    \item $h^0(X,\e\du(nH+K_X))=0$, assuming $\e$ is locally free.
\end{enumerate}
\ma
\begin{proof}
    \gel{1}
    Consider the long exact sequence in cohomology of \eqref{s''}, then $h^n(X,\e(-n))=0$ if and only if $H^{n-1}(X,\id_Y(D-nH))\ra H^n(X,W\du\otimes\os_X(-nH))$ 
    is surjective and $h^n(X,\id_Y(D-nH))=0$.
    Using the cohomology of 
\ses{\id_Y(D-nH)}{\os_X(D-nH)}{\os_Y(D-nH)}
we get 
\[h^n(X,\id_Y(D-nH))=h^n(X,\os_X(D-nH))=h^0(X,\os_X(K_X+nH-D))\]
where in the last step we used Serre duality.

    \gel{2}
    Again by \eqref{s''}, the surjectivity in the previous condition is equivalent to $h^n(X,\e(-n))=h^n(X,\os_X(D-nH))$.
    By Serre duality this is equivalent to $h^0(X,\e\du(K_X+nH))=h^0(X,\os_X(K_X+nH-D))$.
    Dualizing \eqref{s''}, recall \eqref{contoduale'}, and twisting by $\os_X(K_X+nH)$ we get 
    \begin{equation*}
    0\ra \os_X(K_X+nH-D)\ra\e\du(K_X+nH)\ra W\otimes\omega_X(nH)\ra \omega_Y(nH-D)\ra 0
\end{equation*}
We can split this in two short exact sequences.
Using that $h^1(X,\os_X(K_X+nH-D))=0$, we get $0=h^0(X,\e\du(K_X+nH))=h^0(X,\os_X(K_X+nH-D))$ if and only if $d(W)\otimes H^0(X,\omega_X(nH))\ra H^0(Y,\omega_Y(nH-D))$ is injective.

\gel{3}
By Serre duality we immediately get the claim.
\end{proof}

\Co\label{corsurj}
In setting of \thref{surj_totale}, the condition $h^n(X,\e(-n))=0$ is always satisfied if $\e$ is locally free, self U-dual and $h^0(X,\e(-1))=0$.
\io
\begin{proof}
    By the last statement in \thref{surj_totale} we need to verify that $0=h^0(X,\e\du(nH+K_X))=h^0(X,\e\du((n+1)H+K_X)-H))$.
    But, if $\e$ is self U-dual we have $\e\cong \e\du((n+1)H+K_X)$ hence we conclude being $h^0(X,\e(-1))=0$.
\end{proof}

\subsection{Hartshorne-Serre for Ulrich bundles}

Now, we put together the results from the second and the last subsection to obtain an Hartshorne-Serre correspondence for Ulrich bundles.

\Te\label{Ulrichcorrispondenza}
Let $(X,H)$ be a smooth polarised variety of dimension $n\geq 2$ over some algebraically closed field and let $D$ be a divisor.
Suppose that $(X,H)\neq (\p^n,\osn(1))$ and $h^j(X,\os_X(i))=0$ for $-n\leq i\leq -1$ and $1\leq j\leq n-1$. 
There is a bijective correspondence between 
\begin{itemize}
\item isomorphism classes of pairs $(\e,V)$ such that
\begin{enumerate}[1)]
    \item $\e$ is a vector bundle of rank $r>1$ and $det(\e)\cong\os_X(D)$
    \item $V\subset H^0(X,\e)$ has dimension $r-1$, the morphisms $V\otimes\os_X\ra \e$ is injective and drops rank on a locus of pure codimension $2$
    \item $\e$ is Ulrich
\end{enumerate}
\item pairs $(Y,W)$  such that
\begin{enumerate}[i)]
\item\label{II} $Y\subset X$ is a Cohen-Macaulay closed subscheme of pure codimension $2$ 
    \item\label{III} $W\subset Ext^1(\id_Y(D),\os_X)$ has dimension $r-1$ and the evaluation morphism $d(W)\otimes\os_X \ra \omega_Y(-D-K_X)$ is surjective
    \item\label{I}
    \begin{equation}\label{mid1}
        h^j(X,\id_Y(D+iH))=0 \qquad -n\leq i\leq -1,\; 0\leq j\leq n-2
    \end{equation}
    \begin{equation}\label{i1}
        P(\id_Y(D))=rd\binom{t+n}{n}-(r-1)P(\os_X)
    \end{equation}
    \begin{equation}\label{top1}
        h^0(X,\os_X(K_X+nH-D))=0
    \end{equation}
    \begin{equation}\label{surj1}
       H^{n-1}(X,\id_Y(D-nH))\ra H^n(X,V\otimes\os_X(-nH)) 
    \end{equation}
    is surjective or (equivalently) injective
      \end{enumerate}
\end{itemize}
moreover they fit in the exact sequence
\begin{equation}\label{corrispondenza}
    0\ra  V\otimes\os_X\ra \e\ra \id_Y(D)\ra 0.
\end{equation}
\Ma
\begin{proof} 
We know that $\os_X$ is an Ulrich sheaf only if $(X,H)\cong (\p^n,\osn(1))$, see \thref{oulrich} $ii)$.
But if $\e$ is Ulrich and has a trivial summand this latter one must also be Ulrich hence $(X,H)\cong (\p^n,\osn(1))$ contradicting our hypothesis.
Then \thref{vectorcorrispondenza} tells us that 
\[1)+2)+3)\Longrightarrow i)+ii)\Longrightarrow 1)+2)\]
so that, to conclude, we can assume $1)$ and $2)$ and prove that under those assumptions
we have $3)$ if and only if $iii)$.

$\e$ is Ulrich if and only if satisfies the three conditions in \thref{verificaulrich} but the second two coincides exactly with \eqref{i1},\eqref{mid1}.
By the first claim in \thref{surj_totale}, the condition $h^n(X,\e(-n))=0$ in \thref{verificaulrich} is equivalent to \eqref{top1} and surjectivity in \eqref{surj1}.

We are only left to explain why injectivity and surjectivity in \eqref{surj1} are equivalent, if we assume the other facts in $iii)$.
It is enough to show that source and target of this map have the same dimension.
Note that by \eqref{i1} and the other vanishings we have $\chi(\id_Y(D-nH))=-\chi(\os_X^{r-1}(-n))$ but from  we know \[h^n(X,\os_X^{r-1}(-n))=(-1)^n\chi(\os_X^{r-1}(-n))=\]
\[=(-1)^{n-1} \chi(\id_Y(D-nH))=h^{n-1}(X,\id_Y(D-nH)).\]
\end{proof}

The following is a list of remarks about the above construction.

\Oss\label{r1}
\begin{enumerate}[a)]

\item An Ulrich bundle is globally generated by \thref{corollario}, so if we are over $\C$ by \cite{Banica}[Thm. 1] or \cite{Ott}[Thm. 2.8] for a generic $V\subset H^0(X,\e)$ of dimension $r-1$ we get a morphism whose degeneracy locus $Y$ has pure codimension $2$ or is empty.
In particular, each Ulrich bundle on a smooth complex variety such that $c_2(\e)\neq 0$ can be constructed starting from $Y$ as in the hypothesis, recalling \thref{hsc2}.
Note that the Ulrich bundles with $c_2(\e)\sim 0$ are characterised in \thref{c2=0}.

\item Again working over $\C$, if $dim(X)=2,3$ then we can choose $Y$ in the previous theorem to be smooth, this follows again from \cite{Banica}[Thm. 1] or \cite{Ott}[Thm. 2.8].
    Moreover, we can also choose the divisor $D$ containing $Y$ to be smooth: the degeneracy locus of the span of $V$ and a sufficiently general section of $\e$  will be a smooth divisor containing $Y$, again by the above references.
    In higher dimension, we still have some control on the ineluctable singularities, see \cite{Banica}.

\item If $K_X=kH$ for some $k\in\Z$ then $W\cong Ext^1(\id_Y(D),\os_X)$, hence $Y$ alone determines the pair $(\e,V)$, extending the conclusion of \cite{CFK2}[Remark 3.6].
        We already noted, see the line below \eqref{long}, that the map $\phi:H^0(X,\e\du)\ra W$ is zero but, with our new assumption we can argue that
        \[ext^1(\e,\os_X)\cong h^1(X,\e\du)\cong h^{n-1}(X,\e(k))=0\]
        because an Ulrich bundle has no intermediate cohomology by \thref{corollario}.
        Then \eqref{long} gives the thesis.
\end{enumerate}
\one

We will state another version of \thref{Ulrichcorrispondenza}, where we assume some extra but natural hypothesis, to get a simpler criterion for the existence of special Ulrich bundles.
\thref{Ulrichcorrispondenza} depends on a choice of $D$ and the condition that $\e$ is Ulrich tells us the degree $DH^{n-1}$ of this divisor, see \thref{c12}, but when $Pic(X)\neq \Z$ does not identify uniquely $D$ a priori.
Nonetheless, we always have an interesting choice for $D$: we can choose special Ulrich bundles i.e. $D\sim dH$ for $d=\dfrac{r(n+1+k)}{2}$.

\Te\label{Ulrichcorrispondenza'}
Let $(X,H)$ be a smooth polarized variety of dimension $n\geq 3$ over $\C$.
Moreover, suppose that $(X,H)\neq (\p^n,\osn(1))$ and $K_X\sim kH$ for some $k\in\Z$.
There is a bijective correspondence between 
\begin{itemize}
\item isomorphism classes of pairs $(\e,V)$ such that
\begin{enumerate}[1)]
    \item $\e$ is a vector bundle of rank $r>1$
    \item $V\subset H^0(X,\e)$ has dimension $r-1$, the morphisms $V\otimes\os_X\ra \e$ is injective and drops rank on a locus of pure codimension $2$
    \item $\e$ is Ulrich special, i.e. Ulrich such that $det(\e)\cong\os_X(d)$
    with $d=\dfrac{r(n+1+k)}{2}$
\end{enumerate}
\item $Y$ such that
\begin{enumerate}[i)]
\item\label{II'} $Y\subset X$ is a connected, Cohen-Macaulay subscheme of pure codimension $2$ 
\item $\omega_Y(-k-d)$ is globally generated and $h^0(Y,\omega_Y(-k-d))=r-1$
    \item 
    \begin{equation}\label{mid1'}
        h^j(X,\id_Y(D-iH))=0 \qquad 1\leq i\leq n,\; 1\leq j\leq n-2
    \end{equation}
    \begin{equation}\label{i1'}
        P(\id_Y(D))=rd\binom{t+n}{n}-(r-1)P(\os_X)
    \end{equation}
    \begin{equation}\label{top1'}
        h^0(X,\id_Y(D-H))=0
    \end{equation}
   \begin{equation}\label{surj'}
       H^0(Y,\omega_Y(-d-k))\otimes H^0(X,\os_X(n+k))\ra H^0(Y,\omega_Y(n-d))
   \end{equation}
    is surjective or (equivalently) injective.
    \end{enumerate}
\end{itemize}
moreover they fit in the exact sequence
\begin{equation}\label{ucorr'}
    0\ra  V\otimes\os_X\ra \e\ra \id_Y(d)\ra 0.
\end{equation}
\Ma
\begin{proof}
First of all, working over $\C$ and being $H$ ample, Kodaira vanishing gives $h^j(X,\os_X(-i))=0$ for $1\leq i\leq n$ and $1\leq j\leq n-1$, hence we can apply \thref{Ulrichcorrispondenza}.
The hypothesis regarding $(\e,V)$ are unchanged, except that we fix a particular determinant, so we need only to prove that, with our new assumptions, the ones regarding $Y$ are equivalent to those given for $(Y,W)$ in \thref{Ulrichcorrispondenza}.

In both cases $Y$ is Cohen-Macaulay of pure codimension $2$ and  $\omega_Y(-d-k)$ is globally generated.
Being $X$ smooth and $(X,H)\neq (\p^n,\osn(1))$, from \cite{KacKol}[Cor. 2] we get  $k\geq -n$.
It follows that $D\sim dH$, with $d=\dfrac{r(n+1+k)}{2}$, is a positive multiple of $H$, hence an ample divisor.
We can then apply Kodaira vanishing and, being $n\geq 3$, we have $h^i(X,\os_X(-D))=0$ for $i=1,2$ hence as seen in \thref{r00} we obtain 
\[Ext^1(\id_Y(D),\os_X)\cong H^0(Y,\omega_Y(-D-K_X)).\]
We define $W$ to be this $r-1$-dimensional vector space, then we have our pair $(Y,W)$ as in \thref{Ulrichcorrispondenza}.

Notice that condition \eqref{mid1}, \eqref{i1} and the first part of \eqref{top1} are exactly \eqref{mid1'},\eqref{i1'} and \eqref{top1'}.
Moreover, \eqref{top1} is equivalent to $K_X+nH-D=:-A$ being non effective; we will show a bit more: $A$ is ample.
Indeed, it is enough to show that $-A$ is linearly equivalent to a negative multiple of $H$, this means to show
\[0>k+n-d= k+n-\dfrac{r(n+1+k)}{2}\]
which is always satisfied because is true for $r=2$ and the expression on the right is decreasing in $r$ for $r\geq 0$ since $n\geq -k$.
In particular, $K_X+pH-D=-A-(n-p)H$ when $p\leq n$ is anti-ample, so by Kodaira vanishing, being $n\geq 3$, we get $h^1(X,\os_X(k+n-d))=0$.
It follows from \thref{surj_totale} that the condition \eqref{surj1} is equivalent to $d(W)\otimes H^0(X,\os_X(K_X+nH))\ra H^0(Y,\omega_Y(nH-D))$ is surjective.
Moreover, the above isomorphism $W\cong H^0(Y,\omega_Y(-D-K_X))$ implies that this is equivalent to \eqref{surj'}, as desired.

Finally, we are left to prove that $Y$ is connected.
Being $X$ irreducible we have $h^0(X,\os_X)=\C$, so it is enough to prove $h^1(X,\id_Y)=0$ and substitute it in \eqref{subd'}.
The claimed vanishing follows from \eqref{corrispondenza} twisted by $\os_X(-d)$ once we recall that $\e$, being Ulrich has no intermediate cohomology, see \thref{corollario}, and that $h^2(X,\os_X(-d))=0$.
\end{proof}

We discuss another couple of remarks.
\Oss\label{r2}
\begin{itemize}
    \item If, in the context of \thref{Ulrichcorrispondenza'}, we are trying to construct a rank $2$ Ulrich bundle $\e$ then the condition \eqref{surj'} can be ignored.
In fact, by \thref{surj_totale}, this is equivalent to $h^n(X,\e(-n))=0$ which in turn, by the same result, is equivalent to $h^0(X,\e\du(nH+K_X))=0$.
But from \thref{specialdual} a special rank $2$ vector bundle is always U-dual, so by \thref{corsurj} this follows from $h^0(X,\e(-1))=0$, which is a consequence of \eqref{top1'}.

\item 
In the setting of \thref{Ulrichcorrispondenza'} we have
\[h^{n-2}(\os_Y)=h^{n-2}(\os_X)-h^{n-1}(\os_X)+h^{n}(\os_X)+(r-1)h^0(\omega_X(d))+(-1)^{n+1}\binom{n-d}{n}.\]
    Indeed, during the proof of \thref{Ulrichcorrispondenza'} we have seen if $(X,H)\neq (\p^n,\osn(1))$ then $d>0$ so, being $\e$ Ulrich, by \eqref{ucorr'} twisted by $\os_X(-d)$ we get $h^j(X,\id_Y)=0$ for $j\leq n-2$ and 
    \[h^{n-1}(X,\id_Y)=(r-1)h^n(X,\os_X(-d))-h^n(X,\e(-d))+h^{n}(X,\id_Y).\]
    Therefore, being $Y$ of codimension $2$, considering
    \ses{\id_Y}{\os_X}{\os_Y}
    we deduce $h^{n}(X,\id_Y)=h^{n}(X,\os_X)$ and then
    \[h^{n-2}(Y,\os_Y)=h^{n-2}(X,\os_X)+h^{n-1}(X,\id_Y)-h^{n-1}(X,\os_X)=\]
    \[=h^{n-2}(X,\os_X)+(r-1)h^n(X,\os_X(-d))-h^n(X,\e(-d))+h^{n}(X,\os_X)-h^{n-1}(X,\os_X)=\]
    \[=h^{n-2}(X,\os_X)-h^{n-1}(X,\os_X)+h^{n}(X,\os_X)+(r-1)h^0(X,\omega_X(d))+(-1)^{n+1}P(\e)(-d).\]
    where we rearranged and used Serre duality and the fact that $\e$ is Ulrich and $d>0$.
Finally, note that by \thref{hilbertpol}, we know $P(\e)(-d)=\binom{n-d}{n}$ hence we get the statement.
   \end{itemize}
\one

The next step would be to explicitly compute other properties of $Y$ coming from the fact that we know its Hilbert polynomial.
This has been done in \cite{CFK3}[Thm. 3.5] for $dim(X)=3$ and in general in \cite{LopezRaysubv}[Lem. 3.2].

\section{Extensions and deformations}

Here we want to prove some properties regarding extension of vector bundles in order to construct higher-rank Ulrich bundles from low-rank ones.
Such extensions are never stable since by \thref{rhps} an Ulrich sheaf is stable if and only if it has no Ulrich subsheaves.
Therefore, we investigate flat deformations of those extensions to search for stable ones.
In order to control local properties of the base scheme of such families, we will need to work with simple sheaves, as a replacement for stable ones which are not available.
This general strategy dates back at least to \cite{CH}, see also \cite{FaenziK3}, \cite{CFK3} and \cite{CFK3.1} for variations of it.

We try to make this procedure as general as possible, avoiding the hypothesis of the sheaves being Ulrich.
During the preparation of this work, the preprint \cite{ComaschiFaenzi} came out, where similar techniques have been used to deal with instanton bundles on Fano $3$-folds.

\subsection{Extensions}

We will suppose $X$ is a variety over an algebraically closed field $\mathbf{k}$. 
The rank of a sheaf is the rank at the generic point.
In this section, we will adapt some technical results from the literature useful in the next one.

The first and second among the claims below are slight generalisations of \cite{CFK3.1}[Lem. 6.3] while the third extends \cite{CH}[Lem. 4.2].
\Le\label{estensionislope}
Suppose on an $n$-dimensional, smooth polarised variety $(X,H)$ over $\gk=\overline{\gk}$ we have an exact sequence of sheaves
    \sesl{\e}{\f}{\g}{b'} 
    where $\e,\g$ are torsion-free, slope-stable with the same slope $\mu$.
    Then 
    \begin{enumerate}
    \item $\mu(\f)=\mu$ and $\f$ is slope-semistable
    \item for every torsion-free quotient sheaf $\f\twoheadrightarrow\f'$ such that $\mu(\f')=\mu(\f)$ we have $\f'$ isomorphic to one among $\e,\f,\g$ or $(\f')\du\cong\e\du$.
    Moreover, if we are in the case $\f'\cong \e$ then the sequence splits while in the case $(\f')\du\cong\e\du$ we cannot have $\e$ locally free.
    \item If the extension is non-split and there are no injective morphisms $\e\ra\g$ then $\f$ is simple.
    \end{enumerate}
\ma
\begin{proof}
\gel{1}
    First of all, assuming that $\e,\g$ are not both zero, we have
\[\mu(\f)=\dfrac{c_1(\f)H^{n-1}}{rk(\f)}=\dfrac{(c_1(\e)+c_1(\g))\cdot H^{n-1}}{rk(\e)+rk(\g)}=\dfrac{rk(\e)\cdot \mu(\e)+rk(\g)\cdot \mu(\g)}{rk(\e)+rk(\g)}=\mu.\]
Suppose $\f\twoheadrightarrow\f'$ with $0<rk(\f')<rk(\f)$.
By composition we have an induced morphism $\e\ra \f'$ and either this is non-zero or $\f\twoheadrightarrow\f'$ factors through $\g$ and hence we have $\g\twoheadrightarrow \f'$.
In both cases, by slope-semistability of $\e,\g$ we get $\mu(\f')\leq \mu=\mu(\f)$ hence $\f$ is semistable.

\gel{2}
    Let the quotient $\f\twoheadrightarrow\f'$ be torsion-free and such that $\mu(\f')=\mu(\f)$ and $rk(\f')>0$.
    We have the following diagram
    \dia
    0\ar[r] & \e\ar[r]\ar[d,twoheadrightarrow] & \f\ar[r]\ar[d,twoheadrightarrow] & \g\ar[r]\ar[d,twoheadrightarrow] & 0\\
0\ar[r] & \e' \ar[r] & \f' \ar[r] & \g' \ar[r] & 0\\
    \mma
    where $\e'$ is the image of $\e$ in $\f'$ and $\g'$ the cokernel of the inclusion $\e'\ra\f'$.
    Note that by this diagram we have $\e'$ torsion-free, since it is a subsheaf of $\f'$ which is such, and $0\leq rk(\e')\leq rk(\e)$.
    If $rk(\e')=0$ then $\e'=0$ and we have a surjection $\g\twoheadrightarrow\g'\cong\f'\neq 0$ so by slope-stability of $\g$ it must be an isomorphism.
    If $0<rk(\e')<rk(\e)$ then by slope-stability $\mu(\e')>\mu(\e)=\mu=\mu(\f')$ hence we get $0<\mu(\g')<\mu(\f')=\mu$, in particular $\g'\neq 0$, but this would contradict slope-stability of $\g'$.
    
    If $rk(\e')=rk(\e)$ then $\e\cong\e'$, because the kernel of the above morphism $\e\ra\e'$ would have rank $0$ but $\e$ is torsion-free.
    Now, if $\g'=0$ then we get $\f'\cong\e$ and the sequence splits.
    Otherwise, first suppose $rk(\g')>0$ hence $\mu(\g')=\mu$, since $\mu(\f')=\mu=\mu(\e')$, but then slope-stability of $\g$ implies $\g\cong\g'$ so that $\f\cong\f'$. 
    Finally, if $rk(\g')=0$ then $rk(\f')=rk(\e')$ which implies 
    \[c_1(\g')=c_1(\f')-c_1(\e')=\dfrac{rk(\f')}{H^n}(\mu(\f')-\mu(\e'))=0\]
    so that $codim(Supp(\g'))\ge 2$.
    By \cite{HuyLeh}[Prop. 1.1.6 i)] we get $\ext^i(\g',\os_X)=0$ for $i=0,1$ hence dualising the bottom row in the above diagram we get $(\f')\du\cong(\e')\du\cong \e\du$. 
    Note that, if $\e$ is locally free by Serre duality
    \[ext^1(\g',\e')=ext^1(\g',\e)=ext^1(\e\du\otimes\g',\os_X)=h^{n-1}(X,\e\du\otimes\g'(K_X))=0\]
    by dimensional reasons because $\e\du\otimes\g'(K_X)$ is supported on $Supp(\g)$ which has codimension at least $2$.
    Therefore, the lower sequence has to split but this contradicts the fact that $\f'$ is torsion-free.
    The same results follows alternatively by a depth computation.
    
    \gel{3}
Suppose we have a non-zero endomorphism $\alpha:\f\ra \f$, we will show that it is an isomorphism.
Being $\e,\g$ slope-stable and with the same slope, by \thref{stabilitàtrick} the composition 
\[\e\hookrightarrow\f\xra{\alpha}\f\ra \g\] 
must be the zero map or an injection such that $rk(\e)=rk(\g)$. 
But the latter case cannot happen by assumption.
It follows that $\alpha(\e)\subseteq ker(\f\ra\g)=\e$, hence $\alpha$ induces a map $\til{\alpha}:\g\ra\g$ and a diagram 
\dia
0 \ar[r] & \e\ar[r] \ar[d,"\alpha|_{\e}"] & \f\ar[r] \ar[d,"\alpha"] & \g\ar[r] \ar[d,"\til{\alpha}"] & 0 \\
0 \ar[r] & \e\ar[r] & \f\ar[r] & \g\ar[r] & 0 \\
\mma
Our goal is to show that  $\alpha|_{\e}, \Tilde{\alpha}$ are both isomorphisms.
Since $\e,\g$ are slope-stable, hence simple, $\alpha|_{\e}, \Tilde{\alpha}$ must be either zero maps or isomorphisms.
We cannot have $\Tilde{\alpha}=0=\alpha|_\e$ otherwise $\alpha=0$.
If $\Tilde{\alpha}=0$ and $\alpha|_\e$ is an isomorphism then $\alpha:\f\ra \e$ gives a splitting of the sequence.
If $\alpha|_{\e}=0$ and $\Tilde{\alpha}$ is an isomorphism then $\alpha:\g\ra \f$ gives a splitting of the sequence.
Therefore, we are only left with the case in which both $\Tilde{\alpha},\alpha|_\e$ are isomorphism, so that also $\alpha$ is.

We have proved that each element in $End(\f)$ is either $0$ or invertible, hence $End(\f)$ endowed with the operations of sum and composition is a division algebra over the base field $\mathbf{k}$.
Moreover $End(\f)=H^0(\hom(\f,\f))$ is finite dimensional hence, being $\mathbf{k}$ algebraically closed, we must have $End(\f)\cong k$ by \cite{Cohn}[Prop. 5.4.5], i.e. $\f$ is simple. 
\end{proof}

\Oss
In $3)$ we assumed that there is no injective morphism $\e\ra\g$ while in \cite{CH}[Lem. 4.2] there is the weaker hypothesis $\e\neq \g$ because they work with $\e,\g$ locally free.
In our case, we want to deal with sheaves that are not both necessarily locally free, hence we have to assume something in return.
For our applications, we will have either 
$rk(\e)>rk(\g)$ or $rk(\e)=rk(\g)$ and $p(\e)=p(\g)$, so the absence of injective morphisms $\e\ra\g$ is clear.
\one

Our next lemma, and its proof, is parallel to the previous one, we just use Gieseker-stability instead of slope-stability.
The second and third claims are variations of, respectively, \cite{CFK3.1}[Lem. 6.3] and \cite{CH}[Lem. 4.2].
\Le\label{estensionistabili}
Suppose on an $n$-dimensional, smooth polarised variety $(X,H)$ over $\gk=\overline{\gk}$ we have an exact sequence of sheaves
    \sesl{\e}{\f}{\g}{b} 
    where $\e,\g$ are torsion-free, stable with the same reduced Hilbert polynomial $p$.
    Then 
    \begin{enumerate}
    \item $p(\f)=p$ and $\f$ is semistable
    \item for every quotient $\f\twoheadrightarrow\f'$ such that $p(\f')=p(\f)$ we have $\f'$ isomorphic to one among $\e,\f,\g$.
    Moreover, if we are in the case $\f'\cong \e$ then the sequence splits.
    \item If the extension is non-split and $\e,\g$ are non-isomorphic then $\f$ is simple.
    \end{enumerate}
\ma
\begin{proof}
\gel{1}
   First of all, assuming that $\e,\g$ are not both zero, we have
\[p(\f)=\dfrac{P(\f)}{rk(\f)}=\dfrac{P(\e)+P(\g)}{rk(\e)+rk(\g)}=\dfrac{rk(\e)\cdot p(\e)+rk(\g)\cdot p(\g)}{rk(\e)+rk(\g)}=p.\]
Suppose $0\neq\f'\subset\f$ then either we have a non-zero morphism $\f'\ra \g$ or $\f'\subset \e$.
In both cases, by semistability of $\e,\g$ we get $p(\f')\leq p=p(\f)$ hence $\f$ is semistable.

\gel{2}
    Let there be a quotient $\f\twoheadrightarrow\f'\neq 0$ such that $p(\f')=p(\f)$.
    As already seen at the end of the proof of \thref{crit_stab}, we must have $\f'$ torsion-free otherwise projecting to the torsion-free part we would contradict semi-stability of $\f$.
    We have the following diagram
    \dia
    0\ar[r] & \e\ar[r]\ar[d,twoheadrightarrow] & \f\ar[r]\ar[d,twoheadrightarrow] & \g\ar[r]\ar[d,twoheadrightarrow] & 0\\
0\ar[r] & \e' \ar[r] & \f' \ar[r] & \g' \ar[r] & 0\\
    \mma
    where $\e'$ is the image of $\e$ in $\f'$ and $\g'$ the cokernel of the inclusion $\e'\ra\f'$.
    Note that by this diagram we have $\e'$ torsion-free and $0\leq rk(\e')\leq rk(\e)$.
    If $rk(\e')=0$ then $\e'=0$ and we have a surjection $\g\twoheadrightarrow\g'\cong\f'\neq 0$ so by stability of $\g$ it must be an isomorphism since $p(\f')=p=p(\g)$.
    If $0<rk(\e')<rk(\e)$ then by stability $p(\e')>p(\e)=p=p(\f')$ hence we get $0<p(\g')<p(\f')=p$, in particular $\g'\neq 0$, but this would contradict slope-stability of $\g'$.
    
    If $rk(\e')=rk(\e)$ then $\e\cong\e'$, because they are torsion-free.
    Now, $rk(\e')=rk(\f')$ then $p(\e')=p(\e)=p=p(\f')$ would imply $P(\e')=P(\f')$ hence $\g'=0$, $\f'\cong\e$ and the sequence splits.
    If $rk(\e')<rk(\f')$ then, since $p(\f')=p=p(\e')$, we must have 
    \[P(\g')=P(\f')-P(\e')=p(rk(\f')-rk(\e'))=rk(\g)p\]
    hence $p(\g')=p$, but then stability of $\g$ implies $\g\cong\g'$ so that $\f\cong\f'$. 
    
    \gel{3}
Suppose we have a non-zero endomorphism $\alpha:\f\ra \f$, we will show that it is an isomorphism.
Being $\e,\g$ stable and with the same reduced Hilbert polynomial, the composition 
\[\e\hookrightarrow\f\xra{\alpha}\f\ra \g\] 
must be the zero map or an isomorphism, by \autoref{stabilitàtrick}, but the latter case cannot happen by assumption.
It follows that $\alpha(\e)\subseteq ker(\f\ra\g)=\e$, hence $\alpha$ induces a map $\til{\alpha}:\g\ra\g$ and a diagram 
\dia
0 \ar[r] & \e\ar[r] \ar[d,"\alpha|_{\e}"] & \f\ar[r] \ar[d,"\alpha"] & \g\ar[r] \ar[d,"\til{\alpha}"] & 0 \\
0 \ar[r] & \e\ar[r] & \f\ar[r] & \g\ar[r] & 0 \\
\mma
Our goal is to show that  $\alpha|_{\e}, \Tilde{\alpha}$ are both isomorphisms.
Since $\e,\g$ are slope-stable, hence simple, $\alpha|_{\e}, \Tilde{\alpha}$ must be either zero maps or isomorphisms.
We cannot have $\Tilde{\alpha}=0=\alpha|_\e$ otherwise $\alpha=0$.
If $\Tilde{\alpha}=0$ and $\alpha|_\e$ is an isomorphism then $\alpha:\f\ra \e$ gives a splitting of the sequence.
If $\alpha|_{\e}=0$ and $\Tilde{\alpha}$ is an isomorphism then $\alpha:\g\ra \f$ gives a splitting of the sequence.
Therefore, we are only left with the case in which both $\Tilde{\alpha},\alpha|_\e$ are isomorphism, so that also $\alpha$ is.
We have proved that each element in $End(\e)$ is either $0$ or invertible, so we conclude as in \thref{estensionislope}. 
\end{proof}

\Oss
Results analogous to \thref{estensionislope} and \thref{estensionistabili} hold for torsion-free quotients replaced by saturated subobjects.
\one

\subsection{Deformations}

Our goal is to rewrite in general the strategy used in \cite{CH}[Thm. 5.7], \cite{BN}[Thm. 3.2], \cite{CFK3}[Prop. 4.7] and \cite{CFK3.1}[Prop. 6.4] to produce families of stable bundles starting from low-rank ones.
For simplicity, we introduce some terminology to summarise the properties we will need.
Most of the statements are quite intricate, but the generality in which we work will subsume most of the analogous proofs in the literature (from which our arguments are derived).
The main technical tool will be modular families of simple sheaves, see \autoref{moduli}.

We remark that, if we are only interested in showing \textit{Ulrich wildness}, which guarantees the existence of families of arbitrary large dimension of indecomposable non-isomorphic Ulrich bundles, then there is a much easier criterion in \cite{FaenziPons_acm}[Thm. A, Cor. 2.1]. 

\No
Let $X$ be a projective scheme and $\scr{E},\scr{G}$ be two families of coherent sheaves on $X$ flat over integral base schemes $\cb,\cd$.
The families $\scr{E}\boxtimes \os_\cd$ and $\os_\cb\boxtimes\scr{G}$ are flat over $\cb\times\cd$, which is again integral.
Since Euler product is constant in families, see \cite{BanicaPutinarSchumacher}[Satz 3 iii)] we can define the Euler product for families
\[\underline{\chi}(\scr{E},\scr{G}):=\chi(\scr{E}\boxtimes \os_\cd|_{X_{(b,d)}},\os_\cb\boxtimes\scr{G}|_{X_{(b,d)}})=\chi(\scr{E}_b,\scr{G}_d) \quad \text{for some} \; (b,d)\in \cb\times\cd.\]
Being $\cb,\cd$ noetherian, and hence compact, also the following is well defined 
\[\mathfrak{ext}^i(\scr{E},\scr{G}):=\min_{(b,d)\in \cb\times\cd}\{ext^i(\scr{E}_{b},\scr{G}_{d})\}=ext^i(\scr{E}_{b},\scr{G}_{d})\quad (b,d)\in\cb\times\cd\; \text{general enough}\]
where the second equality follows by semicontinuity \cite{BanicaPutinarSchumacher}[Satz 3 i)].
Note that we have
\[\underline{\chi}(\scr{E},\scr{G})=\sum_{i\in\N}\mathfrak{ext}^i(\scr{E},\scr{G}).\]
\ione

Recall that all the sheaves in a flat family share the same Hilbert polynomial and rank, in particular, the same reduced Hilbert polynomial and slope.
Intuitively, we want to deal with well-behaved families of simple sheaves of the maximum possible dimension (in a reasonable sense).

\De\label{defwildext}
Let $(X,H)$ be a polarised scheme over an algebraically closed field $\gk$.
We call a family of sheaves $\scr{E}$ on $X$ flat over $\cb$ \textbf{nice family} if:
\begin{itemize}
    \item $\scr{E}_b$ is torsion free and simple for all $b$, and is stable for general $b$
    \item $ext^2(\scr{E}_b,\scr{E}_b)=0$ for all  $b\in \cb$
    \item $\cb$ is a (smooth) quasi-projective variety which is open in some modular family of simple sheaves on $X$.
\end{itemize}

\noindent We call a pair of families of sheaves $(\scr{E},\scr{G})$ on $X$ flat over $\cb,\cd$ of \textbf{good extension} if:
\begin{itemize}
    \item $\scr{E},\scr{G}$ are nice families and $p(\scr{E}_b)=p(\scr{G}_d)$ for all $b,d$
    \item for $(b,d)\in \cb\times\cd$ general we have $\scr{E}_b\neq\scr{G}_d$.
    \item $\mathfrak{ext}^2(\scr{E},\scr{G})=0=\mathfrak{ext}^2(\scr{G},\scr{E})$ while $\mathfrak{ext}^1(\scr{G},\scr{E})>0$ and $\mathfrak{ext}^1(\scr{E},\scr{G})>0$
\end{itemize}
We call $(\scr{E},\scr{G})$ of \textbf{wild extension} if it is of good extension and
\[\mathfrak{ext}^{1}(\scr{F}_1,\scr{F}_2)>\sum_{i\in\N}\mathfrak{ext}^{2i}(\scr{F}_1,\scr{F}_2) \qquad (\scr{F}_1,\scr{F}_2)=(\scr{G},\scr{G}),(\scr{E},\scr{G}),(\scr{G},\scr{E}).\]
\Ne

\Oss
In \thref{defwildext} notice that, being a pair of good extension is a symmetric condition on the two families while being of wild extension in general is not, since on the second family we ask a stronger condition.
\one

The main point in introducing these notions is the following version of the method used in the proof of \cite{CFK3}[Prop. 4.7]. 
\Le\label{wildext}
Let $(X,H)$ be a polarised scheme over $\gk=\overline{\gk}$.
Let $(\scr{E},\scr{G})$ be a pair of good extension flat over $\cb,\cd$.
\begin{itemize}
    \item For general sheaves $E,G$ in $\scr{E},\scr{G}$, general deformations of $E\oplus G$ are stable hence can be put in a good family $\scr{F}$ flat over $\cc$ and we have
\[dim(\cc)=dim(\cb)+dim(\cd)+\mathfrak{ext}^1(\scr{E},\scr{G})+\mathfrak{ext}^1(\scr{G},\scr{E})-1\]
which is also the dimension of the component of the moduli space of stable sheaves containing the stable elements in $\scr{F}$.
\item We can choose $\scr{F}$ to contain non-split extensions of one among the two forms
\sesl{E}{F_1}{G}{ext11}
\ses{G}{F_2}{E}
\item In addition, if $(\scr{E},\scr{G})$ (resp. $(\scr{G},\scr{E})$) is of wild extension and $\scr{F}$ is defined to contain extensions of type \eqref{ext11} then also $(\scr{F},\scr{G})$ (resp. $(\scr{F},\scr{E})$) is of wild extension.
\end{itemize}
\ma
\begin{proof}
\textbf{Step $1$: Constructing} $\mathbf{\scr{F}}$

    By assumption, we can choose $E,G$ in $\scr{E},\scr{G}$ stable, non-isomorphic, with $p(E)=p(G)$ and such that 
    \[ext^2(E,E)=ext^2(G,G)=ext^2(E,G)=ext^2(G,E)=0.\]
    Therefore, from $\mathfrak{ext}^1(\scr{G},\scr{E})>0$ we have non-split extension as in
    \sesl{E}{F}{G}{est1}
    with $F$ being semi-stable, simple and $p(F)=p(E)=p(G)$ by \thref{estensionistabili} $1),3)$.
    We can find a modular family $\cc$ of simple sheaves containing $F$ and consider some irreducible component where $F$ appears.
    Clearly, the same would hold if we exchange $E$ and $G$ by the symmetry in the definition of family of good extension.
    The smoothness of $\cc$ will follow once we prove $ext^2(\scr{F}_c,\scr{F}_c)=0$ for $c\in\cc$ general, since then, up to shrinking $\cc$, we can suppose this vanishing to hold for all points in it.
    Similarly, as soon as we prove that there is a stable sheaf among them, then it is true for a general one, since this is an open property. 

\textbf{Step $2$: Smoothness and dimension} 

    Next, we show that $ext^2(F,F)=0$. 
    Recalling our choice of $E,G$ and applying $hom(-,E),hom(-,G)$ to \eqref{est1} we get $ext^2(F,E)=0=ext^2(F,G)$.
    Applying $hom(F,-)$ to \eqref{est1} we deduce $ext^2(F,F)=0$ as desired.
    Consider the $\mathcal{C}$-flat family $\scr{F}$ .
    By semicontinuity of cohomology, see \cite{BanicaPutinarSchumacher}[Satz 3 i) p.147], we get $ext^2(\scr{F}_c,\scr{F}_c)=0$ for general $c\in \cc$.

To compute dimensions, we will need to understand the extension of those sheaves.
    Since $E,G$ are stable they are simple hence $hom(E,E)=1=hom(G,G)$.
    Moreover, we have $E\neq G$ hence $hom(E,G)=0=hom(G,E)$ by stability since $p(E)=p(G)$.
    Applying $hom(-,G)$ to \eqref{est1} we get 
    \[hom(F,G)=1 \qquad ext^1(F,G)=ext^1(G,G)+ext^1(E,G).\] 
    Applying $hom(-,E)$ to \eqref{est1}, being the map $Hom(F,E)\ra Hom(E,E)\cong\gk$ equal to $0$ since the sequence is not split, we deduce 
    \[hom(F,E)=0 \qquad ext^1(F,E)=ext^1(E,E)+ext^1(E,G)-1.\]
    Therefore, applying $hom(F,-)$ to \eqref{est1} we conclude 
    \[ext^1(F,F)=ext^1(F,E)+ext^1(F,G)=ext^1(E,E)+ext^1(E,G)+ext^1(G,G)+ext^1(E,G)-1.\]
    The claim on the dimension of $\cc$ follows by deformation theory, since $ext^2(F,F)=0$ and the points in $\cc$ parametrize simple sheaves.
    
\textbf{Step $3$: Stability}\\ 
Next, we need to show the stability of $\scr{F}_c$ for a general $c\in \cc$.
Recall that an extension as in \eqref{est1} is semistable hence, since this is an open property, up to shrinking a bit $\cc$ we can suppose that it holds for any sheaf in $\scr{F}$.
Suppose, by contradiction, that the general element in $\scr{F}$ is not stable so that, by \thref{crit_stab}, we know that there are two families $\scr{K},\scr{Q}$ of sheaves on $X$ flat over an irreducible base $\cq$ such that for all $c\in\cc$ there is some $q\in\cq$ and an exact sequence
\ses{\scr{K}_q}{\scr{F}_c}{\scr{Q}_q}
with $p(\scr{F}_c)=p(\scr{Q}_q)=p(\scr{K}_q)$ and $\scr{K}_q,\scr{Q}_q$ torsion-free.
In particular, this holds with $\scr{F}_c=F$, which fits in the non-split sequence \eqref{est1}, so that we can apply \thref{estensionistabili} $2)$ and deduce that either $\scr{Q}_q\cong F$ or $\scr{Q}_q\cong G$.
But the first case would give $\scr{K}_q=0$ which contradicts $p(\scr{Q}_q)=p(\scr{K}_q)$, hence we deduce $\scr{Q}_q\cong G$ and $\scr{K}_q\cong E$.
Therefore, there is an open set in $U\subset\cq$ such that $\scr{K}_q,\scr{Q}_q$ are simple deformations of $E,G$ if $q\in U$, hence, possibly restricting a bit more $U$, the sheaves  $\scr{K}_q,\scr{Q}_q$ actually belong to the families $\scr{E},\scr{G}$.
In other words, we just proved that the general element in $\scr{F}$ is an extension of sheaves in $\scr{E},\scr{G}$.
    Now, a dimensional count will give us a contradiction.
    
    Since each isomorphism class of sheaves appears at most finitely many times in the modular family, it is enough to parametrize all the  extensions of sheaves in $\scr{E},\scr{G}$ by a scheme of dimension lower than $dim(\cc)$.
    Call $p:X\times \cb\times \cd\ra \cb\times \cd$ the projection.
    We denote by $\ext_p^i(\scr{G},\scr{E})$ the sheaf on $\cb\times \cd$ obtained by sheafifying the presheaf $U\mt Ext^i(\scr{G}|_{p\inv(U)},\scr{E}|_{p\inv(U)})$.
    By \cite{BanicaPutinarSchumacher}[Satz 3, ii) p.147], on the open subset in $\cb\times\cd$ where $ext^1(\scr{G}_b,\scr{E}_b)=\mathfrak{ext}^1(\scr{G},\scr{E})=\rho>0$ the sheaf $\ext_p^1(\scr{G},\scr{E})$ is locally free of rank $\rho$, hence the non-split extensions (up to isomorphism) of a pair of general sheaves in $\scr{E},\scr{G}$ are parametrised by $\p(\ext_p^1(\scr{G},\scr{E}))$.
    There is a sheaf, the \textit{universal extension}, on this space constructed as in \cite{Lange}[Remark 4.3 b)].
    Clearly we have $dim(\p(\ext_p^1(\scr{G},\scr{E})))=dim(\cb)+dim(\cd)+\rho-1<dim(\cc)$ as desired, since $\mathfrak{ext}^1(\scr{E},\scr{G})>0$.

  \textbf{Step 4: Wild extension}\\
  The above computation also proves that $(\scr{F},\scr{G})$ is of good extension, if we can prove that $\mathfrak{ext}^1(\scr{G},\scr{F})>0$ and $\mathfrak{ext}^1(\scr{F},\scr{G})>0$.
  We directly prove the stronger inequalities appearing in the definition of family of wild extension.
  Choose $G'$ a general member of the family $\scr{G}$ and recall that we similarly for $E,F,G$ in \eqref{est1}.
  Applying $hom(G',-),hom(-,G')$ to \eqref{est1} we get 
\[\chi(G',F)=\chi(G',G)+\chi(G',E) \qquad ext^{2i+1}(G',F)\le ext^{2i+1}(G',G)+ext^{2i+1}(G',E)\]
\[\chi(F,G')=\chi(G,G')+\chi(E,G') \qquad ext^{2i+1}(F,G')\le ext^{2i+1}(G,G')+ext^{2i+1}(E,G').\] 
Since Euler characteristic is constant while the $ext^i$-s increase when we specialise our sheaves then we get
\[\mathfrak{ext}^1(\scr{G},\scr{F})-\sum_{i\in\N}\mathfrak{ext}^{2i}(\scr{G},\scr{F})=-\underline{\chi}(\scr{G},\scr{F})-\sum_{i>0}\mathfrak{ext}^{2i+1}(\scr{G},\scr{F})\ge\]
\[\ge-\chi(G',F)-\sum_{i>0}ext^{2i+1}(G',F)\ge\]
\[\ge-\left(\chi(G',G)+\sum_{i>0}ext^{2i+1}(G',G)\right)-\left(\chi(G',E)+\sum_{i>0}ext^{2i+1}(G',E)\right)=\]
\[=-\left(\underline{\chi}(\scr{G},\scr{G})+\sum_{i>0}\mathfrak{ext}^{2i+1}(\scr{G},\scr{G})\right)-\left(\underline{\chi}(\scr{G},\scr{E})+\sum_{i>0}\mathfrak{ext}^{2i+1}(\scr{G},\scr{E})\right)=\]
where in the last step we have equality since we choose $G,G',E$ general in their families.
Continuing the above chain we get
\[=\left(\mathfrak{ext}^1(\scr{G},\scr{G})-\sum_{i\in\N}\mathfrak{ext}^{2i}(\scr{G},\scr{G})\right)+\left(\mathfrak{ext}^1(\scr{G},\scr{E})-\sum_{i\in\N}\mathfrak{ext}^{2i}(\scr{G},\scr{E})\right)>0\]
being $(\scr{E},\scr{G})$ of wild extension.
For $\mathfrak{ext}^1(\scr{F},\scr{G})>\sum_{i\in\N}\mathfrak{ext}^{2i}(\scr{F},\scr{G})$ we argue similarly.

In the case $(\scr{G},\scr{E})$ is of wild extension, we can chose some $E'$ general in $\scr{E}$ and reason similarly: once we prove 
\[\chi(E',F)=\chi(E',G)+\chi(E',E) \qquad ext^{2i+1}(E',F)\le ext^{2i+1}(E',G)+ext^{2i+1}(E',E)\]
\[\chi(F,E')=\chi(G,E')+\chi(E,E') \qquad ext^{2i+1}(F,E')\le ext^{2i+1}(G,E')+ext^{2i+1}(E,E')\] 
the remaining part of the proof is formal.
\end{proof}

The last claim implies that starting only from $\cb$ and $\cd$ of wild extension and applying inductively this lemma, we can obtain extensions of arbitrary high rank.
Under some additional vanishing, we can find a general formula for the dimension of those families.

\Co\label{extpreciso}
In the above setting, suppose $(\scr{E},\scr{G})$ (resp. $(\scr{G},\scr{E})$) is of wild extension. 
Set $\scr{F}_0:=\scr{E}\; ( \text{resp. }\scr{G})$ and $\cc_0:=\cb \;(\text{resp.} \cd)$, then for each $m\in\N$ there is a family $\scr{F}_m$ over a smooth base $\cc_m$ whose elements are Gieseker-stable deformations of extensions of sheaves fitting in
\ses{F_{m-1}}{F_m}{G}
\[(\text{resp.} \;0\to E\to F_m\to F_{m-1}\to 0  \]
Moreover, if
\[\mathfrak{ext}^i(\scr{G},\scr{G})=\mathfrak{ext}^i(\scr{E},\scr{G})=\mathfrak{ext}^i(\scr{G},\scr{E})=0 \qquad \text{for all} \; i\geq 2.\]
then 
\[dim(\cc_m)=dim(\cb)+m\Big( dim(\cd)+\mathfrak{ext}^1(\scr{E},\scr{G})+\mathfrak{ext}^1(\scr{G},\scr{E})-1+(m-1)\mathfrak{ext}^1(\scr{G},\scr{G})\Big)\]
(and similarly exchanging $\scr{G}$ and $\scr{E}$).
\io
\begin{proof}
Applying the second part of \thref{wildext}, we can define inductively $\scr{F}_{m}$ to be the family obtained as an extension of the wild pair $(\scr{F}_{m-1},\scr{G})$ and $\cc_{m}$ its base scheme.
We can compute the above dimension by induction.
The base case $m=0$ follows by the formula in \thref{wildext}.
    To prove the inductive step we can use the same formula, we just need to compute a couple of ext groups.
    A general sheaf $F_m$ in $\scr{F}_m$ is a flat deformation of the direct sum of a sheaf $E$ in $\scr{E}$ and $m$ copies of a sheaf $G'$ in $\scr{G}$.
    Pick another sheaf $G$ in $\scr{G}$.
    We can always suppose $G\neq G'$ otherwise we would have $dim(\cd)=0$ and, being $ext^{2}(\scr{G},\scr{G})=0$, this means $ext^{1}(\scr{G},\scr{G})=0$ hence $\scr{G}$ would not satisfy the properties of wild extension for the pair $(\scr{E},\scr{G})$.
    If we assume $E,G,G'$ to be general, by semicontinuity of cohomology we obtain
\[ 0\leq \mathfrak{ext}^i(\scr{F}_m,\scr{G})=ext^i(F_m,G)\leq ext^i(E,G)+m\cdot ext^i(G',G)=0\]
for all $i\geq 2$.
Similarly, $\mathfrak{ext}^i(\scr{G},\scr{F}_m)=0$ for all $i\geq 2$.
Being $hom(G,E)=0=hom(G,G')$ we have $hom(G,F_m)=0$ and
\[\mathfrak{ext}^1(\scr{G},\scr{F}_m)=ext^1(G,F_m)=-\chi(G,F_m)=-\chi(G,E)-m\chi(G,G')=\]
\[=ext^1(\scr{G},\scr{E})+m\cdot ext^1(\scr{G},\scr{G})=\mathfrak{ext}^1(\scr{G},\scr{E})+m\cdot \mathfrak{ext}^1(\scr{G},\scr{G})).\]
Similarly we deduce $\mathfrak{ext}^1(\scr{F}_m,\scr{G})=\mathfrak{ext}^1(\scr{E},\scr{G})+m\cdot \mathfrak{ext}^1(\scr{G},\scr{G})$.
It follows by the dimension formula in \thref{wildext} that 
\[dim(\cc_{m+1})=dim(\cd)+dim(\cc_m)+\mathfrak{ext}^1(\scr{G},\scr{F}_m)+\mathfrak{ext}^1(\scr{F}_m,\scr{G})-1=\]
\[=dim(\cd)+dim(\cc_m)+\mathfrak{ext}^1(\scr{E},\scr{G})+\mathfrak{ext}^1(\scr{G},\scr{E})+2m\cdot\mathfrak{ext}^1(\scr{G},\scr{G})-1\]
If we now apply the inductive hypothesis, we obtain the claim.
\end{proof}

\chapter{Ulrich sheaves on double covers}

Our next goal would be to study Ulrich bundles on varieties that can be written in a simple way as a finite covering $f:X\ra\p^n$.
Existence results for Ulrich bundles on any smooth cyclic covering of $\p^n$ are already known, see \cite{KuNaPa} for the degree $2$ case and \cite{ParPin} for arbitrary degree, but we will push them further.
Our approach is based on the existence of matrix factorisations; as in the previous ones.
However, adding an idea in \cite{HanselkaKummer}, we will give a different proof which actually works for any $f$ that can be factored through a closed embedding inside a weighted projective space of the form $\p(1^{n+1},m)$, see \thref{ulrichmatrix''}.
It should be noted that the analogous statement for modules over local Cohen-Macaulay rings was proved, with a similar strategy, quite earlier in \cite{BHU}[Thm. 2.5], after the seminal paper \cite{Eis-resolutions}.

Let us briefly explain our reasoning.
We can associate an equation to any divisor $X$ inside $\p(1^{n+1},m)$ and the existence of Ulrich sheaves on $X$ will be equivalent to the possibility of writing some power of this equation as the determinant of some matrix.
Our contribution in \thref{ulrichmatrix''} is to spot a special form for this equation, which gives us an upper bound for the minimal possible rank for an Ulrich bundle on a cyclic covering.
To our knowledge, this is sharpest than the others in the literature and in the case of rank $2$ bundles on double coverings this bound is actually optimal by \thref{branchrango2}.

However, with only algebraic methods we are not able to show that any such Ulrich bundle has to come from this construction.
Therefore, in \thref{section-rk2} we will give a geometric analogue of this algebraic description and show that all Ulrich bundles are actually of this form in \thref{rk2-section}.
The other main result of this chapter is the actual existence of rank $2$, the minimal possible, Ulrich bundles on double coverings of $\p^3$ branched along surfaces of degree $4,6,8$, see \thref{ulrichr2doublecoverexistence}.
In the first two cases $X$ is Fano, so this result was definitely expected, see \cite{CFK3} and \cite{CFK3.1}, while in the third one is Calabi-Yau, a class of varieties where the theory of Ulrich bundle is less developed.

Actually, our proof tells more: if the branch divisor of a double covering $f:X\ra\p^n$ can be written as $p_0^2+p_1p_2+p_3p_4$, where the $p_i$-s have the same degree, then $X$ admits Ulrich bundles of rank $2$ respect to $f^*\osn(1)$.
The geometric way to interpret this is as follows.
The zero locus of a section of such Ulrich bundle is a curve mapped isomorphically by $f$ to the complete intersection of two divisors.
An example of such an image is the complete intersection of $p_1$ and $p_3$, or any other complete intersection of the same degrees on which the branch locus restricts to a square.
A series of geometrical results of this fashion will occupy the last part of this chapter.
Let us end this introduction with another interesting remark.
By \thref{hilbertpol}, any rank $2$ Ulrich bundle $\e$ on a double covering gives a morphism $X\ra Gr(2,4)$, which is a smooth quadric in $\p^5$, such that the $\e$ is the pullback of the universal quotient bundle, which is known to be a spinor bundle hence itself Ulrich.

In the first section, we fix the notation and recall the general theory of cyclic covers and what we call divisorial covers, \thref{defdiv}, which will be our ambient varieties. 
In the second, we survey on the connection between Ulrich sheaves and matrix factorisations, in particular in the case of cyclic coverings, see \thref{ulrichsommeprodotti}, with the aim of proving the main theorems stated above. 
In the last section we study the geometry of Ulrich sheaves on integral double coverings of $\p^n$, in particular their syzygy sheaves in \thref{Discount}, degeneracy loci in \thref{section-rk2} and \thref{rk2-section}, relation with the covering involution in \thref{fisso} and the non-ample loci in \thref{nonampio}.
Most of those ingredients will be used in the study of their moduli spaces in the next chapter, for example the bundles fixed by the involution will give us singular points.

\startcontents[chapters]
\printcontents[chapters]{}{1}{}

\section{Preliminaries on finite coverings}

The following is a recap of well-known facts.
Some of the material on cyclic coverings is, for example, in \cite{BHPV}[§1.16], \cite{EsnaultViehweg}[§3] or \cite{Laz1}[§ 4.1].
Since the connection between Ulrich sheaves and matrix factorisation lives really on the algebraic level, we want to work over arbitrary fields.

\subsection{Divisorial coverings}\label{notazioneproj}

Fix a field $\gk$.
We will always work with a finite surjective morphism $f:X\ra\p^n$ of degree $d$ and $dim(X)=n\ge 1$.
Given such an $f$, we always have an injective morphism $\osn\ra f_*\os_X$; call $E$ its cokernel.
We deduce an embedding $X=Spec_{\osn}(f_*\os_X)\hookrightarrow Spec_{\osn}(Sym(E))$, which in the case $X$ is Gorenstein can also be refined to an embedding of $X$ in $\p(E)$, see \cite{Casnatiekedahl}[Thm. 1.3].
Our goal is to study $f$-Ulrich bundles on $X$ when $f$ factors as $X\subset\p(\osn\oplus\osn(m))$ followed by the standard projection to $\p^n$, we will call those divisorial coverings.
For example, among those are the well-known cyclic coverings.
This total space can be compactified to a weighted projective space in which $X$ sits as a divisor, therefore our work will be guided by the analogy with divisors in $\p^n$.

\De\label{defdiv}
Take a finite surjective map $f:X\ra\p^n$.
We will call $f$ (and by extension $X$) \textbf{divisorial} if $f=\pi\circ i$ where $i$ embeds $X$ in the total space of some line bundle over $\p^n$ and $\pi$ is the standard bundle projection.
\Ne
This definition can be extended to arbitrary base varieties instead of $\p^n$ but we will treat essentially only this case.

\Le\label{divcriterion}
$f$ is divisorial if and only if there is some integer $m\geq 0$ and a surjective morphism of sheaves of $\osn$-algebras $Sym(\osn(-m))\ra f_*\os_X$
\ma
\begin{proof}
There is a controvariant equivalence between the category of $\osn$-algebras and schemes affine over $\p^n$, see \cite{GW1}[Cor. 12.2].
Through this equivalence, a surjective morphism $Sym(\osn(-m))\ra f_*\os_X$ corresponds to an embedding $X\ra Spec (Sym(\osn(-m)))$, where the latter one is the total space of the line bundle $\osn(m)$.
\end{proof}

The case $m=0$ corresponds to disconnected trivial covers, which are contained in $\p^n\times\mathbb{A}^1$, the total space of the trivial line bundle.
So we will always assume $m>0$.

\subsection{Projective bundles and weighted projective spaces}

The variety $Spec (Sym(\osn(-m)))$ is not proper, therefore we prefer to work with some compactification.
The first one is the projective closure.
We are going to review some standard theory of projective bundles: these are all well-known facts which we will repeat just for completeness and to fix the notation once and for all.
We follow the terminology in \cite{EGAII}[Chapter 8 Section 4], for example $\p(-):=Proj(Sym(-))$ will always stand for the projectivisation parametrising $1$-dimensional quotients.

\De
Let $m$ be a positive integer.
We define $P_m:=\p(\os_{\p^n}\oplus\osn(m))$, set $\pi_m:P_m\ra \p^n$ the natural projection, and $\os_{P_m}(1)$ the relative hyperplane bundle.
\Ne    
We will often omit $m$.

    Applying $\p(-)$ associates to the two projections of $\osn\oplus\osn(m)$ onto its factors two section of $\pi$, see \cite{Har}[II Prop. 7.12].
    We will call them, respectively, $H_0:=\p(\osn)$ the \textbf{zero section} and $H_\infty:=\p(\osn(m))$ the \textbf{hyperplane at infinity}.
    Note that $Spec (Sym(\osn(-m)))$ can be identified with $P_m-H_0$, see \cite{EGAII}[Prop. 8.8.4].
    In particular, every divisorial covering of $\p^n$ is naturally embedded in $P_m-H_0$.
    On $P_m$ cohomology computations are easy to handle but sometimes we will prefer a slightly smaller compactification.
    To introduce it let us make a brief detour.
    
    Note that $\os_{P_m}(1)$ is globally generated since $\osn\oplus\osn(m)$ is such, moreover using \cite{GW2} [Prop. 22.86] we have natural isomorphisms 
    \[H^0(P_m,\os_{P_m}(1))\cong H^0(\p^n,\pi_*\os_{P_m}(1))=H^0(\p^n,\osn\oplus\osn(m)).\]
    Let us call $t$ the section of $\os_{P_m}(1)$ corresponding to $1\in H^0(\p^n,\osn)$.
    We claim that the morphism given by $\os_{P_m}(1)$ maps $P_m$ to the cone over the $m$-th Veronese of $\p^n$ contained in the projective space of dimension $\binom{n+m}{n}$: clearly the fibres of $\pi$ are sent to lines, so it is enough to show that $H_0$ is contracted to a point (the vertex) while on $H_{\infty}$ we get the Veronese embedding (the base).
    In fact, since the sections $H_0 \;(resp.\; H_{\infty})$ of $\pi$ corresponds to the projection on the first (second) factor of $\osn\oplus\osn(m)$, the functorial interpretation of the projective bundle as a Grassmannian implies $\os_{P_m}(1)|_{H_0}\cong \os_{H_0}\cong \osn \; (resp. \; \os_{P_m}(1)|_{H_\infty}\cong \osn(m))$, see \cite{Har}[II proof of Prop. 7.12] or \cite{GW1}[§8.6].
    Finally, the restriction map $H^0(P_m,\os_{P_m}(1))\ra H^0(H_{\infty},\os_{P_m}(1)|_{H_{\infty}})\cong H^0(H_{\infty},\osn(m))$ is surjective, hence the induced map is really the $m$-th Veronese.
    In another way, the cone over the $m$-th Veronese of $\p^n$ is, by definition, $Proj\left(\bigoplus_{i\in \N}H^0(P_m,\os_{P_m}(i))\right)$.
Similarly as before we have
\[H^0(P_m,\os_{P_m}(d))\cong H^0(\p^n,\pi_*\os_{P_m}(d))=H^0(\p^n,Sym^d(t\gk\oplus\osn(m)))=\]
\[=H^0\left(\p^n,\bigoplus_{j=0}^dt^j\otimes\osn((d-j)m)\right)=\bigoplus_{j=0}^dt^jH^0(\p^n,\osn((d-j)m))\]
therefore, if we consider $t$ to have degree $m$ and $x_i$ degree $1$, we get that $\bigoplus_{i\in \N}H^0(P_m,\os_{P_m}(i))$ is actually the subring of elements of degree multiple of $m$ in the ring $\gk[x_0,\dots, x_n,t]$.
By \cite{GW1}[Remark 13.7] we have $Proj\left(\bigoplus_{i\in \N}H^0(P_m,\os_{P_m}(i))\right)\cong Proj(\gk[x_0,\dots, x_n,t])$ as claimed.
Note that, for $m=1$ this cone is exactly $\p^{n+1}$ and $P_1\cong Bl_V\p^{n+1}$, hence the corresponding divisorial covers are nothing else than divisors in $\p^{n+1}$.

In the following, we will always confuse vectors in $H^0(P_m,\os_{P_m}(d))$ with homogenous polynomials of degree $dm$ in $\gk[x_0,\dots, x_n,t]$.

\De
We call $\p(1^{n+1},m):=Proj(\gk[x_0,\dots, x_n,t])$ the \textbf{weighted projective space} with weights $(1^{n+1},m)$.
\Ne

Call $V$ the vertex of $\p(1^{n+1},m)$, that is the point of coordinates $[0:\dots 0:1]$.
$\p(1^{n+1},m)-V$ is isomorphic to $P_m-H_0$, hence any divisorial covering $f:X\ra\p^n$ can be embedded in $\p(1^{n+1},m)-V$.
$X$, inside $\p(1^{n+1},m)$, must be defined by an homogeneous polynomial $p=\sum_{i=0}^d t^ib_{d-i}$ of degree $md$, where $b_j$ are homogeneous polynomials of degree $mj$ in the $x_i$-s and $t$ is considered to have degree $m$.

\Le\label{equivalenza}
The following properties of $X$ are equivalent:
\begin{itemize}
    \item there is some divisorial covering $f:X\ra\p^n$ of degree $d$
    \item it is a divisor in $\p(1^{n+1},m)$ not passing through $V$ and defined by a polynomial $p\in\gk[x_0,\dots, x_n,t]$ of degree $dm$
    \item it is defined in $\p(1^{n+1},m)$ by a polynomial with $b_0\neq 0$, hence we can set $p=t^d+\sum_{i=0}^{d-1} t^ib_{d-i}$
    \item is a divisor in $|\os_{P_m}(d)|$ not intersecting $H_0$.
\end{itemize}
\ma
\begin{proof}
We start with the equivalence of the first two and then deduce the rest.
    We have already seen that a divisorial covering $f:X\ra\p^n$ can be embedded inside $\p(1^{n+1},m)$ avoiding $V$.
    Conversely, if $X$ does not pass through $V$ then we can project from this point to the base of the cone, which is $\p^n$.
    The restriction of this map to $X$ is finite since $X$ is proper and intersects each fiber in a finite length subscheme, being $V$ in the closure of all fibers but not contained in $X$.
    If $x\in\p^n$ is a closed point then the fiber of $f\inv(x)$ represents the cycle $H^n$ and, by definition, we have that $deg(H^n)=length(f\inv(x))$.
    The scheme $f\inv(x)$ is defined on $\mathbb{A}^1=\pi\inv(x)$ by the restriction of $p$, hence $length(f\inv(x))$ is the degree of $p$ in the variable $t$.
    
    Note that, $V=[0:\dots 0:1]$ does not satisfy $p$ if and only if $b_0\neq 0$.
    Finally, by definition of $W$ as the image of $P$ under the morphism given by $\os_{P_m}(1)$, its degree $d$ hypersurfaces correspond to divisors in $|\os_{P_m}(d)|$ and the condition of not containing $V$ is equivalent to not intersecting $H_0$.
\end{proof}

Note that we have just parametrised divisorial coverings of degree $d$ by the complement of an hyperplane in $H^0(P_m,\os_{P_m}(d))\cong \gk^{\binom{n+m}{m}+1}$.
In fact, in an equivalent manner, we could take the complement of a hyperplane in $\p(H^0(P_m,\os_{P_m}(d)))$.

\De\label{generaldivisorial}
Fixed $n,d,m$, we say that some property holds for the general divisorial covering if: the set of $p\in H^0(P_m,\os_{P_m}(d))$ such that this property holds on the covering determined by $p$ contains an open subset for the Zariski topology.
\Ne

From now on we will always identify divisorial coverings $f:X\ra\p^n$ with their embeddings in $P_m$ or $\p(1^{n+1},m)$.
Moreover, we set $\os_X(1):=f^*\osn(1)$, so that $(X,\os_X(1))$ is a polarised variety.
    
\Le\label{divisorialcover}
    Let $f:X\ra \p^n$ be a divisorial covering of degree $d$ inside $P_m=\p(\osn\oplus\osn(m)) \;(m\geq 1)$ then 
    \begin{enumerate}[1)] 
        \item $\os_{P_m}(d)|_X\cong \os_X(dm)$.
        \item $X$ is Gorenstein and $\omega_X\cong \os_X(m(d-1)-n-1))$ 
        \item $f$ is flat and $f_*\os_X\cong \bigoplus_{i=0}^{d-1}\osn(-im)$ 
        \item $X$ is connected and $(X,\os_X(1))$ is aCM
        \item Suppose $X$ smooth and irreducible.
        If $\e$ is an $f$-Ulrich bundle of rank $r$ then
        \[deg(\e)=c_1(\e)\cdot c_1(\os_X(1))^{n-1}=rm\dfrac{d(d-1)}{2}.\]
    \end{enumerate}
    \ma
    \begin{proof}
        
        \gel{1}
        By \thref{equivalenza} we have $H_0\cap X=\emptyset$ hence $\os_{P_m}(H_0)|_X\cong \os_X$.
        We know that $\os_{P_m}(H_{0})\cong \os_{P_m}(1)\otimes\pi^*\osn(-m)$, see for example \cite{Sernesi}[Prop. 4.6.2], therefore $\os_{P_m}(1)|_X\cong f^*\osn(m)$ and the claim follows for any $d$ taking tensor products.   

        \gel{2}
        From \cite{GW2}[Thm. 22.86 4)] we know $\omega_{P_m/\p^n}  \cong\os_{P_m}(-2)\otimes\pi^*(\osn(m))$.
        Since $P$ is smooth then the sequence of relative differentials
        \ses{\pi^*\Omega_{\p^n}}{\Omega_{P_m}}{\Omega_{P_m/\p^n}}
        is exact and taking determinants gives \[\omega_{{P_m}}  \cong\omega_{P_m/\p^n}\otimes\pi^*\omega_{\p^n}\cong\os_{P_m}(-2)\otimes\pi^*(\osn(m-n-1)).\]
        Since $X$ is a Cartier divisor in $P_m$, which is smooth, it must be Gorenstein and by the adjunction formula, \cite{GW2}[Cor. 25.130 2)], we have \[\omega_X\cong \omega_{P_m}(d)|_X\cong \os_{P_m}(d-2)\otimes\pi^*(\osn(m-n-1))|_X\cong \os_X(m(d-1)-n-1)\] where we also used $1)$.

        \gel{3}
        By \thref{equivalenza}, the sequence defining $X$ in $P$ is
\ses{\os_{P_m}(-d)}{\os_{P_m}}{\os_X}
If we apply $\pi_*$ to it we obtain
\sesl{\osn}{f_*\os_X}{\bigoplus_{i=1}^{d-1}\osn(-im)}{d''''}
where we used that \cite{GW2}[Prop. 22.86] implies $\pi_*\os_{P_m}(-d)=0$, $\pi_*\os_{P_m}\cong \osn, R^1\pi_*\os_{P_m}=0$ and 
\[R^1\pi_*\os_{P_m}(-d)\cong \left(Sym^{d-2}(\osn\oplus\osn(m))\otimes\osn(m)\right)\du\cong \bigoplus_{i=1}^{d-1}\osn(-im).\]
Note that $f$ is affine and $f_*\os_X$ is locally free so $f$ must also be flat, see \cite{GW1}[Prop. 12.19].
Moreover, the sequence \eqref{d''''} splits being
\[Ext^1\left(\bigoplus_{i=1}^{d-1}\osn(-im),\osn\right)\cong H^1\left(\p^n,\bigoplus_{i=1}^{d-1}\osn(im)\right)=0.\]

\gel{4}
Line bundles on $(\p^n,\osn(1))$ have no intermediate cohomology, therefore the same holds for $(f_*\os_X)(i)$ for all $i\in\Z$.
By projection formula we get
\[f_*(\os_X(i))=f_*(f^*\osn(i))\cong f_*(\os_X\otimes f^*\osn(i))\cong (f_*\os_X)(i).\]
But the map $f$ is finite hence preserves cohomology, so we conclude that $(X,\os_X(1))$ is aCM and $h^0(X,\os_X)=h^0(\p^n,\osn)=1$ hence $X$ has to be connected.

     \gel{5}
    It is enough to apply \thref{c12}, recalling that $\os_X(K_X)\cong \omega_X\cong \os_X(m(d-1)-n-1))$ from part $2)$.
    \end{proof}

\subsection{Cyclic coverings}

As a further specialisation we will get cyclic coverings.
Again, we note that also this definition could be phrased for arbitrary base schemes.
\De\label{defcyclic}
A divisorial covering $f:X\ra\p^n$ is \textbf{cyclic} if the equation of $X$ in $P_m$ is of the form $t^d-b$, with $b\in H^0(\p^n,\osn(dm))$. 

\noindent Fixed $n,d,m$, we say that some property holds for the general cyclic covering if: the set of $b\in H^0(\p^n,\osn(md))$ such that this property holds on the covering determined by $t^d-b$ contains an open subset for the Zariski topology.
\Ne

The degree $2$ case is very special. 

\Le\label{2cyc}
If $char(\gk)\neq 2$ then any flat, finite morphism $f:X\ra\p^n$ of degree $2$ is a cyclic covering.
\ma
\begin{proof}
Being $f$ flat, $f_*\os_X$ is locally free, has rank equal to $deg(f)=2$ and is split by the well-known \thref{trace}.
Since $h^0(X,\os_X)=h^0(\p^n,f_*\os_X)$ we must have $f_*\os_X\cong \osn\oplus\osn(-m)$ for some $m>0$ hence \thref{divcriterion} implies that $X$ is divisorial.
    Inside $W$ our $X$ has equation $p=t^2+b_1t+b_2=0$ but since $2$ is invertible we can always complete the square as $p=(t+\frac{b_1}{2})^2+b_2-\frac{b_1^2}{4}=0$.
    Sending $t\mt t-\frac{b_1}{2}$ and fixing the other variables, we get an automorphism of $W$ that respects the projection to $\p^n$, which just forgets the coordinate $t$, and sends the equation of $X$ to $t^2-\left(b_2+\frac{b_1^2}{4}\right)$ as desired.
\end{proof}

\Le\label{trace}
Suppose $char(\gk)$ is co-prime with $d$.
Let $f:X\ra \p^n$ be a morphism such that $f_*\os_X$ is a locally free $\osn$-module of rank $d$.
The natural map $\osn\hookrightarrow f_*\os_X$ associated with $f$ splits.
\ma
\begin{proof}
    It is enough to define this splitting on some open covering of $\p^n$ trivialising $f_*\os_X$ in a compatible way.
For every such open $U\subset \p^n$, any $s\in H^0(U,f_*\os_X)= H^0(f\inv(U),\os_X)$ acts on $H^0(f\inv(U),\os_X)$ by multiplication, call $m_s$ this endomorphism.
Since $H^0(f\inv(U),\os_X)$ is a free $H^0(U,\osn)$-module of rank $d$ and this action is linear, we have a well-defined notion of trace for $m_s$.
This trace does not depend on the choice of a trivialisation since it does not depend on the choice of a basis of $H^0(f\inv(U),\os_X)$ as an $H^0(U,\osn)$-module, due to the standard identity $tr(MAM\inv)=tr(A)$ for any two square matrices of the same size $A,M$ with $M$ invertible.
All these operations are compatible with restriction maps of sheaves, therefore \[\mathfrak{tr}:f_*\os_X\ra \osn \qquad s\mt \dfrac{tr(m_s)}{d}\]
is well defined, also due to the fact that $d\in \gk^*$.
Clearly, $\mathfrak{tr}|_{\osn}=id_{\osn}$, so we get the desired splitting.
\end{proof}

Embedding a divisorial covering $f:X\ra\p^n$ inside $\p(1^{n+1},m)$ we can recover $f$ as the restriction of the map given by projecting from the vertex, that is removing the variable $t$.
Therefore, we can identify this $\p^n$ with the base of the cone, which is the divisor given by $\{t=0\}$ and corresponds to $H_{\infty}$ in $P_m$.

\De
The \textbf{ramification divisor} $R$ of a cyclic covering $f:X\ra\p^n$ is the divisor given by $\{t=0\}$, since $t$ is seen as a section of $\os_{P}(1)$.
The divisor $B:=f(R)\cong R$ is called \textbf{branch locus}.
\Ne

A local computation shows that, if $char(\gk)$ is coprime with $d$, then $R$ is exactly the locus on which the morphism $f$ is ramified.

\Le\label{ramificazione}
If $f:X\ra\p^n$ is a cyclic covering then $R\in |\os_X(m)|$ and $f^*B=dR$.
\ma
\begin{proof}
    By definition $R\in |\os_{P_m}(1)|_X|$ but by \thref{divisorialcover} we have $\os_{P_m}(1)|_X\cong \os_X(m)$.
    Moreover, being the equation of $X$ of the form $t^d-b=0$, the divisor $f^*B$, which is cut on $X$ by the equation $b=0$, is the same as the one given by $t^d=0$ but this last one is exactly $d$ times $R$.
\end{proof}

\Le\label{cycliciinvolution}
Suppose that $\gk$ has characteristic co-prime with $d$ and is algebraically closed.
On any cyclic covering $f:X\ra\p^n$ of degree $d$ there is an automorphism of order $d$ whose locus of closed fixed points is $R$.
The projection $\pi|_X:X\ra\p^n$ is exactly the quotient by this automorphism, in particular gives an isomorphism $R\cong B$.
\ma
\begin{proof}
Take a primitive $d$-th root of unity $\zeta$.
    We define an automorphism of $\p(1^{n+1},m)$ from the following automorphism of $\gk[x_0,\dots,x_n,t]$:
    \[x_i\mt x_i,\qquad t\mt \zeta t.\]
    The locus of fixed points is exactly the one given by $t=0$, which is $R$ by \thref{ramificazione}.
    Moreover, any degree $d$ cyclic covering is sent to itself since the equation $t^d-b$ is invariant for this action.
    Therefore this descends to an automorphism of $X$ with fixed locus $R$.
    Moreover, by definition, it preserves exactly the fibers of the projection to $\p^n$ thus is the quotient map.
    Since it acts as the identity on $R$ it must map isomorphically this scheme.
\end{proof}

\Le
If $\gk$ contains all $d$-th roots of its elements then a cyclic covering is determined, up to isomorphism, by the branch divisor $B$, hence by the locus $\{b=0\}$ inside $H_{\infty}$.
\ma
\begin{proof}
    Two cyclic covering have the same branch divisor $B$ if and only if they are given by equations of the form $t^d-b=0$ and $t^d-\lambda b=0$ for some $\lambda\in\gk^*$.
    Define $\zeta$ to be a $d$-th root of $\lambda\inv$ then we can define an automorphism of $\p(1^{n+1},m)$ by 
    $\gk[x_0,\dots,x_n,t]$:
    \[x_i\mt x_i,\qquad t\mt \zeta t.\]
    Under this map the locus $\{t^d-\lambda b=0\}$ is sent to $\{\lambda t^d-\lambda b=0\}$ which coincides with $\{t^d-b=0\}$.
\end{proof}

As a consequence, if $\gk$ contains all $d$-th roots of its elements then we can parametrize cyclic coverings of degree $d$ with $|\osn(md)|$ instead of a subspace of $H^0(P_m,\os_{P_m}(d))$.

\subsection{Picard group of divisorial coverings}

The aim of this section is just to recall some  sufficient conditions to have $Pic(X)\cong \Z H$.
As an application, we will see how degeneracy loci of Ulrich bundles behave on them.

\Le\label{Piccover}
Suppose $f:X\ra\p^n$ is a smooth divisorial covering of degree $d$.
If at least one of the following holds
\begin{itemize}
    \item $f$ is cyclic, $char(\gk)$ is coprime with $d$ and $n\geq 4$
    \item $\gk=\C$ and $n>d$
\end{itemize}
then $Pic(X)$ is generated by $\os_X(1)$. 
Under this assumption, for any $f$-Ulrich bundle $\e$ of rank $r$ we have $det(\e)\sim \os_X\left(\dfrac{rm(d-1)}{2}\right)$.
\ma
\begin{proof}
We start with the claim on $Pic(X)$.
In the first case, it follows from \cite{SGA2}[XII Corollaire 3.7] that $Pic(B)\cong Pic(\p^n)$, where $B$ is the branch divisor of $f$.
From the sequence
\ses{\os_X(-R)}{\os_X}{\os_R}
together with the fact that $(X,\os_X(1))$ is aCM, see \thref{divisorialcover}, and $n\ge 4$ we conclude that $H^j(R,\os_R(-i))=0$ for $j=1,2$ and $i>0$.
This implies, together with the fact that we have $B\cong R$ by \thref{cycliciinvolution} and this last one is ample on $X$, that we can apply \cite{SGA2}[XII Corollaire 3.6] and get $Pic(R)\cong Pic(X)$ hence the claim.
    In the second case, we are in a position to apply \cite{Laz}[Prop. 3.1]. 
    
    Therefore, we must have that $det(\e)$ is a power of $\os_X(1)$, hence $\e$ is a special Ulrich bundle.
    The last statement follows from \thref{divisorialcover} $5)$ by recalling that $c_1(\os_X(1))^n=d$.
\end{proof}

From the above, we can deduce some nice properties of degeneracy loci of Ulrich bundles on those special varieties.
We start by recalling some results on degeneracy loci which we will also need later.

\Le\label{tuttedeg}
Suppose $\gk=\overline{\gk}$.
Let $\e$ be an initialised rank $r$ vector bundle on a smooth polarised variety $(X,H)$, i.e. $h^0(X,\e)>0$ but $h^0(X,\e(-1))=0$.
If $Pic(X)\cong \Z H$ then the locus where an injective morphism $\varphi:\os_X^{r-1}\ra\e$ drops rank is empty or has pure codimension $2$ hence we have
\sesl{V\otimes\os_X}{\e}{\id_Y\otimes det(\e)}{degv}
Moreover, if $\e$ has no trivial direct summands and $h^1(X,det(\e)\du)=0$ then only the second possibility takes place.
\ma
\begin{proof}
    Assume that this morphism actually drops rank at some point in $X$, otherwise the thesis is clear.
    If the locus where it happens has a codimension $1$ component, call it $\Delta$, then the restriction morphism $\os_\Delta^{r-1}\ra\e|_\Delta$ is no more injective.
    This means that one among the sections of $\e$ in the image of $\varphi$ vanishes on $\Delta$, that is, $h^0(X,\e(-\Delta)>0$.
    However, being $X$ smooth, $\Delta$ is an effective Cartier divisor, and therefore is linearly equivalent to a non-negative multiple of $H$.
    Now, our assumption that $\e$ is initialised implies $\Delta=\emptyset$, hence $dim(Y)\leq dim(X)-2$.
    Since $coker(\phi)$ has rank $1$ we can apply both parts of \thref{degenerazione'} and conclude that it is isomorphic to $\id_Y\otimes det(\e)$ where $Y\subset X$ is exactly the locus where $\phi$ drops rank.
    Since $\id_Y\otimes det(\e)$ has a resolution of length $1$, by \thref{cmcodim2} we conclude that $Y$ has pure codimension $2$.

    For the last claim, consider the sequence \eqref{degv}.
    We see that $Y\ne \emptyset$ because otherwise $ext^1(det(\e),\os_X)=h^1(X,det(\e)\du)=0$ implies that the above splits contradicting the fact that $\e$ has no trivial summand. 
\end{proof}

\Co\label{tuttedeg2}
In the above setting, if $rk(\e)=2$ then any section of $\e$ gives a sequence as \eqref{degv}.
\io
\begin{proof}
    Being $X$ integral, any non-zero morphism $\os_X\ra\e$ is injective and hence \thref{tuttedeg} applies to it.
\end{proof}

\subsection{Subvarieties lifting to double coverings}

In this section we will face the following problem: given a double covering $f:X\ra\p^n$, when is the preimage of a subvariety $Z\subset \p^n$ made of two components mapping isomorphically to $Z$?
As suggested by \thref{restrizionesezioni}, the zero loci of section of rank $2$ Ulrich bundles on double coverings of $\p^n$ are always of this form, see \thref{section-rk2}.
 
\Le\label{double}
Suppose that $char(\gk)\neq 2$ and $\gk=\overline{\gk}$.
Let $f:X\ra \p^n$ be a double covering such that $f_*\os_X\cong \osn\oplus\osn(-m)$, with $m>0$.
We call $B$ the branch locus, $R\subset X$ the ramification divisor. 
Let $j:Z\hookrightarrow \p^n$ be a closed subvariety and set $\os_Z(i):=j^*\osn(i)$.
The following are equivalent:
\begin{enumerate}[i)]
    \item $f\inv (Z)$ has two irreducible components mapped isomorphically to $Z$ by $f$
    \item $Z\not\subset B$ and over $Z$ we can find a section of $f$
    \item $B|_Z$ is the double of an effective Cartier divisor in $|\os_Z(m)|$ and $Z\not\subset B$. 
    \end{enumerate}
\ma
\begin{proof}
\gel{i)$\Rightarrow$ ii}
    Suppose $f\inv (Z)=Z_1\cup Z_2$, so that $f|_{Z_i}:Z_i\ra Z$ is an isomorphism.
    Its inverse is the section we desire.
    If $Z\subset B$ then $f\inv(Z)\subset f\inv(B)=2R$.
    Topologically we would have $f\inv(Z)\subset R$ but $f|_R$ is an isomorphism hence $f\inv(Z)$ would have only one irreducible component.

\gel{ii)$\Rightarrow$ iii}
    By assumption we have a morphism $s:Z\ra X$ such that $f\circ s=j$.
    Functoriality of pullback implies that 
    \[B|_Z=j^*B=s^*(f^*B)=s^*2R=2s^*R\]
    where we used \thref{ramificazione} to say that $f^*B=2R$.
    Finally notice that
    \[\os_Z(s^*R)\cong s^*\os_X(m)\cong s^*f^*\osn(m)\cong\os_Z(m).\]

\gel{iii)$\Rightarrow$ i}  
    $X$ sits inside $P=\p(\osn\oplus\osn(-m))$ with equation $t^2-b=0$, where $b=0$ is the equation of $B\subset\p^n$.
    If we call $P\times_{\p^n}Z=\p(\os_Z\oplus\os_Z(-m))=:P_Z$ then the scheme $f\inv(Z)$ has equation $t^2-b$ in $P_Z$.
    But $B|_Z$ is the double of a Cartier divisor in $|\os_Z(m)|$, that is $b|_Z$ is the square of some $\beta\in H^0(Z,\os_Z(m))$.
    It follows that $f\inv(Z)$ has two components corresponding to $t\pm\beta=0$; we work out the case $t+\beta$.
    This polynomial defines a morphism
    \[0\ra \os_{P_Z}\xra{t+\beta}\os_{P_Z}(1)\ra \os_{P_Z}(1)\otimes \os_{Z_1}\ra 0\]
    whose pushforward, using \cite{GW2}[Prop. 22.86], is
    \[0\ra \os_{Z}\xra{(1,\beta)}\os_{Z}\oplus\os_Z(m)\ra \os_Z(m)\ra 0\]
    and $Z_1$ is the section of $\pi|_{P_Z}$ over $Z$ corresponding to the surjection on the right.
    It follows that $\pi|_{Z_1}=f|_{Z_1}:Z_1\ra Z$ is an isomorphism.
\end{proof}

\Oss
A similar statement holds in the case of cyclic coverings of arbitrary degree.
More generally, one could ask what is the condition to have two components $Z_1,Z_2$ mapping birationally to $Z$.
In the case $Z$ is a divisor, and we deal with double coverings, a characterisation is given in \cite{CD}.
\one

\Co\label{intersezionebranch}
Suppose any condition in \thref{double} to be satisfied.
If $Z\subset\p^n$ is a variety such that $f\inv(Z)$ has two components $Z_i$ mapped isomorphically on $Z$ by $f$ then 
\[Z_1\cap R=Z_1\cap Z_2=Z_2\cap R,\]
where intersections and equalities are meant scheme-theoretically.
\io
\begin{proof}
We need just to prove that the ideal defined by the intersection of any two among $Z_1,Z_2,R$ is the same defined by the triple intersection.
    $Z_i$ have equations $t\pm\beta=0$ while $R$ has equation $t=0$ by \thref{ramificazione}, hence the claim follows.
\end{proof}

We end by proving a lemma which will be useful later.

\Le\label{idealedoppi}
Let $f:X\ra\p^n$ be a double covering such that $f_*\os_X\cong\osn\oplus\osn(-m)$. If $Y\subset X$ and $f|_Y$ is an isomorphism on the image $Z$, then we have an exact sequence
\ses{\id_Z(m)}{f_*\id_Y(m)}{\osn}
which splits if $h^1(\p^n,\id_Z(m))=0$.
\ma
\begin{proof}
We have an exact sequence
\ses{\id_Y}{\os_X}{\os_Y}
    Since $f$ is finite, $f_*$ is exact and we can construct the diagram 
     \begin{equation}\label{pushideale}
         \begin{tikzcd}
          0 \ar[r] & \id_Z \ar[r] \ar[d] & \osn \ar[r] \ar[d] & \os_Z \ar[r] \ar[d] & 0 \\
          0 \ar[r] & f_*\id_Y \ar[r] \ar[d,twoheadrightarrow] & \osn\oplus\osn(-m) \ar[d,twoheadrightarrow] \ar[r]  & f_*\os_Y\ar[r] & 0 \\
         & \cc \ar[r] & \osn(-m) & &
         \end{tikzcd}
     \end{equation}
     where $Z$ is the schematic image of $Y$, i.e. $\id_Z:=ker(\osn\ra f_*\os_Y)$.
     In particular, by assumption we have $\os_Z\cong f_*\os_Y$.
Snake lemma forces the left column to be short exact and $\cc\cong \osn(-m)$.
Note that we have $ext^1(\osn(-m),\id_Z)=h^{1}(X,\id_Z(m))$, hence if this group is zero then the left column in \eqref{pushideale} is split, as desired.
\end{proof}

\section{Ulrich sheaves and matrix factorizations}\label{sezione0.2}

In this section we will deal with general properties of Ulrich sheaves on divisorial coverings of projective space.
We recap briefly the notation introduced before.
Let $X$ be a scheme and $f:X\ra \p^n$ a degree $d$ divisorial covering of $\p^n$.
Then $X$ is a divisor in $P:=\p(\osn\oplus\osn(m))$ such that $\os_P(X)\cong \os_P(d)$, with $\pi:P\ra\p^n$ the standard projection map and $\os_P(1)$ the relative hyperplane bundle.
We set $\os_X(1):=f^*\osn(1)$, not to be confused with $\os_P(1)|_X\cong \os_X(m)$, by \thref{divisorialcover}.
$\e$ will be an Ulrich sheaf of rank $r$ on $X$ always relative to $f$, or equivalently to $\os_X(1)$.

In the first part of this section we will describe a resolution for $\e$ as a sheaf on $P$ and outline its connection to the theory of matrix factorizations of the equation of $X$.
This will lead to a proof of the existence of Ulrich sheaves on divisorial coverings of $\p^n$.
In the second, we will prove the existence of rank $2$ Ulrich sheaves of some double coverings of $\p^3$.

\subsection{Ulrich sheaves on divisorial coverings} 

Suppose $f:X\ra\p^N$ were an embedding.
We know that $\e$ is Ulrich on $X$ if and only if $f_*\e$ admits a linear resolution as in
\[0\ra \osn(-c)^{r_c}\ra \dots \ra \osn^{r_0}\ra f_*\e\ra 0\] 
where $c=codim(X,\p^N)$, see \cite{EisSch}[Prop. 2.1].
Here we want to derive a similar result in the case when $f$ realizes $X$ as a divisorial covering of $\p^n$.
Since $X$ is a divisor in $P$, we expect a resolution of length $1$ for $\e$ inside that variety.
The next result is essentially contained in \cite{HanselkaKummer}[Prop. 8.1, Remark 8.2], see also \cite{Beadeterminantal}[Thm. A] for the case of divisors in projective space.
We restate it in the projective bundle case and for completeness we reproduce also the proof.
In the following $\gi_l$ is the identity matrix of order $l$.

\Prop\label{Ulrichdivisorial}
Suppose $f:X\ra\p^n$ is a divisorial covering of degree $d$ with equation $p=0$.
If there is a matrix $A$ of order $rd$ whose entries are elements in $H^0(P,\os_P(1))$ such that $p^r=det(A)$ then there is an $f$-Ulrich sheaf $\e$, which seen as a sheaf on $P$ fits in
\begin{equation}\label{resol''}
    0\ra \os_P(-1)^{rd}\xra{A} \os_P^{rd}\ra \e\ra 0.
\end{equation}
If $X$ is integral then $rk(\e)=r$ and also the converse holds, i.e. from an Ulrich sheaf $\e$ we can construct a matrix $A=t\gi_{rd}-A'$, where $A'$ has coefficients in $H^0(\p^n,\osn(m))$.
\One
\begin{proof}
    A matrix as $A$ gives a morphism $\os_P(-1)^{rd}\ra\os_P^{rd}$.
    It is injective since the determinant of the above matrix is $p^r$ hence non-zero on the generic point of $P$.
    Therefore, if we define $\e$ to be the cokernel of the just determined morphism, we obtain the above sequence.
    Being $f:X\ra\p^n$ finite and surjective, it is enough to verify that $f_*\e\cong \pi_*\e$ is a direct sum of trivial bundles, from \thref{proiezione}.
    This can be done by simply taking $\pi_*$ in the above resolution since $\pi_*\os_P\cong\osn$, $\pi_*\os_P(-1)=0= R^1\pi_*\os_P(-1)$ by \cite{GW2}[Thm. 22.86].
    If $X$ is integral then $c_1(\e)\sim rk(\e)[X]$ hence 
    \[d\cdot rk(\e)c_1(\os_P(1))\sim rk(\e)[X]\sim c_1(\e)\sim c_1(\os_P^{rd})-c_1(\os_P(-1)^{rd})\sim rd\cdot c_1(\os_P(1)).\]

For the converse, suppose $\e$ is Ulrich then $f_*\e\cong\osn^{dr}$.
Set $\car:=\bigoplus_{l\in \N}H^0(\p^n,\osn(l))$ with $\car_l:=H^0(\p^n,\osn(l))$ its degree $l$ part.
We know that \[\bigoplus_{i\in \N}H^0(P,\e(i))\cong\bigoplus_{i\in \N}H^0(\p^n,f_*\e(i))\cong\bigoplus_{i\in \N}H^0(\p^n,\osn(i)^{\oplus rd})\cong\car^{dr}\]
where $f_*\e(i)\cong \osn(i)^{rd}$ by projection formula, being $\os_X(i)\cong f^*\osn(i)$.
Moreover, $\bigoplus_{l\in \N}H^0(X,\os_X(l))$, and hence $t\in H^0(P,\os_P(1))$, acts by multiplication on the sections of $\e(i)$.
In particular, this action of $t$ is $\car$-linear hence represented by a matrix $A'$ whose entries are in $\car$.
Being $t|_X\in H^0(X,\os_P(1)|_X)\cong H^0(X,\os_X(m))$ it sends $\car_l$ to $\car_{m+l}$, then the entries of $A'$ must be homogeneous of degree $m$ hence in $H^0(\p^n,\osn(m))$.
Since $p$ vanishes on $X$, it must annihilate $\e$ hence acts like $0$ on $\car^{rd}$.
    This implies that, by formally plugging $A'$ in the polynomial $p$ instead of the variable $t$ we must get $0$.
    It follows that the minimal polynomial of $A'$ divides $p$ and hence they are equal, since $p$ is irreducible being $X$ integral.
    In particular, since all irreducible factors of the characteristic polynomial of $A'$ must divide $p$ they are equal to $p$, so we get that $det(t\gi_{dr}-A')=p^r$.
    Therefore, let us define $A:=t\gi_{dr}-A'$.
    Finally, since $\e$ is globally generated we can choose a presentation $\os_P^{rd}\xra{ev}\e$.
    The matrix $A$ gives a morphism $\os_P(-1)^{rd}\ra\os_P^{rd}$ whose image, by definition, sits in the kernel of $ev$, therefore we get the following commutative diagram 
    \dia
    0 \ar[r] & \os_P(-1)^{rd} \ar[r,"A"] \ar[d] & \os_P^{rd} \ar[r] \ar[d,"id"] & Q  \ar[r] \ar[d,twoheadrightarrow] & 0 \\
    0 \ar[r] & \ck \ar[r] &  \os_P^{rd} \ar[r,"ev"] & \e  \ar[r] & 0. \\
    \mma 
    By the first implication, we know that $Q$ has rank $r$ and is Ulrich on $X$ hence it is torsion-free by \thref{corollario} being $X$ integral.
    Since the kernel of the map $Q\ra\e$ is a torsion subsheaf of $Q$ this must be $0$ hence $Q\cong \e$ and we recover \eqref{resol''}.
\end{proof}

\Oss\label{remark'}
The converse implication in this theorem also follows by a direct cohomological computation.
Indeed, any Ulrich sheaf $\e$ is globally generated by \thref{corollario} and $h^0(X,\e)=dr$ by \thref{hilbertpol}, hence we can construct a surjection $\os_P^{rd}\ra\e$ and we just need to show that the kernel $\ck$ is isomorphic to $\os_P^{rd}(-1)$.
By \cite{Orlov_sod}[Thm. 2.6 and Cor. 2.7] we know that the bounded derived category of coherent sheaves on $P$ has the following semi-orthogonal decomposition:
  \[\mathcal{D}^b(P)=<\os_P(-1)\otimes\pi^*\mathcal{D}^b(\p^n),\pi^*\mathcal{D}^b(\p^n)>=\]
  \[=<\os_P(-1),\os_P(-1)\otimes\pi^*\osn(H), \dots, \os_P(-1)\otimes\pi^*\osn(nH), \os_P, \dots, \pi^*\osn(nH)>\]
  where we used the standard Beilinson decomposition $\mathcal{D}^b(\p^n)=<\os_{\p^n}, \dots, \osn(nH)>$.
  Therefore, to verify $\ck\cong \on{n}(-1)^{dr}$ is enough to prove that $\ck$ is orthogonal to all the above pieces of $\mathcal{D}^b(P)$ except the first one.
  The special form of the matrix $A$ can be recovered up to linear operations on rows and columns, which correspond just to a change of basis in $H^0(P,\os_P^{rd}(-1))$ and $H^0(P,\os_P^{rd})$.
\one

We recover this well-known fact, see  \cite{Beadeterminantal}[Thm. A], \cite{CMRPL}[Thm. 4.2.7].
\Co\label{hypersup}
The datum of an Ulrich sheaf of rank $r$ on an integral hypersurface in $\p^{n+1}$ is equivalent to a linear determinantal representation of the $r$-th power of its equation.
\io
\begin{proof}
    Just consider the case $m=1$ in \thref{divisorialcover}.
\end{proof}

The parameter space for matrices $A$ as in \thref{Ulrichdivisorial} is an affine variety and the condition on the determinant is algebraic hence we get the following.

\Prop\label{constructible}
Fixed $n,m,d,r$, the locus in $H^0(P_m,\os_{P_m}(d))$ corresponding to divisorial coverings having a rank $r$ Ulrich sheaf is constructible (possibly empty); if $r=1$ it is also irreducible.
\One
\begin{proof}
Consider $H^0(P,\os_P(1))^{(rd)^2}$, the space parametrising square matrices of order $rd$ and entries in $H^0(P,\os_P(1))$.
Taking the determinant gives an algebraic map $H^0(P,\os_P(1))^{(rd)^2}\ra H^0(P,\os_P(rd))$ whose image, call it $T$, is clearly irreducible, and is constructible by Chevalley's theorem \cite{GW1}[Thm. 10.20].
If $r=1$ we are done, since the intersection of $T$ and the open set of divisorial coverings is still constructible and irreducible.
Otherwise, we take the intersection of $T$ with the subvariety of $H^0(P,\os_P(rd))$ parametrising $r$-th powers of polynomials in $H^0(P,\os_P(d))$, which is again constructible.
\end{proof}

\Oss
Note that in general the above loci are neither open nor closed.
Indeed, any time the determinant morphism is not dominant this locus is not open, for example a general complex quartic surface in $\p^3$ has Picard group generated by the class of an hyperplane section and hence admits no Ulrich line bundles by \thref{pic1}.

On the other hand, a general cubic surface in $\p^3$ admits Ulrich line bundles with respect to the hyperplane polarisation, for example any smooth one does.
But there are some normal but singular ones which do not, see \cite{Dolgachev}[Thm. 9.3.6, Table 9.2 case XX and also Table 9.1].
\one

\subsection{Matrix factorizations and Ulrich sheaves on divisorial coverings}

We want to focus our attention on the matrix $A$ and on its relation to the polynomial $p$.
This will be particularly important in the case of cyclic coverings.
We start by recalling an important definition.
\De
Let $\mathfrak{R}$ be a commutative unital ring and $\mathfrak{S}\subset \mathfrak{R}$ any subset.
A \textbf{matrix factorization} of an element $x\in \mathfrak{R}$ with coefficients in $\mathfrak{S}$ is the datum of $A_1,\dots , A_l$ matrices with entries in $\mathfrak{S}$ such that $\prod_{i=1}^l A_i=x\gi_n$ for some $n\in \N$.
If all $A_i$ are square matrices of the same order then this order will be the \textbf{size} of the factorization.
Finally, if all $A_i$ are equal we will call this factorization \textbf{cyclic}.
\Ne

Note that if $A_i$ are square matrices then
\[\prod_i det(A_i)=det\left(\prod_iA_i\right)=det(x\gi_n)=x^n\]
i.e any matrix factorization of $X$ produces determinantal representations for some powers of $x$.

From \cite{Eis-resolutions}[Cor. 6.3], the existence of maximally Cohen-Macaulay modules on the quotient of a regular local ring by the ideal generated by a single element $x$ is equivalent to the existence of a matrix factorizations of $x$.
Moreover, the matrices appearing in the factorization of $x$ give a resolution of the above module, similarly to what happens in \eqref{resol''}.

We will extend those results to divisorial covers, which are divisors in special weighted projective spaces.
We start by restating our previous results using the language of matrix factorizations.

\Co\label{factorulrich}
Consider a divisorial covering $f:X\ra\p^n$ of degree $d$ and with equation $p=0$.
If there is a matrix factorization of size $rd$ of the polynomial $p$ such that at least one of the matrices has coefficients in $H^0(\os_P,\os_P(1))$ then there is an Ulrich sheaf $\e$ on $X$.
If $X$ is integral then $rk(\e)=r$ and the converse holds, i.e. if there is an $f$-Ulrich sheaf $\e$ of rank $r$ on $X$ then there is a matrix factorisation of size $rd$ of $p$ with coefficients in $H^0(\os_P,\os_P(1))$.
\io
\begin{proof}
Call $A$ the matrix of the factorization of $p$ having coefficients in $H^0(P,\os_P(1))$, whose existence is granted by hypothesis.
    How already observed, it is clear that $det(A)=p^r$ so we conclude by \thref{Ulrichdivisorial}.
    
Conversely, if there is an Ulrich sheaf $\e$ then we have the resolution \eqref{resol''}.
The matrix $p\gi_{rd}$ induces a morphism $\os_P(-X)^{rd}\cong\os_P(-d)^{rd}\ra\os_P^{rd}$, whose image goes to $0$ when composed with the map $\os_P^{rd}\ra\e$.
    Therefore, we have a commutative diagram as in
    \dia 
    0 \ar[r] & \os_P(-d)^{rd} \ar[r,"id"] \ar[d] & \os_P(-d)^{rd} \ar[d,"p\gi_{rd}"] & &  \\
    0 \ar[r] & \os_P(-1)^{rd} \ar[r,"t\gi_{rd}-A"] &  \os_P^{rd} \ar[r] & \e  \ar[r] & 0. \\
    \mma 
    We conclude that there is a matrix $A'$, determined by the just constructed left vertical map, of size $rd$ such that $(t\gi_{rd}-A)A'=p\gi_{rd}$, hence we have the desired matrix factorization.
\end{proof}

\Co\label{matrixcyclic'''}
Suppose that $X$ is an integral cyclic covering with equation $p=t^d-b$.
There is an $f$-Ulrich sheaf $\e$ of rank $r$ on $X$ if and only if there is a cyclic matrix factorization of the polynomial $b$ of size $dr$ and with coefficients in $H^0(\osn,\osn(m))$.
\io
\begin{proof}
    If there is such an $\e$ then we already constructed in \thref{Ulrichdivisorial} a matrix $A'$ giving $0$ when plugged into $p$, i.e. such that $(A')^d=bId_{dr}$.
    Conversely, a cyclic matrix factorization of $b$ is the datum of a matrix $A'$ of size $rd$, with coefficients in $H^0(\osn,\osn(m))$ such that $(A')^d=bId_{dr}$.
    We deduce that the minimal polynomial of $A'$ must divide, and hence be equal to, $t^d-b$.
    This implies that its characteristic polynomial $det(tId_{dr}-A')$ must be $(t^d-b)^r$, hence we construct an Ulrich sheaf from \thref{Ulrichdivisorial}.
\end{proof}

As shown in the last statement of \cite{BHU}[Thm. 1.2], from any matrix factorization of some polynomial we can obtain a cyclic one, however, this comes at the cost of multiplying the rank of the matrix by $d$.

Now we prove that any divisorial covering admits Ulrich sheaves.
It is possible that this result is already known among experts, since the analogous case for modules over Cohen-Macaulay local rings follows from \cite{BHU}[Thm. 2.5].
However, we could not find any precise reference in the literature except for the case of divisors in $\p^{n+1}$, which is worked out in \cite{CMRPL}[Cor. 4.2.7 and Prop. 4.2.12].

\Te\label{ulrichmatrix''}
Fix some field $\gk$.
Any divisorial covering $f:X\ra\p^n$, in particular any cyclic one, admits an $f$-Ulrich sheaf.
Moreover, if $X$ is integral, $d=deg(f)$ and for some $s\geq 2$\footnote{when $s=1$ we must have $X$ reducible, actually a union of copies of $\p^n$} we have
\[p=\sum_{i=1}^s\prod_{j=1}^dp_{i,j},  \qquad p_{i,j}\in H^0(P_m,\os_{P_m}(1))\] 
then there are Ulrich sheaves of rank $d^{s-2}\varphi(d)$, where $\varphi$ is Euler's totient function.
If $\gk$ contains $d$-th roots of unity then the rank can be taken to be $d^{s-2}$.
\Ma
\begin{proof}
Call $p\in H^0(P,\os_P(d))$ the equation of $X$.
From \thref{factorulrich}, to get an Ulrich sheaf it is enough to construct a matrix factorization of $p$ with coefficients in $H^0(P,\os_P(1))$.
Since $H^0(P,\os_P(1))$ generates $\bigoplus_{i\in\N}H^0(P,\os_P(i))$ as a $\gk$-algebra, we can always write $p$ in the above form for some $s$ and $p_{i,j}\in H^0(P,\os_P(1))$.
But then we just apply \cite{BHU}[Lem. 1.6], or \cite{BHU}[Lem. 1.5] in the case where $\gk$ contains $d$-th roots of unity, and we get the existence and the claimed ranks.
\end{proof}

The following already appeared in \cite{CMRPL}[Porposition 4.3.1] with essentially the same proof.
See \cite{BHU}[Prop. 2.5] for the case of rings.
\Co\label{intcomplete}
Any complete intersection in $\p^N$ admits Ulrich sheaves respect to $\osn(1)$.
\io
\begin{proof}
    Taking $m=1$ in \thref{ulrichmatrix''} we get that any divisor in $\p^N$ admits Ulrich sheaves.
    Then, by applying iteratively \thref{intersezione}, we deduce the same for any complete intersection.
\end{proof}

The next seems to be a new result.
Essentially, we are saying that some divisorial coverings of complete intersections admit Ulrich sheaves.
Actually, in dimension $\ge3$ those are almost all the possible such coverings, see the remark below.
\Co\label{divci}
Suppose $X_2\ra\p^N$ is a divisorial covering and $X_1\subset \p^N$ a complete intersection.
If we denote by $g:X\ra \p^N$ the finite morphism coming from the pullback diagram 
\dia
X:=X_1\times_{\p^N}X_2 \ar[r] \ar[d] & X_2 \ar[d,twoheadrightarrow] \\
X_1 \ar[r, hookrightarrow] & \p^N \\
\mma
then $(X,g^*\on{N}(1))$ admits Ulrich sheaves.
\io
\begin{proof}
    \thref{ulrichmatrix''} and \thref{intcomplete} give the existence of Ulrich sheaves on $X_1,X_2$ and by \thref{pullback} their box product is Ulrich on $X$.
\end{proof}
\Oss
We could define divisorial and cyclic coverings $g:X\ra X_1$ of an arbitrary scheme $X_1$ analogously to what we did for $\p^n$ in \thref{defdiv} and \thref{defcyclic}.
For a complete intersection $X_1$ of dimension $\ge 3$ we know that $Pic(X_1)$ is generated by $\on{N}(1)|_{X_1}$, see \cite{SGA2}[XII Cor. 3.7], and that the restriction morphism $H^0(\p^N,\on{N}(m))\ra H^0(X_1,\os_{X_1}(m))$ is surjective for all $m$.
    Therefore, any divisorial covering of $X_1$ can actually be written as a pullback in the above form.
\one
Recalling \thref{horikawa} we get the next corollary which, as far as we know, is the first example of the existence of Ulrich sheaves on Horikawa surfaces which are not double coverings of $\p^2$. 
\Co\label{horikawa'}
All Horikawa surfaces which are double coverings of quadrics in $\p^3$ admit Ulrich sheaves.
\io
\begin{proof}
    Just take $X_1\subset \p^3$ to be the quadric and $X_2$ a double covering of $\p^3$ branched along a surface $B$ of degree $6$.
\end{proof}

Next we talk about low-degree coverings.
The first is a generalisation of \cite{KuNaPa}.
\Co
Flat double coverings $f:X\ra\p^n$ have $f$-Ulrich sheaf if $char(\gk)\neq 2$.
\io
\begin{proof}
    Since by \thref{2cyc} any such covering is cyclic, hence divisorial, when $char(\gk)\neq 2$, we can apply the previous corollary.
\end{proof}

Similarly, in almost all dimensions, degree $3$ coverings are divisorial.
\Co\label{triple}
For $n\geq 4$, every smooth triple covering $f:X\ra\p^n$ over $\C$ has an $f$-Ulrich bundle.
\io
\begin{proof}
     This $X$ can always be embedded as a divisor in the total space of a line bundle over $\p^n$ by \cite{Laz}[Prop. 3.2] and hence is a divisorial covering by definition.
\end{proof}

We cannot hope for similar results for higher degree coverings because Lazarsfeld gave examples of degree $5$ covering of $\p^n$ of arbitrary high dimension which are not divisorial, see \cite{Laz}[Remark 3.5].
Still, the question of lower-dimensional triple coverings remains open.
\Do
Does every (flat) triple covering of $\p^2$ or $\p^3$ carry an Ulrich sheaf?
\da

\subsection{Cyclic coverings}

Next, we will specialise in the case of cyclic coverings.
If we are interested in the minimal rank of an Ulrich sheaf on a divisorial covering $f:X\ra\p^n$ then, from the previous results, we need to search for matrix factorisation of minimal size for its equation $p$.
This minimal rank highly depends on the specific form of $p$, indeed, in the proof of \thref{ulrichmatrix''} we exploited the fact that $p$ can be written as a sum of $s$ products of forms in $H^0(P,\os_P(1))$.
For cyclic coverings the equation of $X$ is of the form $p=t^d-b$ with $b\in H^0(\p^n,\osn(dm))$, so we can find a more suitable way to express $p$.
For the sake of simplicity we work out only the case in which $\gk$ contains $d$-th roots of unity, but in a similar fashion and using \cite{BHU}[Lem. 1.6], we can also treat the general one.

Our point is that putting the last summand in the equation for $p$ in \thref{ulrichmatrix''} equal to $t^d$, we can further add a term $p_0^d$ without increasing the rank of the matrix needed.

\Prop\label{ulrichsommeprodotti}
Let $f:X\ra\p^n$ be an integral cyclic covering of degree $d$ over a field containing $d$-th roots of unity.
If there is some $s\ge2$ such that the equation $b$ of the branch locus of $f$ can be written as $b=p_0^d+\sum_{i=1}^{s-1}\prod_{j=1}^d p_{i,j}$, where $p_{i,j}\in H^0(\p^m,\osn(m))$, then there exists an $f$-Ulrich sheaf of rank $d^{s-2}$ on $X$.
\One
\begin{proof}
Fix a $d$-th root of unity $\zeta$.
Then $t^d-p_0^d=\prod_{j=1}^d(t-\zeta^jp_0)$ therefore $p=t^d-b$ can be written as a sum of $s$, instead of $s+1$, products of forms in $H^0(P,\os_P(1))$.
By \thref{ulrichmatrix''} there is an Ulrich sheaf of rank $d^{s-2}$.
\end{proof}

As an example we consider the lowest possible values for $s$.
\Es
For $s=2$ then we have $p=t^d-b=t^d-p_0^d-\prod_{j=1}^dp_j$ and there is a rank $1$ Ulrich sheaf on such on $X$.
For $d=2,3$, following the construction in \cite{BHU} we get that $p$ can be expressed, respectively, as the determinant of the following matrices
\begin{equation}
    \begin{pmatrix}\label{det21}
    t-p_0 & p_1\\
    p_2  &  t+p_0\\
\end{pmatrix} \qquad \quad \begin{pmatrix}
    t-p_0 & p_1 & 0 \\
    0 & t-\zeta_3p_0 & p_2 \\
    p_3  & 0 & t-\zeta_3^2p_0\\
\end{pmatrix}.
\end{equation}
In general, we take a matrix having $t-\zeta_d^ip_0$ on the diagonal and the polynomials $p_i$ on the above-diagonal and left down corner.

For $s=3$ and $d=2$ we have $b=p_0^d+\sum_{i=1}^{2}\prod_{j=1}^d p_{i,j}$ and there is a rank $2$ Ulrich sheaf on $X$.
Therefore, we can write $(t^d-b)^2$ as a determinant but even more, as in \cite{Bea}[Cor. 2.4], this polynomial can be seen as a Pfaffian of some skew-symmetric matrix.
Let us illustrate this in case $d=2$. 
We have $t^2-b=t^2-p_0^2-p_1p_2-p_3p_4$ hence from the construction in \cite{BHU} we get $(t^2-b)^2$ equal to
\[det\begin{pmatrix}
    t-p_0 & 0 & p_1 & p_3 \\
    0 & t-p_0 & -p_4 & p_2 \\
    p_2  & -p_3 & t+p_0 & 0 \\
    p_4  & p_1 & 0 & t+p_0 \\
\end{pmatrix}=det\begin{pmatrix}
     0 & t-p_0 & -p_4 & p_2 \\
    -(t-p_0) & 0 & -p_1 & -p_3 \\
    p_4  & p_1 & 0 & t+p_0 \\
    -p_2  & p_3 & -(t+p_0) & 0 \\
\end{pmatrix}\]
where the second matrix is obtained simply by exchanging first with second and third with fourth row, then multiplying by $-1$ the second and fourth.
\Io

Let us illustrate the state of the art for double coverings.
\Oss
Let us suppose $d=2$.
The best upper bound for the minimal rank of an Ulrich bundle on a smooth double covering of $\p^n$ was the one given in \cite{KuNaPa}[Lem. 6.2], their result reads as follows.
If the equation of the branch locus is $b=\sum_{i=1}^{s-1} p_{i,1}p_{i,2}$ then there is an Ulrich sheaf of rank $2^{s-1}$ or $2^{s-2}$.
Our proof implies that to have an Ulrich sheaf of rank $2^{s-2}$ is sufficient that $b=p_0^2+\sum_{i=1}^{s-1} p_{i,1}p_{i,2}$.
Although it is only a slight improvement, it will be optimal in the cases $n=3$ and $m=2,3,4$, see \thref{double3folds}.
\one

\subsection{Minimal rank Ulrich sheaves on some double covering of $\p^3$}

Suppose $char(\gk)\neq 2$.
In this subsection, we will show that the minimal rank of an Ulrich sheaf on the general double covering of $\p^3$ branched along divisors of degree $4,6,8$ is $2$.
The main task is to construct rank $2$ Ulrich sheaves. 
Here we can actually specify what "general" means: we want the equation $b$ of $B$ to be written as $p_0^2+p_1p_2+p_3p_4$, i.e. exactly the form given in \thref{ulrichsommeprodotti} for $s=3$.
We will show that the space of polynomials $b\in H^0(\p^n,\osn(2m))$ satisfying this condition actually contains an open subset.

First, we will prove that a general such $X$ does not admit Ulrich sheaves of rank $1$ for any $n\geq 3$ and $m\geq 2$.
In case $\gk=\C$, this follows from \thref{pic1} since a general $X$ is a smooth double covering of $\p^3$ and hence has cyclic Picard group generated by the pullback of $\on{3}(1)$, due to \thref{Piccover}.

\Le\label{non1}
Suppose $\gk$ is infinite and $P=\p(\osn\oplus\osn(-m))$.
The general $X\in |\os_P(2)|$ has no Ulrich sheaves of rank $1$ respect to $\os_X(1)$ if $n\geq 3$ and $m\geq 2$.
\ma
\begin{proof}
We can reduce to the case $n=3$: indeed, if some $f:X\ra\p^n$ with $n\geq 3$ has an Ulrich sheaf then by \thref{restrizione} on their hyperplane sections, which are double coverings of $\p^3$.
We can suppose $X$ integral, since it is general.
    By \thref{factorulrich} is enough to show that the general polynomial $p\in H^0(P,\os_P(2))$ cannot be written as determinant of a $2\times 2$ matrix with coefficients in $H^0(P,\os_P(1))$.
    The above $2\times 2$ matrices, if $n=3$,  form a vector space of dimension $4h^0(P,\os_P(1))=4(\binom{m+3}{3}+1)$ while the polynomials giving double coverings are $h^0(P,\os_P(2))=\binom{2m+3}{3}+\binom{m+3}{3}+1$ so we need only to prove
    \[\binom{2m+3}{3}+\binom{m+3}{3}+1>4\left(\binom{m+3}{3}+1\right)\]
    which is equivalent to 
    \[\binom{2m+3}{3}-3-3\binom{m+3}{3}>0.\]
    This polynomial in $m$ is easily seen to be always positive for $m\geq 2$.
\end{proof}

Now we turn to positive results.
The idea is to show that for $n=3$ and $m=2,3,4$ a general polynomial in $H^0(P,\os_P(2))$ can be written in the form needed to apply \thref{ulrichsommeprodotti}, that is, the variety parametrising these special polynomials has a dominant map to $H^0(P,\os_P(2))$.
We need some preparation.
Set $V_m:=H^0(\p^3,\osn(m))$.
We will always confuse homogeneous polynomials of degree $m$ with vectors in $V_m$.
\Le\label{diff'}
Consider the algebraic map
\[\phi:V_m^5\ra V_{2m} \qquad (p_0,\dots, p_4)\mt p_0^2+p_1p_2+p_3p_4.\]
Its differential, seen as a linear map $V_m^5\ra V_{2m}$, can be written as
\begin{equation}\label{differential}
    d_{(p_0,\dots, p_4)}\phi: (q_0,q_2,q_1,q_4,q_3)\mt 2p_0q_0+\sum_{i=1}^4p_iq_i. 
\end{equation}
\ma
\begin{proof}
To understand the above differential is enough to understand the differential of the two maps
\[\phi_\mathcal{S}:V_m\times V_m\ra V_{2m} \quad (p_1,p_2)\mt p_1\cdot p_2, \qquad \phi_\mathcal{V}:V_m\ra V_{2m} \quad p\mt p^2.\]
Recall that the tangent space to an affine space can be identified with the affine space itself, see for example \cite{GW1}[Example 6.4].
For any fixed point $(p_1,p_2)\in V_m\times V_m$ the tangent directions in it are parametrized by $V_m\times V_m$ itself, hence by the choice of another point $(q_2,q_1)$.
Varying $t\in \mathbf{k}$ we get a line $(p_1+tq_2,p_2+tq_1)$ which is sent by $\phi_\mathcal{S}$ to the curve $p_1\cdot p_2+t(p_1\cdot q_1+p_2\cdot q_2)+t^2p_1\cdot p_2$.
This implies that 
\[d_{(p_1, p_2)}\phi_\mathcal{S}: (q_2,q_1)\mt p_1\cdot q_1+p_2\cdot q_2 \]
Analogously, 
\[d_{p_0}\phi_\mathcal{V}: q_0\mt 2p_0q_0. \]
hence by summing up the contributions from all the factors we get the claim.
\end{proof}

\Prop\label{join}
Let $\mathbf{k}$ be an infinite field with $char(\mathbf{k})\neq 2$ and suppose $m=2,3,4$.
\begin{enumerate}[i)]
    \item The ideal generated by five general  $p_i\in V_m$ contains all $V_{2m}$.
    \item Given a general homogeneous polynomial $p\in V_{2m}$ of degree $2m$ there are $p_i\in V_m$ such that $p=p_0^2+p_1p_2+p_3p_4$.
\end{enumerate}
\One
\begin{proof}
First we show how to reduce $ii)$ and $i)$ to a single computation.
Consider the algebraic map of affine spaces
\[\phi:V_m^5\ra V_{2m} \qquad (p_0,\dots, p_4)\mt p_0^2+p_1p_2+p_3p_4.\]
If we can show that its differential is surjective in a point in $V_m^5$ then $\phi$ would be smooth in that point by \cite{GW2}[Thm. 18.74].
Being the smooth locus open, this would imply that in some open neighbourhood of this point $\phi$ would be smooth hence open, so $\phi$ would be dominant, giving $ii)$.
But, since we already computed the exact form of $d\phi$ in \thref{diff'}, we infer that the surjectivity of $d\phi$ in a general point is exactly condition $i)$.
Being $char(\gk)\neq 2$ we can ignore the coefficient $2$, so we need to prove that for general $p_i\in V_m$ with $i=0,\dots, 4$ every element in $V_{2m}$ can be written as $\sum_{i=0}^4p_iq_i$.

To achieve this we will use Macaulay2.
I thank Fulvio Gesmundo for suggesting this approach.
To prove the claim over a field of characteristic $0$ is enough to have it over $\Q$.
The cases $m=2,3$ could also be easily carried over by hand.
    \begin{lstlisting}[language=Macaulay2]
--Define our polynomial ring
R=QQ[x,y,z,w]
--Consider three special ideals generated by 5 polynomials of degree 2,3,4
I=ideal{x^2,y^2,z^2,w^2,x*y}
J=ideal{x^3,y^3,z^3,w^3,(x+y+z+w)^3}
K=ideal{x^4,y^4,z^4,w^4,(x+y+z+w)^4}

--Let us compute the dimension of the graded pieces (Hilbert function) of 
--degree 4,6,8 of the quotients of R by those ideals
 hilbertFunction(4,I)
 hilbertFunction(6,J)
 hilbertFunction(8,K)
    \end{lstlisting}
    
To prove the claim over an arbitrary field, we first check the same statement but in the ring $\Z[x,y,z,w]$, then it is enough to tensor with $\gk$ to get the thesis.
In this case we need more complicated polynomials.
\begin{lstlisting}[language=Macaulay2]
    --Define our polynomial ring
R=ZZ[x,y,z,w]
--Consider three special ideals generated by polynomials of degree 2,3,4
I=ideal{x^2,y^2,z^2,w^2,x*y}
J=ideal{x^3+y*z*w,y^3+z*w*x,z^3+w*y*x,w^3+x*y*z,x^2*z}
K=ideal{x^4+x^3*y+x^2*y*z,y^4+y^3*z+y^2*z*w,z^4+z^3*w+z^2*w*x,
w^4+w^3*x+w^2*x*y,x*y*z*w+x^2*y^2+x^2*w^2+z^2*w^2+y^2*z^2+y^2*w^2+x^2*y*z}

--In this case, the Hilbert function is not enough since the quotient of 
--R by those ideals could have torsion, but we can compute a basis for it
 basis(4,R/I)
 basis(6,R/J)
 basis(8,R/K)
 \end{lstlisting}
In the first two cases we are done while in the last the quotient is not $0$. Nevertheless, since we work over a field $\gk$ of characteristic different from $2$, is actually enough to check that $w^8$ is $2$-torsion, so that when we tensor with $\gk$ this quotient becomes $0$.
\begin{lstlisting}[language=Macaulay2]
    (2*w^8) % J==0
    \end{lstlisting}
    \[\]
\end{proof}
\Oss
Part $ii)$ can be restated as a theorem about joins and secant of varieties in $\p(V_{2m})$.
The variety in $\p(V_{2m})$ parametrizing squares of polynomials is the image of the degree $2$ Veronese variety of $\p(V_m)$ under a suitable morphism.
Similarly, the locus in $\p(V_{2m})$ parametrizing classes of product of two polynomials is the image of the Segre variety $\p(V_m)\times \p(V_m)$.
Our claim is that the join of this two varieties is the entire $\p(V_{2m})$.
\one

Now, our main theorem for this section is just a composition of the above lemmas.

\Te\label{ulrichr2doublecoverexistence}
Suppose $char(\gk)\neq 2$ and $\gk$ infinite.
Let $f:X\ra\p^n$ be a general cyclic double covering with branch locus of degree $2m=4,6,8$.
The minimal rank of an $f$-Ulrich sheaf $\e$ is exactly $2$ and those sheaves are stable.
Moreover, if $X$ is smooth then $\e$ is locally free.
\Ma
\begin{proof}
    We have seen in \thref{non1} that a general $X$ has no Ulrich sheaf of rank $1$.
    On the other hand, existence of rank $2$ Ulrich sheaves on any $X$ whose branch divisor can be written as $b=p_0^2+p_1p_2+p_3p_4$ is granted by \thref{ulrichsommeprodotti}.
    Actually, the general polynomial in $b\in H^0(\p^3,\osn(2m))$ for $m=2,3,4$ can be expressed in that way by \thref{join}.
    Moreover, if there are no Ulrich sheaves of rank $1$ on $X$ then the rank $2$ ones must be stable by \thref{stability}.
    If $X$ is smooth then $\e$ must be locally free by \thref{corollario}.
\end{proof}

\subsection{Double coverings of the plane}

This last part is just an appendix on double coverings of $\p^2$ inserted for completeness.
We just review already known examples.

The following was originally proved in \cite{SebastianTripathi}[Thm. 1.2].
\Prop
Over $\C$, a general double covering $f:X\ra\p^2$ admits $f$-Ulrich sheaves of rank $2$.
\One
\begin{proof}
    By \thref{ulrichsommeprodotti} its enough to write any $b\in H^0(\p^2,\on{2}(2m))$ in the form $p_0^2+p_1p_2+p_3p_4$ for $p_i\in H^0(\p^2,\on{2}(m))$; in fact, we will even do it with $p_0=0$.
    This is granted for $X$ general by the first point in \cite{CCG1}[Thm. 5.1].
\end{proof}

For completeness, we also prove the following, which is in \cite{PN_qdp}[Thm. 1.4].
\Prop
The general double covering $f:X\ra\p^2$ branched along a curve of degree $2m$ admits no $f$-Ulrich sheaves of rank $1$ if $m\ge3$.
\One
\begin{proof}
    We can assume $X$ irreducible of equation $t^2-b$.
    If there is a rank $1$ Ulrich sheaf on $X$ then by \thref{Ulrichdivisorial} we can write $t^2-b=det(t\gi_2+A)$ for some $2\times2$ matrix $A$. 
    It follows easily that 
    \[    A'=\begin{pmatrix}
    -p_0 & p_1\\
    p_2  &  p_0\\
\end{pmatrix}, \]
therefore its enough to show that a general polynomial $b\in H^0(\p^2,\on{2}(2m))$ is not of the form $p_0^2+p_1p_2$ for $p_i\in H^0(\p^2,\on{2}(m))$.
Consider the map 
\[\varphi:H^0(\p^2,\on{2}(m))^3\ra H^0(\p^2,\on{2}(2m))\qquad (p_0,p_1,p_2)\mt p_0^2+p_1p_2.\]
By \thref{diff'} we know that its differential is
\[d_{(p_0,p_1,p_2)}\phi: (q_0,q_2,q_1)\mt 2p_0q_0+p_1q_1+p_2q_2\] 
and we see that for general $(p_0,p_1,p_2)$ its kernel has dimension at least $3$, since contains $(p_1,-2p_0,0),(p_2,0,-2p_0),(0,p_2,-p_1)$ and they are independent.
Therefore, we conclude that $\varphi$ is not dominant since
\[dim(Im(\varphi))\le 3\binom{m+2}{2}-3<\binom{2m+2}{2}=h^0(\p^2,\on{2}(2m)) \; \Longleftrightarrow\;  m>2.\]
\end{proof}

\section{Geometry of Ulrich sheaves on double coverings}

In the following, we adopt the previous notation and suppose $\gk$ to be algebraically closed and $char(\gk)\neq 2$, in particular infinite, even though some of the results only need a weaker assumption.
We will work with double coverings $f:X\ra\p^n$.
Recall that, under our assumptions, any flat double covering is cyclic by \thref{2cyc} and the closed points fixed by the involution $\iota$ form the ramification divisor $R$ of $f$.
To simplify the statements, we usually assume $n\ge 3$.

In this section, we study some geometrical properties of an Ulrich sheaf $\e$ on an integral double covering: the zero loci of their sections, the interaction with the involution $\iota$ naturally associated with $X$ and their non-ample loci.
Those results will be helpful later in understanding the moduli spaces of Ulrich bundles on $X$.

\subsection{Syzygy sheaves of Ulrich sheaves on double coverings}

If we consider the embedding $X\subset P_m=p(\osn\oplus\osn(m))$ of a double covering then \thref{Ulrichdivisorial} gives a projective resolution on $P_m$ of an Ulrich sheaf on $X$. 
Restricting it to $X$ it is no longer left exact but it could be extended $2$-periodically as done in \cite{Eis-resolutions}[Thm. 6.1, Cor. 6.3] for the case of modules.
We are mainly interested in the syzygy sheaf of our Ulrich sheaf.

\Prop\label{Discount}
Let $\gk=\overline{\gk}$ and $char(\gk)\neq 2$.
Suppose $f:X\ra\p^n$ is an integral, double covering with $n\geq 1$ and branch locus of degree $2m$ and call $\iota$ the involution of this covering.
If $\e$ is an $f$-Ulrich sheaf of rank $r$ then, taking the evaluation morphism of global sections of $\e$ we get
\sesl{\f(-m)}{\os_X^{2r}}{\e}{discount}
Then $\f\cong \iota^*\e$ is another $f$-Ulrich sheaf of rank $r$.
\One
\begin{proof}
Recall that $X$ can be seen as a divisor of equation $p=t^2-b$ inside $P$ for some $b\in H^0(\p^n,\osn(2m))$, in particular $X\in |\os_P(2)|$.
The equation \eqref{discount} can be put together with the diagram in the proof of \thref{factorulrich} in the commutative diagram below
\dia
0 \ar[r] & \os_P(-2)^{2r} \ar[r,"id"] \ar[d,"M"] & \os_P(-2)^{2r} \ar[d,"p\gi_{2r}"] & &  \\
    0 \ar[r] & \os_P(-1)^{2r} \ar[d] \ar[r,"t\gi_{2r}-A"] &  \os_P^{2r} \ar[d] \ar[r] & \e  \ar[r] \ar[d,"id"] & 0 \\
    0 \ar[r] & \f(-m) \ar[r]  &  \os_X^{2r} \ar[r] & \e  \ar[r] & 0 \\
    \mma
where $M$ is defined to be the matrix with coefficients in $H^0(P,\os_P(1))$ resulting from the diagram.
Twisting the left column by $\os_P(1)$, and recalling that $\os_P(1)|_X\cong \os_X(m)$ by \thref{divisorialcover}, we get immediately that $\f$ is Ulrich by \thref{Ulrichdivisorial}.
We can write $M=tM_1+M_2$ where $M_i$ have coefficients in $H^0(\p^n,\osn(m))$, in particular do not contain $t$.
Commutativity of the higher square implies that
\[(t^2-b)\gi_{2r}=p\gi_{2r}\circ id_{\os_P(-2)^{2r}}=(t\gi_{2r}-A)(tM_1+M_2)=t^2M_1+t(M_2-AM_1)-AM_2\]
therefore we get $M_1=\gi_{2r}$, $M_2=A$ and $A^2=b\gi_{2r}$, hence $A$ gives a cyclic matrix factorization of $b$.
Recall that, by \thref{cycliciinvolution}, $\iota$ is exactly the restriction to $X$ of the automorphism $\iota_P$ of $P$ given by sending $t\mt-t$.
Moreover, we know that $M=t\gi_{2r}+A$ hence 
\[\iota_P(t\gi_{2r}-A)=-t\gi_{2r}-A=-M.\]
But, since the morphisms $\os_P(-2)^{2r}\ra \os_P(-1)^{2r}$ induced by $M$ and $-M$ differ for just a sign, the corresponding cokernels are isomorphic hence we get $\f\cong\iota^*\e$.
\end{proof}

\Oss
For the special case of Ulrich bundles on smooth complex quadrics, Ulrich bundles are the \textit{spinor bundles} and this result was already in \cite{Ottaviani_spinor}[Thm. 2.8].
Being an Ulrich bundle globally generated and $h^0(X,\e)=rd$ by \thref{hilbertpol}, a rank $2$ one determines a morphism $\varphi:X\ra Gr(2,4)$ such that $\varphi^*\mathcal{Q}=\e$, where $\mathcal{Q}$ is a spinor bundle, in particular, under those regularity assumptions, our result could be deduced by Ottaviani's one. 
\one

This result has many nice consequences in the study of moduli spaces of Ulrich sheaves.
As a first corollary we can understand the restriction to hypersurfaces in $|\os_X(m)|$.
For any coherent sheaf $\f$ and any $S\subset X$ we set $\f_S:=\f|_S$. 

\Co\label{restrdoppi}
In the setting of \thref{Discount}, take $S\in|\os_X(m)|$ and let $\e,\g$ be stable $f$-Ulrich sheaves on $X$.
If $\e_S\cong \g_S$ then either $\e\cong \g$ or $\iota^*\e\cong \g$.
\io
\begin{proof}
The strategy of the proof is the same used in \cite{Beafk3}[Prop. 6].
    Applying $Hom(-,\g)$ to \eqref{discount} and $Hom(-,\g_S)$ to the restriction of \eqref{discount} to $S$ we get the commutative diagram 
    \dia
    0 \ar[r] & Hom(\e,\g) \ar[d] \ar[r] &  H^0(X,\g)^{2r} \ar[d] \ar[r] & Hom(\iota^*\e(-m),\g)   \ar[d]  \\
    0 \ar[r] & Hom(\e_S,\g_S) \ar[r]  &  H^0(S,\g_S)^{2r} \ar[r] & Hom(\iota^*\e(-m)_S,\g_S)   \\
    \mma
    We need to prove that either $\iota^*\e\cong \g$ or the restriction map $Hom(\e,\g)\ra Hom(\e_S,\g_S)$ is surjective, since in this case the isomorphism $\e_S\cong\g_S$ lifts to a non-zero morphism from $\e$ to $\g$ which must be an isomorphism by stability, see \thref{stabilitàtrick}, since the two share the same reduced Hilbert polynomial, see \thref{hilbertpol}.
    The central vertical map comes from the cohomology of 
    \sesl{\g(-m)}{\g}{\g_S}{g} 
    and hence is an isomorphism being $\g$ Ulrich.   
    By diagram chasing $Hom(\e,\g)\ra Hom(\e_S,\g_S)$ is surjective if $Hom(\iota^*\e(-m),\g)\ra Hom(\iota^*\e(-m)_S,\g_S)$ is injective.
    This last morphism comes from applying $Hom(\iota^*\e(-m),-)$ to \eqref{g}, in particular its kernel is
    \[Hom(\iota^*\e(-m),\g(-m))\cong Hom(\iota^*\e,\g)\]
    which is either $0$ or, by stability, $\iota^*\e\cong\g$.
\end{proof}

\Co
In the setting of \thref{Discount}, if $\e$ is locally free and stable then either $\e\cong \iota^*\e\du(m)$ or the multiplication map $H^0(X,\e)\otimes H^0(X,\e)\ra H^0(X,\e\otimes\e)$ is injective.
\io
\begin{proof}
    Tensoring \eqref{discount} by $\e$ and recalling that the term $\os_X^{2r}$ is exactly $H^0(X,\e)\otimes\os_X$, we get
    \ses{\e\otimes\iota^*\e(-m)}{\e\otimes H^0(X,\e)}{\e\otimes\e}
    Taking global sections we get that the kernel of the multiplication map involved in the statement is exactly $H^0(\e\otimes\iota^*\e(-m))\cong Hom(\iota^*\e\du(m),\e)$. 
    Since $\omega_X\cong \os_X(m-n-1)$ by \thref{divisorialcover}, the vector bundle $\iota^*\e\du(m)$ is exactly the U-dual of $\iota^*\e$, in particular is Ulrich by \thref{dual}.
    From the stability of $\e$ it follows that either $\e\cong \iota^*\e\du(-m)$ or there are no homomorphism between them, as desired.
\end{proof}

In the very special case of quartic double solids, the following will be helpful to understand the smoothness of points representing Ulrich bundles in the corresponding moduli space.
\Co\label{extquartic}
In the setting of \thref{Discount}, if $n=3, m=1,2,3$ and $\e$ is locally free then $ext^2(\e,\e)\cong \hom(\iota^*\e(4-2m),\e)$.
\io
\begin{proof}
Note that $\iota^*\e$ is locally free being $\e$ such.
From \thref{divisorialcover} we know that $\omega_X\cong \os_X(m-4)$, so by Serre duality and Ulrich condition we get 
\[ext^i(\e,\os_X)=ext^i(\e(m-4),\omega_X)=h^{3-i}(\e(m-4))=0 \qquad  i=2,3.\]
    Appling $Hom(\e,-)$ to \eqref{discount} gives $ext^2(\e,\e)=ext^3(\e,\iota^*\e(-m))$ but with some standard manipulation and Serre duality we get
    \[ext^3(\e,\iota^*\e(-m))=ext^3(\os_X,\e\du\otimes\iota^*\e(-m))=h^3(X,\e\du\otimes\iota^*\e(-m))=\]
    \[=\hom(\e\du\otimes\iota^*\e(-m),\os_X(m-4))=\hom(\iota^*\e(4-2m)
    ,\e).\]
\end{proof}

\subsection{Zero loci of Ulrich bundles on double coverings}

We are going to describe the zero loci of sections of Ulrich sheaves, so that we can construct them using some version of Hartshorne-Serre's correspondence.

We begin with some standard computations that will be useful later.
Recall that we denote by $R\subset X$ the ramification locus of $f$.
As usual, we call $Z\subset \p^n$ a \textbf{complete intersection of type} $\mathbf{(m,m)}$ if we have a short exact sequence of the form
\ses{\osn(-m)}{\osn^2}{\id_Z(m)}
In particular, $h^0(\p^n,\id_Z(m))=2$ and for every $l\geq m$ this induces  surjective multiplication maps 
\[H^0(\p^n,\osn(l-m))\otimes H^0(\p^n,\id_Z(m))\ra H^0(\p^n,\id_Z(l)),\]
that is, the equation of any hypersurface of degree $l$ containing $Z$ (as a subscheme) is generated by the two of degree $m$.

\Le\label{numericsci}
Let $Z\subset \p^n\; (n\geq 3)$ be an $(m,m)$ complete intersection, then
\begin{itemize}
\begin{multicols}{2}
    \item $\omega_Z\cong\os_Z(2m-n-1)$
    \item $h^0(\p^n,\id_Z(m))=2$
    \item $h^0(Z,\os_Z(m))=\binom{m+n}{n}-2$
    \item $h^j(\p^n,\id_Z(i))=0,\;\; i\in \Z,\; 1\leq j\leq n-2$
    \item $h^0(\p^n,\id_Z(2m))=2\binom{m+n}{n}-1$
    \item $h^0(Z,\os_Z(2m))=\binom{2m+n}{n}-2\binom{m+n}{n}+1$
    \end{multicols}
    \item $P(\id_Z(m))(t)=2\binom{t+n}{n}-\binom{t+n-m}{n}$
\end{itemize}
\ma
\begin{proof}
    The first claim follows by adjunction formula and the fact that $\omega_{\p^n}\cong\osn(-n-1)$.
    Being $Z$ an $(m,m)$ complete intersection, we have an exact sequence
    \sesl{\osn(-2m)}{\osn(-m)^2}{\id_Z}{ci}
    Since $n\geq 3$ we have $h^j(\p^n,\osn(i))=0$ for $i\in \Z, 1\leq j\leq n-1$ hence, from the above sequence it immediately follows that $h^j(\p^n,\id_Z(i))=0$ for all $i\in \Z$ and $1\leq j\leq n-2$.
    Similarly, we can easily compute the global sections of $\id_Z(m),\id_Z(2m)$.
    The rest of the claims follows using the previous results and the cohomology of the adequate twists of
    \ses{\id_Z}{\osn}{\os_Z}
    Finally, the last claim follows by additivity of Hilbert polynomial in the sequence \eqref{ci} twisted by $\osn(m)$ and the fact that $P(\osn)(t)=\binom{t+n}{n}$.
\end{proof}

Recall that an Ulrich sheaf $\e$ of rank $2$ on a double covering $X$ is special if $c_1(\e)\sim \omega_X(n+1)$ by \thref{divisorialcover}.
To justify the assumptions in the next proposition, remember that if we have a finite surjective $f:X\ra\p^n$, in \thref{restrizionesezioni} we proved that the zero locus $Y$ of a section of an $f$-Ulrich bundle cannot contain schematic fibers of $f$, in particular, if $f$ has degree $2$ then $f|_Y$ is injective.

\Prop\label{section-rk2}
Suppose $\gk=\overline{\gk}$ and $char(\gk)\neq 2$.
Let $f:X\ra\p^n$ be an integral, double covering with $n\geq 3$ and consider $Y\subset X$ closed of codimension $2$ such that $f|_Y$ is an isomorphism on the image and $Z:=f(Y)$ is an $(m,m)$ complete intersection.
We have $f_*\id_Y\cong \id_Z\oplus\osn(-m)$ and any such $Y$ is the zero locus of a unique (up to scalar) section of a rank $2$ $f$-Ulrich sheaf $\e$, which is special, reflexive, satisfies $\e\du\cong \e(-m)$ and sits in 
\sesl{\os_X}{\e}{\id_Y(m)}{ext}
\One
\begin{proof} 
By \thref{numericsci} and \thref{divisorialcover} we have that 
\[\omega_Y(-mH)\otimes\omega_X\du\cong \os_Y(m-n-1)\otimes\omega_X\du\cong \os_Y\]
hence is globally generated by just one section.
By \thref{divisorialcover} $3)$ we deduce that $(X,\os_X(1))$ has no intermediate cohomology.
Therefore, as done in \thref{r00}, we can see that $h^0(Y,\omega_Y(-D)\otimes\omega_X\du)=ext^1(\id_Y(D),\os_X)$. 
Being $Y\cong Z$ connected, any such $Y$ gives rise to a unique (up to isomorphism) non-split extension \eqref{ext}.
By \thref{numericsci} we have $ext^1(\osn(-m),\id_Z)=h^{1}(X,\id_Z(m))=0$, hence by \thref{idealedoppi} we have $f_*\id_Y\cong \id_Z\oplus\osn(-m)$.

To prove that $\e$ is Ulrich we want to apply \thref{verificaulrich}.
Since $f$ is affine, hence preserves cohomology, and we have $f_*\os_X\cong \osn\oplus\osn(-m)$, it follows that
\[P(\id_Y(m))=P(f_*\id_Y(m))=P(\osn)+P(\id_Z(m))=\quad (\text{by \thref{numericsci}})\]
\[=3\binom{t+n}{n}-\binom{t+n-m}{n}=4\binom{t+n}{n}-P(\os_X).\]
In addition, for $1\leq j\leq n-2$ and for all $i\in\Z$ or for $j=0$ and $i<m$ we have
\[h^j(X,\id_Y(i))=h^j(\p^n,f_*\id_Y(i))=0.\]
Therefore, we are left to prove that $h^n(X,\e(-n))=0$, which will be a consequence of the non-splitting of \eqref{ext}.
By Serre duality
\[h^n(X,\e(-n))\cong hom(\e(-n),\omega_X)\cong h^0(X,\e\du(n)\otimes\omega_X\du)\cong h^0(X,\e\du(-m-1)).\]
Apply $\Hom(-,\os_X)$ to \eqref{ext}: we get an exact sequence which can be split in
\sesl{\os_X(-m)}{\Hom(\e,\os_X)}{\ck}{10}
\sesl{\ck}{\os_X}{\ext^1(\id_Y(m),\os_X)}{11}
which are everywhere exact except, possibly, for the right side of the second.
Since the extension class of \eqref{ext} is non-zero by definition, it follows that the right map in \eqref{11} is non-zero and hence surjective, being $Y$ connected and
\[\ext^1(\id_Y(m),\os_X)\cong \ext^2(\os_Y(m),\os_X)\cong \ext^2(\os_Y,\omega_X)\otimes\omega_X\du(-m)\cong\]
\[\cong\omega_Y(-m)\otimes\omega_X\du\cong \os_Y.\]
It follows that actually $\ck\cong \id_Y$ and hence \eqref{10} is a non-split extension of $\id_Y$ and $\os_X(-m)$, i.e. isomorphic to \eqref{ext} up to twisting by $\os_X(m)$.
In particular $\e\du\cong \e(-m)$ so that $h^0(X,\e\du(-m-1))=h^0(X,\e(-1))=0$.
From this we also deduce immediately that $\e\dd\cong (\e(-m))\du\cong \e$. 
\end{proof}

\Do
In case $X$ is not smooth, are the above sheaves $\e$ still locally free?
If not, are they at least in the subcategory of $\mathscr{D}^b(X)$ generated by perfect complexes?
This is especially interesting in view of the connection between matrix factorisations (and hence Ulrich sheaves) and the singularity category of $X$, see \cite{Orlovmatrix}[Thm. 3.9].
\da

Under some additional assumption, we have a converse to the previous statement.
The fact that $c_2(\e)\neq 0$ is expected, since this is true over $\C$ by \thref{c2=0}, hence that for any Ulrich $\e$ we can find a corresponding $Y$.

\Prop\label{rk2-section}
Suppose $f:X\ra\p^n$ is a smooth double covering over a field $\gk=\overline{\gk}$ with $char(\gk)\neq 2$ and $n\ge 3$.
Let $\e$ be a rank $2$, special $f$-Ulrich bundle, then a general global section of $\e$ gives an extension as in \eqref{ext}
\sesl{\os_X}{\e}{\id_Y(m)}{ext'}
for some $Y\subset X$ such that $f|_Y$ an isomorphism onto a complete intersection $Z$ of type $(m,m)$.
If $Pic(X)\cong \Z H$ then the same happens for all sections in $H^0(X,\e)$.
\One
\begin{proof}
     By \cite{Ful}[Example 12.1.11], a general section of $\e$ vanishes along a codimension $2$ subscheme $Y$ as in \eqref{ext'} since our base field is infinite and $\e$, being Ulrich on a smooth scheme, is globally generated and locally free.
If $Pic(X)\cong \Z H$ then this happens for any section by \thref{tuttedeg2}, because an Ulrich bundle has no trivial direct summand except that on $\p^n$, see \thref{2-3} and \thref{oulrich}.
     
     We can draw again the diagram \eqref{pushideale}, where we recall that $Z$ is defined to be the schematic image of $Y$, hence $\id_Z=ker(\osn\ra f_*\os_Y)$.
     \begin{equation}\label{pushideale'}
         \begin{tikzcd}
          0 \ar[r] & \id_Z \ar[r] \ar[d] & \osn \ar[r] \ar[d] & \os_Z \ar[r] \ar[d] & 0 \\
          0 \ar[r] & f_*\id_Y \ar[r] \ar[d] & \osn\oplus\osn(-m) \ar[d] \ar[r]  & f_*\os_Y\ar[r] & 0 \\
         & \cc \ar[r] & \osn(-m) & &
         \end{tikzcd}
     \end{equation}
     Everything commutes and we define $\mathcal{C}$ to be the cokernel of the induced map $\id_Z\ra f_*\id_Y$.
     From diagram chasing it follows that $\mathcal{C}\subset \osn(-m)$.
     $f_*\e\cong\osn^4$, see \thref{hilbertpol}, is globally generated but this property propagates to quotients hence, being $f_*$ exact, also $f_*\id_Y(m)$ is globally generated and therefore also $\mathcal{C}(m)$ is.
     This implies that we have a non-zero map $\osn\ra \mathcal{C}(m)$.
     Composing with the inclusion $\mathcal{C}(m)\subset \osn$ we get a non-zero endomorphism of $\osn$ which must therefore be an isomorphism.
     This implies that $\mathcal{C}\cong \osn(-m)$ and, looking at \eqref{pushideale'}, also that $\os_Z\cong f_*\os_Y$.
     Since $f|_Y$ is finite, being $f$ such, we conclude that $f|_Y$ is an isomorphism on the image.
     
     We have $h^0(X,\e)=4$ hence, by \eqref{ext'}, $h^0(X,\id_Y(m))=3$.
     Recall that being $f$ finite, we have $h^0(\p^n,f_*\id_Y(m))=h^0(X,\id_Y(m))=3$.
     From the left vertical sequence in \eqref{pushideale'} we infer that
     \[3=h^0(\p^n,f_*\id_Y(m))\ge h^0(\p^n,\id_Z(m))\geq h^0(\p^n,f_*\id_Y(m))-h^0(\p^n,\osn)=2.\]
     We claim that $ h^0(\p^n,\id_Z(m))=2$.
     Indeed, $f_*\id_Y(m)$ is globally generated, therefore, $h^0(\p^n,\id_Z(m))=3$ would imply $\id_Z(m)\ra f_*\id_Y(m)$ to be surjective and hence an isomorphism, which is not.
     
     Applying $Hom(\os_X,-)$ we deduce that the extension class of the sequence
     \ses{\id_Z(m)}{f_*\id_Y(m)}{\osn}
     vanishes, hence it is split.
     Now, applying $f_*$ to \eqref{ext'} we obtain 
     \ses{\osn\oplus\osn(-m)}{\osn^4}{\id_Z(m)\oplus\osn}
     but, since $X$ is integral any non-zero endomorphism of $\osn$ is an automorphism, hence, after simplification, we get 
     \ses{\osn(-m)}{\osn^2}{\id_Z(m)}
     i.e. $Z$ is a complete intersection of type $(m,m)$.
\end{proof}

\Co\label{unique}
For a smooth double covering $f:X\ra\p^n$ any rank $2$, special $f$-Ulrich bundle must be as the ones constructed in \thref{section-rk2}.
\io

\subsection{Non existence results}

We apply the previous criterion to prove that not only having branch divisor of equation $b=p_0^2+p_1p_2+p_3p_4$ is sufficient to have Ulrich sheaves, as found in \thref{ulrichsommeprodotti}, but it is also necessary, at least in the smooth case.
This will lead us to some non-existence results.

\Le\label{branchrango2}
Let $char(\gk)\neq 2$ and $\gk=\overline{\gk}$.
Let $f:X\ra \p^n$ be a double covering with $n\ge 3$ and branch locus defined by $b=0$ for some polynomial $b$.
Suppose $p_i\in H^0(\p^n,\osn(m))$ for $i=0,\dots 4$.
\begin{enumerate}\label{1)}
    \item We can write $b=p_1p_2+p_3p_4$ with $p_1,p_3\not\equiv 0$ and without common factor if and only if there is some $Y\subset R$ such that $Z:=f(Y)$ is a complete intersection of type $(m,m)$.
    
    \item If $X$ is smooth and there is a special Ulrich bundle $\e$ of rank $2$, then we have $b=p_0^2+p_1p_2+p_3p_4$.
\end{enumerate}
\ma
\begin{proof}
\gel{1}
Suppose $Y\subset R$ such that $Z=f(Y)\subset B\subset\p^n$ is a complete intersection of two divisors of degree $m$ in $\p^n$.
Call $p_1,p_3$ two polynomials defining these divisors, in particular they are not identically $0$ and have no component in common.
Since $Z\subset B$, the complete intersection assumption implies that $b$ is in the ideal generated by $p_1,p_3$, therefore, has the above form.

On the other direction, we can define $Z$ to be the complete intersection $p_1=0=p_3$.
Since $Z\subset B$, its set-theoretic (but not scheme-theoretic!) preimage is $Y\subset R$ and $f|_Y$ is an isomorphism on the image since already $f|_R$ is.

\gel{2}
  We use \thref{rk2-section} to construct from an Ulrich bundle of rank $2$ our $Y$ such that $f|_Y$ is an isomorphism on some $Z$ as above.
  Call $p_1$ and $p_3$ two polynomials defining $Z$, so that $b=p_1p_2+p_3p_4+q'$ and $B|_Z$ is given by $q'|_Z$.
  We can suppose $Y\not\subset R$, otherwise we just apply the previous statement, then by \thref{double} we know that $B|_Z$ is the square of some divisor in $|\os_Z(m)|$.
Since $H^0(\p^n,\osn(m))\ra H^0(Z,\osn(m))$ is surjective, see \thref{numericsci}, it follows that there is another polynomial $q\in H^0(\p^n,\osn(m))$ such that $q'|_Z=q^2|_Z$.
Therefore, $q'-q^2=p_1q_1+p_3q_3$ and we can write $b=p_1(p_2+q_1)+p_3(p_4+q_3)+q^2$ as desired.
\end{proof}

We have seen in \thref{constructible} that the locus in $|\os_P(d)|$ of divisorial coverings having Ulrich sheaves of some fixed rank $r$ is always constructible, but in general is not clear if it is irreducible when $r>1$.

\Co
Let $char(\gk)\neq 2$ and $\gk=\overline{\gk}$.
The locus of smooth, double coverings having an Ulrich bundle of rank $2$ is irreducible, if not empty.
\io
\begin{proof}
Double coverings are parametrised by the equation $b$ of their branch locus, that is by $H^0(\p^n,\osn(2m))$. 
    We can apply the second part of \thref{branchrango2}, hence any smooth $X$ bearing an Ulrich bundle of rank $2$ has branch locus of the form $b=p_0^2+p_1p_2+p_3p_4=0$.
    Consider the image of the morphism 
    \[V_{m}^5\ra V_{2m} \qquad (p_0,\dots, p_4)\mt  p_0^2+p_1p_2+p_3p_4,\]
    which is an irreducible subspace.
    The locus we are interested in is irreducible, coinciding with the open subset of this image that corresponds to the equations of smooth divisors.
\end{proof}

Another application of the previous criterion is non-existence of Ulrich bundles.
\Prop\label{nonexistence}
Suppose $char(\mathbf{k})\neq 2$, $\gk= \overline{\gk}$ and either $n\geq 5$ or $n=4$ and $m\geq 2$ or $n=3$ and $m\geq 5$.
For each $n,m$ satisfying the above conditions, there are no \text{special} rank $2$ Ulrich bundles on the general double solid. 
\One
\begin{proof}
    By \thref{branchrango2}, it is enough to show that the general polynomial $b$ defining the branch locus, with the above restrictions on $n$ and $m$, cannot be written in the form $b=p_0^2+p_1p_2+p_3p_4$.
    If $n\geq 5$ then clearly a $b$ written like this is non-generic: it is singular at least on the non-empty intersection of the hypersurfaces defined by $p_i$.
    
    In general, the possible $b$ depend on $h^0(\p^n,\osn(2m))=\binom{2m+n}{n}$ parameters while the polynomials on the right side depend on at most $5\cdot h^0(\p^n,\osn(m))=5\binom{m+n}{n}$ parameters.
    We want to prove that $\binom{2m+n}{n}>5\binom{m+n}{n}$ which is equivalent to
    \[\prod_{i=1}^n(2m+i)-5\prod_{i=1}^n(m+i)>0.\]
    Putting $n=3,4$ we get two polynomials in $m$ and we can easily verify that they are positive in the above range, except for $n=4$ and $m=2$.

    This last case deserves a better analysis.
    Set $V_l$ to be the vector space $H^0(\p^4,\os_{\p^4}(l))$ and consider
    \[\phi:V_2^5\ra V_{4} \qquad (p_0,\dots, p_4)\mt p_0^2+p_1p_2+p_3p_4.\]
    We will prove that for any point $(p_0,\dots ,p_4)\in V_{2}^5$ the corank of $d_{(p_0,\dots ,p_4)}\phi$ is at least $5$.
    It will follow that the variety parametrizing the desired polynomials has codimension at least $5$, thus proving our claim.
    
    We can easily compute as above that $dim(V_2^5)=5dim(V_2)=75$ and $dim(V_4)=70$.
    Moreover, in \thref{diff'} we showed that
    \begin{equation*}
    d_{(p_0,\dots, p_4)}\phi: (q_0,q_2,q_1,q_4,q_3)\mt 2p_0q_0+\sum_{i=1}^4p_iq_i. 
\end{equation*}
For any choice of two among the $p_i$-s, we present a vector in the kernel of the above differential: if we choose $p_1,p_2$ just take $(0,-p_1,p_2,0,0)$.
    Those choices are $\binom{5}{2}=10$ and give independent vectors for a general choice of $p_i$-s, hence the kernel of $d_{(p_0,\dots ,p_4)}\phi$ has dimension at least $10$ and its corank is at least $5$, as desired.
    Furthermore, taking a general $(p_0,\dots ,p_4)$ we can easily see that this corank is actually $5$: for example, if we use $x,y,z,v,w$ as variables is enough to take $(p_0,\dots ,p_4)=(x^2,y^2,z^2,v^2,w^2)$.
    Indeed, the only monomials of degree $4$ which are not in the ideal generated by $x^2,y^2,z^2,v^2,w^2$ are $xyzv,xyzw,xyvw,xzvw,yzvw$.
    \[\]
\end{proof}

\subsection{The involution}

We will focus on understanding when $\iota^*\e\cong \e$, at least under some assumptions.

\Prop\label{fisso}
Suppose $f:X\ra \p^n$ is a smooth, double covering with $n\geq 2$ and $\gk=\overline{\gk}$ with $char(\gk)\neq 2$.
Let $\e$ be a rank $2$ $f$-Ulrich bundle on $X$ and $\iota$ the covering involution.
Consider the following statements:
\begin{enumerate}[i)]
    \item we have an exact sequence
    \sesl{\os_X}{\e}{\id_Y(D)}{ext''}
    with $Y=\iota^* Y$ and $\os_X(D)\cong det(\e)$ (in particular $Y$ has codimension $2$)
    \item $\e\cong \iota^*\e$
\end{enumerate}
We have $i)\Rightarrow ii)$, while, if we also have $Pic(X)\cong \Z H$, then $ii)\Rightarrow i)$.
\One
\begin{proof}
\gel{$\mathbf{i)\Rightarrow ii}$}  
$Y=\iota^* Y$ implies $\iota^*\id_Y\cong \id_{\iota^*Y}$ hence applying $\iota^*$ to $\eqref{ext''}$  we get 
\sesl{\os_X}{\iota^*\e}{\id_Y(D)}{ext'''}
Applying $Hom(-,\os_X)$ to \eqref{ext''} and recalling \thref{grado0} we get
\[0=Hom(\e,\os_X)\ra Hom(\os_X,\os_X)\ra Ext^1(\id_Y(D),\os_X)\ra Ext^1(\e,\os_X)\ra\dots\] 
but by Serre duality $Ext^1(\e,\os_X)\cong H^{n-1}(\e(K_X))=0$, being $\omega_X\cong \os_X(m-n-1)$ by \thref{divisorialcover}.
We conclude that $ext^1(\id_Y(D),\os_X)=1$ which implies the desired isomorphism, since $\e$ and $\iota^*\e$ are both locally free and hence \eqref{ext''} and \eqref{ext'''} are two non-split extensions.

\gel{$\mathbf{ii)\Rightarrow i}$}  
   The automorphism $\iota$ acts linearly on $H^0(X,\e)$.
   Recall that, being $\iota$ an involution and $char(\gk)\neq 2$, the eigenvalues can only be $\pm 1$ and it is diagonalisable. 
   Therefore, at least one of the eigenspaces is non-empty.
   Let us choose an eigensection $\os_X\ra\e$.
   Being $Pic(X)\cong \Z H$, by \thref{tuttedeg2} the above morphism has to drop rank on a codimension $2$ locus, hence we have the sequence \eqref{ext''}.
Moreover, by construction, $\iota^*$ preserves the left morphism in this sequence, hence also its cokernel. We just proved that $\id_Y(D)\cong \iota^*\id_Y(D)\cong \id_{\iota^*Y}(D)$ which implies $\iota^*Y=Y$.
\end{proof}

\Co\label{fisso'}
With the above notation, suppose $char(\gk)\neq 2$, $\gk=\overline{\gk}$ and $Pic(X)$ is generated by $\os_X(1)$.
Then, there exists an $f$-Ulrich bundle $\e$ of rank $2$ such that $\e\cong \iota^*\e$ if and only if there is some $Y$ contained in the ramification divisor of $f$ such that \eqref{ext''} holds.
\io
\begin{proof}
    We can apply both conclusions in \thref{fisso} hence, $\e\cong \iota^*\e$ if and only if there exists $Y\subset X$ of codimension $2$ such that $Y=\iota^* Y$ is the zero locus of a section of $\e$, as in \eqref{ext''}.
    But from \thref{rk2-section}, or even from \thref{restrizionesezioni}, $f|_Y$ is injective hence $Y$ must be pointwise fixed by $\iota$.
    We conclude that $Y\subset R$ by \thref{cycliciinvolution}.
    Conversely, if $Y\subset R$ then clearly $Y=\iota^* Y$.
\end{proof}

\subsection{Positivity properties}

Suppose $\gk=\C$.
Here we want to understand the positivity properties of a special $f$-Ulrich bundle $\e$ of rank $2$ on a smooth double covering $f:X\ra\p^n$ with $n\ge 3$.
We already know that $\e$ is globally generated by \thref{corollario}.
Here we will show that the codimension $2$ subschemes in $X$ on which $\e$ restrict to a non-ample bundle are exactly the zero loci of sections of $\iota^*\e$, where $\iota$ is the involution on $X$ associated to $f$.

Consider the projective bundle $q:\p(\e)\ra X$. Define $\xi:=c_1(\os_{\p(\e)}(1))$ and $L:=c_1(q^*\os_X(1))$.
Recall that, being the base field $\C$ and $n>2$, we have $Pic(X)\cong \Z H$ by \thref{Piccover}.
Therefore, $Pic(\p(\e))=<\xi,L>$ by \cite{GW2}[Prop. 24.69], in particular, linear and numerical equivalence coincide on those varieties.
Recall that a divisor $D$ is \textbf{Nef} if $D\cdot C=deg(D|_C)\ge 0$ for all curves $C\subset X$; they form a convex cone in $Pic(X)\otimes\Q$.
Also effective divisors on $X$ form a cone and its closure is called the \textbf{Pseudo-effective} cone.
Divisors which are in the interior of this cone are called \textbf{big}; this is not the usual definition but they are equivalent by \cite{Laz1}[Thm. 2.2.26]. 

\Le\label{coni}
In the above setting, for any special rank $2$ Ulrich bundle $\e$ on $X$ we have $\xi^4=0$.
Both Nef and pseudo-effective cone of divisors of $\p(\e)$ are spanned by $L,\xi$, in particular $\e$ is not even big\footnote{Recall that a vector bundle $\e$ is called big if $\os_{\p(\e)}(1)$ is a big line bundle}.
\ma
\begin{proof}
Since $\os_X(1),\e$ are globally generated also $L,\xi$ are.
It follows that they are both Nef.
By \cite{Laz1}[Thm. 2.2.16], they are not big if we verify that their top self-intersection is $0$.
Clearly 
\[L^{n+1}\sim (q^*H)^{n+1}\sim q^*(H^{n+1})\sim 0.\]
Since $\xi$ is globally generated and $h^0(\p(\e),\os_{\p(\e)}(1))=h^0(X,\e)=4$ by \thref{hilbertpol}, its linear system induces a morphism $\phi:\p(\e)\ra\p^3$ such that $\phi^*\on{3}(1)\cong \os_{\p(\e)}(1)$.
It follows that $\xi^4\sim 0$.
This implies that both $L,\xi$ are Nef not big hence are extremal in both the Nef and the pseudo-effective cone of $X$.
\end{proof}

\Prop\label{nonampio}
For $n\ge 3$, let $f:X\ra\p^n$ be a smooth, complex double covering, $\e$ a (special) $f$-Ulrich bundle of rank $2$ and call $q:\p(\e)\ra X$ the standard projection.
Then:
\begin{enumerate}
    \item the morphism $\phi:\p(\e)\ra\p^3$ determined by $|\os_{\p(\e)}(1)|$ is surjective.
    \item The morphism $\phi:\p(\e)\ra\p^3$ is isomorphic to the family $\psi:\mathcal{U}\ra\p(H^0(X,\iota^*\e)\du)\cong \p^3$ of zero loci of sections of $\iota^*\e$, in particular $\phi$ is flat.
\end{enumerate}
\One
\begin{proof}
\gel{1}
    Before continuing, we will need some more notation.
    Remember that, being $Pic(X)\cong\Z H$, $\e$ has to be special hence by \thref{Piccover} we deduce $c_1(\e)\sim mH$ hence $q^*c_1(\e)\sim mL$.
Call $\gamma=q^*c_2(\e)$ and $\iota^*\gamma:=q^*\iota^*c_2(\e)$.
Let us recall that the Chow ring $A^{\bullet}(\p(\e))$ is isomorphic to $\Z[\xi,L]/(\xi^2-mL\cdot\xi+\gamma)$ by \cite{Ful}[Example 8.3.4]\footnote{our $\p(\e)$ is actually $P(\e\du)$ with Fulton's notation}.
Identifying a $0$-cycle with its degree we have 
\[H^{n-2}\cdot c_2(\e)=m^2 \qquad H^n=2,\]
the first by \thref{rk2-section}, and the second since $|H|$ induces the double covering $f$.
From those we also deduce 
\[\xi\cdot L^n=\xi\cdot q^*(H^n)=2 \qquad \xi\cdot \gamma\cdot L^{n-2}=\xi\cdot q^*(c_2(\e)\cdot H^{n-2})=m^2.\]
We have
    \begin{equation}\label{csi3}
        \xi^3\sim \xi\cdot \xi^2=\xi (mL\cdot\xi-\gamma)=mL(mL\xi-\gamma)-\xi\gamma=\xi(m^2L^2-\gamma)-mL\gamma\neq 0
    \end{equation}
    by the structure of $A^{\bullet}(\p(\e))$ or, for example, because we can compute
    \[\xi^3\cdot L^{n-2}=(\xi(m^2L^2-\gamma)-mL\gamma)\cdot L^{n-2}=m^2\xi\cdot L^n-\xi\cdot\gamma\cdot L^{n-2}-mL^{n-1}\cdot \gamma=\]
    \begin{equation}\label{contocubo}
        =2m^2-m^2-m\cdot q^*(H^{n-1}\cdot Y)=m^2.
    \end{equation}
    Then $\phi$ is surjective, since otherwise its image would have dimension at most $2$ and hence $\xi^3=\phi^*c_1(\on{3}(1))^3=0$.

\gel{2}
Choose a section of $\iota^*\e$ whose zero locus, by \thref{rk2-section}, is a codimension $2$ subscheme $\iota^*Y\subset X$.
The proof is quite technical but the idea is simple: we can show that the restriction of $\e$ on $\iota^*Y$ has a trivial quotient, implying that we have a lift $\iota^*Y\cong \p(\os_{\iota^*Y})\subset \p(\e)$ which is contained in some fiber of $\phi$.
Then we repeat the previous construction in families, in order to get a morphism $\cu\ra\p(\e)$ which we prove is actually an isomorphism.

Let $Y$ be the zero locus of a section of $\e$ so that we have
    \sesl{\os_X}{\e}{\id_Y(m)}{solita}
    Restricting 
    \ses{\id_Y(m)}{\os_X(m)}{\os_Y(m)}
    to $\iota^*Y$ we get the right exact
    \ses{\id_Y(m)|_{\iota^*Y}}{\os_{\iota^*Y}(m)}{\os_{Y\cap \iota^*Y}(m)}
    hence we have a surjection from $\id_Y(m)|_{\iota^*Y}$ to 
    \[Ker\left(\os_{\iota^*Y}(m)\twoheadrightarrow\os_{Y\cap \iota^*Y}(m)\right)=\id_{Y\cap \iota^*Y/\iota^*Y}(m)\cong \os_{\iota^*Y}(-R)(m)\cong \os_{\iota^*Y}\]
    since we have $Y\cap \iota^*Y=R\cap\iota^*Y$, by \thref{intersezionebranch}, and $\os_X(R)\cong\os_X(m)$.
    We derive a surjection \[\e\twoheadrightarrow\e|_{\iota^*Y}\twoheadrightarrow\id_Y(m)|_{\iota^*Y}\twoheadrightarrow\os_{\iota^*Y}.\]
    The kernel of the morphism $\e|_{\iota^*Y}\twoheadrightarrow \os_{\iota^*Y}$ has to be locally free of rank $1$ hence isomorphic to $det(\e|_{\iota^*Y})\cong\os_{\iota^*Y}(m)$, so we deduce the sequence
    \sesl{\os_{\iota^*Y}(m)}{\e|_{\iota^*Y}}{\os_{\iota^*Y}}{zeri}   

    Now consider the flat family $\psi:\cu\ra\p^3\cong\p(H^0(X,\iota^*\e)\du)$ whose fibers are the zero loci $\psi_t:=\iota^*Y_t$ of sections of $\iota^*\e$, where $t\in\p^3$.
    Note that $\cu$ is irreducible of dimension $n+1$ since $\iota^*Y_t$ have codimension $2$ in $X$ and the general one is irreducible being connected, since isomorphic to a complete intersection in $\p^n$, and smooth, being $\e$ globally generated.
    Call $e:\cu\ra X$ the evaluation morphism sending each curve to its realisation in $X$.
    \dia
    & \p(\e) \ar[rd, "q"]\ar[ld, "\phi"] & &  \cu \ar[rd, "\psi"]\ar[ld, "e"] \ar[ll, "g"]&  \\
    \p^3 & & X & & \p^3 \\
    \mma
    Then, $\psi_*e^*(\e(-m))\cong \cl$ is a line bundle by Grauert's theorem, see \cite{GW2}[Thm. 23.140], since $h^0(\psi_t,(e^*\e(-m))_t)=h^0(\iota^*Y_t,\e|_{\iota^*Y_t}(-m))=1$ for all $t\in\p^3$.
    By projection formula we have 
    \[\psi_*(e^*\e(-m)\otimes\psi^*\cl\du)\cong \cl\otimes\cl\du\cong\on{3}\]
    hence $h^0(\cu,e^*\e(-m)\otimes\psi^*\cl\du)=1$.
    In other words, we have a morphism $e^*\os_X(m)\otimes\psi^*\cl\du\ra e^*\e$ whose cokernel is a line bundle since restricted to all the fibers of $\psi$ coincides with $\os_{\psi_t}$ by \eqref{zeri}.
    A determinant computation implies we have an exact sequence
    \sesl{e^*\os_X(m)\otimes\psi^*\cl\du}{e^*\e}{\psi^*\cl}{magico}
    
    By universal property of $\p(\e)\ra X$, this gives us a morphism $g:\cu\ra\p(\e)$ over $X$ such that $g^*\os_{\p(\e)}(1)\cong \psi^*\cl$.
    Our next goal is to show that $g$ is birational.
    Since $g^*\phi^*\on{3}(1)\cong g^*\os_{\p(\e)}(1)\cong \psi^*\cl$, then the only curves that the morphism $\psi\circ g$ could contract are the fibers of $\psi$.
    But $g|_{\psi_t}$ is an isomorphism on the image, since $e=q\circ g$ and $e|_{\psi_t}$ is an isomorphism on the image, by definition of $e$.
    Therefore, $g$ is finite and the fibers of $\phi\circ g$ are finite unions of $\psi_{t_i}$-s.
    To show that it is generically injective it is enough to show that the general fiber of $\phi$ contains just one $g(\psi_{t})$.
    But, if $s\neq t$ then we have $e(\psi_s)\neq e(\psi_t)$ hence $g(\psi_s)\neq g(\psi_t)$.
    Suppose, by contradiction, that the general fiber of $\phi$ contains both $g(\psi_t), g(\psi_s)$ with $g(\psi_t)\neq g(\psi_s)$.
    Being $\phi$ surjective, its general fiber has dimension $n-2$ and as a cycle is rationally equivalent to $\xi^3$.
    By dimensional reasons, $g(\psi_t)$ and $g(\psi_s)$ are both irreducible components of this fiber therefore, we have $\xi^3\sim g(\psi_t)+g(\psi_s)+\eta$ for some effective cycle $\eta$.
   But then by \eqref{contocubo}
    \[m^2=\xi^3\cdot L^{n-2}=(g(\psi_t)+g(\psi_s)+\eta)\cdot L^{n-2}\ge (g(\psi_t)+g(\psi_s))\cdot L^{n-2}=\]
    so by projection formula it becomes
    \[=g_*(\psi_t+\psi_s)\cdot L^{n-2}=(\iota^*Y_t+\iota^*Y_s)\cdot H^{n-2}=2m^2,\]
    a contradiction.
    Therefore, $g$ is a generically injective morphism between two varieties hence is birational, since it induces a degree $1$ field extension on their generic points.
    But then, being $\p(\e)$ smooth, the map $g$ must be an isomorphism by Zariski's main theorem, \cite{GW1}[Cor. 12.88]. 
\end{proof}

\Oss
In the quadric case we have already seen that $\e$ is a spinor bundle and the $Y$-s are isomorphic to $\p^{n-2}$.
Actually, this second contraction $\phi$ is again a projective bundle and in \cite{SzurekWisniewski}[(3.4)], using pretty different methods, it has been shown that it is the projectivization of a \emph{null-correlation bundle}.
\one

We end by an interesting but unrelated remark.
\Oss
By the universal property of Grassmannians, see \cite{GW1}[Chapter 8 §4], any rank $2$ Ulrich bundle on a double covering $f:X\ra\p^n$ gives us a morphism $\Phi_\e:X\ra Gr(2,4)$ such that $\Phi_\e^*\mathcal{Q}\cong \e$, where $\mathcal{Q}$ is the universal quotient bundle on it.
Note that $Gr(2,4)$ is actually the smooth quadric in $\p^5$ and $Q$ is a spinor bundle, see \cite{Ottaviani_spinor}[Def. 1.3], which is Ulrich respect to the restriction of $\on{5}(1)$. 

The map $\Phi_\e$ cannot be an embedding for $n\ge 3 $.
This is clear for $n\ge 4$, while for $n=3$ follows from the fact that, by Lefschetz hyperplane theorem, see \cite{SGA2}[XII Cor. 3.7], any divisor in this quadric has the Picard group generated by the restriction of $\on{5}(1)$, which is very ample.
However, for $n=3$ we have $det(\e)\cong \os_X(m)$ very ample, so that the claims in the second part of \cite{LopezSierra_geom}[Thm. 1] are not equivalent if we only assume the polarisation to be ample and globally generated.

\one

\chapter{Some moduli spaces of Ulrich bundles}

Our goal is to understand as much as possible the components of those moduli spaces containing Ulrich bundles on double coverings of $f:X\ra\p^3$ with branch locus of degree $2m=4,6,8$.
Existence in the rank $2$ case have been shown in \thref{ulrichr2doublecoverexistence}.
Local properties of those spaces can be studied with deformation theory, often translating the issue on the zero loci of sections of those bundles, in the spirit of Hartshorne-Serre correspondence.
Those results, together with the machinery developed in \thref{wildext}, will also give us a way to completely solve the existence problem for arbitrary high rank Ulrich bundles for $m=2,3$, following \cite{CH} and \cite{CFK3}.
In particular, we conclude that those varieties are \textit{Ulrich wild}, meaning that there are arbitrary large families of non-isomorphic indecomposable Ulrich bundles.

Note that the case $m=4$ is really different, since we get a $0$-dimensional moduli space and, for $X$ general, those bundles have no non-split extensions, but still very interesting since we can prove that those are spherical objects in $\mathcal{D}^b(X)$.

More problematic is the situation regarding global properties.
Even though we are able to relate moduli spaces on $X$ with other moduli spaces on surfaces contained in it, which are better understood, we could not address the question of irreducibility of the former ones.
Nonetheless, in the case $m=2$ we prove in \thref{tyurin}, following Tyurin's theorem in \cite{Beafk3}[§2], that the space parametrising Ulrich bundles of rank $2$ sits as (an open subset of) a Lagrangian subvariety in an Hyperkhaler one.  
Furthermore, still for $m=2$, we can prove that the general Ulrich bundle of rank $r\ge2$ on a smooth divisor $\Sigma\in |\os_X(1)|$ extends to the whole $X$, see \thref{restrizdomi}.
  
In the first section, we compute dimension and analyse smoothness of the components of those moduli spaces containing Ulrich bundles of rank $2$, see \thref{double3folds},
Moreover, we study the locus of vector bundles fixed by the involution $\iota$ in \thref{finitefixedsheaves}.
The second section is entirely devoted to the case $m=2$: the \textit{quartic double solid}.
In addition to the results already cited, let us mention that, by adapting an argument in \cite{Faenzi11}[Thm. D], we obtain that \thref{ulrichr2doublecoverexistence} holds on any smooth quartic double solid.
Finally, we notice that in \thref{qdpwild} we give a complete classification of special Ulrich bundles on a smooth Del Pezzo surface of degree $2$.

\startcontents[chapters]
\printcontents[chapters]{}{1}{}

\section{Moduli spaces of rank $2$ Ulrich bundles on some double solids}

In the following we will always assume $k=\C$.
Fix $f:X\ra\p^3$ a double covering branched along a divisor $B$ of degree $2m=4,6,8$.
Recall that $f_*\os_X\cong \osn\oplus\osn(-m)$.
We will suppose our $X$ to be smooth, even though some partial result holds without that hypothesis.
In this setting, from \thref{Piccover} we have that $Pic(X)$ is generated by $\os_X(1)=f^*\osn(1)$, hence all Ulrich bundles are special.
In particular, there are no Ulrich line bundles on $X$ from \thref{pic1} so all rank $2$ ones are stable by \thref{stability}.
Moreover, Noether-Lefschetz's theorem implies that for a general branch locus $B\subset \p^3$ also $Pic(B)$ is cyclic generated by the hyperplane section.
Under our assumption, fixed $(X,\os_X(1))$ and some rank $r$, we can form the Gieseker-Maruyama moduli space of semi-stable sheaves which contains all Ulrich bundles of rank $r$, since they are always semi-stable by \thref{stability} and have a well determined Hilbert polynomial by \thref{hilbertpol}.
Since being Ulrich is open in families, stable Ulrich bundles form an open subset in the stable locus of those moduli spaces by \cite{HuyLeh}[Cor. 4.3.5].

\subsection{Dimension of moduli spaces}

We will use deformation theory to study the Hilbert scheme of the curves $Y$ corresponding to Ulrich bundles $\e$ as in \thref{section-rk2}, and from this deduce properties of the moduli spaces of those sheaves.
To do this, we need to know the cohomology of $N_{Y/X}$ and $\e\otimes\e\du$.

The key observation is that double coverings $f:X\ra\p^3$ sit inside $P_m=\p(\on{3}\oplus\on{3}(m))$ as divisors in $|\os_P(2)|$.
Recall that we denote $\pi:P_m\ra\p^3$ the projection.
Suppose that the equation of $X$ inside $P_m$ is given by $t^2-b=0$ where $b=p_0^2+p_1p_2+p_3p_4$ with $p_i\in H^0(\p^3,\osn(m))$, we have seen in  \thref{join} that this is the case for the general $X$.
Inside $P_m$ we can view $Y$ as the complete intersection of the divisors in $|\os_{P_m}(1)|$ given by $t-p_0,p_1,p_3$.
Indeed, clearly $Y\subset X$ and, since $t-p_0=0$ defines a section of $\pi$, our $Y$ is isomorphic to the complete intersection $Z$ in $\p^3$ given by $p_1=0=p_3$.
Observe that, by definition, we have $\n_{Y/P}\cong \os_{P_m}(1)|_Y^3$ 
As before, we set $\os_Y(1):=\os_X(1)|_Y=f^*\on{3}(1)|_Y$, which should not be confused with $\os_{P_m}(1)|_Y\cong \os_Y(m)$, since already $\os_{P_m}(1)|_X\cong \os_X(m)$ as seen in \thref{divisorialcover}.

With the above notation, we have a chain of closed embeddings $Y\subset X\subset P$ such that $X\in |\os_{P_m}(2)|$.
We conclude that $\n_{X/P_m}\cong \os_X(2)$
Moreover, being $Y$ isomorphic to a complete intersection in $\p^3$ and $X$ smooth, by \cite{GW2}[Cor. 19.44] $Y$ is locally complete intersection in $X$, hence also $\n_{Y/X}\du$ is locally free.
By \cite[\href{https://stacks.math.columbia.edu/tag/06CC}{Tag 06CC}]{Stacks}, we have the exact sequence 
\sesl{(\n_{X/P_m}|_Y)\du}{\n_{Y/P_m}\du}{\n_{Y/X}\du}{conormali}
Since all the sheaves involved are locally free, the dual of the above sequence stays exact and its cohomology gives a long exact sequence in which appears
\[\dots\ra H^0(Y,\n_{Y/P_m})\xra{\rho}H^0(Y,\n_{X/P_m}|_Y)\xra{\delta} H^1(Y,\n_{Y/X})\ra H^1(Y,\n_{Y/P_m})\ra\dots\]
Therefore, to understand $h^1(Y,\n_{Y/X})$ we analyse $\delta$.
Let us choose coordinates $x_0,x_1,x_2,x_3$ on $\p^3$.
\Le\label{contonormali}
With the above notation, the dual of \eqref{conormali} reads
\sesl{\n_{Y/X}}{\os_Y(m)^3}{\os_Y(2m)}{normali}
where the right map is given by taking scalar product with $(t+p_0,p_2,p_4)$ and then
\[Im(\delta)=Im\Big(H^0(Y,\n_{X/P_m}|_Y)\ra H^1(Y,\n_{Y/X})\Big)\cong \Big(\C[x_0,x_1,x_2,x_3]/(p_0,p_1,p_2,p_3,p_4)\Big)_{2m},\]
where the subscript stands for the homogeneous part of that degree.
\ma
\begin{proof}
Let us first understand the map $H^0(Y,\n_{Y/P_m})\ra H^0(Y,\n_{X/P_m}|_Y)$.
We can define a morphism $\os_{P_m}(-2)\ra \os_{P_m}(-1)^3$ by the vector $(t+p_0,p_2,p_4)$ and a morphism $\os_{P_m}(-1)^3\ra \os_{P_m}$ by scalar product with $(t-p_0,p_1,p_3)$.
Recall that the polynomial defining $X$ can be written as $p=(t-p_0)(t+p_0)-p_1p_2-p_3p_4$ and the multiplication by $p$ gives $\os_{P_m}(-2)\cong \id_{X/P_m}\subset\os_{P_m}$.
Those maps fits in the left diagram 
\dia
\os_{P_m}(-2) \ar[rr, twoheadrightarrow, "\cdot p"] \ar[d,"\ensuremath{(t+p_0,p_2,p_4)}"'] & &  \id_{X/P_m} \ar[d,hookrightarrow] & &  \os_Y(2m)   & 
& \n_{X/P_m}|_Y \ar[ll,"\cong"]  \\
\os_{P_m}(-1)^3\ar[rr, twoheadrightarrow, "\ensuremath{\cdot (t-p_0,p_1,p_3)}"] & & \id_{Y/P_m} & &  \os_Y(m)^3\ar[u,"\ensuremath{\cdot(t+p_0,p_2,p_4)}"]  & & \n_{Y/P_m} \ar[u,twoheadrightarrow]  \ar[ll,"\cong"]  \\
\mma
which becomes the right one we apply $Hom(-,\os_Y)$, in particular the vertical maps are scalar product with $(t-p_0,p_1,p_3)$.
Now recall that $Y$ is mapped isomorphically by $f$, or equivalently by $\pi:P_m\ra\p^3$, to a curve $Z\subset \p^3$ which is a complete intersection of the two surfaces $p_1,p_3$.
We obtain 
\dia
H^0(\n_{X/P_m}|_Y) & H^0(\os_Y(2m)) \ar[l,"\cong"] &  H^0(\on{3}(2m)) \ar[l,twoheadrightarrow] & H^0(\id_Z(2m)) \ar[l,hookrightarrow] \\
H^0(\n_{Y/P_m}) \ar[u,"\rho"]  & H^0(\os_Z(m)^3) \ar[l,"\cong"] \ar[u] & H^0(\on{3}(m)^3) \ar[l,twoheadrightarrow] \ar[u,"\rho'"] \\
\mma
where we used that the restriction maps from line bundles on $\p^3$ to $Z$ are surjective on global sections, see \thref{numericsci}.
By our construction, the ideal of $Z$ is generated by the two polynomials $p_1,p_3$; the image of the rightmost horizontal map is exactly $(p_1,p_3)_{2m}$.
So, it follows that 
\[H^0(Y,\n_{X/P}|_Y)\cong H^0(\p^3,\osn(2m))/H^0(\p^3,\id_Z(2m))\cong \left(\C[x_0,x_1,x_2,x_3]/(p_1,p_3)\right)_{2m}\]
Since, using $\pi$, $\p^3$ is identified with the locus $t=0$ inside $P$, the map $\rho'$ is identified with $\rho$, the scalar product with $(p_0,p_2,p_4)$.
It follows that the image of $\rho'$ is $(p_0,p_2,p_4)_{2m}$.
Since $ker(\delta)=Im(\rho)$ and by commutativity $Im(\rho')$ projects into $Im(\rho)$, we conclude \[Im(\delta)=H^0(Y,\n_{X/P_m}|_Y)/ker(\delta)\cong \left(\C[x_0,x_1,x_2,x_3]/(p_1,p_3)\right)_{2m}/Im(\rho')\cong\]
\[\cong \left(\C[x_0,x_1,x_2,x_3]/(p_1,p_3)\right)_{2m}/(p_0,p_2,p_4)_{2m}\cong \Big(\C[x_0,x_1,x_2,x_3]/(p_0,p_1,p_2,p_3,p_4)\Big)_{2m}.\]
\end{proof}

\Le\label{HSrk2}
Let $(X,H)$ be a smooth, $n$-dimensional polarised variety and $\e$ a stable, special rank $2$ Ulrich bundle on it fitting in 
\sesl{\os_X}{\e}{\id_Y(l)}{solita'''}
then 
\begin{itemize}
    \item $hom(\e,\e)=1$
    \item $ext^1(\e,\e)=h^{0}(X,\n_Y)-h^0(X,\e)+1$
    \item $ext^j(\e,\e)=h^{j-1}(X,\n_Y)$ for all $1<j<n-1$.
\end{itemize} 
Moreover, if $(X,H)$ is a Fano pair, that is $K_X\sim kH$ with $k<0$, then $ext^j(\e,\e)=h^{j-1}(X,\n_Y)$ also holds true also for $j=n-1,n$.
\ma
\begin{proof}
Recall that $\e|_Y\cong \n_Y$ by \thref{r01} then, tensoring by $\e$ the standard sequence associated to $Y$ we get
   \ses{\e\otimes\id_Y}{\e}{\n_{Y/X}}
   Recalling that Ulrich bundles are aCM, \thref{defUlrich}, we can compute 
   \[ext^1(\e,\e)=h^{0}(X,\n_Y)-h^0(X,\e)+h^0(X,\e\otimes\id_Y)\]
   and $h^{j}(X,\e\otimes\id_Y)=h^{j-1}(X,\n_{Y/X})$ for $i>1$.
   
Being $\e$ stable it is also simple and then $hom(\e,\e)=1$.
  Tensoring \eqref{solita'''} by $\e\du\cong \e(-l)$ becomes
    \sesl{\e(-l)}{\e\otimes\e\du}{\e\otimes\id_Y}{utd}
    then simplicity of $\e$ and the Ulrich condition imply that $h^0(X,\e\otimes\id_Y)=1$ and 
    \[ext^j(\e,\e)=h^j(X,\e\otimes\e\du)= h^i(X,\e\otimes\id_Y)\]
    for $0<j<n-1$.
    Being $(X,H)$ a Fano pair and $\e$ special and of rank $2$, by \thref{c12} we have 
    \[lH\sim c_1(\e) \sim (n+1)H+K_X\sim (n+1+k)H\]
    hence $1\le l\le n$ and then also $h^n(X,\e(-l))=0$ and the above claim also holds for $j=n-1,n$.
\end{proof}

\Prop\label{dimensionideformazione}
Suppose $f:X\ra\p^3$ is a smooth double covering with branch divisor of degree $2m$ and equation $b:=p_0^2+p_1p_2+p_3p_4=0$.
Moreover, call $Y\subset X\subset P_m$ the complete intersection of the divisors in $|\os_{P_m}(1)|$ given by $t-p_0,p_1,p_3$.
If the ideal generated by the $p_i$-s contains all the degree $2m$ polynomials then we get
\[\begin{tabular}{|c | c c c c c c|}
\hline
   m  & $h^0(\n_Y)$ & $h^1(\n_Y)$ & $hom(\e,\e)$ & $ext^1(\e,\e)$ & $ext^2(\e,\e)$ & $ext^3(\e,\e)$ \\
     \hline
    2 & 8 & 0 & 1 & 5 & 0 & 0 \\ 
    \hline
    3 & 9 & 0 & 1 & 6 & 0 & 0 \\
    \hline
    4 & 3 & 3 & 1 & 0 & 0 & 1 \\
    \hline
\end{tabular}\]
    where $\e$ is the rank $2$ Ulrich bundle corresponding to $Y$, as in \thref{section-rk2}.
\One
\begin{proof}
    Since $Y$ is isomorphic to $Z$ then by \thref{numericsci} we have $\omega_Y\cong \os_Y(2m-4)$.
    It follows that $h^1(Y,\n_{Y/P_m})=h^1(Y,\os_Y(m)^3)$ equals $0$ if $m=2,3$ and equals $3$ if $m=4$.
    In any case $h^1(Y,\n_{X/P_m}|_Y)=h^1(Y,\os_Y(2m))=0$.
    Gathering all the information, the cohomology sequence obtained from 
    \eqref{normali} gives 
    \[0\ra H^0(\n_{Y/X})\ra H^0(\os_Y(m))^3 \ra H^0(\os_Y(2m)) \xra{\delta} H^1(\n_{Y/X})\ra H^1(\os_Y(m))^3 \ra 0 \]
    From \thref{contonormali} we know that $Im(\delta)=0$ hence it follows immediately that $h^1(\n_{Y/X})=3 h^1(\os_Y(m))$ and by \thref{numericsci}
\[h^0(\n_{Y/X})=3h^0(\os_Y(m)) -h^0(\os_Y(2m))=3h^0(\os_Z(m)) -h^0(\os_Z(2m))=\]
\[=5\binom{m+3}{3}-\binom{2m+3}{3}-7\]
    hence we deduce the first two columns of the table.
   Since $X$ is smooth, $\e$ is locally free.
   Moreover, $Pic(X)\cong \Z H$ by \thref{Piccover}, therefore $\e$ is stable since a destabilising line bundle would be Ulrich by \thref{stability} and they cannot exist by \thref{pic1}.
   In the Fano cases $m=2,3$ we conclude by \thref{HSrk2}.
   If $m=4$ we have $\omega_X\cong\os_X$ so that by Serre duality 
   \[ext^i(\e,\e)=ext^i(\e\otimes\e\du,\os_X)=h^{3-i}(X,\e\otimes\e\du)=ext^{3-i}(\e,\e),\]
   hence again we conclude by \thref{HSrk2}.
\end{proof}

\Oss
Call $\Gamma:=Y\cup \iota(Y)$ and suppose $Y$ to be smooth and $Y\neq \iota(Y)$.
By \cite{HarHir}[Cor. 3.2] we have the exact sequence 
\ses{\n_{Y/X}}{(\n_{\Gamma/X})\vert_Y}{T^1_{\Gamma}}
where $T^1$ is Schlessinger's functor so, being $\Gamma$ a nodal curve having as only nodes $Y\cap \iota(Y)$, we have $T^1_{\Gamma}\cong \os_{Y\cap \iota(Y)}$.
Moreover, recall that $\Gamma=f^*(f(Y))$ is a complete intersection given by a section of $\os_X(m)^2$, hence $\n_{\Gamma/X}\cong \os_{\Gamma}(m)^2$.
In particular, $\n_{\Gamma/X}|_Y\cong \os_Y(m)^2$.
The image $H^0(X,\n_{Y/X})\ra H^0(X,\os_Y(m)^2)$ consists of deformations of $\Gamma$ that extend deformations of $Y$, hence in which $\Gamma$ stays reducible.
In case $m=2,3$, the fact that $H^1(Y,\n_{Y/X})=0$ implies that $H^0(X,\os_Y(m)^2)\ra H^0(X,T^1_{\Gamma})$ is surjective and hence the nodes of $\Gamma$ can be smoothed independently, see \cite{Sernesi84}[Prop. 1.6].
For $m=4$ this is no longer true: indeed we have $h^1(Y,\n_{Y/X})=3$ while $h^1(Y,\os_Y(m)^2)= 2h^1(Y,\omega_Y)=2$ therefore the morphism $H^0(X,\os_Y(m)^2)\ra H^0(X,T^1_{\Gamma})$ has $1$-dimensional cokernel.
\one

\Oss\label{maxh1}
Note that, choosing a smooth surface $S\in |\id_Y(m)|$ we have the exact sequence
\sesl{\n_{Y/S}}{\n_{Y/X}}{\n_{S/X}|_Y}{norm}
which actually reads
\ses{\os_Y}{\n_{Y/X}}{\os_Y(m)}
since we have seen in \thref{r01} that $\n_{Y/X}\cong\e|_Y$, in particular $det(\n_{Y/X})\cong\os_Y(m)$.
It follows that 
\[h^1(Y,\n_{Y/X})\leq h^1(Y,\os_Y)+h^1(Y,\os_Y(m))=\begin{cases}
    p_a(Y) \qquad m=2,3
\\
p_a(Y)+1\qquad m=4. \end{cases}\]
This bound is attained if $Y\subset R$.
In fact, $T_X|_R\cong T_R\oplus\n_{R/X}\cong T_R\oplus\os_R(R)$ so that the inclusion $T_Y\ra T_X|_Y$ factors through $T_R|_Y$.
The cokernel of the former map, which by definition is the normal bundle of $Y$ in $X$, therefore equals $\n_{Y/R}\oplus\os_Y(R)$.
By the sequence \eqref{norm} with $S=R$ we get $\n_{Y/R}\cong\os_Y$, so the desired claim.
\one

\Te\label{double3folds}
Let $f:X\ra\p^3$ be a general smooth complex double covering branched along a surface of degree $2m=4,6,8$.
There are generically smooth components of the moduli space of sheaves on $X$ whose general point $[\e]$ represents a stable Ulrich bundle of rank $2$ with $ext^2(\e,\e)=0$, and which have dimensions $5,6,0$ for, respectively, $m=2,3,4$.
\Ma
\begin{proof}
By \thref{ulrichr2doublecoverexistence} we know that a general $X$ admits rank $2$ Ulrich bundles. By \thref{join} the equation of the branch locus satisfies the hypothesis in \thref{dimensionideformazione}.
It follows that, for a general $X$ there is some Ulrich bundle $\e$ with $ext^2(\e,\e)=0$ and, by standard deformation theory, this must be a smooth point of the corresponding moduli space, which has dimension equal to $ext^1(\e,\e)$ around $[\e]$.
All Ulrich bundles of rank $2$ on $X$ are stable therefore, in each component containing one Ulrich bundle those sheaves actually form an open subset.
\end{proof}

We end with an important remark.

\Oss
Note that for $X$ general and $m=4$ the above rank $2$ Ulrich bundles are \textbf{spherical}, that is $ext^i(\e,\e)=0$ for $i\neq 0,3$ and $1$ otherwise.
Such sheaves are of fundamental importance to understand the bounded derived categories of coherent sheaves on $X$, since they give auto-equivalences of $\mathscr{D}^b(X)$ see \cite{SeidelThomas}[Thm. 1.2].
From our study we can say that on a general $X$ there are at least $2$ different such bundles: $\e$ and $\iota^*\e$.
Indeed, $\iota^*\e\neq\e$ by \thref{fisso'} since for general $X$ there are no $Y$-s as above in the ramification locus being $Pic(R)$ cyclic.
See \thref{noic} and \thref{curvefisse} for the precise codimension of the locus of branch divisors for which this becomes an isomorphism.
\one

We have seen that in the case $m=4$ the above moduli spaces have to be reducible (for $X$ general) but we are not able to compute the number of irreducible components.
For $m=2,3$ irreducibility is still an open question.

\Do
Is the moduli space of sheaves with the same Hilbert polynomial as rank $2$ Ulrich bundles irreducible?
If not, how many components do contain Ulrich bundles?
\da

\subsection{Fixed loci of the involution on $\mathscr{M}$}

The involution $\iota$ of $X$ lifts to an involution on $\mathscr{M}$, which we will still denote by $\iota$.
\thref{fisso} tells us that its fixed points are the sheaves which have a section vanishing on some $Y\subset R$, where $R$ is the ramification divisor of $f$.
The curves $Y$, and the sheaves $\e$, for which this happens are somehow extremal; see \thref{maxh1}.
In this subsection, we will work with any integer $m\geq 2$.

First, we can compute the codimension of the locus of surfaces in $|\on{3}(2m)|$ that contain a complete intersection as $Z$, which coincides with the locus of surfaces whose equation can be written as $p_1p_2+p_3p_4=0$ by \thref{branchrango2}.
\Prop\label{noic}
For any $m\geq 2$, the locus of surfaces in $|\on{3}(2m)|$ whose equation can be written as $p_1p_2+p_3p_4=0$ with $p_i\in H^0(\p^3,\on{3}(m))$ is irreducible of dimension $4\binom{m+3}{3}-5$.
\One
\begin{proof}
Any $(m,m)$ complete intersection has the same Hilbert polynomial and they form an open subset, which we will denote by $\h_{Z}$, of the adequate Hilbert scheme.
Correspondingly, there is an open subscheme $\h_{Z,B}$ of a flag Hilbert scheme parametrizing pairs $(Z,B)$ with $Z\subset B\in |\osn(2m)|$, see for example \cite{Sernesi}[Chapter 4.5], with two projections as in
    \dia
    & \h_{Z,B}    \ar[ld,"\pi_1"] \ar[rd,"\pi_2"]  \\
    \h_{Z}  & & \vert\on{3}(2m)\vert\cong\p^{\binom{2m+3}{3}-1}\\
    \mma
     $\h_{Z}$ is isomorphic to the open subscheme of $Gr(2,H^0(\p^3,\on{3}(m)))$ parametrizing subspaces with base locus of dimension $1$, see \cite{Benoist}[Section 2.2, Prop. 2.2.5], in particular it is smooth and irreducible of dimension $2\left(\binom{m+3}{3}-2\right)$.
    The fiber of $\pi_1$ over a fixed $Z$ is a closed subscheme of $\vert\on{3}(2m)\vert$ supported on $\vert\id_Z(2m)\vert\cong \p^{2\binom{m+3}{3}-2}$, see \thref{numericsci}.
    It follows that $\h_{Z,B}$ is irreducible of dimension $4\binom{m+3}{3}-6$.
    The general point in $\h_{Z,B}$ correspond to a pair $(Z,B)$ of smooth varieties hence to compute the dimension of the general fiber of $\pi_2$ we can reduce to this case.
    Finally, a computation shows that any fibre of $\pi_2$ over a smooth surface has dimension $1$ (this follows also from the more general \thref{rango2su3fold}) hence the dimension of the image of $\pi_2$ is $4\binom{m+3}{3}-5$.   
\end{proof}

\Oss\label{curvefisse}
We can easily compute that the codimension of the image of $\pi_2$ is $1,10,31$, respectively, for $m=2,3,4$. In the first two cases, this coincides with $p_a(Z)$, which is also the maximum possible value of $h^1(\n_{Y/X})$ from \thref{maxh1}.
While, for $m=4$ we have $p_a(Z)=33$ and $h^1(\n_{Y/X})=34$.

This suggests that, at least for $m=2,3$, there could be a stratification of $\vert\on{3}(2m)\vert$ on which, if $Y$ is the zero locus of a section of an Ulrich sheaf,  $h^1(\n_{Y/X})$ assumes all the possible values between $0$ (general case) and $p_a(Z)$.
\one

Next we show that if $X$ is such that $\iota$ admits some fixed Ulrich bundles, then actually the number of those sheaves is finite and strongly linked to the geometry of the branch locus $B$ of $f$.
We will need a general lemma.

\Le\label{rango2su3fold}
Suppose $X$ is a smooth complex $3$-fold with $h^1(X,\os_X)=0$ and $Pic(X)\cong\Z H$.
Let $\e$ be a rank $2$ vector bundle without trivial summands sitting in
\sesl{\os_X}{\e}{\id_Y(D)}{exy}
with $Y$ of codimension $2$.
Let $S$ be any smooth, irreducible, ample surface such that $S\sim D$ and containing $Y$. 
Then $Y^2=0$ on $S$, $\os_S(Y)$ is globally generated by $2$ sections and any divisor in the pencil $|\os_S(Y)|$ is the zero locus of a section of $\e$.
\ma
\begin{proof}
Since $S$ is smooth, $Y$ is Cartier on it hence the sequence
    \ses{\n_{Y/S}}{\n_{Y/X}}{\os_X(S)|_Y}
    is exact.
    By \thref{r01} we know that $\n_{Y/X}\cong \e|_Y$, in particular $c_1(\n_{Y/X})\sim c_1(\e)|_Y\sim S|_Y$, so that $\n_{Y/S}\cong \os_Y$.
    It follows that $\os_S(Y)|_Y\cong \n_{Y/S}\cong\os_Y$ hence $Y^2=0$.
    Putting $h^1(X,\os_X)=0$ in
    \ses{\os_X(-S)}{\os_X}{\os_S}
    and using Kodaira vanishing we get $h^1(X,\os_S)=0$.
    Therefore by
    \ses{\os_S}{\os_S(Y)}{\os_S(Y)|_Y\cong\os_Y}
    we conclude that $\os_S(Y)$ is globally generated by $2$ sections.
    
    The surface $S$ gives a non-zero element in $h^0(X,\id_Y(D))$ then, being $h^1(X,\os_X)=0$, from \thref{surjdeg} there is a section $s'\in H^0(X,\e)$ such that $S=\{s\wedge s'=0\}$.
    From \thref{tuttedeg2} and our assumptions, we get that the zero locus of $s'$ must be some curve $Y'\subseteq X$.
    Clearly $Y'\subset S$ and, being $S$ smooth, we find $Y'\cap Y=\emptyset$, hence $Y\cdot Y'=0$.
    Moreover, $S$ contains the zero loci of all sections of $\e$ in the span $<s,s'>$ and they are all disjoint, otherwise $Y'\cap Y\ne\emptyset$.
    Call $Y''$ any such curve.
    Then, being $\os_S(Y)$ globally generated, $Y''$ must be a component of some divisor linearly equivalent to $Y$.
    Since they both represent $c_2(\e)$ by \thref{hsc2}, they have the same degree respect to $H$ hence $Y\sim Y''$.
    We conclude that the zero loci of all sections of $\e$ in the span $<s,s'>$ form the entire $|Y|$. 
\end{proof} 

\Prop\label{finitefixedsheaves}
Let $f:X\ra\p^3$ be a smooth double covering branched along a surface of degree $2m\ge4$.
Isomorphism classes of Ulrich bundles $\e$ of rank $2$ on $X$ such that $\e\cong \iota^*\e$ form a finite set with cardinality equal to half the number of pencils of complete intersection curves of type $(m,m)$ on $B$. In particular, if $X$ is general there is none.
\One
\begin{proof}
    Since $f|_R:R\ra B$ is an isomorphism, any such pencil in $B$ correspond to one in $R$.
    Those curves $Y$ give a Ulrich bundle $\e$ through Hartshorne-Serre construction, see \thref{section-rk2}, and from \thref{fisso'} we know $\e\cong \iota^*\e$ being $Y\subset R$.
    Moreover, by \thref{divisorialcover} $det(\e)\cong\os_X(m)$ is ample, so any such pencil gives a unique $\e$ by \thref{rango2su3fold}.
    Conversely, given an Ulrich bundle $\e$ fixed by $\iota$ we get an action of $\iota$ on $H^0(X,\e)\cong\C^4$.
    The condition that the zero locus of a section is contained in $R$ implies that $\iota$ acts as a scalar on it, hence the section is an eigenvector.
    
    Being $\iota$ an involution, it is diagonalisable and there are at most two eigenspaces, corresponding to the possible eigenvalues $\pm 1$.
    We will prove that both have dimension $2$.  
    It is enough to prove that both eigenspaces have dimension at least $2$, since $H^0(X,\e)\cong\C^4$.
    First of all, we show that the action of $\iota$ is not a multiple of the identity.
    By \thref{rk2-section} we have a sequence
    \ses{\os_X}{\e}{\id_Y(2)}
    in particular, $h^0(X,\id_Y(2))=3$ hence we can choose a smooth surface $S\neq R$ containing $Y$. 
    Then \thref{rango2su3fold} implies that on $S$ we have $Y^2=0$ and all curves in the pencil $|Y|$ come from sections of $\e$, in particular some of them are not in $R$ hence are not fixed by $\iota$ hence do not come from eigensections. 
    Now suppose $Y\subset R$ is the zero locus of some eigensection of $\e$, this is a curve since we can apply \thref{tuttedeg2}, being $Pic(X)\cong \Z H$ by \thref{Piccover}.
    Now \thref{rango2su3fold} gives $Y^2=0$ and all curves in the pencil $|Y|$ come from sections of $\e$.
    Being inside $R$, these curves are fixed by $\iota$, hence correspond to eigensections of $\e$ which must be in the same eigenspace of the one giving $Y$.
    Therefore, this eigenspace of $H^0(X,\e)$ has dimension at least $2$.
    
    Therefore, both eigenspaces have dimension $2$ and give pencils of curves in $R$ whose projection on $B$ gives a pencil of complete intersection curves of type $(m,m)$ by \thref{rk2-section}.
    We conclude that each Ulrich bundles furnishes two pencils hence the number of bundles is half the number of pencils as desired.
    The final claim follows from \thref{noic}.
\end{proof}

\subsection{Higher rank Ulrich bundles}
In case $m=2,3$, that is the \textit{quartic} and \textit{sextic double solid}, there are slope-stable Ulrich bundles of any even rank on a general $X$.
Since the dimension of the base scheme of those sheaves becomes arbitrary large we obtain that general sextic double solids are Ulrich wild.

\Prop\label{rango>sextic}
Let $f:X\ra\p^3$ be a smooth complex sextic (resp. quartic) double solid admitting $f$-Ulrich bundles $\e$ of rank $2$ with $ext^2(\e,\e)=0$\footnote{Recall that these exists for a general such variety by \thref{double3folds}}.
There are slope-stable $f$-Ulrich bundles of any even rank $r=2\rho$ obtained as deformations of extension as in
\sesl{\e}{\f_r}{\f_{r-2}}{extr}
Moreover, for every $\rho$ their moduli space contains a generically smooth component of dimension $5\rho^2+1$ (reps. $4\rho^2+1$) whose general element $\f'$ satisfies $ext^2(\f',\f')=0$.
\One
\begin{proof}
We deal with the sextic case $m=3$, the other one is formally identical, one just needs to change the numbers appearing.
Indeed, for the case $\rho=1$, define $\cc_2$ to be a smooth, open subset in some irreducible component of a modular family of simple sheaves on $X$ parametrizing rank $2$ Ulrich bundles such that $ext^2(\e,\e)=0$, in particular there is a flat family $\scr{F}_2$ of such sheaves on it.
We claim that the pair $(\scr{F}_2,\scr{F}_2)$ is of wild extension.
Then by \thref{extpreciso} we have the existence of families $\scr{F}_{r+2}$ of stable rank $r+2=2(\rho+1)$ bundles, over a smooth irreducible base, obtained as deformations of extensions of bundles in $\scr{F}_2$ and $\scr{F}_{r}$.
Since being Ulrich is an open property and extension of Ulrich bundles is Ulrich by \thref{2-3}, the general element in those families are again Ulrich bundles.
Moreover, for an Ulrich bundle stability implies slope-stability by \thref{stability}.

    Let us prove the above claim.
    By construction, $\scr{F}_2$ is a nice family, as defined in \thref{defwildext}, in particular $ext^2(\e,\e')=0$ for general $\e,\e'$ in $\scr{F}_2$.
    To prove that $(\scr{F}_2,\scr{F}_2)$ is of wild extension is enough, due to semicontinuity of cohomology, to show that for some stable Ulrich bundle $\e$ of rank $2$ we have $ext^1(\e,\e')=5$ for general $\e,\e'$ in $\scr{F}_2$.
    Recall that, by \thref{divisorialcover}, $X$ is Fano of index $1$, i.e. $\omega_X\cong \os_X(-1)$. 
    By Serre duality
    \[ext^3(\e,\e)=h^3(X,\e\otimes\e\du)=h^0(X,\e\otimes\e\du(-1))=0\]
    since $\e$ is simple. 
    Now, by \thref{double3folds} we have $ext^1(\e,\e)=h^1(X,\e\otimes\e\du)=6$ so that
    \[\chi(\e\otimes\e\du)=\sum_{i=0}^3h^i(X,\e\otimes\e\du)=h^0(X,\e\otimes\e\du)-h^1(X,\e\otimes\e\du)=1-6=-5.\]
    Then, since for $\e'\neq \e$ we have $hom(\e,\e')=0$, for a general $\e'$ we have
    \[ext^1(\e,\e')=h^1(X,\e\otimes\e\du)=-\chi(\e\otimes\e\du)=5\]
    by semicontinuity of cohomology.

    Finally, by \thref{extpreciso} the base of the family $\scr{F}_{r}$ of rank $2\rho$ Ulrich bundles we have constructed has dimension
    \[dim(\cc_2)+(\rho-1)\Big(dim(\cc_2)-1+\rho\cdot\mathfrak{ext}^1(\scr{E}_2,\scr{E}_2)\Big)=6+5\rho(\rho-1)+5(\rho-1)=5\rho^2+1.\]
\end{proof}

\Oss
This result is optimal for $m=3$: there cannot be Ulrich bundles of odd rank on a smooth sextic double solid by the formula for the first Chern class in \thref{c12}.
However, in the case $m=2$ there exists also odd rank Ulrich bundle as soon as it is at least $3$, we will show this in \thref{qdswild}.

Finally, for $m=4$ we have $ext^1(\e,\e)=0$ for $\e$ a general Ulrich bundle we constructed on $X$, see \thref{dimensionideformazione}.
Nevertheless, we do not know any obstruction to the existence of higher rank ones.
\one

\section{The quartic double solid case}

Suppose $\gk=\C$.
We will work out more the case $X$ smooth and $m=2$.
Note that such an $X$ is Fano of index $2$, by \thref{divisorialcover}, and degree $\left(\frac{-K_X}{2}\right)^3=2$.
Those variety are usually called \textit{quartic double solids}. 
Recall that any rank $2$ Ulrich bundle on $X$ is automatically stable and fits in a sequence
\ses{\os_X}{\e}{\id_Y(2)}
with $Y$ a connected curve of degree $4$ and arithmetic genus $1$, see \thref{rk2-section}.

\subsection{Existence for all smooth QDS}
First, the existence of Ulrich bundles can be deduced from arguments in the existing literature and without genericity assumptions.
We start with a construction.
It essentially coincides with the one given in \cite{Faenzi11}[Thm. D Step 1], see also \cite{Bea}[Lem. 6.2] for another possibility.
\Prop\label{deformazione} 
Every smooth quartic double solid contains smooth irreducible curves of genus $1$ and degree $4$ not contained in any divisor in $|H|$.
The Hilbert scheme $\mathscr{H}^{4}_1$ has a generically smooth component of dimension $8$.
\One
\begin{proof}
By adjunction, smooth divisors $\Sigma \in|H|$ are smooth Del Pezzo surfaces of degree $2$, hence they are isomorphic to the blow up of $\mathbb{P}^{2}$ in $7$ sufficiently general distinct points $P_i$, see \cite{Dolgachev}[Prop. 6.3.9].
In particular, $\operatorname{Pic}(\Sigma)$ is generated by $L$, the pullback for this blow-up of the class of a line in $\mathbb{P}^{2}$, and the exceptional divisors $E_{i}$ with $i=1, \ldots 7$. Moreover we know
\[
L^{2}=1 \quad E_{i}^{2}=-1 \quad L \cdot E_{i}=0 .\]
Note that by adjunction
\[H_{\Sigma} \sim-(K_{X}+\Sigma)|_{\Sigma} \sim-K_{\Sigma}=3 L-\sum_{i=1}^{7} E_{i}\]

Set $C\sim 3 L-\sum_{i=1}^{5} E_{i}$.
Using the above formulas, we can verify that $C\cdot H=-K_{\Sigma} C=4$ and $C^{2}=4$, so by adjunction $C$ has arithmetic genus $1$. 
To show that linear system $\left|3 L-\sum_{i=1}^{5} E_{i}\right|$ has no base points we can apply \cite{Har}[V Prop. 4.3] being the $P_i$-s general.
Moreover, $\left|3 L-\sum_{i=1}^{5} E_{i}\right|$ has dimension $2$ by \cite{Har}[V Cor. 4.4] and is not composed with a pencil since $(3 L-\sum_{i=1}^{5} E_{i})\cdot E_1=1$.
This implies, by Bertini's theorem, that the general curve in $\left|3 L-\sum_{i=1}^{5} E_{i}\right|$ will be smooth and irreducible.

Now we have curves satisfying all properties but always contained in divisors in $|H|$, the idea is to deform them and show that the generic deformation is contained in no hyperplane. 
Let us call $\mathscr{H}_{1}^{4}$ the Hilbert scheme of subscheme having Hilbert polynomial $4t$ with respect to $H$. We can consider the incidence variety
\[\mathcal{I}:=\left\{([C], [\Sigma]) \;|\; [C] \in \mathscr{H}_{1}^{4}, [\Sigma] \in|H| \text { and } C \subset \Sigma\right\}\]

and the diagram

\dia
& & \mathcal{I} \ar[dl, "p"] \ar[dr, "q"]\\
& \mathcal{H}^{4}_1 & & \vert H\vert \cong \p^{3} \\
\mma

Let $C \subset \Sigma$ be a smooth curve in a smooth surface such that $[C] \in \mathscr{H}_{1}^{4}$.
There is an exact sequence
\[\left.0 \rightarrow \mathcal{N}_{C / \Sigma} \rightarrow \mathcal{N}_{C / X} \rightarrow \mathcal{N}_{\Sigma / X}\right|_{C} \rightarrow 0\]
where $\mathcal{N}_{C / \Sigma}\cong \os_\Sigma(C)|_{C}$ and $\mathcal{N}_{\Sigma / X}|_{C}\cong \os_{C}(1)$.
We can easily compute that

\[
h^{0}\left(C, \mathcal{N}_{C / \Sigma}\right)=C^2=4=CH=h^{0}\left(C,\left.\mathcal{N}_{\Sigma / X}\right|_{C}\right) \]
\[h^{1}\left(C, \mathcal{N}_{C / \Sigma}\right)=0=h^{1}\left(C,\left.\mathcal{N}_{\Sigma / X}\right|_{C}\right)
\]

so we get $h^{1}\left(C, \mathcal{N}_{C / X}\right)=0$ and $h^{0}\left(C, \mathcal{N}_{C / X}\right)=8$ hence $\mathscr{H}_{1}^{4}$ is smooth of dimension $8$ near $C$. 
Take an open subset $U$ of the component of $\mathscr{H}_{1}^{4}$ containing $C$ such that all curves in $U$ are smooth. Now consider the map $q$ restricted to $p^{-1}(U) \subset \mathcal{I}$. 
The image of this map is contained in $|H|$ hence it has dimension at most $3$. 
A fiber of this map over $S \in H$ consists of smooth curves in $S$ of genus 1 and degree $4$. 
Each of those curves is contained in a linear system of dimension $4$ on $S$, because we already saw that $h^{0}\left(C, \mathcal{N}_{C/\Sigma}\right)=4$, in particular the fiber of $q|_{p^{-1}(U)}$ must be union of open subsets in those linear systems so must have dimension $4$. 
Therefore, each irreducible component of $p^{-1}(U)$ can have dimension at most $7$ hence no component of $p^{-1}(U)$ can dominate $U$, which has dimension $8$.
Since $p$ is not dominant, there are curves in a neighbourhood of $C$ satisfying our claim.
\end{proof}

\Te\label{qdsrk2}
Over a smooth complex quartic double solid $X$ there always exist stable Ulrich bundles $\e$ of rank $2$ and they all satisfy $ext^3(\e,\e)=0$.
The open subscheme of the moduli space parametrising Ulrich bundles on $X$ is smooth of dimension $5$.
In particular we have $ext^2(\e,\e)=0$ except at the finite number of points representing sheaves such that $\iota^*\e\cong \e$, where the tangent space has dimension $6$.
Those points do not appear on a general $X$.
\Ma
\begin{proof}
    The existence proof goes, as usual, through the Hartshorne-Serre construction.
    We can show that there is a degree $4$, genus $1$ smooth irreducible curve $Y\subset X$ which is not contained in any hyperplane via \thref{deformazione}.
    In particular $\omega_Y\cong \os_Y$ is globally generated by just one section.
    Then, we can apply Hartshorne-Serre's construction in the form \thref{Ulrichcorrispondenza'}: indeed, \eqref{top1'} is exactly the fact that $h^0(X,\id_Y(1))=0$, \eqref{surj'} can be ignored by \thref{r2} and next we show \eqref{i1'} and \eqref{mid1'}.
    For the first, we already know that $\e$ has the correct Hilbert polynomial by \thref{rk2-section}, since sits in
    \sesl{\os_X}{\e}{\id_Y(2)}{idyseq}
    and $Y$ has genus $1$ and degree $4$ as a complete intersection of type $(2,2)$ in $\p^3$.
    For the second, from the sequences 
    \sesl{\id_Y(2-i)}{\os_X(2-i)}{\os_Y(2-i)}{idyseq1}
    we get $h^1(X,\id_Y(2-i))=0$ for $i=1,2,3$: $i=3$ is obvious, $i=2$ follow from $Y$ being connected and $h^1(X,\os_X)=0$, while $i=1$ by using $h^1(X,\os_X(1))=0=h^0(X,\id_Y(1))=0$.
    Moreover, similarly to what done in \thref{dimensionideformazione}, we can easily see from \eqref{idyseq}, \eqref{idyseq1} and the Ulrich condition that
    \[ext^3(\e,\e)=ext^3(\id_Y(2),\e)=0.\]

    Any Ulrich bundle is stable so if $\iota^*\e\neq \e$ we have $ext^2(\e,\e)=hom(\e,\iota^*\e)=0$ by \thref{extquartic} and hence $[\e]$ is a smooth point.
    A stable sheaf is always simple so for any Ulrich bundle $\e$ we have $hom(\e,\e)=1$.
    By \thref{dimensionideformazione} we get $ext^3(\e,\e)=0$, $ext^1(\e,\e)=5$, hence $\chi(\e\otimes\e\du)=-4$.

    If $\iota^*\e\cong \e$, then repeating the above computation we get $ext^2(\e,\e)=hom(\e,\iota^*\e)=1$ so this time we have $ext^1(\e,\e)=6$.
    Those sheaves form a finite set due to \thref{finitefixedsheaves} but the Ulrich locus is open hence they live inside a smooth component and are actual singular points.
    Finally, by \thref{finitefixedsheaves} we know that sheaves fixed by $\iota$ only appear if the Picard rank of $B$ is greater than $1$, in particular $B$ is non-general.
\end{proof}

\Oss
For smooth complex Fano $3$-fold of index $2$ and degree $\left(\frac{-K_X}{2}\right)^3>2$ any stable Ulrich bundle of rank $2$ is always a smooth point in its moduli space, see \cite{CFK3}[Prop. 4.2].
We have seen that the general smooth complex quartic double solids also enjoys the above property, but there are some special ones, still smooth, for which it fails.
However, this phenomenon is not new when considering moduli spaces of Bridgeland stable objects: in fact the Kuznetsov component of $\mathscr{D}^b(X)$ is an Enriques category, since its Serre functor is given by $\iota^*(-)[2]$, therefore \cite{PerryPertusiZhao}[Thm. 1.8 1)] applies predicting exactly that objects fixed by $\iota$ give non-smooth points.
Similar results also hold moduli spaces of Gieseker-stable sheaves on standard Enriques surfaces, see \cite{Kim_Enriques}[§2 p.88 Thm.]. 
\one

\subsection{Restriction to $K3$-s}

We want to study the relation between the moduli space of sheaves on $X$ and the ones living on anti-canonical divisors.
This will lead to interesting subvarieties of moduli spaces on $K3$-surfaces in $|\os_X(2)|$.
Since $f_*\os_X(2)\cong\osn\oplus\osn(2)$, the general surface in this linear system is not a pullback from $\p^3$, in particular is isomorphic to a quartic in $\p^3$, see also \cite{CD}[Prop. (2.2)], and its Picard is cyclic generated by the restriction of $\on{3}(1)$ by Noether-Lefschetz's theorem.

First, we deal with the problem of when the restrictions of two Ulrich sheaves are isomorphic.
\Le\label{injrest}
If $X$ is a smooth quartic double solid and $\e$ is a rank $2$ Ulrich bundle then $\e_S\cong\e'|_S$ implies $\e'\cong \e$ if $S\in |\os_X(2)|$ is a sufficiently general surface.
\ma
\begin{proof}
Suppose we have a stable, rank $2$ Ulrich bundle $\e$ on $X$ fitting in 
    \ses{\os_X}{\e}{\id_Y(2)}
    From \thref{restrdoppi} if there is another Ulrich bundle $\e'$ such that $\e_S\cong\e'|_S$ then $\e'\cong\e$ or $\iota^*\e$.
    It remains to show that the former holds for a general $S$.
    
    Restricting the above sequence to $S$ it remains exact, since $Y\not\subset S$ by genericity, hence we deduce that $\e'|_S$ has a section vanishing on the schematic intersection $Y\cap S$.
    Being $\e'$ Ulrich, from the sequence
    \ses{\e'(-2)}{\e'}{\e'|_S}
    the restriction map $H^0(X,\e')\ra H^0(S,\e'|_S)$ is an isomorphism hence there is a section of $\e'$ whose zero locus, call it $Y'$, is such that $Y\cap S=Y'\cap S$.
    Put it in another way, $S$ and $Y'$ cut the same divisor on $Y$.
    But, being $h^1(X,\id_Y(2))=0$, the surfaces in $|\os_X(2)|$ cut the complete linear system $|\os_Y(2)|\cong\p^7$ while the curves $Y'$ zeros of sections of $\iota^*\e$ move in $\p(H^0(X,\iota^*\e)\du)\cong \p^3$.
    We conclude that for $S$ and $\e$ general we do not have $Y\cap S=Y'\cap S$ hence $\e'\not\cong \iota^*\e$.
\end{proof}

Call $\cu_2$ an irreducible open subset of the moduli space parametrising rank $2$ (stable) Ulrich bundles with $h^2(\e\otimes\e\du)=ext^2(\e,\e)=0$; note that for $X$ general all rank $2$ Ulrich bundles satisfy this property by \thref{finitefixedsheaves}.

\Prop\label{tyurin}
Let $X$ be a smooth quartic double solid and $S\in |\os_X(2)|$ be a sufficiently general surface.
Restriction to $S$ gives an morphism $\rho_S:\cu_2\ra \mathcal{M}^s_S$, where $\mathcal{M}^s_S$ is the stable locus in the adequate moduli space of semi-stable sheaves on $S$.
The map $\rho_S$ is an isomorphism onto its image, which is an open subset of a Lagrangian submanifold\footnote{Indeed, since $\mathcal{M}^s_S$ is isomorphic to an open in the moduli space of sheaves containing $\e(-1)|_S$, it is birational to an $OG10$-type manifold.}.
\One
\begin{proof}
    From a theorem of Tyurin, see \cite{Beafk3}[§1], we know that restriction to $S$ maps $\cu_2$ onto a Lagrangian subvariety of the moduli (algebraic) space of simple sheaves and is étale on its target.
    This essentially follows by considering 
    \ses{\e\otimes\e\du(-2)}{\e\otimes\e\du}{\e\otimes\e\du|_S} 
    and using $h^1(\e\otimes\e\du(-2))=h^2(\e\otimes\e\du)=0$.
    We can prove that, under our assumptions, $\rho_S(\e)$ is always stable hence, since the moduli space of stable bundles forms an open subscheme in the moduli of simple ones, the map $\rho_S$ in our statement is actually well posed.
    We prove that $\e_S$ is even slope-stable.
    Having rank $2$, it can only be destabilised by a line bundle.
 Since $det(\e)\cong\os_X(2)$, it shares the slope with $\os_X(1)$.
    But, being $S$ general, we can assume that it is smooth and its Picard group is generated by $\os_X(1)|_H$ so from 
    \ses{\e(-3)}{\e(-1)}{\e(-1)|_S}
    we see that $0=h^0(S,\e|_S(-1))=hom(\os_S(1),\e)$ hence $\e$ is slope-stable.
    
    By \thref{injrest} we deduce that, for $S$ general, the morphism $\rho_S$ has at least one fiber made up of just one point. 
    Being $\rho_S$ etale on the image we conclude that it is an actual isomorphism on the image.
\end{proof}

\Oss
Clearly, if $R$ is the ramification locus of $f$, we have $\e|_R\cong \iota^*\e|_R$ for all $\e$ hence $R$ is never "general" in the sense of the previous result.
\one

The first part of this proof works also in the case $m=3$ but we are not able to compute the degree in this case.
\Prop
Let $X$ be a smooth sextic double solid and $S\in |\os_X(1)|$ be a general surface.
Restriction to $S$ gives an morphism $\rho_S:\cu_2\ra \mathcal{M}^s_S$, where $\mathcal{M}^s_S$ is the stable locus in the adequate moduli space of semi-stable sheaves on $S$.
The map $\rho_S$ is étale onto its image, which is an open subset of a Lagrangian submanifold.
\One

\subsection{Lines and rank $2$ Ulrich bundles}

This subsection is just a collection of some well-known facts on lines on $X$ and some computations regarding them that are needed in the following.
Here by lines we mean curves $l$ in $X$ such that $l\cdot H=1$.
Those are automatically smooth and rational, hence isomorphic to $\p^1$, since they cover birationally lines in $\p^3$.
Such curves are of extreme importance to understand the geometry of varieties, as for example are used in the classification of Fano $3$-folds, see \cite{IskPro}.
Furthermore, they are fundamental in our study of Ulrich bundles on Fano $3$-folds of index $2$ as $\id_l(1)$ has the same reduced Hilbert polynomial of an Ulrich bundle, even not being a sheaf of this kind. 

The following has been proven in \cite{Isk79}[III Prop.1.3, Remark 1.7], \cite{Wel}[(1.1), (1.2)] and \cite{KuzProShr}[Lem. 2.6].
\Le\label{normalerette}
Given a line $l\subset X$ on a smooth quartic double solid we have $3$ possibilities for the normal bundle:
\[\n_l=\os_l^2,\; \os_l(-1)\oplus\os_l(1),\; \os_l(-2)\oplus\os_l(2).\]
The last one happens only for lines inside the ramification locus, whose number is finite and actually $0$ for a general $X$ hence the Hilbert scheme of lines is a surface whose only non-smooth points, if any, corresponds to those lines.
Moreover, this Hilbert scheme is irreducible and its general point corresponds to a line whose normal bundle is $\os_l^2$.
\ma

We gather some information on ideal sheaves of lines.
We will give a proof, but note that from \cite{ComaschiFaenzi}[Lem. 2.4], see also \cite{KuzProShr}[Lem. B.5.6], these are essentially equivalent to knowing the normal bundle of the line.

\Le\label{calcoliretta}
If $l\subset X$ is a line then 
\begin{itemize}
    \item $\id_l(1)$ is slope-stable with Hilbert polynomial $2\binom{t+3}{3}$
    \item $hom(\id_l,\id_l)=1$ and $ext^3(\id_l,\id_l)=0$
    \item if $l\not\subset R$ then $ext^1(\id_l,\id_l)=2$ and $ext^2(\id_l,\id_l)=0$ otherwise they are, respectively, $3$ and $1$.
\end{itemize}
\ma
\begin{proof}
    $\id_l$ is slope-stable being torsion-free of rank $1$, in particular $hom(\id_l,\id_l)=1$.
    From the sequence
    \sesl{\id_l(1)}{\os_X(1)}{\os_l(1)}{idretta}
    recalling that $f_*\os_X(1)\cong \osn(1)\oplus\osn(-1)$ and Using Riemann--Roch on $l$ we get
    \[P(\id_l(1))(t)=P(\os_X(1))(t)-P(\os_l(1))(t)=\chi(\os_X(1+t))-\chi(\os_l(1+t))=\]
    \[=\binom{t+4}{3}+\binom{t+2}{3}-(t+2)=2\binom{t+3}{3}.\]
    Applying $Hom(-,\id_l(1))$ to \eqref{idretta} we get the long exact sequence
    \[\dots \ra Ext^1(\os_X,\id_l)\ra Ext^1(\id_l,\id_l)\ra Ext^2(\os_l,\id_l)\ra Ext^2(\os_X,\id_l)\ra\]
    \[\ra Ext^2(\id_l,\id_l)\ra Ext^3(\os_l,\id_l)\ra Ext^3(\os_X,\id_l)\ra Ext^3(\id_l,\os_l)\ra 0\]
    Note that $ext^i(\os_X,\id_l)=h^i(X,\id_l)=0$ for all $i$ hence we get $ext^i(\id_l,\id_l)=ext^{i+1}(\os_l,\id_l)$ which for $i=3$ is clearly $0$.

    Moreover, applying $Hom(\id_l(1),-)$ to \eqref{idretta} and using $ext^1(\os_X,\os_l)=h^1(X,\os_l)=0$ together with $hom(\id_l,\id_l)=hom(\id_l,\os_X)$ we get 
    \[ext^1(\id_l,\id_l)=hom(\id_l,\os_l)=h^0(l,\n_l)\]
    which has been already computed during the proofs of \cite{Wel}(1.1) and (1.2)].
    Finally, since $ext^j(\os_X(i),\id_l)=0$ for all $j$ and $i=0,1$, by definition the sheaf $\id_l$ sits in the Kuznetsov component $\mathscr{K}$ of the semi-orthogonal decomposition $\scr{D}^b(X)=<\scr{K},\os_X,\os_X(1)>$ see \cite{Kuznetsov_Perry}[§8.1].
    We denote by $[n]$ the functor which consists in shifting of $n$ places to the left an element in $\scr{D}^b(X)$.
    Therefore, using the description of the Serre functor $S_\scr{K}(-)=\iota^*(-)[2]$ in the above reference, we have
    \[ext^2_{\os_X}(\id_l,\id_l)=hom_{\scr{D}^b(X)}(\id_l,\id_l[2])=hom_{\scr{K}}(\id_l,\id_l[2])=hom(\id_l,\iota^*\id_l)\]
    and, by stability, this last group is $0$ except when $\iota^*\id_l=\id_l$, that is $l\subset R$.
\end{proof}

We will also need to understand the restriction of an Ulrich bundle $\e$ on lines on $X$.
We will show that fixed a line, the general Ulrich bundle restricts in a balanced way to it.

\Le\label{restrette}
Fixed a rank $2$ Ulrich bundle $\e$, for any line $l\subset X$ we have $\e|_l\cong \os_l(1)^2$ or $\e|_l\cong \os_l\oplus\os_l(2)$.
If $l$ is a fixed line, then a general $\e$ among the ones constructed in \thref{deformazione} has a section whose vanishing locus is smooth and meets (schematically) $l$ in exactly one point, so that $\e|_l\cong \os_l(1)^2$.
\ma
\begin{proof}
    Since $\e$ is globally generated, also $\e|_l$ is, but being $det(\e)\cong\os_X(2)$ and any vector bundle on $\p^1$ split we have only the two above possibilities.
    
    Call $\scr{H}$ the universal family on the Hilbert scheme $\mathcal{H}^4_1$ of curves of degree $4$ and arithmetic genus $1$ on $X$; denote by $q_1,q_2$ the two projections on $X,\mathcal{H}^4_1$.
    Given any line $l\subset X$, the scheme $q_1\inv(l)$ has codimension $2$ in $\scr{H}$ but intersects the general fiber of $q_2$, which is $1$-dimensional, in a finite number of points, hence $q_2(q_1\inv(l))\subset \mathcal{H}^4_1$ is a divisor; call it $\mathcal{D}$.
    \dia
    & \scr{H} \ar[rd, "q_1"]\ar[ld, "q_2"] & &  \mathcal{J} \ar[rd, "\pi_1"]\ar[ld, "\pi_2"] &  \\
    X & & \mathcal{H}^4_1 & & \vert\os_X(2)\vert \\
    \mma
    Now we consider the incidence variety $\mathcal{J}\subset |\os_X(2)|\times \mathcal{H}^4_1$ of pairs $(S,C)$ such that $C\subset S$ and $\pi_1,\pi_2$ call the projections onto the two factors.
    
    Take some irreducible component $\mathcal{D}_0$ of $\mathcal{D}$ meeting the locus of smooth curves $C$ contained in a hyperplane $\Sigma$ such that $l\not\subset\Sigma$.
    This exists because in the linear system of $C\subset \Sigma$ there are smooth curves passing through the unique point $l\cap \Sigma$.
    The general fiber of $\pi_2$ are all isomorphic to $\p^2\cong \p(H^0(X,\id_Y(2))\du)$, then $\pi_2\inv(\cd_0)$ is an irreducible divisor in $\mathcal{J}$.
    Its general point is of the form $(S,C)$ with $C\cap l$ just one point and $C$ smooth, by semicontinuity of cohomology and openness of the smooth locus, since $\pi_2\inv(\cd_0)$ contains pairs with $C\subset \Sigma$ a smooth curve.
    The general fiber of $\pi_1$ is a finite union of $\p^1$-s, since such is the linear system of an elliptic curve on a $K3$ surface.
    Then $\pi_2\inv(\cd_0)$ intersects the general fiber of $\pi_1$ in a finite number of points, since for $(S,C)$ general we have $|\os_S(C)|\cong \p^1$ at most two of those curves can intersect $l$, hence $\pi_1(\pi_2\inv(\cd_0))$ is dominant on the image of $\pi_1$. 
    We have just proved that, if $[S]$ is general in the image of $\pi_1$, meaning that is general among the surfaces in $|\os_X(2)|$ containing an elliptic curve $C$ of degree $4$, then we can assume that $[S]$ is in $\pi_1(\pi_2\inv(\cd_0))$, hence contains an elliptic curve $C$ of degree $4$ and $C\cap l=Q$ is just one point.
    Therefore, a general Ulrich bundle constructed in \thref{deformazione} can be written as
    \ses{\os_X}{\e}{\id_C(2)}
    with $C\cap l=Q$.
    Restricting 
    \sesl{\id_C}{\os_X}{\os_C}{idc} to $l$ we get a surjective map $\id_C|_l\twoheadrightarrow \id_{Q/l}\cong\os_l(-1)$.
    Then, by composition we get a surjective morphism $\e|_l\twoheadrightarrow\id_C(2)|_l\ra\os_l(1)$ whose kernel, by degree reasons, must be $\os_l(1)$ and we conclude.
    \end{proof}

\Co\label{restrizioni_sezioni_rette}
Call $\iota$ the involution of the covering $f:X\ra \p^3$.
For any line $l\subset X$ not contained in the ramification locus the following are equivalent
\begin{enumerate}[1)]
    \item $\e|_l\cong \os_l(1)^2$
    \item the natural restriction $H^0(X,\e)\ra H^0(\iota^*l,\e|_{\iota^*l})$ is surjective
\end{enumerate}
\io
\begin{proof}
Note that by \thref{hilbertpol} we have $h^0(X,\e)=4$ while $h^0(\iota^*l,\e|_{\iota^*l})=4$ by \thref{restrette}.
Therefore, $H^0(X,\e)\ra H^0(\iota^*l,\e|_{\iota^*l})$ is surjective if and only if it is injective if and only if $H^0(X,\e\otimes\id_{\iota^*l})=0$ if and only if there is no section of $\e$ vanishing on $\iota^*l$.

\gel{1)$\Rightarrow$ 2}
Any section of $\e$ vanishing on $\iota^*l$ either vanishes also on $l$ or, once restricted to $l$, gives a section of $\e|_l$ vanishing on $l\cap \iota^*l$.
The first possibility is excluded by \thref{restrizione} since $l\cup \iota^* l$ is the complete intersection of two hyperplane sections on $X$ while the second cannot happen being $\e|_l\cong \os_l(1)^2$.

\gel{2)$\Rightarrow$ 1}
By \thref{restrette}, the only other possibility except $\e|_l\cong \os_l(1)^2$ would be $\e|_l\cong \os_l\oplus\os_l(2)$ in which case we would get a trivial quotient $\e\twoheadrightarrow\os_l$.
Hence $l$ can be lifted to a curve in $\p(\e)$ on which $\os_{\p(\e)}(1)$ is trivial, by universal property of $\p(\e)$.
But all such curves project to sections of $\iota^*\e$ by \thref{nonampio} hence, by applying $\iota^*$, there is a section of $\e$ whose zero locus contains $\iota^*l$ hence we reach a contradiction.
\end{proof}

We need to know some ext groups to apply \thref{wildext}, but in the final section we will recur back again to similar computations, so we get ahead with the work in the next lemma.

\Le\label{contiextqds}
Suppose $l$ is a line on $(X,H)$ and $\e$ a rank $2$ Ulrich bundle, then 
\begin{enumerate}[i)]
\item $hom(\id_l(1),\e)=0=hom(\e,\id_l(1))$
\item $ext^1(\e,\id_l(1))=2=ext^1(\id_l(1),\e)$
\item $ext^j(\id_l(1),\e)=0=ext^j(\e,\id_l(1))$ for $j=2,3$. 
\end{enumerate}
Moreover, if $\e$ is general respect to $l$ then also $ext^1(\id_l(1),\e(-1))=0=ext^1(\e,\id_l)$.
\ma
\begin{proof}
\gel{i}
 We have seen in \thref{calcoliretta} that $\e,\id_l(1)$ share the same reduced Hilbert polynomial.
   Since they are also stable but not isomorphic then any morphism among them is zero.
   
   \gel{ii), iii}
Recall that by \thref{restrette} we have $\e|_l\cong \os_l(1)^2,\os_l\oplus\os_l(2)$.
    Note that for a rank $2$ Ulrich bundle $\e$ we have $\e\du\cong\e(-2)$ hence
    \[ext^j(\id_l(1),\e)=ext^j(\id_l(1)\otimes\e\du,\os_X)=ext^j(\id_l\otimes\e(-3),\omega_X)=h^{3-j}(X,\id_l\otimes\e(-3))\]
    by Serre duality, while \[ext^j(\e,\id_l(1))=ext^j(\os_X,\e\otimes\id_l(-1))=h^j(X,\e\otimes\id_l(-1)).\]
    Twisting by $\e(l)$ for $l=-1,-3$ the sequence
    \sesl{\id_l}{\os_X}{\os_l}{seqretta'}
    we get $h^{k}(X,\id_l\otimes\e(l))=h^{k-1}(X,\e|_l(l))$ for all $k$ since $h^k(X,\e(l))=0$ for all $k$ being $l=-1,-3$ and $\e$ Ulrich.
    Therefore, for $j=1,2,3$ we have  
    \[ext^j(\id_l(1),\e)=h^{2-j}(X,\e|_l(-3))=2,0,0 \qquad ext^j(\e,\id_l(1))=h^{j-1}(X,\e|_l(-1))=2,0,0.\]

    For the last claim, reasoning as above we get 
    \[ext^1(\id_l(1),\e(-1))=h^{2}(X,\id_l\otimes\e(-2))\qquad ext^1(\e,\id_l)=h^1(X,\e\otimes\id_l(-2))\]
    so that twisting by $\e(-2)$ the sequence \eqref{seqretta'} and using that $\e$ is Ulrich we get 
    \[ext^1(\id_l(1),\e)=h^{1}(X,\e|_l(-2))=0 \qquad ext^1(\e,\id_l(1))=h^{0}(X,\e|_l(-2))=0,\]
    by generality of $\e$ and $\thref{calcoliretta}$.
\end{proof}

\subsection{Odd rank Ulrich bundles}

We are going to show that on a general quartic double solid there are Ulrich bundles of any rank $r\geq 2$.
The even rank case could be treated as in \thref{rango>sextic} but we will follow the idea in \cite{CFK3}, which holds for all Fano $3$-folds of index $2$ with $-\frac{K_X}{2}$ globally generated, and construct arbitrary rank Ulrich bundles by extensions and deformations of rank $2$ ones and ideal sheaves of lines.
In practice, the strategy to construct such bundles consists of applying the standard machinery of \thref{wildext}, while proving that they are actually Ulrich will definitely require more care.
Indeed, in our proof, we will need to consider some auxiliary sheaves obtained as elementary transformations of Ulrich ones along line bundle on lines.
We start with those.

\Le\label{inattesi}
On a smooth, complex quartic double solid $X$ consider a line $l$ with normal bundle $\n_l\cong\os_l^2$ and a sufficiently general stable rank $2$ Ulrich bundle $\e$ such that $ext^2(\e,\e)=0$ and $\e|_l\cong \os_l(1)^2$.
Define $\e_l$ to be a kernel of the form
\sesl{\e_l}{\e}{\os_l(1)}{tnat}
Then $\e_l$ is a rank $2$ slope-stable sheaf such that $ext^2(\e_l,\e_l)=0$ hence $[\e_l]$ is a smooth point of its moduli space which has dimension $9$. 
Moreover, a general deformation $\e_l'$ of $\e_l$ is locally free and we have $ext^1(\os_X(1),\e_l')=1$.
\ma
\begin{proof}
    We have $\mu(\e)=\mu(\e_l)$ and any slope-destabilising subsheaf for $\e_l$ would be the same for $\e$, which has no slope-destabilising subsheaves being slope-stable hence also $\e_l$ is slope-stable.
    By \eqref{tnat} and the Ulrich condition on $\e$ we have
   \[ext^j(\os_X(1),\e_l)=h^j(X,\e_l(-1))=h^{j-1}(X,\os_l)=0\]
   for all $j\ne 1$ so the same holds for a general deformation $\e_l'$.
   Moreover, for $j=1$ we get $ext^j(\os_X(1),\e_l)=h^{0}(X,\os_l)=1$ and again the same must hold for $\e_l'$, by semicontinuity of cohomology and being $\chi(\e_l)=\chi(\e_l')$.

   By \thref{restrette} a sufficiently general $\e$ has a section $\sigma$ whose zero locus $Y$ is smooth and meets (schematically) $l$ in just one point.
   This implies that $\sigma$ is sent by $H^0(X,\e)\ra H^0(l,\e|_l)\cong H^0(l,\os_l(1)^2)$ to a global section of one of the two summands of $\os_l(1)^2$.
   Therefore, considering the projecting $\e\ra\e|_l\ra \os_l(1)$ to the other summand we conclude that $\sigma$ comes from a section of $ker(\e\ra\os_l(1))$, therefore we can form the following commutative diagram
    \dia
       0 \ar[r] & \os_X \ar[r] \ar[d, "id"] &  \e_l \ar[r] \ar[d] & \id_C(2) \ar[r] \ar[d] & 0 \\
       0 \ar[r] & \os_X \ar[r]  &  \e \ar[r] \ar[d] & \id_Y(2) \ar[r] \ar[d] & 0 \\
       & & \os_l(1) \ar[r, "id"]  & \os_l(1) &  \\
   \mma
   From this we can construct
   \dia
       &  &  &  \os_l(-1) \ar[d] & \\
       0 \ar[r] & \id_C \ar[r] \ar[d] & \os_X \ar[d, "id"] \ar[r] & \os_C \ar[r] \ar[d] & 0 \\
       0 \ar[r] & \id_Y \ar[r] \ar[d] & \os_X \ar[r]  & \os_Y \ar[r] & 0 \\
       &  \os_l(-1)  &  &  & \\
   \mma
   hence $C=Y\cup l$ and $\os_l(-1)\cong \id_{Y/C}\cong \id_{Y\cap l/l}$ so $length(Y\cap l)=1$.
   It follows that $C$ is a nodal curve, in particular locally complete intersection, and, being $Y$ of degree $4$ and arithmetic genus $1$ by \thref{rk2-section}, we deduce that $C$ has degree $5$ and arithmetic genus $1$.
    By \thref{HSrk2} we have $h^1(Y,\n_Y)=ext^2(\e,\e)=0$ while $h^1(Y,\n_Y(-1))=0=h^1(Y,\n_Y)$ being $l$ general enough, hence as in \cite{HarHir}[Thm. 4.1] we conclude that $h^1(C,\n_C)=0$ and $C$ can be smoothed.

    By Riemann--Roch on $C$ it follows that 
    \[h^0(C,\n_C)=\chi(\n_C)=deg(det(\n_C))=2H\cdot C=10.\]
    Note that \cite{ComaschiFaenzi}[Lem. 2.4] implies that $ext^1(\id_C,\id_C)=h^0(C,\n_C)=10$ and $ext^2(\id_C,\id_C)=h^1(C,\n_C)=0$.
    We have $hom(\id_C(2),\os_X)=hom(\id_C,\os_X(-2))=0$ by \thref{stabilitàtrick} since they are both slope-stable by definition.
    By Serre duality
    \[ext^i(\id_C(2),\os_X)=ext^i(\id_C,\os_X(-2))=h^{3-i}(X,\id_C),\]
    which is $0$ for $i=2$ being $C$ connected.
    Therefore, applying $Hom(\id_C(2),-)$ to
    \sesl{\os_X}{\e_l}{\id_C(2)}{solitac}
    we get $ext^2(\id_C(2),\e_l)=0$ and 
    \[ext^1(\id_C(2),\e_l)=ext^1(\id_C(2),\id_C(2))+ext^1(\id_C(2),\os_X)-hom(\id_C,\id_C)=10\]
    being, using also the previous computations, 
    \[ext^1(\id_C(2),\os_X)=h^{2}(X,\id_C)=h^{1}(X,\os_C)=1.\]
    Furthermore, applying $Hom(-,\e_l)$ to \eqref{solitac} we conclude $ext^2(\e_l,\e_l)=0$, being $ext^2(\os_X,\e_l)=h^2(X,\e_l)=0$, and 
    \[ext^1(\e_l,\e_l)=ext^1(\id_C(2),\e_l)+ext^1(\os_X,\e_l)+hom(\e_l,\e_l)-hom(\os_X,\e_l)=\]
    \[=11+h^1(X,\e_l)-h^0(X,\e_l)=11+h^0(X,\os_l(1))-h^0(X,\e)=9\]
    where in the last step we used the cohomology sequence of \eqref{tnat}.
    We just proved smoothness and computed the dimension of the moduli space of $\e_l$.

    Finally, by Serre duality 
    \[ext^i(\id_C(2),\os_X)=ext^i(\id_C,\os_X(-2))=h^{3-i}(X,\id_C)=1,0\]
    depending whether $i=1$ or not, where in the last step we used the sequence
    \ses{\id_C}{\os_X}{\os_C}
    and the fact that $C$ is connected and of genus $1$.
    so being cohomology semicontinuous and Euler characteristics constant in flat families we get $ext^1(\id_{C'}(2),\os_X)=1$ for a general deformation $C'$ of $C$.
    Then, the universal family over some open neighbourhood of $[C]$ in its Hilbert scheme satisfies the assumptions in \thref{modularcorrispondenza} hence we obtain a flat family of sheaves corresponding to those curves through Hartshorne-Serre.
    Since $C$ is smoothable, the general element in the Hilbert scheme will be a smooth genus $1$ curve $C'$, in particular $\omega_{C'}\cong \os_{C'}$ is globally generated by just one section.
    Then, being $h^1(X,\os_X(-2))=0=h^2(X,\os_X(-2))$, by \thref{r00} and \thref{vectorcorrispondenza} the sheaf corresponding to such a $C'$ is a vector bundle which specialises to $\e_l$.
\end{proof}

\Oss
\begin{itemize}
    \item The locally free deformations of sheaves of the form $\e_l$ as in \eqref{tnat} have already been studied: indeed, since they come through Hartshorne-Serre by elliptic quintics, they coincide with $\e_C(1)$ where $\e_C$ is as in \cite{MarTik_qds}[Eq. (6.1)].
    In particular, their argument not only computes the dimension of the moduli space containing $\e_l$ but also gives its irreducibility.
    \item Once we prove the existence of Ulrich bundles of rank higher than $2$ and that they restrictions to lines behave well, see \thref{restrette'}, then we could repeat the construction of the elementary transformations in \eqref{tnat} with $\e$ replaced by one of those sheaves.
    However, it is not clear to us whether the corresponding curve $C$ can be smoothed in this general situation.
\end{itemize}
\one

Now we are ready for our main result. 
Note that, by \thref{calcoliretta}, $\id_l(1)$ has the same reduced Hilbert polynomial of an Ulrich sheaf moreover it has almost all the vanishings to make it an Ulrich sheaf.
The following theorem has roots in this similarity.
This observation has already been used for an analogous purpose in \cite{LMS} and in \cite{CFK3}.

\Te\label{qdswild}
Let $f:X\ra\p^n$ be a smooth complex quartic double solid.
Then there are slope-stable Ulrich bundles of any odd rank $r\geq 3$ that are deformations of sheaves $\f$ fitting in
\sesl{\e}{\f}{\g}{extul}
with $\e$ a general Ulrich bundle of rank $2$ and $\g$ either a general Ulrich bundle of odd rank or $\id_l(1)$ for some general line $l\subset X$.
The components of the moduli space containing those Ulrich bundles have dimension $r^2+1$ and are reduced since for a general such sheaf $\f'$ we have $ext^2(\f',\f')=0$.
\Ma
\begin{proof}
\textbf{Step 1: Construction of the families}\\
    We want to apply \thref{wildext} to the pair $(\scr{F}_1,\scr{F}_2)$, where $\scr{F}_1$ is an irreducible component of the family of twisted ideals of lines $\id_l(1)$ for $l\not\in R$ and $\scr{F}_2$ is an open subset of an irreducible component of the family of rank $2$ Ulrich bundles found in \thref{deformazione}.
    Thanks to \thref{qdsrk2}, we can choose $\scr{F}_2$ such that its members are not fixed by the involution of the covering.
    By \thref{calcoliretta} $\scr{F}_1$ is a nice family.
    Since there are no Ulrich line bundles on $X$, the sheaves in $\scr{F}_2$ are stable Ulrich bundles by \thref{stability}.
    By \thref{qdsrk2}, we conclude that also $\scr{F}_2$ is a nice family. 
    The sheaves in those families share the same reduced Hilbert polynomial $2\binom{t+3}{3}$.
    
    By \thref{qdsrk2}, for any $\e$ in $\scr{F}_2$ we have $ext^1(\e,\e)=5$, in particular $ext^1(\e,\e)>hom^1(\e,\e)+ext^2(\e,\e)$.
    Adding the computations in \thref{contiextqds}, we just verified that $(\scr{F}_1,\scr{F}_2)$ is a pair of wild extension.
    Hence, by the first part of \thref{extpreciso}, for any $r\geq 3$ we get smooth families $\scr{F}_r$, over some base $\cu_r$, of rank $r$ Gieseker stable torsion-free sheaves on $X$ presented as deformations of extensions of the form \eqref{extul}.
    Specifically, we want to use the form of \thref{extpreciso} written in parentheses; hence, for $r=3$ we have an extension of a rank $2$ Ulrich bundle and $\id_l(1)$ while for the subsequent ones we have an extension of two Ulrich bundles, where the one on the left has rank $2$.
    Moreover, by \thref{qdsrk2} we have $ext^3(\e,\e)=0$, so, using other vanishings in \thref{contiextqds}, we can apply the second part of \thref{extpreciso} and easily compute $dim(\cu_r)=r^2+1$ .
    Note that being $det(\e)\cong \os_X(2)$ and $det(\id_l(1))\cong\os_X(1)$, we get that a rank $r$ sheaf constructed in the above way has determinant isomorphic to $\os_X(r)$.

\textbf{Step 2: Reduction to the case}$\mathbf{r=3}$\\
We prove by induction that for any $r\ge 3$ the general sheaf in $\mathscr{F}_{r}$ is a slope-stable Ulrich bundle.
Essentially, we need only to treat the base case $r=3$ since, as described before, in the inductive step we are taking extensions of two Ulrich bundle so it is again Ulrich by \thref{2-3} and it is also slope-stable because any stable Ulrich bundle is such by \thref{stability}.
Therefore, we are only left to prove the base case $r=3$.

Note that any slope-stable sheaf in $\mathscr{F}_{3}$ is automatically Ulrich: indeed, from \eqref{extul} we get $h^j(X,\f(i))=0$ for all $j$ if $i=-1,-2$ and for $j=0,1$ if $i=-3$ but by slope-stability we also have
\[h^3(X,\f(-3))=hom(\f(-3),\omega_X)=hom(\f,\os_X(1))=0\]
so, being $\chi(\f(-3))=0$, we also deduce $h^2(X,\f(-3))=0$, i.e. $\f$ is Ulrich.
Being $X$ smooth an Ulrich sheaf is locally free by \thref{corollario}.\footnote{Actually, as shown in \cite{CFK3}[Lem. 4.9], it is a direct consequence of \thref{restrette} that already non-split extensions as in \eqref{extul} are locally free.}

\textbf{Step 3: There is a flat family over an irreducible base containing the general sheaf in $\scr{F}_3$ and some} $\mathbf{\widetilde{\f}}$ \textbf{as in}
\sesl{\e_l'}{\widetilde{\f}}{\os_X(1)}{cattivi'}
\textbf{where} $\mathbf{\e_l'}$ \textbf{is a general deformation of } $\mathbf{\e_l}.$\\
Consider all the semistable sheaves which have the same Hilbert polynomial as $\e\oplus\id_l(1)$.
Those form a bounded family by \cite{HuyLeh}[Thm. 3.3.7] hence by \cite{HuyLeh}[Lem. 1.7.6 ii)], the Castelnuovo-Mumford regularity of those sheaves is bounded.
Then, by \cite{Mum}[Lecture 14 Proposition] we can suppose they are quotients of $\os_X(-l)^\rho$ for some $l,\rho>>0$, that is, they appear as members of some Quot-scheme $\cc$.
Clearly among them we have $\e\oplus\id_l(1)$, we show that they are smooth points in $\cc$.
Indeed, consider an exact sequence 
\sesl{\ck}{\os_X(-l)^\rho}{\e\oplus \id_l(1)}{q}
so that by \cite{HuyLeh}[Prop. 2.2.8] it is enough to check $ext^1(\ck,\e\oplus \id_l(1))=0$.
This holds true by taking $Hom(-,\e\oplus \id_l(1))$ in \eqref{q} because
\[ext^1(\os_X(-l)^\rho,\e\oplus \id_l(1))=\rho \cdot h^1(X,\e(l)\oplus \id_l(1+l))=0\]
by Serre vanishing, being $l>>0$, and $ext^2(\e\oplus \id_l(1),\e\oplus \id_l(1))=0$ by the computations in Step $1$ and \thref{calcoliretta}.
We deduce that there is only one irreducible component of $\cc$ containing $\e\oplus \id_l(1)$, so we can suppose $\cc$ irreducible. 
It follows that $\cc$ also contains all semistable deformations of $\e\oplus \id_l(1)$, such as the general element in the family $\scr{F}_3$.

Consider a sheaf $\e_l$ fitting in \eqref{tnat} and a split extension in \eqref{extul} with $\g=\id_l(1)$, that is $\f=\e\oplus\id_l(1)$.
We can form the commutative diagram
\dia
 & & & \os_l(1) \ar[d] \\
       0 \ar[r] & \e_l \ar[r] \ar[d] & \e\oplus\id_l(1) \ar[d, "id"] \ar[r] & \id_l(1)\oplus\os_l(1) \ar[r] \ar[d] & 0 \\
       0 \ar[r] &  \e \ar[r] \ar[d] & \e\oplus\id_l(1) \ar[r]  & \id_l(1) \ar[r]  & 0 \\
       & \os_l(1)  & &  \\
\mma
Note that, $\os_X(1)$ can be regarded as a non-split extension
\ses{\id_l(1)}{\os_X(1)}{\os_l(1)} hence $\f_l$ as in 
\sesl{\e_l}{\f_l}{\os_X(1)}{cattivi}
specialises to $\e\oplus\id_l(1)$.
But $\widetilde{\f}$ as in \eqref{cattivi'} specialises to $\f_l$ hence to $\e\oplus\id_l(1)$, so that a general $\widetilde{\f}$ is also semistable.
Therefore, the universal family of $\cc$ is a flat family $\scr{C}$ over an irreducible base and contains the general sheaf of the form $\widetilde{\f}$ and the general sheaf in $\scr{F}_3$.
 
\textbf{Step 4: Slope-stability}\\
The family $\scr{F}_3$ over the irreducible base $\cu_3$ gives a morphism $\psi:\cu_3\ra \cm$, where $\cm$ is an irreducible component of the corresponding moduli space of semistable sheaves.
Being $\cu_3$ an irreducible component of a modular family by construction, the general isomorphism class of simple, in particular stable, sheaves in $\cm$ appears in $\scr{F}_3$ and it does at most finitely many times, so that $\psi$ is dominant and generically finite.
In particular, $dim(\cm)=dim(\cu_3)=10$.
Moreover, also the family $\scr{C}$ constructed in the previous step gives a morphism $\psi':\cc\ra \cm$, which is again dominant since $\scr{C}$ contains the general element in $\scr{F}_3$.
We assume by contradiction that the general sheaf in $\cm$ is not slope-stable, then the same applies to the general element in $\scr{C}$.
Then, we can apply \thref{crit_stab} to the family $\scr{C}$ which gives, up to base change to another irreducible Quot-scheme, a slope-destabilising quotient $\scr{C}\twoheadrightarrow\scr{Q}$, with the general element of $\scr{Q}$ torsion-free.

Let us specialise this surjection to a quotient $\widetilde{\f}\twoheadrightarrow \widetilde{Q}$, where $\widetilde{\f}$ is as in \eqref{cattivi'}.
Recall that by \thref{inattesi} we can assume $\e_l'$ in \eqref{cattivi'} to be locally free and the extension non-split.
Let us consider the case in which $\widetilde{Q}$ has torsion.
It follows by \thref{estensionislope} and \eqref{cattivi'} that the torsion-free part of $\widetilde{Q}$ must be $\os_X(1)$, in particular $P(\widetilde{Q})>P(\os_X(1))$.
Moreover, a general element $Q'$ of $\scr{Q}$ would be a torsion-free deformation of $\widetilde{Q}$, hence a torsion-free rank $1$ sheaf with $P(Q')=P(\widetilde{Q})$, in particular $det(Q')=det(\os_X(1))$ being $Pic(X)$ discrete.
But then we would have an embedding of $Q'$ in its double dual, which must be $\os_X(1)$ and hence $P(Q')\le P(\os_X(1))<P(\widetilde{Q})$, a contradiction. 

Therefore, $\widetilde{Q}$ is torsion-free and again by \thref{estensionislope} and \eqref{cattivi'} $\widetilde{Q}=\os_X(1)$ and $Ker(\widetilde{\f}\twoheadrightarrow \widetilde{Q})=\e_l'$ as in \eqref{cattivi'}.
The general element in $\scr{Q}$ is a locally free deformation of $\os_X(1)$ hence coincides with $\os_X(1)$ while the general element in $ker(\scr{C}\twoheadrightarrow\scr{Q})$ must be a deformation of $\e_l'$.
Since $ext^1(\os_X(1),\e_l')=1$ by \thref{inattesi}, the isomorphism class of a general element in $\scr{C}$ is the unique non-split extension of some $\e_l'$ and $\os_X(1)$, meaning that the moduli space containing $\e_l'$ dominates the image of $\psi'$ in $\cm$.
But the former has dimension $9$ by \thref{inattesi} while we have already seen that $Im(\psi)$ is dense in $\cm$ which has dimension $10$, giving our contradiction.
\end{proof}

\Oss\label{R}
\begin{itemize}
    \item From Step $3$, we know that the closure of the locus of extensions $\widetilde{\f}$ meets the locus of extensions $\f$ at least in the points representing sheaves of the form $\e\oplus\id_l(1)$.
    We do not know if the closure of the former locus contains the latter, a fact that could simplify the above proof.
    \item The above argument is highly inspired by the one given in \cite{CFK3}[§4] to prove the analogous statement on other Fano $3$-folds of index $2$ with $-\frac{K_X}{2}$ very ample.
    We believe that there is a small gap in the presentation given in \cite{CFK3}[Lem. 4.13], since it is not clear to us why $\cq^*$ should specialise to $\cq_0^*\cong \e_r$, being the sheaves $\cq,\cq_0$ not locally free.
    A posteriori, the existence of sheaves $\widetilde{\f}$ as in \eqref{cattivi'} suggests that more work needs to be done.
    However, our proof actually carry over even in their case.
\end{itemize}
\one

\Co\label{restrette'}
For every integer $r\geq 2$ the moduli space of rank $r$ Ulrich bundles has a reduced component of dimension $r^2+1$ whose general element $\f'$ satisfies $ext^2(\f',\f')=0$ and $\f'|_l\cong \os_l(1)^r$ for a general line $l$.
\io
\begin{proof}
   We already saw in \thref{rango>sextic} that for $r$ even the desired moduli spaces have dimension $r^2+1$ and such bundles are obtained as deformations of subsequent extensions as in
\sesl{\e}{\f}{\g}{extr'}
where $\e,\g$ are general rank $2,r-2$ Ulrich bundles.
Analogously, \thref{qdswild} computes the dimension of the moduli space of Ulrich bundles of odd rank obtained as deformations of $\f$ in \eqref{ext'}, where $\e$ is again Ulrich of rank $2$ and $\g$ could be either $\id_{l'}(1)$ for some line $l'$ or some odd rank Ulrich bundle.
The condition on the second ext-group has also been verified in the above results, so we need only to show the one regarding restriction to lines.

Recall that by \thref{restrette} the thesis is clear for a general $\e$ as above while, for a sufficiently general line $l'\subset X$ we can suppose $l\cap l'=\emptyset$, in particular $\id_l(1)|_{l'}\cong \os_{l'}(1)$ as desired.
We deduce the claim for any $\f$ as in \eqref{ext'}, where $\g$ is either $\id_{l'}$ or Ulrich of rank $2$.
But then, for a general $\f'$ which is deformation of $\f$ we have $\f'|_{l'}\cong \os_{l'}(1)^{r+1}$ since, being $det(\f')\cong \os_X(r+1)$, the above condition on the restriction is equivalent to $h^0(l',\f|_{l'}(-2))=0$, which is open by semicontinuity of cohomology.
Therefore, by induction on the number of extensions we conclude.
\end{proof}

\subsection{Ulrich bundles on quartic double planes}

In this section, we will show that smooth hyperplane sections of our $X$ are Ulrich wild, having families of Ulrich bundles of any positive rank and increasing dimension.
Explicitly, first we classify all Ulrich line bundles and then focus on higher rank ones with the maximum possible $c_2$.

Take $\Sigma\in |\os_X(1)|$ a smooth surface, we follow the notation in \thref{deformazione}, where we have seen that $\Sigma$ is also a blow-up of $\p^2$ in $7$ points.
By abuse of notation we denote by $H:=-K_\Sigma\sim 3L-\sum_{i=1}^7E_i$; note that, by adjunction, this is actually our previous $H$ restricted to $\Sigma$.
We remark that our polarisation will always be $H$.

The fact that Ulrich line bundles on del Pezzo surfaces come from rational normal curves has been noted in \cite{PL_T} and \cite{CKM_delpezzo}[Prop. 2.9], see also \cite{CFK2}[Thm. 2.1].
\Prop
There are exactly $126$ isomorphism classes of Ulrich line bundles on any smooth $\Sigma\in |\os_X(1)|$ corresponding to conics on it, i.e. rational curves $C$ such that $CH=2$ and $C^2=0$.
\One
\begin{proof}
    Suppose a line bundle $\os_\Sigma(C)$ is Ulrich.
    In view of the fact that $c_2(\os_\Sigma(C))=0$ and by Riemann--Roch
    \[1=\chi(\os_\Sigma)=\dfrac{K_\Sigma^2+c_2(\Sigma)}{12}=\dfrac{2+deg(c_2(\Sigma))}{12} \; \Rightarrow \; deg(c_2(\Sigma))=10,\]
    by \thref{c12} we deduce that $CH=2$ and $C^2=0$.
    By adjunction $C$ must be rational since we have
    \[-2=C^2+CK_\Sigma=2p_a(C)-2.\]

    We can write down explicitly all the classes of curves in $Pic(\Sigma)$ satisfying the two equations $CH=2$ and $C^2=0$.
    Using the shortcut $(a; a_1,\dots a_7)$ for $C\sim aL-\sum_{i=1}^7a_iE_{i}$ and exponents to represent repeated values, the desired classes are of the form
    \[(1;1,0^6)\quad (2;1^4,0^3)\quad (3;2,1^5,0)\quad (4;2^3,1^4)\quad (5;2^6,1).\]
   Note that taking U-duals, i.e. replacing $C$ with $2H-C$, sends the coefficient $a$ to $6-a$ hence acts on the above five types.
   We can see that there are respectively
   \[7 \quad 35 \quad 42 \quad 35 \quad 7\]
   isomorphism classes of line bundles of the above types, for a total amount of $126$.
\end{proof}

For future use we study the extensions of those line bundles.

\Le\label{ext1''}
Suppose $C_1,C_2$ are two curves giving Ulrich line bundles on $\Sigma$.
Then 
\begin{itemize}
    \item $\chi(\os_\Sigma(C_1-C_2))=1-C_1\cdot C_2$
    \item $ext^1(\os_\Sigma(C_2),\os_\Sigma(C_1))=C_1\cdot C_2-1$ if $C_1\not\sim C_2$ and $0$ otherwise
    \item $ext^2(\os_\Sigma(C_2),\os_\Sigma(C_1))=h^2(\Sigma,\os_\Sigma(C_1-C_2))=0$
    \item $0\leq C_1\cdot C_2\leq 4$ and the extremal values are achieved when $C_1\sim C_2$ and $C_1\sim 2H-C_2$ respectively.
\end{itemize}
\ma
\begin{proof}
By Riemann--Roch we have
    \[\chi(\os_\Sigma(C_1-C_2))=1+\dfrac{(C_1-C_2)^2+(C_1-C_2)H}{2}=1-C_1\cdot C_2.\]
    Note that $(C_1-C_2)H=0$ hence $h^0(\Sigma,\os_\Sigma(C_1-C_2))=1$ if and only if $C_1\sim C_2$ and $0$ otherwise.
    By Serre duality $h^2(\Sigma,\os_\Sigma(C_1-C_2))=h^0(\Sigma,\os_\Sigma(C_2-C_1-H))=0$ for all $C_1,C_2$.
    It follows that 
    \[ext^1(\os_\Sigma(C_2),\os_\Sigma(C_1))=h^1(\Sigma,\os_\Sigma(C_1-C_2))=\]
    \[=-\chi(\os_\Sigma(C_1-C_2))+h^0(\Sigma,\os_\Sigma(C_1-C_2))=-1+C_1\cdot C_2+h^0(\Sigma,\os_\Sigma(C_1-C_2)).\]
    If $C_1\sim C_2$ this quantity is $0$ otherwise is exactly $C_1\cdot C_2-1$.

    We can easily estimate $C_1\cdot C_2$.
    Indeed, Ulrich bundles are globally generated by \thref{corollario} hence the curves corresponding to Ulrich line bundles are effective and base-point free divisors. 
    It follows that $C_1\cdot C_2\geq 0$ and we have equality if and only if $C_1\sim C_2$.
    Similarly, replacing $C_2$ with its U-dual, $C_1\cdot (2H-C_2)\geq 0$ and we have equality if and only if $C_1\sim 2H-C_2$, that is $C_1\cdot C_2\leq 2H\cdot C_1=4$ and we have equality if and only if $C_1\sim 2H-C_2$.
\end{proof}

In \cite{MRPL}[Thm. 4.10] it was proved, among other things, that for any $r\geq 2$ there is an irreducible family of dimension $r^2+1$ of simple, semistable, special Ulrich bundles with the maximal possible $c_2$, which is $r^2$ from \cite{MRPL}[Remark 4.11].
Let us remark that in \cite{MRPL}[Thm. 4.10] it is proved that the general such bundles are stable but with respect to a polarisation of the form $H+(n-3)L$ for some unknown $n>>0$.
We are particularly interested in these since they are the ones which extend to $X$, as we will see at the end of this work, \thref{restrizdomi}.
In \cite{CFK2}[Thm. 4.8] it is proven that on blow-ups of general points in $\p^2$ there are stable Ulrich bundles of any rank and they also compute the dimension of some irreducible and reduced components of the moduli spaces.

Nevertheless, we will reprove this fact.
Recall that a rank $r$ Ulrich bundle $\e$ on $\Sigma$ is special if $c_1(\e)\sim rH$.
Moreover, by \thref{c12}, if $\e$ is special then the degree of $c_2(\e)$ is determined and we can check its $r^2$, the maximum possible.
We will need a lemma.
\Le\label{lemmafine}
Let $(\Sigma,-K_\Sigma)$ be a smooth polarised del Pezzo surface of degree $2$.
\begin{enumerate}
    \item If $T$ is a $0$-dimensional subscheme of length $r^2$ with $h^0(\Sigma,\id_T(r-1))=0$ then $ext^1(\id_T(r-1),\os_\Sigma)=r-1$.
    \item If $\e$ is a special rank $r$ Ulrich bundle on $(\Sigma,-K_\Sigma)$ then it sits in a sequence
   \sesl{Ext^1(\id_{T}(r),\os_{\Sigma})\otimes\os_\Sigma}{\e}{\id_T(r)}{ulrichr}
   where $T$ is a $0$-dimensional subscheme of length $r^2$ with $h^0(\Sigma,\id_T(r-1))=0$.
\end{enumerate}
\ma
\begin{proof}
\gel{1}
It follows by 
   \ses{\id_T(r-1)}{\os_\Sigma(r-1)}{\os_T(r-1)}
   that 
   \[h^1(\Sigma,\id_T(r-1))=h^0(T,\os_T(r-1))-h^0(\Sigma,\os_\Sigma(r-1))=r^2-1-r(r-1)=r-1.\]
   Serre duality tells us that $ext^1(\id_T(r),\os_\Sigma)=h^1(\Sigma,\id_T(r-1))=r-1$. 
   
\gel{2}
    Since $\e$ is special, that is $c_1(\e)\sim rH$, we have $c_1(\e)^2=2r^2$ and by \thref{c12} we get $c_2(\e)=r^2$ (recall that $12=c_2(\Sigma)+K_{\Sigma}^2$ by Noether's formula).
   By \cite{Ful}[Example 14.4.3 i)] choosing $r-1$ general sections we get an injective morphism $\os_\Sigma^{r-1}\ra\e$, so by \thref{Ulrichcorrispondenza} we get an exact sequence 
   \ses{W\otimes\os_\Sigma}{\e}{\id_T(r)}
   with $W\subset Ext^1(\id_{T}(r),\os_{\Sigma})$ of dimension $r-1$.
   Since $T$ represents $c_2(\e)$ by \thref{hsc2}, it is a subscheme of length $r^2$.
   Moreover, $h^0(X,\e(-1))=0$ and $h^1(\Sigma,\os_\Sigma(-1))=0$ imply $h^0(\Sigma,\id_T(r-1))=0$.
   Therefore, by the first part we get $ext^1(\id_{T}(r),\os_{\Sigma})=r-1$ hence $W=Ext^1(\id_{T}(r),\os_{\Sigma})$.
\end{proof}

\Te\label{qdpwild}
On a smooth, complex, polarised quartic double plane $(\Sigma,-K_\Sigma)$ there are slope-stable special Ulrich bundles of rank $r$ for any $r\geq 2$ which fit in 
\ses{\os_\Sigma^{r-1}}{\e}{\id_T(r)}
with $T\subset \Sigma$ a general $0$-dimensional subscheme of length $r^2$.
Moreover, all Ulrich bundles of fixed rank are in the same irreducible component of the moduli space, which has dimension $r^2+1$ and is reduced.
\Ma
\begin{proof}
\textbf{Step 1: Construction of the families}\\
Choose an Ulrich line bundle $\cl$ on $\Sigma$.
   Then $(\cl,\cl\du(K_\Sigma))$ is a pair of good extension by \thref{ext1''} hence by \thref{wildext} we have a $5$-dimensional good family $\scr{F}_2$ of stable rank $2$ Ulrich bundles.
   Note that the above pair is not of wild extension because $\chi(\cl,\cl)=\chi(\os_\Sigma)=1\not <0$.
   Therefore, we need a new pair.
   We can easily check that $(\scr{F}_2,\scr{F}_2)$ would give us special, slope-stable even rank Ulrich bundles as in \eqref{rango>sextic}.
   Furthermore, by taking subsequent extensions of rank $1$ and rank $2$ Ulrich bundles and then deformations we could obtain stable Ulrich bundles of any rank, but those won't be special.
   Instead, for odd ranks we follow the proof given on the quartic double solid.
   
   We define $\scr{F}_1$ to be the family of ideal sheaves of points in $\Sigma$ twisted by $\os_\Sigma(1)$.
   This is constructed as follows: take the ideal sheaf of the diagonal in $\Sigma\times\Sigma$ and twist it by the pullback of $\os_\Sigma(1)$.
   For any point $Q\in \Sigma$ clearly $\id_Q(1)$ is slope stable and $\Sigma$ itself is the moduli space of those sheaves.
   In particular, being a smooth variety of dimension $2$ we get $Ext^1(\id_Q,\id_Q)\cong T_{Q}\Sigma\cong \C^2$.
   From the sequence
   \sesl{\id_Q}{\os_\Sigma}{\os_Q}{idw}
   we obtain $0=h^2(\Sigma,\id_Q)=ext^2(\os_\Sigma,\id_Q)$ since by Serre duality $h^2(\Sigma,\os_\Sigma)=h^0(\Sigma,\omega_\Sigma)=0$.
   By applying $Hom(-,\id_Q)$ to \eqref{idw} we conclude that $ext^2(\id_Q,\id_Q)=0$, hence $\scr{F}_1$ is a good family.
   Moreover, we get $\chi(\id_Q,\id_Q)=-1$ hence, being $hom(ext^1(\id_Q,\id_{Q'})=0$ by stability, we get $ext^1(\id_Q,\id_{Q'})=1$ for $Q\neq Q'$.
   
   Now consider the new pair $(\scr{F}_2,\scr{F}_1)$. 
   In order to show that it is of wild extension, being $\Sigma$ a surface and $\chi(\id_Q,\id_Q)<0$, is enough to show that for any Ulrich bundles $\e$ of rank $2$ and any point $Q\in \Sigma$
   \[ext^1(\id_Q(1),\e)=2=ext^1(\e,\id_Q(1)) \qquad ext^2(\id_Q(1),\e)=0=ext^2(\e,\id_Q(1)).\]
   For any rank $2$ Ulrich bundle $\e$ we have $\e\du\cong\e(-2)$. 
   Therefore, for $i=1,2$
   \[ext^i(\e,\id_Q(1))=h^i(\Sigma,\e\du\otimes\id_Q(1)=h^i(\Sigma,\e\otimes\id_Q(-1)),\]
   which can be computed by twisting \eqref{idw}.
   Similarly we calculate
   \[ext^i(\id_Q(1),\e)=ext^i(\e\du\otimes\id_Q(1),\os_\Sigma)=ext^i(\e\otimes\id_Q(-2),\omega_\Sigma)=h^{2-i}(\Sigma,\e\otimes\id_Q(-2)),\]
   where in the last step we used Serre duality.
   By \thref{extpreciso}, for any odd $r$ we get smooth $r^2+1$-dimensional families $\scr{U}_r$ of stable, special sheaves as deformations of extensions.
   Clearly, among them we have restrictions of Ulrich bundles on $X$ obtained in \thref{qdswild}, those remain Ulrich bundles hence the general element in $\scr{U}_r$ must be locally free and Ulrich, in particular slope-stable.

\textbf{Step 2:Irreducibility}\\   
   Since $h^i(\Sigma,\os_\Sigma(r-1))=0$ for $i>0$, by Riemann--Roch 
   \[h^0(\Sigma,\os_\Sigma(r-1))=\chi(\os_\Sigma(r-1))=1+r(r-1)=r^2-r+1.\]
   Therefore, if we choose $T\subset \Sigma$ a general subscheme of dimension $0$ and length $r^2$ then $h^0(\Sigma,\id_T(r-1))=0$ and by \thref{lemmafine} we get $ext^1(\id_T(r-1),\os_\Sigma)=r-1$.
   Consider the Hilbert scheme $\Sigma^{[r^2]}$ parametrising all the dimension $0$ and length $r^2$ subschemes of $\Sigma$; it is smooth and irreducible of dimension $2r^2$.
   We can apply \thref{modularcorrispondenza} to the universal family on $\Sigma^{[r^2]}$ and get a family $\scr{E}$ of rank $r$ sheaves, corresponding to the pairs $(\id_T(r),Ext^1(\id_{T}(r),\os_{\Sigma}))$ through Hartshorne-Serre, over some open cover of $\Sigma^{[r^2]}$.
   By \thref{lemmafine} any special Ulrich bundle belongs to $\scr{E}$, hence the general sheaf in this family is semistable.
   Then \thref{mappamodulare} implies that $\Sigma^{[r^2]}$, which is irreducible, dominates the locus of Ulrich bundles in the moduli space as desired.
\end{proof}

We remark that the existence of Ulrich bundles on $\Sigma$ can be easily deduced just by restricting Ulrich bundles on a quartic double solid having $\Sigma$ as an hyperplane section, see \thref{restrizione}.
The point of Step $1$ of the above proof is to ensure that they are stable.

\subsection{Restriction to hyperplane sections}

Here, we want to study restrictions of Ulrich bundles from $X$ to divisors $\Sigma\in|\os_X(1)|$.
Recall that the restriction of an Ulrich bundle to an hyperplane section stays Ulrich by \thref{restrizione}, in particular semistable.
Our main result is that, for a general $X$, the moduli space of Ulrich bundles on $X$ dominates the corresponding one for $\Sigma$.
As far as we know, the only other example of such a phenomenon happens for $X$ a cubic $3$-fold, see \cite{CH}[Cor. 5.12], and relies on some computer-aided computations.
If we restrict a stable bundle on $X$ to $\Sigma$ we are not sure it remains stable, so we prefer to work with simple ones in order to make some local computation.

\Le\label{restretale}
Let $\e$ be a stable rank $2$ Ulrich bundle on a smooth $X$ and fitting in 
\sesl{\os_X}{\e}{\id_Y(2)}{sezione}
then the following are equivalent:
\begin{enumerate}[1)]
    \item $\e|_\Sigma$ is simple and restriction induces an isomorphism $Ext^1(\e,\e)\cong Ext^1(\e|_\Sigma,\e|_\Sigma)$
    \item $H^1(X,\e\otimes\e\du(-1))\cong H^2(X,\e\otimes\e\du(-1))=0$
    \item $\e|_Y\cong \n_{Y}$ is Ulrich on $(Y,\os_Y(1))$ where $\os_Y(1):=\os_X(1)|_Y$
    \item $H^1(Y,\e|_Y(-1))=0$. 
    \end{enumerate}
\ma
\begin{proof} 
\gel{1} $\Longleftrightarrow$ \gel{2}
Since $\e$ is locally free we have $Ext^1(\e,\e)\cong H^1(X,\e\otimes\e\du)$ and the same applies to $\e|_\Sigma$.
Moreover, by Serre duality we have $h^1(X,\e\otimes\e\du(-1))=h^2(X,\e\otimes\e\du(-1))$.
    Therefore, the cohomology sequence of
    \sesl{\e\otimes\e\du(-1)}{\e\otimes\e\du}{\e\otimes\e\du|_\Sigma}{rsigma}
     gives the claim, recalling that $\e$ is stable hence simple.
    
    \gel{2} $\Longleftrightarrow$ \gel{3}
    Twisting \eqref{sezione} by $\e(-1)\du\cong \e(-3)$ we get 
    \ses{\e(-3)}{\e\otimes\e\du(-1)}{\e\otimes\id_Y(-1)}
    moreover, we have also the sequence 
    \ses{\e\otimes\id_Y(-1)}{\e(-1)}{\e(-1)|_Y}
    Since $\e$ is Ulrich, both $\e(-1)$ and $\e(-3)$ have no cohomology hence it follows that for $i=1,2$ we have 
    \[H^i(X,\e\otimes\e\du(-1))\cong H^i(X,\e\otimes\id_Y(-1))\cong H^{i-1}(X,\e|_Y(-1)).\]

    \gel{3} $\Longrightarrow$ \gel{4}
    This is clear by definition.

    \gel{4} $\Longrightarrow$ \gel{3}
   We just need to prove that $H^0(Y,\e|_Y(-1))=0$ which is equivalent to $\chi(\e|_Y(-1))=0$.
   $Y$ has arithmetic genus $1$, hence the thesis follows by Hirzebruch--Riemann--Roch and the fact that $c_1(\e(-1))=0$.
\end{proof}

In particular, we conclude that properties $3)$ and $4)$, which seem to depend on the particular curve $Y$, actually depend only on the vector bundle $\e$.

Let $\Sigma$ be any divisor in $|\os_X(1)|$.
Let $\cu_r$ be an irreducible component of the open subset of the moduli space parametrising rank $r$ stable Ulrich bundles on $X$ constructed in \thref{rango>sextic} or \thref{qdswild}.
Note that all those sheaves are special, being $Pic(X)\cong\Z$.
We have a restriction morphism $\rho_\Sigma:\cu_r\ra \cm^{ss}_r$, where $\cm^{ss}_r$ is the unique irreducible component of the moduli space of semistable sheaves of rank $r$ on $\Sigma$ which contains special Ulrich bundles, see \thref{qdpwild}.

\Te\label{restrizdomi}
Let $X$ be a general complex quartic double solid and $\Sigma$ be a smooth divisor in $|\os_X(1)|$.
For any $r\ge 2$, the restriction map $\rho_\Sigma:\cu_r\ra \cm^{ss}_r$, which sends Ulrich bundles of rank $r$ on $X$ to Ulrich bundles of rank $r$ on $S$ is generically étale.
\Ma
\begin{proof}
Let $\e=\e_r$ be a stable Ulrich bundle in $\cu_r$.
Call $\cu^\circ$ the open subset of $\cu_r$ where $H^1(X,\e\otimes\e\du(-1))=0$.
We will prove the statement supposing that $\cu^\circ\neq \emptyset$ and then show that such a claim holds for a general $X$.

\textbf{Step $1$: Reduction to cohomology vanishing}\\
From \thref{restretale} we know that if $\e$ gives a point in $\cu^\circ$ then $\e|_\Sigma$ is a simple Ulrich bundle, hence on $\cu^\circ$ we have a family of simple, special Ulrich bundles on $\Sigma$.
Therefore, up to replacing $\cu^\circ$ with an étale base change, we have a morphism $\rho_\Sigma':\cu^\circ\ra\scr{S}$, where $\scr{S}$ is any modular family of simple sheaves on $\Sigma$ containing special Ulrich bundles of rank $r$.

Recall that, the tangent space to $\cu_r$, hence also to $\cu^\circ$, in $[\e]$ is $Ext^1(\e,\e)$ by \cite{HuyLeh}[Cor. 5.4.2] while the tangent spaces to the modular family in $[\e|_\Sigma]$ is $Ext^1(\e|_\Sigma,\e|_\Sigma)$ by definition of this family, see \autoref{defstability}.
If $j:\Sigma\ra X$, then the differential of $\rho_\Sigma'$ agrees with the map 
\[H^1(X,\e\otimes\e\du)\ra H^1(\Sigma,\e_\Sigma\otimes\e_\Sigma\du)\cong H^1(X,j_*(\e_\Sigma\otimes\e_\Sigma\du))\]
induced by restriction to $\Sigma$, but this is exactly the one coming from the sequence \eqref{rsigma} and we already know that, under our assumption, it is an isomorphism by \thref{restretale}.
Hence $\rho_\Sigma'$ is étale around $\e$, in particular this map hits the stable locus in $\scr{S}$ because its image is dense and the latter is open. 
But, if both $\e$ and $\e|_\Sigma$ are stable then there are some neighbourhoods, containing only stable bundles, of those points on which both $\rho_\Sigma$ and $\rho_\Sigma'$ are defined and agree.
This proves the claim, under our initial assumption.

\noindent\textbf{Step $2$: Rank $2$ case}\\
We are only left to prove that the locus of rank $r$ Ulrich bundles such that $h^1(X,\e_r\otimes\e_r\du(-1))=0$ is non-empty for a general $X$.
We start with $r=2$.
We can assume $X$ smooth hence, from \thref{rk2-section} any Ulrich bundle $\e$ fits in a sequence like \eqref{sezione} so it is enough to prove that $h^1(Y,\n_Y(-1))=0$.
     Once we fixed $Y$ we also have the equation of the branch locus of $f:X\ra\p^3$ in the form $p_0^2+p_1p_2+p_3p_4$ by \thref{branchrango2}.
     Then, we get the exact sequence form \thref{contonormali} which, twisted by $\os_X(-1)$ reads 
     \ses{\n_{Y/X}(-1)}{\os_Y(1)^3}{\os_Y(3)}
where the right map is given by taking scalar product with $(t+p_0,p_2,p_4)$ hence
\[Im\Big(H^0(Y,\os_Y(1)^3)\ra H^1(Y,\os_Y(3))\Big)\cong \Big(\C[x_0,x_1,x_2,x_3]/(p_0,p_1,p_2,p_3,p_4)\Big)_{3},\]
Being $h^1(Y,\os_Y(1))=0$ we deduce that $h^1(Y,\n_Y(-1))=0$ if and only if the homogenous degree $3$ part of $\C[x,y,z,w]/(p_0,p_1,p_2,p_3,p_4)$ is $0$.
    The same strategy used in \thref{join} holds: if we can prove this claim for a special choice of  $p_0,p_1,p_2,p_3,p_4$ then it holds for a general one.
    Clearly the ideal $(x^2,y^2,z^2,w^2,xy+zw)$ contains all monomials of the types $x^3$ and $x^2y$.
    Moreover, $xyz=z(xy+zw)-w(z^2)$ and similarly for the other triple products.
    Those monomials are a basis of the space of degree $3$ polynomials hence our claim follows.

\noindent\textbf{Step $3$: Rank $3$ case}\\
By \thref{qdswild} Ulrich bundles of rank $3$ are obtained as deformation of extensions 
\ses{\e}{\f}{\id_l(1)}
where $l\subset X$ is a line and $\e$ is a rank $2$ Ulrich bundle. 
We prove that if $l$ and $\e$ are general enough then
\begin{equation}\label{rango3}
ext^1(\f,\f(-1))=ext^1(\f,\id_l(1))=ext^1(\id_l(1),\f(-1))=0
\end{equation}
so by semicontinuity the first vanishing will hold for a general rank $3$ Ulrich bundle, as desired.
Again by semicontinuity, it is enough to show the thesis on a split extension $\f=\e\oplus\id_l(1)$.
In this case $ext^1(\f,\f(-1))$ is the sum of four terms but 
\[ext^1(\e,\e(-1))=ext^1(\e,\id_l)=ext^1(\id_l(1),\e(-1))=0\]
by the previous Step and the last claim of \thref{calcoliretta}, which applies by generality of $\e$.
Therefore, we only need to check $ext^1(\id_l(1),\id_l)=0$.
Applying $Hom(\id_l(1),-)$ to the sequence
\ses{\id_l}{\os_X}{\os_l}
we get 
\[ext^1(\id_l(1),\id_l)=hom(\id_l(1),\os_l)=h^0(X,\hom(\id_l(1),\os_l))=h^0(X,\n_l(-1))=0\]
by \thref{normalerette} since $l$ can be chosen general.
The other two claims follow immediately by the same reasoning since 
\[ext^1(\f,\id_l)=ext^1(\id_l(1),\id_l)+ext^1(\e,\id_l)=0\]
\[ext^1(\id_l(1),\f)=ext^1(\id_l(1),\id_l(1))+ext^1(\e,\id_l(1))=0.\]

\noindent\textbf{Step $4$: General case}\\ 
We argue by induction on $r$ and we want to prove the following: for a general rank $r$ Ulrich bundle $\e_r$ we have 
\[ext^1(\e_r,\e_r(-1))=ext^1(\e_r,\e(-1))=ext^1(\e,\e_r(-1))=0\]
where $\e$ is a general rank $2$ Ulrich bundle in $\cu_2$.
The base cases are $r=2,3$ done in the previous steps.
By \thref{rango>sextic} and \thref{qdswild} Ulrich bundles $\e_r$ of rank $r>3$ are obtained as deformation of extensions 
\ses{\e}{\f}{\g}
where $\e,\g$ are rank $2,r-2$ Ulrich bundles corresponding to a point in $\cu_{2},\cu_{r-2}$.
In particular, if we prove 
\[ext^1(\f,\f(-1))=ext^1(\f,\e(-1))=ext^1(\e,\f(-1))=0\]
for a split extension $\f\cong \e\oplus\g$ with $\e,\g$ general then, by semicontinuity, those vanishings will hold for a general Ulrich bundle of rank $r$ and the first one gives us the thesis.
Now $ext^1(\f,\f(-1))$ is the sum of four numbers but $ext^1(\e,\e(-1))=0$ by Step $2$, being $\e$ general, while
\[ext^1(\g,\g(-1))=ext^1(\g,\e(-1))=ext^1(\e,\g(-1))=0\]
by inductive assumption.
Similarly, we compute $ext^1(\g,\e(-1))=0=ext^1(\e,\g(-1))$.
\end{proof}

\Oss
We are not able to compute the generic degree of $\rho_\Sigma$.
Nevertheless, reasoning as in \thref{restrdoppi}, it can be seen that if $\e|_\Sigma\cong\g|_\Sigma$ and $\g\neq \e$ then we must have a non-zero element in $Hom(\iota^*\e(-1),\g)$.
Clearly this degree is at least $2$ since $\e|_\Sigma\cong\iota^*\e|_{\Sigma}$.
\one

We already noticed that $\id_l(1)$ somehow plays the role of Ulrich bundles of rank $1$.
We know that $X$ is covered by lines, that is, for any point in $X$ we can find a line that passes through it.

\Co\label{interpolazione}
Let $X$ be a general complex quartic double solid.
For any $r\ge 2$, given $r^2$ general points on some smooth $\Sigma\in |\os_X(1)|$ there is a smooth degree $r^2$ and genus $\frac{2}{3}r(r^2-3r+2)+(r-1)^2$ projectively normal curve passing through them.
\io
\begin{proof}
    By \thref{qdpwild}, for a general $T\subset \Sigma$ of dimension $0$ and degree $r^2$ we find a sequence
    \sesl{V\otimes\os_\Sigma}{\e_\Sigma}{\id_T(r)}{solitasigma}
    with $\e_\Sigma$ a general Ulrich bundle of rank $r$ on $\Sigma$ and $V\subset H^0(\Sigma,\e_\Sigma)$ general.
    From \thref{restrizdomi} we know that there is an Ulrich bundle $\e$ on $X$ such that $\e|_\Sigma\cong \e_\Sigma$.
    Since the restriction morphism $H^0(X,\e)\ra H^0(\Sigma,\e_\Sigma)$ is an isomorphism by \thref{restrizionesezioni}, $V$ gives a general $(r-1)$-dimensional subspace of $H^0(X,\e)$ hence we find an exact sequence 
    \sesl{V\otimes\os_X}{\e}{\id_Y(r)}{solitax}
    which restricts to \eqref{solitasigma}.
    In particular, $\id_Y|_{\Sigma}\cong \id_T$ hence $Y\cap\Sigma=T$ schematically.
    A general such $Y$ is smooth, see for example \thref{r1}, and its degree and genus can be computed from \cite{CFK3}[Thm. 3.5].
    Projective normality follows from the cohomology of twists of \eqref{solitax} once we recall that $(X,\os_X(1))$ is aCM and an Ulrich bundle is also aCM by \thref{definizioneUlrich}.
\end{proof}

\newpage

\printbibliography

\end{document}